\documentclass[10pt,reqno]{amsart}

\usepackage{amsmath}
\usepackage{amsthm}
\usepackage{amsfonts}
\usepackage{amssymb}
\usepackage{mathrsfs}
\usepackage{graphicx}
\usepackage[T1]{fontenc}
\usepackage[OT2,T1]{fontenc}
\usepackage[latin1]{inputenc}
\usepackage{lmodern}
\usepackage[all,cmtip]{xy}
\usepackage{cancel}
\usepackage[normalem]{ulem}
\usepackage{bbm}
\usepackage[vcentermath]{youngtab}
\usepackage{ stmaryrd }
\usepackage{calc}
\usepackage{empheq}
\usepackage{xcolor}
\usepackage{setspace}
\usepackage{epigraph}
\usepackage{mathdots}
\usepackage{braket}
\usepackage{longtable}
\usepackage{float}
\usepackage{changepage}
\usepackage{pdflscape}
\usepackage{enumerate}

\usepackage[pdftex,     
          plainpages=false,   
           breaklinks=true,    
           colorlinks=true,
           pdftitle=My Document
           pdfauthor=My Good Self
           colorlinks=true,
	    urlcolor=blue,
	    citecolor=red,
	    linkcolor=blue
          ]{hyperref}
          
\usepackage[a4paper,top=4cm,bottom=4cm,left=2.8cm,right=2.8cm]{geometry}

\theoremstyle{plain}
\newtheorem{teorema}{Theorem}[section]

\theoremstyle{definition}
\newtheorem{defi}[teorema]{Definition}

\newtheorem{es}[teorema]{Example}
\theoremstyle{plain}
\newtheorem{conj}[teorema]{Conjecture}
\newtheorem{cor}[teorema]{Corollary}
\newtheorem{lemma}[teorema]{Lemma}
\newtheorem{prop}[teorema]{Proposition}
\theoremstyle{remark}
\newtheorem{oss}[teorema]{Remark}

\newenvironment{sistema}
{\left\lbrace\begin{array}{@{}l@{}}}
{\end{array}\right.}
\newlanguage\fakelanguage
\newcommand\cyr{\fontencoding{OT2}\fontfamily{wncyr}\selectfont
   \language\fakelanguage}
\DeclareTextFontCommand{\textcyr}{\cyr}

\numberwithin{equation}{section}

\DeclareMathOperator{\Hom}{Hom}
\DeclareMathOperator{\End}{End}
\DeclareMathOperator{\Ob}{Ob}

\usepackage{calligra}

\DeclareMathOperator{\HOM}{\mathscr{H}\text{\kern -3pt {\calligra\large om}}\,}

\makeatletter
\newsavebox{\@brx}
\newcommand{\llangle}[1][]{\savebox{\@brx}{\(\m@th{#1\langle}\)}%
  \mathopen{\copy\@brx\kern-0.5\wd\@brx\usebox{\@brx}}}
\newcommand{\rrangle}[1][]{\savebox{\@brx}{\(\m@th{#1\rangle}\)}%
  \mathclose{\copy\@brx\kern-0.5\wd\@brx\usebox{\@brx}}}
\makeatother

\usepackage{framed,color}
\definecolor{shadecolor}{rgb}{0.92, 0.92, 0.92}
\definecolor{ambra}{rgb}{1.0, 0.75, 0.0}
\definecolor{ametista}{rgb}{0.6, 0.4, 0.8}
\definecolor{auburn}{rgb}{0.43, 0.21, 0.1}
\definecolor{ballblue}{rgb}{0.13, 0.67, 0.8}
\definecolor{cadmiumgreen}{rgb}{0.0, 0.42, 0.24}
\definecolor{candypink}{rgb}{0.89, 0.44, 0.48}
\definecolor{caribbeangreen}{rgb}{0.0, 0.8, 0.6}
\definecolor{airforceblue}{rgb}{0.36, 0.54, 0.66}

\newcommand{\bsh}{\begin{shaded}}
\newcommand{\esh}{\end{shaded}}

\begin{document}
\title{Helix structures in Quantum Cohomology of Fano Varieties}

\author{G. Cotti$^{(\ddag)}$, B. Dubrovin$^{(\dag)}$, D. Guzzetti$^{(\dag)}$}

\begin{abstract}
In this paper we consider a conjecture formulated by the second author in occasion of the 1998 ICM in Berlin (\cite{dubro0}). This conjecture states the equivalence, for a Fano variety $X$, of the semisimplicity condition for the quantum cohomology $QH^\bullet(X)$ with the existence condition of full exceptional collections in the derived category of coherent sheaves $\mathcal D^b(X)$. Furthermore, in its quantitative formulation, the conjecture also prescribes an explicit relationship between the monodromy data of $QH^\bullet(X)$ and characteristic classes of both $X$ and objects of the exceptional collections. In this paper we reformulate a refinement of \cite{dubro0}, which corrects the \emph{ansatz} of \cite{dubro4} for what concerns the conjectural expression of the central connection matrix. We clarify the precise relationship between the  refined conjecture presented in this paper and $\Gamma$-conjecture II of S. Galkin, V. Golyshev and H. Iritani (\cite{gamma1,gamma2}). Through an explicit computation of the monodromy data and a detailed analysis of the action of the braid group on both the monodromy data and the set of exceptional collections, we prove the validity of our refined conjecture for all complex Grassmannians $\mathbb G(r,k)$. From these results, it is outlined an explicit description of the ``geography'' of the exceptional collections realizable at points of the small quantum cohomology of Grassmannians, i.e. corresponding to the monodromy data at these points. In particular, it is proved that Kapranov's exceptional collection appears at points of the small quantum cohomology only for Grassmannians of small dimension (namely, less or equal than 2). Finally, a property of \emph{quasi-periodicity} of the Stokes matrices of complex Grassmannians, along the locus of the small quantum cohomology, is described.
\end{abstract}

\address[$\ddag$]{Max-Planck-Institut für Mathematik, Vivatsgasse 7 - 53111 Bonn, GERMANY}
\address[$\dag$]{SISSA, Via Bonomea 265 - 34136 Trieste ITALY. 
\begin{enumerate}
\item[] Giordano Cotti's ORCID ID: 0000-0002-9171-2569
\item[] Boris Dubrovin's ORCID ID: 0000-0001-9856-5883
\item[] Davide Guzzetti's ORCID ID: 0000-0002-6103-6563
\end{enumerate}}
\email[G. Cotti]{gcotti@sissa.it, gcotti@mpim-bonn.mpg.de}
\email[B. Dubrovin]{dubrovin@sissa.it}
\email[D. Guzzetti]{guzzetti@sissa.it}

\maketitle
\tableofcontents


\section{Introduction}
In this paper we consider the conjecture formulated by the second author in occasion of the 1998 ICM in Berlin \cite{dubro0}.  The primary and genuine aim of the conjecture was a characterization of smooth projective Fano varieties having semisimple quantum cohomology and the computation of the monodromy data of the corresponding semisimple Frobenius manifolds in algebro-geometric terms. 

The original version of the conjecture can be described in two different parts, a \emph{qualitative} and a \emph{quantitative} one. In the first qualitative part, for a smooth projective Fano variety $X$,  the semisimplicity condition of the quantum cohomology $QH^\bullet(X)$ is conjectured to be equivalent to existence of a full exceptional collection in $\mathcal D^b(X)$, the derived category of coherent sheaves on $X$. It consists of an ordered collection $\frak E=(E_1,\dots, E_n)$ of objects in $\mathcal D^b(X)$ satisfying the semi-orthogonality conditions 
\[\Hom^\bullet(E_i,E_i)\cong \mathbb C,
\]
\[\Hom^\bullet(E_j,E_i)\cong 0,\quad \text{if }j>i.
\]
Furthermore, in order to be \emph{full}, the collection $\frak E$ must generate the category $\mathcal D^b(X)$ as a triangulated one. The second quantitative (and maybe the most important) part of the conjecture predicates an \emph{analytic} and \emph{explicit} relationship between the monodromy data of the quantum cohomology $QH^\bullet(X)$ and algebro-geometric data of the objects of $\frak E$. Such a relationship is based on the analytic theory of Frobenius manifolds, as introduced and extensively developed in \cite{dubro0, dubro1, dubro2} and further refined in \cite{CDG} in order to include semisimple coalescence points.  

The core of this analytic theory consists of a local identification of the semisimple part of a Frobenius manifold (in this paper $QH^\bullet(X)$) with the space of parameters of isomonodromic deformations of differential operators with rational coefficients of the form
\begin{equation}
\label{opdiffu}
\Lambda(z,u):=\frac{d}{dz}-\left(U+\frac{1}{z}V(u)\right),
\end{equation}
the matrix $U={\rm diag}(u_1,\dots, u_n)$ being a diagonal matrix and $V(u)$ antisymmetric. The entries of $U$, which define a system of canonical coordinates at semisimple points of the Frobenius manifold (the \emph{canonical coordinates}, see Section \ref{psigauge}), are the parameters of deformation of the family of differential operators. The dependence of $V$ on the parameters $u$ is such that the system  $\Lambda Y=0$ admits fundamental matrix solutions $Y(z,u)$ with monodromy data independent of small deformations of $u$ (this is the \emph{Isomonodromicity Property} described in Theorems \ref{06.07.17-4}, \ref{iso2}, \ref{iso2.1}). The geometrical condition which implies the isomonodromicity of the deformation \eqref{opdiffu} is the flatness of the \emph{extended deformed  connection} $\widehat\nabla$ associated with  the Frobenius manifold $QH^\bullet(X)$ (see Section \ref{extdefconn}).

Through such a description, some local invariants of the Frobenius structure called \emph{monodromy data} are introduced (and defined up to some non-canonical choices, as summarized in Section \ref{monlocalmod}), namely the \emph{Stokes} and \emph{central connection} matrices $(S,C)$. Vice versa, knowing the monodromy data $(S,C)$ one can reconstruct, by solving certain Riemann--Hilbert problem, the local structure of the semisimple Frobenius manifold (see below for the details). The analytic continuation of the local Frobenius structure can be completely described in terms of an action of the braid group on the monodromy data.

Let us briefly recall the definition of the monodromy data. In correspondence with a point $p$ in $QH^\bullet(X)$, with canonical coordinates $ u=( u_1,\dots,  u_n)$, we consider three solutions $Y_{\rm left/right/0}(z,u)$ of the differential system $\Lambda(z,u)Y(z,u)=0$.  The solutions $Y_{\rm left/right}(z,u)$ are uniquely determined by their asymptotic expansion 
\[Y_{\rm left/right}(z, u)\sim F(z, u)\exp(z U),
\]
for $|z|\to \infty$ in sectors of  $\mathcal R:=\widetilde{\mathbb C\setminus{0}}$ (the universal cover) with angular opening greater than $\pi$, respectively  containing the two half-planes 
\[\Pi_{\operatorname{right}} (\phi):=\left\{z\in\mathcal R\colon \phi-\pi< \arg\ z< \phi\right\},\quad \Pi_{\operatorname{left}} (\phi):=\left\{z\in\mathcal R\colon \phi< \arg\ z< \phi+\pi\right\},
\quad  \phi\in[0;2\pi[.
\] 
Here $F(z, u)$ is a formal series satisfying
\[
F(z,u)=\sum_{m\geq 0}F_m(u)z^{-m},\quad F_0\equiv\mathbbm 1, \quad F(-z,u)^T\cdot F(z,u)=\mathbbm 1.
\]

The solution  $Y_0(z,u)$ is a fundamental matrix in {\it Levelt form} 
\begin{equation}
\label{16luglio2018-2}
Y_0(z,u)= \Psi(u)(1+O(z))z^\mu z^R,
\end{equation}
where $ \Psi(u)$is a matrix describing a change of frames in the holomorphic tangent bundle of $QH^\bullet(X)$, namely
\begin{equation}\label{03.05.18-1}\frac{\partial}{\partial t^\alpha}=\sum_i\Psi_{i\alpha}(u)\cdot \eta\left(\frac{\partial}{\partial u_i},\frac{\partial}{\partial u_i}\right)^{-\frac{1}{2}}\frac{\partial}{\partial u_i},
\end{equation} 
so that it takes $V(u)$ to the diagonal form $\mu$;  $(1+O(z))$ stands for a convergent Taylor expansion at $z=0$, and   $R$ is a nilpotent matrix. Among the fundamental matrices in Levelt form, we select the so called   \emph{enumerative-topological solution} $Y_{\rm top}(z,u)$, see below, because of its important geometrical interpretation. Indeed,  $Y_{\rm top}(z,u)$   admits   an expansion given by a generating function of genus 0 Gromov--Witten invariants of $X$ with their gravitational descendants. If $(T_\alpha)_\alpha$ is a fixed homogeneous basis of $H^\bullet(X,\mathbb C)$, and $(t^\alpha)_{\alpha=1}^n$ are the  flat coordinates  with respect to the Poincaré metric $\eta$,
\[\eta(v_1,v_2):=\int_Xv_1\cup v_2,\quad v_1,v_2\in H^\bullet(X,\mathbb C),
\]then 
\begin{align*}
& Y_{\rm top}(z,u)=\Psi(u)\cdot Z_{\rm top}(z,t(u)),
\\
& Z_{\rm top}(z,t):=\Theta_{\rm top}(z,t)\cdot z^\mu z^{R},
 \end{align*}
where 
\begin{equation}
\label{16agosto2018-3}
 \mu\in{\rm End}_{\mathbb C}(H^\bullet(X,\mathbb C)),\quad \mu(T_\alpha):=\frac{\deg T_\alpha-\dim_{\mathbb C}X}{2}\cdot T_\alpha,\quad 
 R=c_1(X)\cup.
\end{equation}
Moreover,  $\Theta_{\rm top}(z,t)$ is  a convergent Taylor series in variable $z$
$$
\Theta_{\rm top}(z,t)_\lambda^\gamma:=\delta_\lambda^\gamma+\sum_{n=0}^\infty \left(\sum_{k=0}^\infty\sum_{\beta\in\text{Eff}(X)}\sum_{\alpha_1,\dots,\alpha_k}\frac{h_{\lambda,k,n,\beta,\underline{\alpha}}^\gamma}{k!}\cdot t^{\alpha_1}\dots t^{\alpha_k}\right)\cdot z^{n+1},
$$
with $t$-dependent coefficients, defined by 
$$
h_{\lambda,k,n,\beta,\underline{\alpha}}^\gamma:=\sum_{\delta}\eta^{\delta\gamma}\int_{[\overline{\mathcal M}_{0,k+2}(X,\beta)]^\text{virt}}c_1(\mathcal L_1)^n\cup\text{ev}_1^*T_\lambda\cup\text{ev}_2^*T_{\delta}\cup\prod_{j=1}^k\text{ev}^*_{j+2}T_{\alpha_j},
$$
where $\mathcal L_1$ is the $1$-st tautological cotangent bundle on the Deligne--Mumford stack
 $\overline{\mathcal M}_{0,k+2}(X,\beta)$ of stable maps of genus 0, with degree
  $\beta$ and $(k+2)$-punctures and target space $X$. As explained above, the topological-enumerative
   solution is an element of a class of solutions of $\Lambda(z,u)Y=0$ in Levelt form: 
   all other such solutions are parametrized by the points of an algebraic group, in this
    paper\footnote{In the notations of \cite{CDG}, such a group would be denoted by $\mathcal C_0(\eta,\mu, R)$, where $\eta$ denotes the Poincaré metric on $H^\bullet(X,\mathbb C)$, $\mu$ is the
     grading operator, and $R$ denotes the endomorphism of $H^\bullet(X,\mathbb C)$ of $\cup$-multiplication by $c_1(X)$.} denoted by $\mathcal C_0(X)$ (see Sections \ref{extdefconn} for precise definitions), 
     and they are given by
\begin{equation}
\label{16agosto2018-1}
Y_0(z,u)A,\quad A\in\mathcal C_0(X),
\end{equation}
where $Y_0$ has been fixed, for example choosing $Y_0=Y_{\rm top}$. 

The Stokes and central connection matrices $(S,C)$ are defined through the equations
\[Y_{\rm right}(z, u)=Y_{0}(z, u)C,\quad Y_{\rm left}(z, u)=Y_{\rm right}(z, u)S,\quad \text{for all }z\in\mathcal R.\]
The Isomonodromicity property guarantees that the matrices $(S,C)$, seen as functions of the parameters $u$ (i.e. of the point of $QH^\bullet(X)$), are constant in open regions of the Frobenius manifold $QH^\bullet(X)$, called \emph{$\ell$-chambers} (here $\ell$ stands for an oriented line in the complex plane with slope $\phi\in[0;2\pi[$; see Definition \ref{chamberdefi}). The values of these numerical invariants in different chambers can be related through the action of the braid group $\mathcal B_n$ on the pair $(S,C)$ (see Section \ref{monlocalmod}).

 The quantitative part of the conjecture mentioned above gives an exact prescription of these invariants in terms of characteristic classes of the objects of an exceptional collection $\frak E$.

In \cite{dubro0} the Stokes matrix $S$, computed at any point $p\in QH^\bullet(X)$ (with respect to any choice of orthonormalized idempotents of Section \ref{psigauge}, any oriented line $\ell$ in the complex plane of slope $\phi\in[0;2\pi[$, and a suitable order of canonical coordinates, the \emph{$\ell$-lexicographical one}; see Section \ref{monlocalmod} for details) was conjectured to be equal to the Gram matrix of the Grothendieck--Euler--Poincaré product
\[\chi(E,F):=\sum_{i}(-1)^i\dim_\mathbb C\Hom^i(E,F),\quad E,F\in \Ob(\mathcal D^b(X)),
\]associated with some exceptional collection $\frak E$. For what concerns the central connection matrix $C$,  the original formulation of the conjecture did not completely identify its geometrical counterpart in  $\mathcal D^b(X)$. The only observation appearing in \cite{dubro0} is an \emph{ansatz} for the general structure of the central connection matrix, namely 
\[C=C'\cdot C'',
\]where $C''$ is a matrix whose column entries are the components of the {graded} Chern character ${\rm Ch}(E_\ell)$ of the objects of $\frak E$, namely
\[{\rm Ch}(E_\ell):=\sum_{h=1}^{{\rm rk} E_\ell}\exp(2\pi i \delta_{\ell,h}),\quad\delta_{\ell,j}\text{'s being the Chern roots of }E_\ell,
\]and where $C'$ is a matrix  only required to commute with the operator of classical $\cup$-multiplication
\[c_1(X)\cup (-)\colon H^\bullet(X,\mathbb C)\to H^\bullet(X,\mathbb C).
\]

\subsection{Refinement of the conjecture} After the partial and incomplete formulation of the conjecture in \cite{dubro0}, a crucial insight for its refinement was suggested by L. Katzarkov, M. Kontsevich and T. Pantev in \cite{kkp}. In the context of a non-commutative Hodge theoretical extension of Homological Mirror Symmetry, the authors of \cite{kkp} firstly recognized in the entries of the central connection matrix for Projective Spaces $\mathbb P^n_{\mathbb C}$ (briefly shown in \cite{dubro0}, without details about their computation) the components of a characteristic class of $\mathbb P^n_{\mathbb C}$, obtained by applying the F. Hirzebruch's construction of multiplicative genera (\cite{hirz}) to the Taylor expansion of a $\Gamma$-function (see below for more details). The same characteristic class was the main ingredient used by Katzarkov, Kontsevich and Pantev for defining a \emph{rational} structure on the $A$-model $\bf{nc}$-Hodge structure associated with a compact symplectic manifold $(X,\omega)$.

Recently, the  necessity for a deeper understanding and refinement of the conjecture of \cite{dubro0} has increased. In such a direction, two main contributions require to be mentioned\footnote{In the very recent and interesting paper \cite{tar-var}, it is proved that the $\Gamma$-classes prescribe the leading terms of the asymptotic expansions of solutions of the \emph{equivariant} quantum differential equations for the cotangent bundles of partial flag varieties. We believe that further investigations are needed for connecting the results of the present work with the ones presented in \cite{tar-var} in the equivariant case. See also Remark \ref{sasha}.}. 
\begin{enumerate}
\item In \cite{dubro4}, the second author suggested that the column entries of the central connection matrix $C$ should be equal to the components of the characteristic classes  
\begin{equation}\label{classdubr}\frac{1}{(2\pi)^\frac{d}{2}}\widehat\Gamma^-_X\cup{\rm Ch}(E_i),\quad d=\dim_\mathbb CX,
\end{equation}$E_i$ being objects of an exceptional collection. 

\item Almost contemporarily to \cite{dubro4}, in the papers \cite{gamma1} and \cite{gamma2} S. Galkin, V. Golyshev and H. Iritani proposed a set of conjectures, called $\Gamma$-conjectures (I and II) describing the exponential asymptotic behaviour of flat sections of the \emph{quantum connection} (namely, the extended deformed connection $\widehat\nabla$ defined on $QH^\bullet(X)$ mentioned above). It is claimed that $\Gamma$-conjecture II  refines the conjecture of \cite{dubro0}, and it identifies the column entries of the central connection matrix, defined as above, with the components of the characteristic classes
\begin{equation}\label{classggi}
\frac{1}{(2\pi)^\frac{d}{2}}\widehat\Gamma^+_X\cup{\rm Ch}(E_i),\quad d=\dim_\mathbb CX,
\end{equation}
$E_i$ being objects of an exceptional collection.
\end{enumerate}
Here, $\widehat\Gamma^\pm_X$ denote the characteristic classes\footnote{Curiously enough, when we started to address the problem of looking for the compatibility of both proposals in \cite{dubro4},\cite{gamma1,gamma2}, we were further confused by an evident typo in \cite{kkp}. In the last sentence of the proof of Proposition 3.1, the series defining the $\Gamma$-class is said to be the Taylor series of $\Gamma(1+t)$, although in the rhs of the subsequent formula it appears the Taylor series of $\Gamma(1-t)$.} of $X$ associated with the Taylor power series, centered at $t=0$, of the functions $\Gamma(1\pm t)$ through the construction described by F. Hirzebruch (\cite{hirz}):
\[\widehat\Gamma^\pm_X:=\prod_{j=1}^{\dim_{\mathbb C}X}\Gamma(1\pm \delta_j),\quad\delta_j\text{ are the Chern roots of }TX,
\]
\[\Gamma(1- t)=\exp\left\{\gamma t+\sum_{n=2}^\infty\frac{\zeta(n)}{n}t^n\right\}=1+\gamma t+\cdots ,
\]where $\gamma$ denotes the Euler--Mascheroni constant, and $\zeta$ the Riemann zeta function. As observed in \cite{kkp} (see Remark 3.3 of \emph{loc. cit.}), the $\Gamma^\pm$-classes (and/or their $\cup$-multiplicative inverses) previously appeared in the literature in several contexts such as in the work \cite{kont99} on deformation quantization, in mirror symmetry for Calabi-Yau varieties \cite{libg}, in the description of the {\bf{nc}}-motives of the Landau--Ginzburg model of a toric Fano \cite{goly1,goly2}. Furthermore, the $\Gamma$-characteristic class of $X$ naturally appears in the study of the $\mathbb S^1$-equivariant geometry of the free loop space $\mathcal LX$, and consequently in Givental's equivariant Floer Theory (see \cite{rong}, \cite{iri9}, \cite{givehomgeo} and \cite{gamma2}).

Our explicit computations for the simple case of $\mathbb G(2,4)$, described in Section 6 of \cite{CDG}, suggested that both proposals of the conjecture formulated in \cite{dubro4} and \cite{gamma1, gamma2}  require some refinements and/or clarifications, at least as far as the central connection matrix $C$ is concerned. Indeed,  Theorem 6.2 of \cite{CDG} proves that the central connection matrix for $\mathbb G(2,4)$ can be of both forms \eqref{classdubr} and \eqref{classggi}, which belong to the same ${\mathcal C}_0(\mathbb G(2,4))$-orbit,  if computed with respect to two different solutions $Y_0(z,u)$ in Levelt form \eqref{16luglio2018-2} at $z=0$, related by the action  \eqref{16agosto2018-1},  \emph{none of the two coinciding with the topological-enumerative one}\footnote{This will be proved to hold true for any even dimensional smooth projective variety $X$. For odd dimensional varieties one has an extra factor $\sqrt{-1}$. See point (1) of Theorem \ref{teoequivintro} below.}. 


 In this paper, strong evidences are given in favor of the following general conjecture, which refines the last part of the conjecture in \cite{dubro0}, and which we prove for $\mathbb{G}(r,k)$.

\begin{shaded}
\begin{conj}[cf. Conjecture \ref{congettura}]
\label{conjeintro}
Let $X$ be a smooth Fano variety of Hodge--Tate type, i.e. for which\footnotemark
$$h^{p,q}(X)=0,\quad \text{if }p\neq q.$$  
\begin{enumerate}
\item The quantum cohomology $QH^\bullet(X)$ is semisimple if and only if there exists a full exceptional collection in the derived category of coherent sheaves $\mathcal D^b(X)$. 

\item If $QH^\bullet(X)$ is semisimple, then for any oriented line $\ell$ (of slope $\phi\in[0;2\pi[$) in the complex plane there is a correspondence between $\ell$-chambers of  $QH^\bullet(X)$ and founded helices, i.e. helices with a marked foundation $\frak E_\ell=(E_1,\dots, E_n)$ in the derived category $\mathcal D^b(X)$. 

\item The monodromy data of the system \eqref{opdiffu}, computed in a $\ell$-chamber $\Omega_\ell$, being $(u_1,...,u_n)$ in the lexicographical order, are related to the following geometric data of the corresponding exceptional collection $\frak E_\ell=(E_1,\dots, E_n)$ (the marked foundation):
\begin{enumerate}
\item the Stokes matrix is equal to the \emph{inverse} of the Gram matrix of the Grothendieck--Poincaré--Euler product on $K_0(X)_{\mathbb C}:=K_0(X)\otimes_\mathbb Z\mathbb C$, computed with respect to the exceptional basis $([E_i])_{i=1}^n$
\[ S^{-1}_{ij}=\chi(E_i, E_j);
\]
\item the Central Connection matrix $C$, defined by $Y_{\rm right}(z,u)=Y_{\rm top}(z,u)C$, coincides with the matrix associated with the $\mathbb C$-linear morphism 
\[\textnormal{\textcyr{D}}^-_X\colon K_0(X)_{\mathbb C}\to H^\bullet(X,\mathbb C)\colon E\mapsto \frac{i^{\bar d}}{(2\pi)^\frac{d}{2}}\widehat\Gamma^-_X\cup\exp(-\pi i c_1(X))\cup{\rm Ch}(E),
\]where $d=\dim_\mathbb C X$, and $\bar d$ is the residue class $d $ {\rm (mod 2)}. The matrix is computed with respect to the exceptional basis $([E_i])_{i=1}^n$ and any pre-fixed basis $(T_\alpha)_{\alpha=1}^n$ in cohomology (see Section \ref{notations}).
\end{enumerate}
\end{enumerate}
\end{conj}
\end{shaded}
\footnotetext[5]{Here $h^{p,q}(X):=\dim_{\mathbb C} H^q(X,\Omega_X^p)$, with $\Omega_X^p$ the sheaf of holomorphic $p$-forms on $X$, denotes the $(p,q)$-Hodge number of $X$.}

In Section \ref{chmainconj} we also show that the identifications between the monodromy data and the geometry of the derived category can be further enriched, according to the following result.

\begin{shaded}
\begin{teorema}[cf. Theorem \ref{03.08.17-1}]\label{theoconjintro} Let $X$ be a smooth Fano variety of Hodge--Tate type for which Conjecture \ref{conjeintro} holds true. Then, all admissible operations with the monodromy data have a geometrical counterpart in the derived category $\mathcal D^b(X)$, as summarized in Table \ref{enrichconjintro} at the end of this Introduction. In particular, we have the following:
\begin{enumerate}
\item Mutations of the monodromy data $(S,C)$ (see Def. \ref{10novembre2018-2}) correspond to mutations of the exceptional basis (see Def. \ref{10novembre2018-1}).

\item Different choices of branches of the square roots in \eqref{03.05.18-1}, and hence of the $\Psi$-matrix, correspond to shifts of objects of the exceptional collections. This reflects on an action of $(\mathbb Z/2\mathbb Z)^n$ on the objects of the exceptional basis.

\item The monodromy data $(S, C^{(k)})$ computed with respect to other fundamental matrix solutions $Y^{(k)}_{\rm left/right}$ of system \eqref{opdiffu}, having the prescribed asymptotic expansion in rotated sectors 
\[Y^{(k)}_{\rm left/right}(z,t)\sim Y_{\rm formal}(z,t),\quad z\in e^{2\pi i k}\Pi_{\rm left/right}(\phi), \quad |z|\to \infty,\quad k\in\mathbb Z,
\]uniformly in $t$, are associated as in points {\rm (3a)-(3b)} of Conjecture \ref{conjeintro}, with different foundations of the helix, related to the marked one by an iterated application of the Serre functor $(\omega_X\otimes - )[\dim_{\mathbb C}X]\colon \mathcal D^b(X)\to  \mathcal D^b(X)$. 

\item The group ${\mathcal C}_0(X)$ is isomorphic to a subgroup of the identity component of the isometry group ${\rm Isom}_{\mathbb C}(K_0(X)_\mathbb C, \chi)$: more precisely, the morphism
\[
{\mathcal C}_0(X)\to {\rm Isom}_{\mathbb C}(K_0(X)_\mathbb C,\chi)_0
\colon\quad
A\mapsto \left(\textnormal{\textcyr{D}}^-_X\right)^{-1}\circ A\circ \textnormal{\textcyr{D}}^-_X
\]
defines a monomorphism. In particular, ${\mathcal C}_0(X)$ is abelian.
\end{enumerate}  
\end{teorema}
\end{shaded}

Let us briefly summarize and clarify the exact relationships of $\Gamma$-conjecture II of S. Galkin, V. Golyshev and H. Iritani (\cite{gamma1,gamma2}) with our Conjecture \ref{conjeintro}. As mentioned above, the proposed formula \eqref{classggi} cannot be a compatible refinement of the original conjecture of \cite{dubro0}, essentially due to several different choices of normalizations, done in \cite{gamma1}, not completely standard in the theory of Frobenius manifolds. Let us underline the main differences:
\begin{enumerate}
\item in \cite{gamma1} another (flat) extended deformed connection, that we denote $\widehat\nabla^{\rm GGI}$, is considered on the Frobenius structure $QH^\bullet(X)$ (see Section \ref{relggi} for the precise definition). This connection can be identified with $\widehat\nabla$ only up to an identification of the spectral parameters\footnote{ We denote by $z$ and $\lambda$ our spectral parameter and the one in \cite{gamma1}, respectively.} $\lambda=z^{-1}$.

\item Despite this possible identification, the differential-geometrical meaning of the isomonodromic problem attached to $QH^\bullet(X)$ in \cite{gamma1,gamma2} is different from ours: it is defined as a flatness condition for a \emph{vector field} 
rather than for a differential form 
(the differential of a deformed flat coordinate). This implies that the isomonodromic system discussed in \cite{gamma1,gamma2} can be identified with the differential system  \eqref{opdiffu}  above only up to an identification of the spectral parameters given by 
\begin{equation}
\label{identsp}
\lambda=e^{\pm\pi i}z^{-1}.
\end{equation}
More precisely, the equation \eqref{ggi}
\footnote{For the convenience of the reader, when we discuss the relationship of our results with those of \cite{gamma1,gamma2}, if we want to keep the notations of S. Galkin, V. Golyshev and H. Iritani, not always coinciding with ours, we use an upper-script GGI. So, e.g. the matrix $Y^{\rm GGI}$ has as columns the components of vector fields with respect to the coordinated vectors $\frac{\partial}{\partial t^\alpha}$, whereas the column-entries of our matrix $Y$ are always intended to be components with respect to the orthonormalized idempotent vectors (see Section \ref{psigauge})
\[\eta\left(\frac{\partial}{\partial u_i}, \frac{\partial}{\partial u_i}\right)^{-\frac{1}{2}}\frac{\partial}{\partial u_i}.
\]}
\[\widehat\nabla^{\rm GGI}_{\frac{\partial}{\partial\lambda}}Y^{\rm GGI}=0,\]
is identified via \eqref{identsp} with the equation
\[\Lambda (\Psi\cdot Y^{\rm GGI})=0,
\]
with
\[ Y^{\rm GGI}=(Y^{\rm GGI,\alpha}_i)_{\alpha,i},\quad Y^{\rm GGI}_i=\sum_\alpha Y_i^{\rm GGI,\alpha}\frac{\partial}{\partial t^\alpha}\in\Gamma(\pi^*TQH^\bullet(X)),\] where $\pi\colon\mathbb C^*\times QH^\bullet(X)\to  QH^\bullet(X)$.
\item Furthermore, the solution $Y^{\rm GGI}_{\rm top}(\lambda)$ with respect to which the central connection matrix is computed in \cite{gamma1,gamma2} is analogous to our topological-enumerative solution, but it does not coincide with its specialization under the identification \eqref{identsp}. Namely, the solution $Y^{\rm GGI}_{\rm top}(\lambda)$ is in a \emph{different Levelt form}, at $0\in QH^\bullet(X)$ being equal to
\begin{align*}
Y^{\rm GGI}_{\rm top}(\lambda)&=\Theta_{\rm top}(-\lambda^{-1},0)\lambda^{-\mu}\lambda^{c_1(X)\cup -}
\\
&=\Theta_{\rm top}(-\lambda^{-1},0)e^{\mp i\pi\mu}\underbrace{z^{\mu}e^{\pm i\pi c_1(X)}z^{-\mu}}_{\text{ polynomial in }z}\cdot z^\mu z^{-(c_1(X)\cup -)}.
\end{align*}
In particular, notice that the exponent $-(c_1(X)\cup)$  has exactly  the opposite sign of the ``natural'' one appearing in \eqref{16agosto2018-3}, usually considered for any \emph{good} Frobenius manifolds and which, in the case of quantum cohomologies, comes from the classical limit point (see Proposition 2.2 and Corollary 2.1 of \cite{dubro2} and Remark \ref{natlev}).
\end{enumerate} 

By keeping track of all these discrepancies, in Section \ref{relggi} we give detailed proves of the following results.
\bsh
\begin{teorema}[cf. Proposition \ref{propgammapiu} and Theorem \ref{teoequivgamma2}]
\label{teoequivintro}
Let $X$ be a smooth projective variety of Hodge--Tate type for which point \emph{(3.b)} of Conjecture \ref{conjeintro} holds true. Namely, let the central connection matrix $C$ (computed with respect to some choice of the $\Psi$-matrix and of an oriented line $\ell$) be the matrix associated with the morphism $\textnormal{\textcyr{D}}^-_X$ and some exceptional collection $\frak E=(E_1,\dots, E_n)$. Then:
\begin{enumerate}
\item There exist matrices $A_{\pm}\in\mathcal C_0(X)$ such that the central connection matrix computed with respect to the solution $Y_{\rm top}A_{\pm}$ (and with respect to the same choices of $\Psi$ and $\ell$) has as columns the components of the characteristic classes
\[\frac{i^{\bar d}}{(2\pi)^\frac{d}{2}}\widehat\Gamma^\pm_X\cup {\rm Ch}(E_i).
\] In particular, if $X$ has even dimension $d$, the central connection matrix can be put in both of the forms \eqref{classdubr}-\eqref{classggi}, i.e. the ones predicted by \cite{dubro4} and by $\Gamma$-conjecture \emph{II}.

\item The validity of point \emph{(3.b)} of Conjecture \ref{conjeintro} is equivalent to $\Gamma$-conjecture \emph{II} for the system \[\widehat\nabla^{\rm GGI}_{\frac{\partial}{\partial\lambda}}Y=0,
\]considered in \cite{gamma1}. More precisely, through the identification \eqref{identsp}, the central connection matrix $C^{\rm GGI}$, computed with respect to the same choice of the $\Psi$-matrix and of the oriented line $\ell$, and with respect to the analogue of the topological solution $Y^{\rm GGI}_{\rm top}(\lambda)$ (not in \emph{natural} Levelt form) has as columns the components of the characteristic classes
\[\frac{1}{(2\pi)^\frac{d}{2}}\widehat\Gamma^+_X\cup {\rm Ch}(E_i'),\quad E_i':=\begin{sistema}
E_i^*,\quad \text{if }\lambda=e^{\pi i}z^{-1}\\
\\
\kappa^{-1}(E_i^*),\quad \text{if }\lambda=e^{-\pi i}z^{-1}
\end{sistema}
\] 
where $\frak E^*=(E_n^*,\dots, E_1^*)$ denotes the geometrical dual exceptional collection,
\[E_i^*:={\bf R}\HOM_{\mathcal O_X}^{\bullet}(E_i,\mathcal O_X),
\]and $\kappa:=(-\otimes \omega_X)[\dim_\mathbb CX]$ denotes the Serre functor.
\end{enumerate}
\end{teorema}
\esh

\subsection{Results for complex Projective Spaces} In Section \ref{computationsproj}, we focus on the case of complex Projective Spaces $\mathbb P^{k-1}_\mathbb C$. There we prove the validity of Conjecture \ref{conjeintro}, we explicitly compute the central connection matrix at points of the small quantum cohomology, and we carry out a detailed analysis of the braid group action on the monodromy data and on the corresponding exceptional collections. In particular we complete the study initiated by the third author in \cite{guzzetti1}, where point (3a) of Conjecture \ref{conjeintro} was proved. Let us summarize  the main results obtained.

\begin{shaded}
\begin{teorema}[cf. Theorem \ref{risultato1}, Corollary \ref{corexceptional1}]
\label{teocpn1intro}
Conjecture \ref{conjeintro} is true for all complex Projective Spaces $\mathbb P^{k-1}_{\mathbb C}$, $k\geq 2$. More precisely, the central connection matrix computed at $0\in QH^\bullet(\mathbb P^{k-1}_{\mathbb C})$ with respect to an oriented line $\ell$ of slope $\phi\in]0;\frac{\pi}{k}[$ coincides with the matrix associated with the morphism
\[\textnormal{\textcyr{D}}^-_{\mathbb P^{k-1}_\mathbb C}\colon K_0(\mathbb P^{k-1}_{\mathbb C})_\mathbb C\to H^\bullet(\mathbb P^{k-1}_{\mathbb C},\mathbb C)
\]computed with respect to the exceptional bases obtained by projecting on the $K_0$-group suitable shifts of the following exceptional collections:
\vskip 3mm
{\bf CASE $k$ EVEN:}\[\left(\mathcal O\left(\frac{k}{2}\right),\bigwedge\nolimits^1\mathcal T\left(\frac{k}{2}-1\right),\mathcal O\left(\frac{k}{2}+1\right),\bigwedge\nolimits^3\mathcal T\left(\frac{k}{2}-2\right),\dots,\mathcal O(k-1),\bigwedge\nolimits^{k-1}\mathcal T\right);
\]
\vskip 3mm
{\bf CASE $k$ ODD:}
\[\left(\mathcal O\left(\frac{k-1}{2}\right),\mathcal O\left(\frac{k+1}{2}\right),\bigwedge\nolimits^2\mathcal T\left(\frac{k-3}{2}\right),\mathcal O\left(\frac{k+3}{2}\right),\bigwedge\nolimits^4\mathcal T\left(\frac{k-5}{2}\right),\dots,\mathcal O\left(k-1\right),\bigwedge\nolimits^{k-1}\mathcal T\right).
\]
Here, we denote by $\mathcal O$ and $\mathcal T$ the structural and the tangent sheaf of $\mathbb P^{k-1}_\mathbb C$ respectively, and more in general by $\bigwedge\nolimits^p\mathcal T(q)$ the tensor product
\[\left(\bigwedge\nolimits ^p \mathcal T\right)\otimes\mathcal O(q).
\]
\end{teorema}
\end{shaded}

To the best of our knowledge, the above result is the first \emph{explicit} description of the exceptional collections that actually arise from the numerical values of the monodromy data as described by the original conjecture of \cite{dubro0}. We remark that the exceptional collections appearing in Theorem \ref{teocpn1intro} are in the same $\mathcal B_k$-orbit of the Beilinson exceptional collection $\frak B:=(\mathcal O,\dots, \mathcal O(k-1))$. Hence, it is worthwhile to understand for which Projective Spaces there exists suitable choices of signs for the $\Psi$-matrix, and oriented lines $\ell$ for which the monodromy data computed along the small quantum locus $H^2(\mathbb P^{k-1}_{\mathbb C},\mathbb C)$ are associated with the Beilinson exceptional collection $\frak B$. The following result gives us the answer.

\begin{shaded}
\begin{teorema}[cf. Theorem \ref{teotreccebeilinson}, Corollary \ref{corbeilip12}]\label{beilinsonthintro}
The Beilinson exceptional collection $\frak B$ is associated with the monodromy data computed in a chamber of the small quantum cohomology if and only if $k=2,3$.
\end{teorema}
\end{shaded}

Potentially, Theorem \ref{teocpn1intro}  can give us information about some region of the \emph{big} quantum cohomology of complex Projective Spaces: if we were able to solve  the Riemann--Hilbert problem associated with the monodromy data corresponding to $\frak B$, this could lead us to an explicit representation of the analytic continuation of the genus 0 Gromov--Witten potential of $\mathbb P^{k-1}_{\mathbb C}$.

In order to prove Theorem \ref{beilinsonthintro}, a careful analysis of the hidden symmetries of the Stokes phenomenon is carried on. By using symmetries of the regular polygons (which represent the spectrum of the operator $\mathcal U$ along the small quantum locus), and studying properties of all Stokes factors computed in \cite{guzzetti1}, we obtain the following result.

\begin{shaded}
\begin{teorema}[cf. Theorem \ref{quasi-period.th1}]
\label{thcpn2intr0}
The monodromy data computed at any  point of the small quantum cohomology of $\mathbb P^{k-1}_\mathbb C$ with $k\geq 2$, with respect to any  choice of an oriented line $\ell$, are obtained from those computed at $0\in QH^\bullet(\mathbb P^{k-1}_\mathbb C)$ with respect to a line of slope $\phi\in]0;\frac{\pi}{k}[$ by acting with a braid of the form
\[\omega_{1,k}\omega_{2,k}\omega_{1,k}\omega_{2,k}\dots,
\]where\footnotemark
\begin{itemize}
\item {\bf if $k$ is even} we set
\[\omega_{1,k}:=\prod_{\substack{i=2\\ i\ \text{even}}}^k\beta_{i-1,i},\quad \omega_{2,k}:=\prod_{\substack{i=3\\ i\ \text{odd}}}^{k-1}\beta_{i-1,i};
\]
\item {\bf if $k$ is odd} we set
\[\omega_{1,k}:=\prod_{\substack{i=3\\ i\ \text{odd}}}^k\beta_{i-1,i},\quad \omega_{2,k}:=\prod_{\substack{i=2\\ i\ \text{even}}}^{k-1}\beta_{i-1,i}.
\]
\end{itemize}
The corresponding exceptional collections are obtained (up to shifts) by acting with the above braids on the collections of Theorem \ref{teocpn1intro}.

Moreover, if we denote by $S_{\mathbb P^{k-1}_\mathbb C}(p,\phi)$ the Stokes matrix computed at a point $p\in H^2(\mathbb P^{k-1}_\mathbb C,\mathbb C)$, with respect to a line $\ell(\phi)$ of slope $\phi\in\mathbb R$, and in the $\ell$-lexicographical order, then the following facts hold.
\begin{enumerate}
\item If $\sigma$ denotes the generator of $H^2(\mathbb P^{k-1}_\mathbb C,\mathbb C)$, then the Stokes matrix $S_{\mathbb P^{k-1}_\mathbb C}(t\sigma,\phi)$, with $t\in\mathbb C$, is a function of ${\rm Im}(t)+k\phi$. We will write this as \[S_{\mathbb P^{k-1}_\mathbb C}(t\sigma,\phi)=S_{\mathbb P^{k-1}_\mathbb C}({\rm Im}(t)+k\phi),\quad t\in\mathbb C.
\] 
\item The Stokes matrix satisfies the \emph{quasi-periodicity} condition
\[S_{\mathbb P^{k-1}_\mathbb C}(p,\phi)\sim S_{\mathbb P^{k-1}_\mathbb C}\left(p,\phi+\frac{2\pi i}{k}\right),
\]where $A\sim B$ means that the matrices $A,B$ are in the same $(\mathbb Z/2\mathbb Z)^k$-orbit.
\item The entries 
\[
S_{\mathbb P^{k-1}_\mathbb C}(p,\phi)_{j,j+1}\text{ and }S_{\mathbb P^{k-1}_\mathbb C}\left(p,\phi+\frac{\pi i}{k}\right)_{j,j+1}
\]
differ by some signs for all $p\in H^2(\mathbb P^{k-1}_\mathbb C,\mathbb C)$, $\phi\in\mathbb R$ and for any $j=1,\dots k-1$. In particular, the $(k-1)$-tuple
\[
\left(\left|S_{\mathbb P^{k-1}_\mathbb C}(p,\phi)_{1,2}\right|, \left|S_{\mathbb P^{k-1}_\mathbb C}(p,\phi)_{2,3}\right|,\dots, \left|S_{\mathbb P^{k-1}_\mathbb C}(p,\phi)_{k-1,k}\right|\right)
\]does not depend on $p$ and $\phi$, and it is equal to
\[\left(\binom{k}{1},\dots,\binom{k}{k-1}\right).
\]
\end{enumerate}
\end{teorema}  
\end{shaded}
\footnotetext[8]{Here $\beta_{i-1,i}$ denotes one of the generator of the braid group $\mathcal B_n$. By identifying $\mathcal B_n$ with the mapping class group of a punctured disk (the punctures being the canonical coordinates $u_i$'s), the element $\beta_{i-1,i}$ corresponds to the elementary transformation given by a counterclockwise rotations of $u_{i-1}$ wrt $u_i$.}

Finally, we also obtained some results concerning the group ${\mathcal C}_0(\mathbb P^{k-1}_\mathbb C)$, refining point (3) of Theorem \ref{theoconjintro}. 

\begin{shaded}
\begin{teorema}[cf. Theorem \ref{teoctildecpn}, Corollary \ref{corctildecpn}]
Let for brevity $\mathbb P:=\mathbb P^{k-1}_\mathbb C$. The group ${\mathcal C}_0(\mathbb P)$ is an abelian unipotent algebraic group of dimension $[\frac{k}{2}]$. In particular, the exponential map defines an isomorphism
\[{\mathcal C}_0(\mathbb P)\cong \underbrace{\mathbb C\oplus\dots\oplus\mathbb C}_{[\frac{k}{2}]\text{ copies}}.
\]With respect to the basis $(1,\sigma,\dots, \sigma^{k-1})$ of $H^\bullet(\mathbb P,\mathbb C)$, the group ${\mathcal C}_0(\mathbb P)$ is described as follows
\[{\mathcal C}_0(\mathbb P)=\left\{C\in GL(k,\mathbb C)\colon C=\sum_{i=0}^{k-1}\alpha_iJ_i,\quad\alpha_0=1,\quad 2\alpha_{2n}+\sum_{\substack{i+j=2n\\1\leq i,j}}(-1)^i\alpha_i\alpha_{j}=0, \quad 2\leq 2n\leq k-1\right\},
\]
where the matrix $J_i$ is defined by
\[(J_i)_{ab}=\delta_{i,a-b}.
\]
In particular,  ${\mathcal C}_0(\mathbb P)$ is isomorphic to the identity component of the isometry group ${\rm Isom}_\mathbb C(K_0(\mathbb P)_\mathbb C,\chi)$.
\end{teorema}
\end{shaded}

\subsection{Results for complex Grassmannians} As an application of the \emph{abelian-nonabelian correspondence}, described for the specific case of complex Grassmannians in Sections \ref{secclasscohgrass} for classical cohomology and in Section \ref{secquantcohgrass} for the quantum cohomology, in Section \ref{chconjgrass} we explicitly compute the monodromy data of $QH^\bullet(\mathbb G(r,k))$ at points of the small quantum cohomology, deducing them from the corresponding monodromy data for the Projective Space $\mathbb P^{k-1}_\mathbb C$. Notice that for almost complex Grassmannians $\mathbb G(r,k)$, the points of their small quantum cohomologies are \emph{semisimple coalescence points} (i.e. points at which the Frobenius algebra is semisimple but with some coalescing canonical coordinates $u_i$'s; see \cite{cotti0}). At these points of the Frobenius manifold, the monodromy data are still well defined, and locally constant, thanks to the main  results of \cite{CDG0} and \cite{CDG}.  

In the following statement, we denote by $\bigwedge\nolimits^rA$ the $r$-th exterior power of a matrix $A\in M_{k}(\mathbb C)$ (also called $r$-th \emph{compound matrix} of $A$), namely the matrix of all $r\times r$ minors of $A$, ordered in the lexicographical order. Let us summarize the main results.

\begin{shaded}
\begin{teorema}[cf. Theorem \ref{teoexcollgrass}, Corollary \ref{08.02.18-1}, Theorem \ref{quasiperiodgrass}]
Let $\ell$ be an oriented line of slope $\phi\in]0;\frac{\pi}{k}[$, admissible\footnotemark  at both points 
\[p=t\sigma_1\in H^2(\mathbb G(r,k),\mathbb C)\quad\text{and}\quad \hat p:=(t+(r-1)\pi i)\sigma\in H^2(\mathbb P^{k-1}_\mathbb C,\mathbb C),
\]$\sigma$ and $\sigma_1$ being the Schubert classes generating the second cohomology groups of $\mathbb P^{k-1}_\mathbb C$ and $\mathbb G(r,k)$ respectively. For a suitable choice of signs of the $\Psi$-matrices, the monodromy data of $\mathbb G(r,k)$ are given by
\[S_{\mathbb G(r,k)}(p,\phi)=\bigwedge\nolimits^rS_{\mathbb P^{k-1}_\mathbb C}(\hat p,\phi),\quad C_{\mathbb G(r,k)}:=i^{-\binom{k}{r}}\left(\bigwedge\nolimits^r C_{\mathbb P^{k-1}_\mathbb C}(\hat p,\phi)\right)e^{\pi i(r-1)\sigma_1\cup(-)}.
\]
In particular, Conjecture \ref{conjeintro} holds true for the Grassmannian $\mathbb G(r,k)$. The exceptional collections associated with the above monodromy data are (modulo shifts) in the $\mathcal B_n$-orbit of the twisted Kapranov exceptional collection
\[\left(\mathbb S^\lambda\mathcal S^\vee\otimes\mathscr L\right)_\lambda,\quad \mathscr L:=\det\left(\bigwedge\nolimits^2\mathcal S^\vee\right),
\]where $\mathbb S^\lambda$ denotes the $\lambda$-th Schur functor and $\mathcal S$ the tautological bundle on $\mathbb G(r,k)$. Furthermore, the Stokes matrix satisfies the following conditions:
\begin{enumerate}
\item it has the following functional form
\[S_{\mathbb G(r,k)}(t\sigma_1,\phi)=S_{\mathbb G(r,k)}({\rm Im }\,t+k\phi);
\]
\item it is \emph{quasi-periodic} along the small quantum locus, in the sense that
\[S_{\mathbb G(r,k)}(p,\phi)\sim S_{\mathbb G(r,k)}\left(p,\phi+\frac{2\pi i}{k}\right),
\]where $A\sim B$ means that the matrices $A$ and $B$ are in the same orbit under the action of $(\mathbb Z/2\mathbb Z)^{\binom{k}{r}}$;

\item the upper-diagonal entries 
\[S_{\mathbb G(r,k)}(p,\phi)_{j,j+1},\quad S_{\mathbb G(r,k)}\left(p,\phi+\frac{\pi i}{k}\right)_{j,j+1}
\]differ by some signs, and 
\[|S_{\mathbb G(r,k)}(p,\phi)_{j,j+1}|\in\left\{\binom{k}{1},\dots,\binom{k}{k-1}\right\}\cup\left\{0\right\}.
\]
\end{enumerate}
\end{teorema}
\end{shaded}
\footnotetext[9]{An oriented line $\ell$ will be said to be \emph{admissible} at a point $p\in QH^\bullet(X)$ if it does not contain any of the Stokes rays
\[R_{ij}:=\left\{\sqrt{-1}\rho(\overline{u_i(p)}-\overline{u_j(p)})\colon\ \rho\in\mathbb R_+\right\},\quad i, j=1,\dots, n,\quad i\neq j.
\]}
\begin{shaded}
\begin{cor}[cf. Corollary \ref{corkapranovgrass}]
The Kapranov exceptional collection $(\mathbb S^\lambda\mathcal S^\vee)_\lambda$, twisted by a suitable line bundle, is associated with the monodromy data of $\mathbb G(r,k)$ at points of the small quantum locus if and only if $(r,k)=(1,2),(1,3),(2,3)$. In this cases the line bundle is trivial and the Kapranov collection coincides with the Beilinson one\footnotemark.
\end{cor}
\end{shaded}
\footnotetext[10]{Notice that $\mathbb G(2,3)\cong\mathbb P((\mathbb C^3)^\vee)\cong\mathbb P^2_\mathbb C$ by duality. }

\subsection{Plan of the paper} In Section \ref{GWQCOH} we recall basic notions in Gromov--Witten and quantum cohomology theories for smooth projective varieties. We describe the Frobenius structure naturally defined on $QH^\bullet(X)$, and we briefly summarize the analytic theory of semisimple Frobenius manifolds, their isomonodromic description and the main properties of their monodromy local invariants.

In Section \ref{Helixsec}, for completeness and convenience of the reader, we review the general theory of Helices in triangulated categories as developed by the Moscow School of Algebraic Geometry (see \cite{ru}, \cite{helix}). We recall basic notions and properties of \emph{exceptional objects}, \emph{exceptional collections} and more general \emph{semiorthogonal decomposition} in a $\mathbb K$-linear triangulated category $\mathscr D$, and we define their mutations under the action of the braid group. After recalling the properties of \emph{admissibility} of full triangulated subcategories and of \emph{saturatedness} of $\mathscr D$, the problem of existence of Serre functors is also discussed. We finally introduce the notions of \emph{dual} exceptional collections and of \emph{helix} generated by an exceptional collection.

In Section \ref{Mukaisec} we focus on \emph{unimodular  Mukai lattice} structures, introduced and studied by A.L. Gorodentsev (\cite{Go2.2, Go3}). Particular attention is given to the case of \emph{exceptional Mukai lattices}, i.e. those admitting an exceptional basis, an important example being furnished by the Grothendieck group $K_0(\mathscr D)$ of a $\mathbb K$-linear triangulated category $\mathscr D$ admitting a full exceptional collection. The mutations of exceptional bases under the action of the braid group, the canonical operator and the isometry group ${\rm Isom}(V,\langle\cdot,\cdot\rangle)$ are introduced and described. Furthemore, the complete isometric classification of Mukai spaces due to Gorodentsev is outlined (Theorem \ref{13.07.17-5} and Theorem \ref{classcompleta}). We also consider the geometrical case of the Grothendieck group of a smooth projective variety $X$ admitting a full exceptional collection in $\mathcal D^b(X)$. We briefly recall that the existence of such a collection implies a \emph{motivic decomposition} of $X$, and hence strong constraints are deduced on its geometry and topology. Finally, results on the isometric classification of the Grothendieck groups $K_0(X)\otimes_\mathbb Z\mathbb C$ with non-degenerate Euler--Poincaré form are presented. 

In Section \ref{chmainconj} we review the original version of the conjecture of \cite{dubro0}. This conjecture states the equivalence of the condition of semisimplicity of the quantum cohomology of a Fano variety $X$ with the condition of existence of a full exceptional collection  in $\mathcal D^b(X)$, and it also prescribes the monodromy data $(S,C)$ in  geometric terms with respect to the objects of the exceptional collection. After reviewing the results available in the literature partially confirming the conjecture, we  formulate a refined and complete version of the conjecture (Conjecture \ref{conjeintro}/\ref{congettura}), including a prescription also for the central connection matrix $C$. We also explain how heuristically the conjecture  should follow from M. Kontsevich's proposal of Homological Mirror Symmetry. Finally, we describe the precise relationships between Conjecture \eqref{congettura} with $\Gamma$-conjecture II of S. Galkin, V. Golyshev and H. Iritani.

In Section \ref{computationsproj} we prove validity of Conjecture \ref{congettura} for all complex Projective Spaces $\mathbb P^{k-1}_\mathbb C$. After computing the topological-enumerative solution for the system of deformed flat coordinates, we show that the group ${{\mathcal  C}}_0(\mathbb P^{k-1}_\mathbb C)$, which describes the ambiguity in the choice of a solution in the Levelt  form at $z=0$, is isomorphic to the identity component of the isometry group ${\rm Isom}_\mathbb C(K_0(\mathbb P^{k-1}_\mathbb C), \chi(\cdot,\cdot))_0$. Hence, we compute the central connection matrix at the point $0\in QH^\bullet(\mathbb P^{k-1}_\mathbb C)$ w.r.t a line $\ell$ of slope $0<\phi<\frac{\pi}{k}$. By completing the braid analysis developed in \cite{guzzetti1}, we recognize in the computed monodromy data the geometric information, as prescribed by the Conjecture \ref{congettura}, associated with an explicit mutations of the Beilinson exceptional collection (see Theorem \ref{risultato1} and Corollary \ref{corexceptional1}). After studying in detail  the $\ell$-chamber decomposition along the small quantum locus, a property of \emph{quasi-periodicity} of Stokes matrix is shown (Theorem \ref{quasi-period.th1}). From this property, we deduce that the only Projective Space for which the monodromy data are the ones associated with the Beilinson collection are $\mathbb P^1_\mathbb C$ and $\mathbb P^2_\mathbb C$ (Corollary \ref{corbeilip12}). For all other Projective Spaces the data corresponding to the Beilinson exceptional collection can be computed in chambers of the big quantum cohomology, by means of the action of the braid group. 

In Section \ref{chcoalgrass}, using the \emph{(Quantum) Abelian/Non-abelian correspondence} of \cite{BCFK, BCFK2, CFKS}, we prove validity of Conjecture \ref{conjeintro}/\ref{congettura} for Grassmannians  by using the results of the previous Section. In particular, we show that the monodromy data computed at the points of the small quantum cohomology, with respect to an oriented line $\ell$ in the complex plane, are the prescribed geometric data associated with an exceptional collection which can be mutated into the Kapranov exceptional collection twisted by a line bundle (Theorem \ref{teoexcollgrass}).\newline

\addtocontents{toc}{\protect\setcounter{tocdepth}{1}}
\subsection{Acknowledgements} We would like to thank P. Belmans, M. Bertola, U. Bruzzo, S. Cecotti, B. Fantechi, S. Galkin, V. Golyshev, V. Gorbounov, C. Hertling, H. Iritani, A. Its, M. Kontsevich, C. Korff, C. Sabbah, F. Sala, M. Smirnov, J. Stoppa, I. Strachan, A. Varchenko, D. Yang for several comments and very useful discussions. The first author is grateful to the Max-Planck-Institut für Mathematik in Bonn, where a significant part of this work was developed, for support and providing excellent working conditions. 
\addtocontents{toc}{\protect\setcounter{tocdepth}{2}}

\landscape
$\quad$\newline

$\quad$\newline

$\quad$\newline
\begin{table}[h]
\centering
\caption{Identifications enriching the Conjecture of \cite{dubro0}}
\label{enrichconjintro}
\begin{tabular}{|c|l|c|l|c|l|}
\hline
\multicolumn{2}{|c|}{Frobenius Manifolds $QH^\bullet(X)$}                                                                                                                                                                                                               & \multicolumn{2}{c|}{Derived category $\mathcal D^b(X)$}                                                                                                                 & \multicolumn{2}{c|}{Grothendieck group $K_0(X)_\mathbb C$}                                                                                                                             \\ \hline\hline
\multicolumn{2}{|c|}{Stokes matrix $S$}                                                                                                                                                                                                                                 & \multicolumn{2}{l|}{}                                                                                                                                                   & \multicolumn{2}{c|}{\begin{tabular}[c]{@{}c@{}}inverse of the Gram matrix \\ $G_{ij}:=\chi(E_i,E_j)$\\ for an exceptional basis $([E_i])_i$\end{tabular}}                              \\ \hline
\multicolumn{2}{|c|}{Central connection matrix $C$}                                                                                                                                                                                                                     & \multicolumn{2}{l|}{}                                                                                                                                                   & \multicolumn{2}{c|}{\begin{tabular}[c]{@{}c@{}}matrix associated with the morphism \\ \textcyr{D}$^-_X\colon K_0(X)_\mathbb C\to H^\bullet(X,\mathbb C)$\end{tabular}}               \\ \hline
\multicolumn{2}{|c|}{\begin{tabular}[c]{@{}c@{}}action of the braid group \\ $\mathcal B_n$ on the monodromy data  \end{tabular}}                                                                                                                                                               & \multicolumn{2}{c|}{\begin{tabular}[c]{@{}c@{}}action of $\mathcal B_n$ on \\ the set of exceptional collections\end{tabular}}                                          & \multicolumn{2}{c|}{\begin{tabular}[c]{@{}c@{}}action of $\mathcal B_n$ on \\ exceptional bases\end{tabular}}                                                                          \\ \hline
\multicolumn{2}{|c|}{\begin{tabular}[c]{@{}c@{}}action of the group \\ $(\mathbb Z/2\mathbb Z)^{\times n}$ on the monodromy data\end{tabular}}                                                                                                                                                & \multicolumn{2}{c|}{shifts of exceptional collections}                                                                                                                  & \multicolumn{2}{c|}{\begin{tabular}[c]{@{}c@{}}projected shifts, i.e. change of signs, \\ of exceptional bases\end{tabular}}                                                           \\ \hline
\multicolumn{2}{|c|}{\begin{tabular}[c]{@{}c@{}}action of the group \\ ${\mathcal C}_0(X)$ on the monodromy data\end{tabular}}                                                                                                                                                      & \multicolumn{2}{c|}{\begin{tabular}[c]{@{}c@{}}action of a subgroup \\ of autoequivalences\\ ${\rm Aut}(\mathcal D^b(X))$\\ on the set of exceptional collections\end{tabular}}                                  & \multicolumn{2}{c|}{\begin{tabular}[c]{@{}c@{}}action of a subgroup of\\ the identity component of  \\ ${\rm Isom}_\mathbb C(K_0(X)_\mathbb C,\chi)$\\ on the set of exceptional bases\end{tabular}}                     \\ \hline
\multicolumn{2}{|c|}{\begin{tabular}[c]{@{}c@{}}complete ccw rotation of the line $\ell$, \\ action of the generator of the center $Z(\mathcal B_n)$,\\ action of the matrix $M_0\in{\mathcal C}_0(X)$,\\ $M_0:=\exp(2\pi i \mu)\exp(2\pi i R)$,\\ on the monodromy data\end{tabular}} & \multicolumn{2}{c|}{\begin{tabular}[c]{@{}c@{}}action of the Serre functor \\ $(\omega_X\otimes-)[\dim_\mathbb CX]$\\ on the set of exceptional collections\end{tabular}} & \multicolumn{2}{c|}{\begin{tabular}[c]{@{}c@{}}action of the canonical operator \\ $\kappa\colon K_0(X)_\mathbb C\to K_0(X)_\mathbb C$\\ on the set of exceptional bases\end{tabular}} \\ \hline
\end{tabular}
\end{table}
\endlandscape

\section{Gromov--Witten Theory and Quantum Cohomology}\label{GWQCOH}
\subsection{Notations and preliminaries}\label{notations}
Let $X$ be a smooth projective complex variety. In order not to deal with Frobenius superstructures, we will suppose\footnote{To our purposes, this is not a restrictive condition, as the reader can see in Theorem \ref{hmt}.} that the variety $X$ has vanishing odd cohomology, i.e. $H^{2k+1}(X,\mathbb C)\cong 0$ for $0\leq k$. Let us fix a homogeneous basis $(T_1,\dots, T_N)$ of $H^\bullet(X,\mathbb C)=\bigoplus_{k}H^{2k}(X,\mathbb C)$ such that
\begin{itemize}
\item $T_1=1$ is the unity of the cohomology ring;
\item $T_2,\dots, T_r$ span $H^2(X,\mathbb C)$.
\end{itemize}
We will denote by $\eta\colon H^\bullet(X,\mathbb C)\times H^\bullet(X,\mathbb C)\to\mathbb C$ the Poincaré metric
\[\eta(\xi,\zeta):=\int_X\xi\cup\zeta,
\]and in particular 
\[\eta_{\alpha\beta}:=\int_XT_\alpha\cup T_\beta.
\]

If $\beta\in H_2(X;\mathbb Z)/\text{torsion}$, we denote by $\overline{\mathcal M}_{g,n}(X,\beta)$ the Kontsevich-Manin moduli stack of $n$-pointed, genus $g$ stable maps with target $X$ of degree $\beta$, which parametrizes equivalence classes of pairs $((C_g,\bold x); f)$, where:
\begin{itemize}
\item  $(C_g,\bold x)$ is an $n$-pointed algebraic curve of genus $g$, with at most nodal singularities and with $n$ marked points $\bold x=(x_1,\dots, x_n)$, and $f\colon C_g\to X$ is a morphism such that $f_*[C_g]\equiv \beta$. Two pairs $((C_g,\bold x);f)$ and $((C'_g,\bold {x'});f')$ are defined to be equivalent if there exists a bianalytic map $\varphi\colon C_g\to C'_g$ such that $\varphi(x_i)=x'_i$, for all $i=1,\dots, n$, and $f'=\varphi\circ f$.

\item The morphisms $f$ are required to be \emph{stable}: if $f$ is constant on some irreducible component of $C_g$, then that component as a pointed curve should admit only a finite number of automorphisms (in other words, it must have at least 3 distinguished points, i.e. points that are either nodes or marked ones). 
\end{itemize}

We will denote by $\operatorname{ev}_i\colon \overline{\mathcal M}_{g,n}(X,\beta)\to X\colon ((C_g,\bold x); f)\mapsto f(x_i)$ the naturally defined evaluation map, and by $\psi_i\in H^2(\overline{\mathcal M}_{g,n}(X,\beta);\mathbb Q)$ the Chern classes of tautological cotangent line bundles $$\mathcal L_i\to \overline{\mathcal M}_{g,n}(X,\beta),\quad \mathcal L_i|_{((C_g,\bold x); f)}=T^*_{x_i}C_g,\quad \psi_i:=c_1(\mathcal L_i).$$ 
Using the construction of \cite{fantechi} of a \emph{virtual fundamental class} $[\overline{\mathcal M}_{g,n}(X,\beta)]^{\text{virt}}$ in the Chow ring $CH_\bullet(\overline{\mathcal M}_{g,n}(X,\beta))$, and of degree equal to the expected dimension
\[[\overline{\mathcal M}_{g,n}(X,\beta)]^{\text{virt}}\in CH_D(\overline{\mathcal M}_{g,n}(X,\beta)), \quad D=(1-g)(\dim_{\mathbb C}X-3)+n+\int_\beta c_1(X),
\]a good theory of intersection is allowed on the Kontsevich-Manin moduli stack. 

We can thus define the \emph{Gromov--Witten invariants (with descendants) of genus $g$, with $n$ marked points and of degree $\beta$} of $X$ as the integrals (whose values are rational numbers)
\begin{equation}\label{gw1}\langle\tau_{d_1}\gamma_1,\dots,\tau_{d_n}\gamma_n\rangle_{g,n,\beta}^X:=\int_{[\overline{\mathcal M}_{g,n}(X,\beta)]^{\text{virt}}}\prod_{i=1}^n\operatorname{ev}_i^*(\gamma_i)\cup\psi_i^{d_i},
\end{equation}
\[ \gamma_i\in H^\bullet(X,\mathbb C),\quad d_i\in\mathbb N,\quad i=1,\dots, n.
\]Since by \emph{effectiveness} (for an axiomatic treatment of the Gromov--Witten invariants we follow \cite{manin}, \cite{kon} and \cite{cox}) the integral is non-vanishing only for effective classes $\beta\in\operatorname{Eff}(X)\subseteq H_2(X;\mathbb Z)$, the generating function of rational numbers \eqref{gw1}, called \emph{total descendent potential} (or also \emph{gravitational Gromov--Witten potential}, or even \emph{Free Energy}) of \emph{genus $g$} is defined as the formal series
\begin{equation}\label{21.06.17}\mathcal F^X_g(\gamma,\bold Q):=\sum_{n=0}^\infty\sum_{\beta\in \text{Eff}(X)}\frac{\bold Q^\beta}{n!}\langle\underbrace{\gamma.\dots,\gamma}_{n\text{ times}}\rangle_{g,n,\beta}^X,
\end{equation}
where we have introduced (infinitely many) coordinates $\bold{t}:=(t^{\alpha,p})_{\alpha,p}$
\[\gamma=\sum_{\alpha,p}t^{\alpha,p}\tau_pT_{\alpha},\quad \alpha=1,\dots, N,\ p\in\mathbb N,
\] and formal parameters
\[\bold Q^\beta:=Q_2^{\int_\beta T_2}\cdot\dots\cdot Q_r^{\int_\beta T_r},\quad Q_i \text{'s are elements of the Novikov ring }\Lambda:=\mathbb C[\![Q_2,\dots,Q_r]\!].
\]
The free energy $\mathcal F^X_g\in\Lambda[\![{\bold t}]\!]$ can be seen as a function on the \emph{large phase-space}, and restricting the free energy to the \emph{small phase space} (naturally identified with $H^\bullet(X,\mathbb C)$), 
\[F^X_g(t^{1,0},\dots, t^{N,0}):=\mathcal F^X_g(\bold t)|_{t^{\alpha,p}=0,\ p>0},
\]
one obtains the generating function of the Gromov--Witten invariants of genus $g$.

\subsection{Quantum cohomology and its semisimplicity}
By the Divisor Axiom, the genus $0$ Gromov--Witten potential $F^X_0(t)$, can be seen as an element of the ring $\mathbb C[\![t^1,Q_2e^{t^2},\dots,Q_re^{t^r},t^{r+1},\dots, t^N]\!]$: in what follows we will be interested in cases in which $F^X_0$ is the expansion of an analytic function, i.e.
\[F^X_0\in\mathbb C\left\{t^1,Q_2e^{t^2},\dots,Q_re^{t^r},t^{r+1},\dots, t^N\right\}.
\]Without loss of generality, we can put $Q_2=\dots=Q_r=1$, and $F^X_0(t)$ defines an analytic function in an open neighborhood $\mathcal D\subseteq H^\bullet(X,\mathbb C)$ of the point\footnote{This means that there exist two positive real numbers $\varepsilon, C$ such that $F^X_0(t)$ is convergent and analytic on the open set 
\[|t^i|<\varepsilon,\quad i=1,r+1,\dots, N,
\]
\[{\rm{Re}}\ t^i<-C,\quad i=2,\dots, r.
\]}
\begin{align}\label{classical1}
t^i&=0,\quad i=1,r+1,\dots, N,\\
\label{classical2}
\operatorname{Re} t^i &\to -\infty,\quad i=2,\dots, r.
\end{align}
The function $F^X_0$ is a solution of WDVV equations (for a proof see \cite{kon}, \cite{manin}, \cite{cox}), and thus it defines an analytic Frobenius manifold structure on $\mathcal D$ (\cite{dubronapoli, dubro1, dubro0, dubro2, CDG}), characterised by the following objects:
\begin{itemize}
\item the flat metric is given by the Poincaré metric $\eta$;
\item the unity vector field is $T_0=1$, using the canonical identifications of tangent spaces with cohomology
$$
T_p\mathcal D\cong H^\bullet(X,\mathbb C)\colon \partial_{t^\alpha}\mapsto T_\alpha;
$$
\item the Euler vector field is
\begin{equation}\label{eulerfield}E:=c_1(X)+\sum_{\alpha=1}^N\left(1-\frac{1}{2}\deg T_\alpha\right)t^\alpha T_\alpha.
\end{equation}
\end{itemize}
More precisely, by the Point Mapping Axiom, the Gromov--Witten potential can be decomposed into a \emph{classical} term and a \emph{quantum} correction as follows
\begin{align}\label{gwpot1}F^X_0(\gamma)&=F_{\text{classical}}+F_{\text{quantum}}\\
\label{gwpot2}
&=\frac{1}{6}\int_X\gamma^3+\sum_{k=0}^\infty\sum_{\beta\in\operatorname{Eff}(X)\setminus\left\{0\right\}}\frac{1}{k!}\langle\underbrace{\gamma,\dots,\gamma}_{k\text{ times}}\rangle^X_{0,k,\beta},\quad\text{where }\gamma=\sum_{\alpha=1}^Nt^\alpha T_\alpha.
\end{align}
Consequently, the product on the algebra $(T_p\mathcal D, \circ_p)$ at a point $p$ defined by
\begin{equation}\label{23.06.17-1}T_\alpha\circ_p T_\beta:=\sum_{\gamma,\delta}\left.\frac{\partial^3F^X_0}{\partial t^\alpha\partial t^\beta\partial t^\gamma}\right|_p\eta^{\gamma\delta}T_\delta,\quad p\in\mathcal D,
\end{equation}defines a deformation of the classical cohomological $\cup$-product. The associativity of $\circ_p$ is equivalent to validity of the WDVV equations
\[\frac{\partial^3F^X_0}{\partial t^\alpha\partial t^\beta\partial t^\gamma}\eta^{\gamma\delta}\frac{\partial^3F^X_0}{\partial t^\delta\partial t^\epsilon\partial t^\mu}=\frac{\partial^3F^X_0}{\partial t^\mu\partial t^\beta\partial t^\gamma}\eta^{\gamma\delta}\frac{\partial^3F^X_0}{\partial t^\delta\partial t^\epsilon\partial t^\alpha},
\]and it is easily seen that
\[\frac{\partial^3F^X_0}{\partial t^1\partial t^\alpha\partial t^\beta}=\eta_{\alpha\beta},
\]the variable $t^1$ appearing only in the classical term of $F^X_0$. Hence, the resulting algebras $(T_p\mathcal D,\circ_p)$ are \emph{Frobenius}, in the sense of the following
\begin{defi}
Let $\mathbb K$ be a field\footnote{We are interested in the case $\mathbb K=\mathbb C$.}, and let $(A,*,1)$ be a finite dimensional associative, commutative and unital $\mathbb K$-algebra, endowed with a non-degenerate symmetric bilinear form $\eta\colon A\times A\to\mathbb K$. We say that $(A,*,1,\eta)$ defines a \emph{Frobenius algebra} if the following  conditions are satisfied
\[\eta(v*u,w)=\eta(v,u*w),\quad u,v,w\in A.
\]
\end{defi} Furthermore, the Gromov--Witten potential $F^X_0$ satisfies also the quasi-homogeneity condition 
\[\frak L_E F^X_0=(3-\dim_\mathbb C X)\cdot F^X_0~\text{plus at most quadratic  terms}.
\]

\begin{defi}
The Frobenius manifold structure defined on the domain of convergence $\mathcal D$ of the Gromov--Witten potential $ F^X_0$, solution of the WDVV problem, is called \emph{Quantum Cohomology} of $X$, and denoted by $QH^\bullet(X)$. 
By the expression \emph{small quantum cohomology of $X$} (or \emph{small quantum locus}) we denote the Frobenius structure attached to points in $\mathcal D\cap H^2(X,\mathbb C)$. In case of convergence, the potential $F^X_0$ (and hence the whole Frobenius structure) can be maximally analytically continued to an unramified covering of an open subdomain of $H^\bullet(X,\mathbb C)$. We refer to this global Frobenius structure as the \emph{big quantum cohomology of $X$}, and it will be still denoted by $QH^\bullet(X)$.
\end{defi}

Although no general results guarantee the convergence of the Gromov--Witten potential $F^X_0$ for a generic smooth projective variety $X$, for some classes of varieties it is known that the sum defining $F^X_0$ at points of the small quantum cohomology (at which $t^1=t^{r+1}=\dots=t^N=0$) is finite. This is the case for
\begin{itemize}
\item Fano varieties,
\item varieties admitting a transitive action of a semisimple Lie group. 
\end{itemize}
For the proof see \cite{cox}. Notice that for these varieties the small quantum locus coincides with the whole space $H^2(X,\mathbb C)$. Conjecturally, for Calabi--Yau manifolds the series defining $F^X_0$ is convergent in a neighborhood of the classical limit point (see \cite{cox}, \cite{kon}).

\begin{oss}
In literature, the name \emph{small quantum cohomology} of $X$ is usually reserved to the family of Frobenius algebras parametrized either by points of $H^2(X,\mathbb C)/(2\pi i H^2(X,\mathbb Z))$, or by points of  $H^2(X,\mathbb C)$, through the uniformization $q_i=\exp(t^i)$. In this paper we adopt the second convention, since we want to develop a detailed analysis on the dependence of a set of local invariants of the Frobenius manifold structure (the so called \emph{monodromy data}) on the parameters $t^i$'s. See for example Section \ref{symcpn} and Section \ref{symgrass}.
\end{oss}

In this paper we will focus on smooth projective varieties $X$ whose (big) quantum cohomology is a \emph{semisimple} Frobenius manifold.

\begin{defi} A smooth projective variety $X$ admits semisimple quantum cohomology if there exists an open dense subset of points $p\in QH^\bullet(X)$ at which the associated algebra $(T_pQH^\bullet(X),\circ_p)$ satisfies the following equivalent conditions:
\begin{enumerate}
\item it is semisimple,
\item it has vanishing Jacobson ideal,
\item it is without nilpotents,
\item it is isomorphic as $\mathbb C$-algebra to $\mathbb C^{N}$ (with component-wise multiplication).
\end{enumerate}
A point $p\in QH^\bullet(X)$ whose associated Frobenius algebra is semisimple will be called a \emph{semisimple point}, for short.
\end{defi}

Notice that the classical Frobenius cohomological algebra $(H^\bullet(X,\mathbb C),\cup)$, corresponding to the limit point \eqref{classical1}-\eqref{classical2}, is not semisimple, if $\dim X\neq 0$, since it clearly contains nilpotent elements. By quantum deformation of the $\cup$-product, it may happen that the semisimplicity condition is satisfied. The problem of characterizing smooth projective varieties with semisimple quantum cohomology is far from being solved.  The following result shows that the assumption on $X$ considered above, of having odd-vanishing cohomology $H^{\rm{odd}}(X,\mathbb C)\cong 0$, is a necessary condition in order to have semisimplicity of the quantum cohomology $QH^\bullet(X)$.

\begin{teorema}[\cite{hmt}]{\label{hmt}}
If $X$ is a smooth projective variety whose quantum cohomology $QH^\bullet(X)$ is a semisimple analytic Frobenius manifold, then $X$ is of Hodge--Tate type\footnote{Here, $h^{p,q}(X)$ denotes the $(p,q)$-Hodge number of $X$, namely $$h^{p,q}(X):=\dim_{\mathbb C}H^q(X,\Omega^p_X),$$ where $\Omega^p_X$ is the sheaf of holomorphic $p$-forms on $X$.}, i.e.
\[h^{p,q}(X)=0,\quad\text{if }p\neq q.
\]
In particular, $X$ is with odd-vanishing cohomology.
\end{teorema}

For some classes of varieties, such as some Fano threefolds \cite{ciolli}, toric varieties \cite{iri7}, and some homogeneous spaces \cite{cmp10}, it has been proved that points of the small quantum cohomology are all semisimple. Ordinary complex Grassmannians are among these varieties. More general homogeneous spaces may have non-semisimple small quantum cohomology (\cite{cmp10}, \cite{chaput-perrin}, \cite{gms}). Some sufficient conditions for other Fano varieties are given in \cite{perrin}. 

\begin{oss}
Remarkably, under the assumptions of convergence of the genus 0 Gromov--Witten potential $F^X_0$ and semisimplicity of the quantum cohomology $QH^\bullet(X)$, it can be shown (\cite{coat-iri}) that there exist two real positive constants $C,\varepsilon$ such that, for any $g\geq 0$, the power series \eqref{21.06.17} defining the genus $g$ total descendant potential  $\mathcal F^X_g$ is convergent on the infinite-dimensional polydisc
\begin{align*}
 |t^{\alpha, p}|&<\varepsilon\frac{p!}{C^p}\quad \text{for $\alpha=1,\dots, N$, and $p\in\mathbb N$,}\\
 |Q_i|&<\varepsilon\quad \text{for }i=2,\dots, r.
\end{align*}
\end{oss}

\subsection{Extended deformed connection and topological solution}\label{extdefconn} The study and the classification of semisimple Frobenius manifolds have been introduced and extensively developed in \cite{dubro1, dubro0, dubro2} and further refined in \cite{CDG}. One of the main objects catching the main properties of the general theory of Frobenius manifolds is the so called \emph{extended deformed connection} $\widehat\nabla$. Starting from the flat Levi-Civita connection $\nabla$ of the metric $\eta$ defined on the holomorphic tangent bundle $TQH^\bullet(X)$, whose flat coordinates are $(t^\alpha)_{\alpha=1}^N$, let us define the family of connections $\nabla^{(z)}$, with $z\in\mathbb C$, as follows
\begin{equation}\label{14.11.18-1}\nabla^{(z)}_ZY:=\nabla_ZY+z Z\circ Y,
\end{equation}for any vector fields $Y,Z$. The associativity of the quantum product together with the fact that the quantum product is induced by the potential $F^X_0$ through equation \eqref{23.06.17-1} are equivalent to the \emph{flatness} of the connections $\nabla^{(z)}$ for any $z\in \mathbb C$. The family $\left\{\nabla^{(z)}\right\}_{z\in\mathbb C^*}$ ca be further rigidified into a unique connection $\widehat\nabla$ defined on the pull-backed bundle $\pi^*TQH^\bullet(X)$, where $\pi\colon \mathbb C^*\times QH^\bullet(X)\to QH^\bullet(X)$ is the canonical projection. Let us introduce the $(1,1)$-tensors $\mathcal U,\mu\in\Gamma(\operatorname{End}(TQH^\bullet(X)))$ defined by
\begin{equation}\label{06.07.17-3}\mathcal U(Y):=E\circ Y,\quad \mu(Y):=\frac{2-\dim_\mathbb CX}{2}Y-\nabla_YE,\quad Y\in\Gamma(TQH^\bullet(X)),
\end{equation}
 and let us pull-back on $\pi^*TQH^\bullet(X)$ all the tensors $\eta,\mathcal U,\mu$, and denote them by the same symbols. The Levi-Civita connection $\nabla$ can be lifted on $\pi^*TQH^\bullet(X)$ in such a way that
\[\nabla_{\partial _z}Y=0
\]for any section $Y$ of the inverse image sheaf $\pi^{-1}\mathscr T$, where $\mathscr T:=\Gamma(-, TQH^\bullet(X))$ is the tangent sheaf of $QH^\bullet(X)$. Thus, we can define the extended deformed connection $\widehat\nabla$ as follows
\begin{align*}
\widehat\nabla_ZY:&=\nabla_ZY+z Z\circ Y,\\
\widehat\nabla_{\partial _z}Y:&=\nabla_{\partial_z}Y+\mathcal U(Y)-\frac{1}{z}\mu(Y),
\end{align*}
for any section $Y$ of the pull-back sheaf $\pi^*\mathscr T$. Remarkably, the whole connection $\widehat\nabla$ is \emph{flat} ( see \cite{dubro1, dubro2, CDG}). This implies existence of functions $\tilde t=\tilde t(t,z)$, called \emph{deformed flat coordinates}, such that
\[\widehat\nabla d\tilde t=0,\quad d:=\sum_\alpha \frac{\partial}{\partial t^\alpha}dt^\alpha.
\]Given independent functions $(\tilde t^1,\dots, \tilde t^N)$ as above, we have a system of $\widehat \nabla$-flat coordinates $(z, \tilde t^1,\dots, \tilde t^N)$ on the product $\mathbb C^*\times QH^\bullet(X) $. The $\widehat \nabla$-flatness condition can be rewritten in terms of the $\eta$-gradients as follows 
\begin{align}\label{06.07.17-1}
\partial_\alpha\zeta&=z\mathcal C_\alpha\zeta,\\
\label{06.07.17-2}
\partial_z\zeta&=\left(\mathcal U+\frac{1}{z}\mu\right)\zeta,
\end{align}
where we set
\[\left(\mathcal C_\alpha\right)_\beta^\gamma=c_{\alpha\beta}^\gamma,\quad \zeta=(\zeta^1,\dots,\zeta^{N})^T,\quad \zeta^\alpha=\eta^{\alpha\beta}\frac{\partial \tilde t}{\partial t^\beta}.\]
Notice that, because of the commutativity of the Frobenius algebras, and the compatibility of the product with the metric $\eta$, the following properties hold for the operators $\mathcal U,\mu$ appearing in equation \eqref{06.07.17-2}:
\[\mathcal U^T\eta=\eta\mathcal U,\quad \mu^T\eta+\eta\mu=0.
\]Furthermore, in flat coordinates   $(t^\alpha)_{\alpha=1}^N$ of $\nabla$,  the (1,1)-tensor $\mu$ is in diagonal form
\[\mu={\rm diag}(\mu_1,\dots,\mu_N),\quad \mu_\alpha=\frac{1}{2}\left(\deg T_\alpha-\dim_\mathbb CX\right).
\]
In order to describe fundamental solutions of the differential system \eqref{06.07.17-1}-\eqref{06.07.17-2} let us introduce the following definitions. Let $(V,\eta,\mu)$ be the datum of 
\begin{enumerate}
\item a $N$-dimensional complex vector space $V$,
\item a symmetric non-degenerate bilinear form $\eta\colon V\times V\to\mathbb C$,
\item a diagonalizable endomorphism $\mu\in\End(V)$ which is $\eta$-antisymmetric
\[\eta(\mu a,b)+\eta(a,\mu b)=0,\quad\text{ for all }a,b\in V.
\]
\end{enumerate}

\begin{defi}[\cite{dubro2,dubro3}] Let $(V,\eta,\mu)$ as above. An endomorphism $A\in\End(V)$ is $\mu$-\emph{nilpotent} if
\[AV_{\mu_\alpha}\subseteq\bigoplus_{m\geq 1}V_{\mu_\alpha+m}\quad\text{for any }\mu_\alpha\in\operatorname{spec}(\mu),
\]and where we introduced the family $(V_\lambda)_{\lambda\in\mathbb C}$ of subspaces of $V$ defined by
\[V_\lambda:=\left\{v\in {V}\colon\mu(v)=\lambda v\right\}.
\] 
In particular such an operator is nilpotent in the usual sense. A $\mu$-nilpotent operator $A$ can be uniquely decomposed in components $A_k\in\End( V)$, $k\geq 1$, such that
\[A_k V_{\mu_\alpha}\subseteq V_{\mu_\alpha+k}\quad\text{for any }\mu_\alpha\in\operatorname{spec}(\mu),\quad A=\sum_{k\geq 1}A_k.\]
\end{defi}

\begin{defi}[\cite{dubro2, CDG}] Let $(V,\eta,\mu)$ as above. Let us define on $V$ a new non-degenerate bilinear form $\left\{\cdot,\cdot\right\}$ by the equation
\[
\left\{a,b\right\}:=\eta\left(e^{i\pi\mu}a,b\right),\quad\text{for all }a,b\in V.
\]
We define the $(\eta,\mu)$-\emph{parabolic orthogonal group}, denoted by $\mathcal G(\eta,\mu)$, as the complex Lie group of all $\left\{\cdot,\cdot\right\}$-isometries $G\in\operatorname{GL}(V)$ of the form
\[
G=\mathbbm 1_V+\Delta
\]with $\Delta$ a $\mu$-nilpotent operator. Its Lie algebra $\frak g(\eta,\mu)$ coincides with the set of all $\mu$-nilpotent operators $R$ which are also $\left\{\cdot,\cdot\right\}$-\emph{skew-symmetric} in  the sense that
\[\left\{Rx,y\right\}+\left\{x,Ry\right\}=0.
\]In particular, any such operator $R$ commutes with the operator $e^{2\pi i\mu}$. 
\end{defi}

\begin{defi}If $(V_1,\eta_1,\mu_1)$ and $(V_2,\eta_2,\mu_2)$ are two triples satisfying the properties (1),(2),(3) above, we call a morphism of triples $f\colon (V_1,\eta_1,\mu_1)\to (V_2,\eta_2,\mu_2)$ the datum of a linear morphism $f\colon V_1\to V_2$ compatible with both the metrics and the operators $\mu$'s in the following sense:
\[\eta_1(x,y)=\eta_2(f(x),f(y)),\quad x,y\in V_1,
\]
\[f\circ\mu_1=\mu_2\circ f.
\]The notion of isomorphism of triples naturally follows.
\end{defi}
Given a smooth projective complex variety $X$, we can canonically associate an isomorphism class of  triples $[(V,\eta,\mu)]$ to its quantum cohomology $QH^\bullet(X)$. Such a class is called \emph{spectrum} of $QH^\bullet(X)$ in the terminology of Frobenius manifolds theory (see \cite{dubro2,dubro3,CDG}). At each point $p\in QH^\bullet(X)$, we have a triple $(T_pQH^\bullet(X),\eta_p,\mu_p)$ satisfying all the properties above. All these triples are (non-canonically) isomorphic. Indeed, using the Levi-Civita connection $\nabla$, we can identify all tangent spaces $T_pQH^\bullet(X)$ by parallel transport. Such an identification is not canonical, since it depends on the paths connecting two points. 

\begin{defi}[\cite{dubro2,CDG}] Let $X$ be a smooth projective complex variety, and let $[(V,\eta,\mu)]$ be the spectrum of $QH^\bullet(X)$.
By abuse of notation, we will denote the group of orthogonal parabolic operators associated with  $(V,\eta,\mu)$  (and its Lie algebra) simply by $\mathcal G(X)$ (by $\frak g(X)$, respectively).
\end{defi}

\bsh
\begin{teorema}[\cite{dubro1,dubro2,CDG}]\label{06.07.17-4} Let $X$ be a smooth projective variety, for which the Gromov--Witten potential $F^X_0$ is convergent, and whose quantum cohomology $QH^\bullet(X)$ is not necessarily semisimple.
\begin{enumerate}
\item The differential system \eqref{06.07.17-1}-\eqref{06.07.17-2} admits analytic fundamental matrix solutions of the form
\[Z(t,z)=\Phi(t,z)z^\mu z^R,
\]where
\[\Phi(t,z)=\sum_{k\in\mathbb N}\Phi_k(t)z^k,\quad \Phi_0\equiv\mathbbm 1,\quad \Phi(t,-z)^T\eta\Phi(t,z)=\eta,
\]and where the matrix $R$ {\rm independent} of $t$ is the matrix associated with an endomorphism in $\frak g(X)$ and computed with respect to the basis $(T_\alpha)_{\alpha=1}^N$ of $H^\bullet(X,\mathbb C)$. A solution of such a form will be said to be in {\rm Levelt  form} at $z=0$. The series $\Phi$  converges in $\mathbb{C}$, because $z=0$ is a Fuchsian singularity of the system \eqref{06.07.17-2}. 

\item Solutions of \eqref{06.07.17-1}-\eqref{06.07.17-2} in Levelt  form at $z=0$ are not unique. Given two of them
\[Z(t,z)=\Phi(t,z)z^\mu z^R,\quad \widetilde Z(t,z)=\widetilde \Phi(t,z)z^\mu z^{\widetilde R},
\]there exists a unique matrix $G$ associated with an endomorphism in $\mathcal G(X)$ (and with respect to the fixed basis $(T_\alpha)_{\alpha=1}^N$) such that
\[\widetilde Z(t,z)=Z(t,z)\cdot G,
\]
\[\widetilde R= G^{-1}\cdot R\cdot G,\quad \widetilde\Phi(t,z)=\Phi(t,z)\cdot P_G(z),
\]where 
\begin{align*}P_G(z):&=z^\mu\cdot G\cdot z^{-\mu}\\
&=\mathbbm 1+z\Delta_1+z^2\Delta_2+\dots   \quad\hbox{(finite sum)},
\end{align*}
the matrices $(\Delta_k)_{k\geq 1}$ being the components of $(G-\mathbbm 1)$.
\item The operator of classical $\cup$-multiplication 
\[c_1(X)\cup (-)\colon H^\bullet(X,\mathbb C)\to H^\bullet(X,\mathbb C)
\]is a $\mu$-nilpotent, $\left\{\cdot,\cdot\right\}$-skew-symmetric endomorphism, i.e. an element of $\frak g(X)$.
\end{enumerate}
\end{teorema}
\esh

\begin{oss}\label{natlev}
    The choice of the  exponent $R$  as in point (3) of Theorem \ref{06.07.17-4} is a \emph{canonical} one, naturally induced by the reduction of system \eqref{06.07.17-1}-\eqref{06.07.17-2} at the classical limit point \eqref{classical1}-\eqref{classical2}. It is a general feature of all \emph{good} Frobenius manifolds, i.e. with good analytical properties in the sense of \cite{dubro1,dubro2}, namely whose potential $F(t)$ is an analytic perturbation of a cubic term 
\[F(t)=\frac{1}{6}c_{\alpha\beta\gamma}t^\alpha t^\beta t^\gamma+\sum_{k,\ell\geq 0}A_{k,\ell}(t'')^\ell\exp(k\cdot t'),
\]
where $k,\ell$ are multi-indices, and the coordinates $(t^\alpha)_{\alpha=1}^N$  are subdivided in two classes $(t',t'')$ with $\deg t'= 0$ and $\deg t''\neq 0$, the degree of the coordinate $t^\beta$ being defined as the constant $(1-q_\beta)$ in the expression for the Euler vector field
\[E=\sum_\alpha \left( (1-q_\alpha)t^\alpha+r_\alpha\right)\frac{\partial}{\partial t^\alpha},\quad r_\alpha\neq 0 \text{ only if }q_\alpha=1.
\]  Under the assumption of convergence of the Gromov-Witten potential, quantum cohomologies of smooth projective varieties are within this class of Frobenius manifolds, as it is manifested from the structure of the Gromov--Witten potential (see equations \eqref{eulerfield}, \eqref{gwpot1}-\eqref{gwpot2}). For this class of Frobenius manifolds it can be shown that the exponent $R$  of point (3) in Theorem \ref{06.07.17-4} can be chosen to be the limit $\mathcal U_{\rm cl}$ of the tensor $\mathcal U$ at the classical point 
\[{\rm Re}(t')\to -\infty,\quad t''=0.
\]
 Sometimes, we will refer to such a choice of Levelt form as the \emph{natural} Levelt  form.
\end{oss}

\begin{defi}[\cite{CDG}]Let $X$ be a smooth projective variety. We define the group $ {\mathcal C}_0(X)$ as the isotropy subgroup of the operator $c_1(X)\cup(-)\in {\rm End}(H^\bullet(X,\mathbb C))$ under the adjoint action
\[{\rm Ad}\colon\mathcal G(X)\to {\rm Aut}(\frak g(X))\colon G\mapsto G\cdot(-)\cdot G^{-1}.
\]In the notations of the paper \cite{CDG}, such a group is denoted by ${\mathcal C}_0(\eta,\mu, c_1(X)\cup(-))$.
\end{defi}

It is clear from Theorem \ref{06.07.17-4} that, even for a fixed $R$, solutions in Levelt  form at $z=0$ are not unique, having a freedom in choice of the series $\Phi$. 

\bsh
\begin{cor}[\cite{CDG}] Under the same assumptions of Theorem \ref{06.07.17-4}, solutions of the differential system \eqref{06.07.17-1}-\eqref{06.07.17-2} in Levelt  form at $z=0$, and with specified exponent $R\equiv c_1(X)\cup(-)$, are not unique. They are parametrized by the group ${\mathcal C}_0(X)$: given two of them, 
\[Z(t,z)=\Phi(t,z)z^\mu z^R,\quad \widetilde Z(t,z)=\widetilde\Phi(t,z)z^\mu z^ R,
\]there exists a unique element $G\in{\mathcal C}_0(X)$ such that 
$$
\widetilde{Z}=ZG,
$$
or equivalently
\[\widetilde\Phi(t,z)=\Phi(t,z)\cdot P_G(z),
\]where $P_G(z):=z^\mu\cdot G\cdot z^{-\mu}$.
\end{cor}
\esh

This freedom in the choice of solutions in Levelt  forms at $z=0$ is a typical phenomenon of all \emph{resonant }Frobenius manifolds\footnote{A Frobenius manifold is called \emph{resonant} if there exist at least two eigenvalues of $\mu$ whose difference is a non-zero integer.} (see \cite{CDG}). Although quantum cohomologies are Frobenius manifolds of this type, the following result shows that in this enumerative-geometrical case a \emph{canonical} choice can be done.

\bsh
\begin{prop}[\cite{dubro1, dubro2, CDG}, \cite{gamma1}]\label{18.07.17-2}
For any smooth projective variety $X$, the system of differential equation \eqref{06.07.17-1}-\eqref{06.07.17-2}  admits the following solution
\begin{align*}
&Z_{\rm top}(z,t):=\Theta_{\rm top}(z,t)\cdot z^\mu z^{c_1(X)\cup(-)},\\
\Theta_{\rm top}(z,t)_\lambda^\gamma:&=\delta_\lambda^\gamma+\sum_{k,n=0}^\infty\sum_{\beta\in{\rm Eff}(X)}\sum_{\alpha_1,\dots,\alpha_k}\frac{h_{\lambda,k,n,\beta,\underline\alpha}^\gamma}{k!}\cdot t^{\alpha_1}\dots t^{\alpha_k}\cdot z^{n+1},\\
h_{\lambda,k,n,\beta,\underline\alpha}^\gamma:=\sum_{\epsilon}&\eta^{\epsilon\gamma}\int_{[\overline{\mathcal M}_{0,k+2}(X,\beta)]^\text{virt}}c_1(\mathcal L_1)^n\cup\text{ev}_1^*T_\lambda\cup\text{ev}_2^*T_{\epsilon}\cup\prod_{j=1}^k\text{ev}^*_{j+2}T_{\alpha_j},
\end{align*}
whose coefficients are Gromov--Witten invariants of $X$ with gravitational descendants. Furthermore, if $X$ is a Fano manifold, among all solutions $$\Phi(t,z)z^\mu z^{c_1(X)\cup(-)}$$ in Levelt form at $z=0$, the {\rm topological-enumerative} solution above is the unique one for which the product
\[H(t,z)=z^{-\mu}\Phi(t,z)z^\mu,
\]computed at points of the small quantum cohomology (i.e. $t^i=0$ unless $i=2,\dots r$), 
\begin{enumerate}
\item is holomorphic at $z=0$,
\item and moreover \[H(t,0)=\exp\left(\left(\sum_{i=2}^rt^iT_i\right)\cup(-)\right).
\]
\end{enumerate}
\end{prop}
\esh

\subsection{Idempotent vielbein and $\Psi$-matrix}
\label{psigauge}
 
 Let us introduce the following definition.

\begin{defi}\label{17agosto2018-1}
The \emph{bifurcation set} $\mathcal B_X\subseteq QH^\bullet(X)$ is defined as the set of points at which the spectrum of $\mathcal U(p)$, i.e. the operator of quantum multiplication by the Euler vector field (see equation \eqref{06.07.17-3}), is not simple. The \emph{caustic} $\mathcal K_X$ is defined as the set of points $p\in QH^\bullet(X)$ at which the associated Frobenius algebra is not semisimple.
\end{defi}

\begin{teorema}[\cite{dubronapoli, dubro1, dubro2, CDG}]
Let $\Omega\subseteq QH^\bullet(X)\setminus \mathcal K_X$ be a simply connected open set.
 It is possible to label the eigenvalues $\left\{u_i\right\}_{i=1}^{N}$ of $\mathcal U(p)$ at any $p\in\Omega$ in 
 such a way that  an $N$-tuple $(u_1,\dots, u_N)$ of single-valued and holomorphic functions is well-defined 
 on $\Omega$. Furthermore, such functions can be used as a system of local coordinates
  (called \emph{canonical}) on $\Omega$: at any $p\in\Omega$ the coordinate vector fields 
  $\left\{\frac{\partial}{\partial u_i}\right\}_{i=1}^{N}$ coincide with the idempotents vector fields 
  of the Frobenius algebra $(T_p\Omega,\circ_p,\eta_p)$, i.e.
\[\frac{\partial}{\partial u_i}\circ\frac{\partial}{\partial u_i}=\delta_{ij}\frac{\partial}{\partial u_i},\quad \eta\left(\frac{\partial}{\partial u_i},\frac{\partial}{\partial u_j}\right)=0\text{ if }i\neq j.
\] 
\end{teorema}

\begin{defi}Let $\Omega\subseteq QH^\bullet(X)\setminus\mathcal K_X$ be a simply connected open subset on which an ordering of canonical coordinates $(u_1,\dots, u_{N})$ has been fixed. Let us define at any point $p\in\Omega$ a system of normalized idempotent vector fields 
\begin{equation}\label{5luglio2017-1}f_i|_p:=\frac{1}{\eta\left(\frac{\partial}{\partial u_i},\frac{\partial}{\partial u_i}\right)^\frac{1}{2}}\cdot \left.\frac{\partial}{\partial u_i}\right|_p,\quad i=1,\dots, N,
\end{equation}where a determination of the square roots has been chosen in such a way that $f_i$ are (necessarily single-valued and) holomorphic on $\Omega$. We define on $\Omega$ the matrix $\Psi$ through the equation
\[\frac{\partial}{\partial t^\alpha}=\sum_i\Psi_{i\alpha}f_i.
\]

\end{defi}

\subsection{Monodromy data as local moduli} \label{monlocalmod}  Once an ordering of canonical coordinates and a determination of the matrix $\Psi$ has been fixed on a simply-connected open subset $\Omega\subseteq QH^\bullet(X)\setminus \mathcal K_X$, the system \eqref{06.07.17-1}-\eqref{06.07.17-2} can be rewritten in the idempotent vielbein. If $Y:=\Psi Z$, then we have 
\begin{align}
\label{16.07.17-1}
\partial_iY&=(zE_i+V_i)Y,\\
\label{16.07.17-2}
\partial_z Y&=\left(U+\frac{1}{z}V\right)Y,
\end{align}
where $\partial_i:=\frac{\partial}{\partial u_i}$, $(E_i)_{hk}=\delta_{ih}\delta_{ik}$, $V_i:=\partial_i\Psi\cdot\Psi^{-1}$, $U:={\rm diag}(u_1,\dots, u_{N})$ with $u_i:=u_i(p)$, $p\in\Omega$, and $V:=\Psi\cdot\mu\cdot\Psi^{-1}$. The general theory developed in \cite{dubro1, dubro0, dubro2}, extended and refined in \cite{CDG0, CDG}, provides local invariants of the Frobenius structure defined on $QH^\bullet(X)$ through the study of the monodromy of \eqref{06.07.17-2}, or equivalently \eqref{16.07.17-2}.

Since at $z=\infty$ the system \eqref{06.07.17-2} admits an irregular singularity (a singularity of second kind), in order to describe the Stokes phenomenon let us fix an oriented line $\ell$ in the complex plane. We recall the following definitions, following the general description given in \cite{CDG}. 

\begin{defi}\label{chamberdefi}
An oriented line $\ell$ in the complex plane will be said to be \emph{admissible at a point} $p\in QH^\bullet(X)$ if $\ell$ does not contain any Stokes ray at $p$, i.e. any ray
\[R_{ij}(p):=\left\{z\in\mathbb C\colon z=\-\sqrt{-1}\rho(\overline{u_i}-\overline{u_j}),\ \rho\in\mathbb R_+\right\},\quad i,j=1,\dots, N,\quad i\neq j.
\] 
An open connected component of  the set of points  $p\in QH^\bullet(X)$  satisfying 
\begin{enumerate}
\item the canonical coordinates $(u_i(p))_{i=1}^{N}$ are pairwise distinct,
\item the line $\ell$ is admissible at $p$.
\end{enumerate}
will be called an $\ell$-chamber, and will be denoted by  $\Omega_\ell$.
\end{defi}
\begin{oss}
Note that the definition above is well-posed since it does not depend on the choice either of the order of canonical coordinates or of a branch of the $\Psi$-matrix. The topology of an $\ell$-chamber $\Omega_\ell$ can be non trivial, e.g. it can be non-simply-connected. Despite of this, it can be shown (see \cite{CDG}) that on any $\ell$-chamber canonical coordinates can be coherently labeled so that they define a system of holomorphic and single-valued functions on $\Omega_\ell$, as well as a single-valued determination of a branch of the $\Psi$-matrix is possible. This follows from the fact that for any $z\in \Omega_\ell$ the inclusions 
\[\xymatrix{
\Omega_\ell\ \ar@{^{(}->}[rr]^{i\quad\quad}&& QH^\bullet(X)\setminus\mathcal K_X\ \ar@{^{(}->}[rr]^{j}&&QH^\bullet(X),
}
\]induce morphisms in homotopy 
\[\xymatrix{
\pi_1(\Omega_\ell,z)\ \ar@{^{(}->}[rr]^{i_*\quad\quad}&& \pi_1(QH^\bullet(X)\setminus\mathcal K_X,z)\ \ar@{^{(}->}[rr]^{j_*}&&\pi_1(QH^\bullet(X),z),
}
\]such that ${\rm im}(i_*)\cap{\rm ker}(j_*)=\left\{0\right\}$.
\end{oss}

\bsh
\begin{teorema}[\cite{dubro2, CDG0, CDG}]\label{16.07.17-3}
Let $\Omega\subseteq QH^\bullet(X)\setminus\mathcal K_X$ be a simply-connected subset with a fixed ordering of canonical coordinates $u\colon \Omega\to\mathbb C^{N}\colon p\mapsto (u_1(p),\dots, u_{N}(p))$, and a fixed holomorphic branch of the $\Psi$-matrix. 
 The following facts hold true on  a suitable restriction of $\Omega$,
\begin{enumerate}
\item The system of differential equations \eqref{16.07.17-1}-\eqref{16.07.17-2} admits a unique formal solution of the form
\begin{align*}Y_{\rm formal}(z,u)&=F(z,u)\exp(zU)\\
F(z,u)=\sum_{m\geq 0}\frac{F_m(u)}{z^m},\quad F_0&\equiv\mathbbm 1,\quad F(-z,u)^T\cdot F(z,u)=\mathbbm 1,
\end{align*}
where the functions $F_m$'s are holomorphic on $u(\Omega)$.
\item Let $\ell$ be an oriented line of slope $\phi\in[0;2\pi[$ in the complex plane. For any $k\in\mathbb Z$, there exist two solutions $Y^{(k)}_{\rm left/right}$of the system \eqref{16.07.17-1}-\eqref{16.07.17-2}, analytic and single-valued on $\mathcal R\times u(\Omega)$ and uniquely characterized by the asymptotic expansion 
\[Y^{(k)}_{\rm left/right}(z,u)\sim Y_{\rm formal }(z,u),\quad |z|\to\infty,\quad z\in e^{2\pi i k}\Pi_{\rm left/right}(\phi),
\]uniformly on any compact subset of $u(\Omega)$, and where 
\begin{align*}\Pi_{\operatorname{right}} (\phi):&=\left\{z\in\mathcal R\colon \phi-\pi< \arg\ z< \phi\right\},\\
\Pi_{\operatorname{left}} (\phi):&=\left\{z\in\mathcal R\colon \phi< \arg\ z< \phi+\pi\right\}.
\end{align*}

\item For any $u\in u(\Omega)$, for any $k\in\mathbb Z$ we have that
\[Y^{(k)}_{\operatorname{right/left}}(e^{2\pi ki}z,u)=Y^{(0)}_{\operatorname{right/left}}(z,u),\quad z\in\mathcal R.
\]
\end{enumerate}
\end{teorema}
\esh

\begin{oss}
The precise meaning of the asymptotic relation in (2) of Theorem \ref{16.07.17-3} is the following:  
\[
\forall K\Subset u(\Omega),
\
 \forall h\in\mathbb N,
 \
  \forall \overline{\mathcal S}\subsetneq e^{2\pi k i }\Pi_{\operatorname{right/left}} (\phi),
  \
   \exists C_{K,h,\overline{\mathcal S}}>0\colon \text{ if }z\in\overline{\mathcal S}\setminus\left\{0 \right\}\text{ then }\]
   \[
    \sup_{u\in K}\left\| Y^{(k)}_{\operatorname{right/left}}(z,u)\cdot\exp(-z U(u))-\sum_{m=0}^{h-1}\frac{F_m(u)}{z^m}\right\|
    <\frac{C_{K,h,\overline{\mathcal S}}}{|z|^h}.
\]
Here $\overline{\mathcal S}$ denotes any unbounded closed sector of $\mathcal R$ with vertex at $0$.  
\end{oss}

\begin{defi}\label{defistokcenconn}
Let $p\in QH^\bullet(X)\setminus\mathcal K_X$ be a semisimple point and fix an ordering of canonical coordinates $u(p)=(u_1(p),\dots, u_{N}(p))$ and a branch of the matrix $\Psi(p)$. If the oriented line $\ell$ of slope $\phi\in[0;2\pi[$ is admissible at $p$, we define the matrices $(S^{(k)},S^{(k)}_-, C^{(k)})_{k\in\mathbb Z}$ through the equations
\begin{align*}
Y^{(k)}_{\rm left}(z,u(p))&=Y^{(k)}_{\rm right}(z,u(p))S^{(k)},\quad z\in\mathcal R,\\
Y^{(k)}_{\rm left}(e^{2\pi i}z,u(p))&=Y^{(k)}_{\rm right}(z,u(p))S^{(k)}_-,\quad z\in\mathcal R,\\
Y^{(k)}_{\rm right}(z,u(p))&=Y_{0}(z,u(p))C^{(k)},\quad z\in\mathcal R.
\end{align*}
Here,  
\begin{itemize}
\item the functions $Y^{(k)}_{\rm left/right}(z,u)$ denote the solutions of \eqref{16.07.17-1}-\eqref{16.07.17-2} described in Theorem \ref{16.07.17-3}, and defined on $\mathcal R\times u(\Omega)$ for a sufficiently small open simply-connected neighborhood of $p$;
\item the function $Y_0(z,u)$ is a solution of \eqref{16.07.17-1}-\eqref{16.07.17-2} of the form
\[Y_0(z,u(p))=\Psi(p)\cdot Z(z,t(p)),
\]where $Z(z,t)=\Phi(z,t)z^\mu z^R$ is a solution of \eqref{06.07.17-1}-\eqref{06.07.17-2} in Levelt form at $z=0$. 
We will say that also the solution $Y_0$ is in Levelt  form at $z=0$. 
\end{itemize}
The matrices $$S:=S^{(0)},\quad S_-:=S^{(0)}_-,\quad C:=C^{(0)}$$ are respectively called \emph{Stokes and Central Connection matrices of $QH^\bullet(X)$} at the point $p$ (with respect to the fixed ordering of canonical coordinates, the fixed branch of the $\Psi$-matrix, the fixed oriented line $\ell$, the fixed solution $Y_0$ in Levelt  form at $z=0$).  
\end{defi}

The following Theorems summarize the main properties of $(S, S_-,C)$.

\bsh
\begin{teorema}\cite{dubro2, CDG}\label{19.03.18-1} Let $p\in QH^\bullet(X)\setminus\mathcal K_X$ be a semisimple point. For  
\begin{itemize}
\item any fixed ordering of canonical coordinates at $p$, 
\item any choice of branch of the $\Psi$-matrix at $p$, 
\item any oriented line $\ell$ of slope $\phi\in[0;2\pi[$ admissible at $p$, 
\item and any solution $Y_0$ in Levelt  form at $z=0$, 
\end{itemize} the Stokes matrices $S,S_-$ and the central connection matrix $C$ at $p$ satisfy the following properties:
\begin{enumerate}
\item the whole family of matrices $(S^{(k)},S^{(k)}_-, C^{(k)})_{k\in\mathbb Z}$ can be reconstructed from the triple $(S,S_-,C)$. Namely, for all $k\in\mathbb Z$ 

$$
S^{(k)}=S,\quad S_-^{(k)}=S_-,\quad C^{(k)}=M_0^{-k}~ C,
$$

where $M_0=\exp(2\pi i \mu)\exp(2\pi i R)$.
\item For all $k\in\mathbb Z$ and all $z\in\mathcal R$
\begin{align*}
Y^{(k)}_{\operatorname{right}}(e^{2\pi i}z,u(p))&=Y^{(k)}_{\operatorname{right}}(z,u(p))~ S_-~ S^{-1},\\
Y^{(k)}_{\operatorname{left}}(e^{2\pi i}z,u(p))&=Y^{(k)}_{\operatorname{right}}(z,u(p))~ S^{-1}~ S_-.
\end{align*}
\item We have that
\begin{align*}
 S_-&= S^T,\\
  S_{ii}=1,&\quad i=1,\dots, N,\\
      S_{ij}~{\rm can ~ be}~\neq 0\text{ with }i\neq j\text{ only if }&u_i\neq u_j \text{ and }R_{ij}\subset \Pi_{\operatorname{left}}(\phi).
\end{align*}
 \item $CS^TS^{-1}C^{-1}=M_0=e^{2\pi i\mu}e^{2\pi iR};$
\item $S=C^{-1}e^{-\pi i R}e^{-\pi i \mu}\eta^{-1}(C^T)^{-1};$
\item $S^T=C^{-1}e^{\pi i R}e^{\pi i \mu}\eta^{-1}(C^T)^{-1}.$
\end{enumerate}
\end{teorema}
\esh

\bsh
\begin{teorema}\cite{dubro2, CDG}\label{iso2}
Let $\ell$ be an oriented line in the complex plane, of slope $\phi\in[0;2\pi[$, and let $\Omega_\ell$ be an $\ell$-chamber of $QH^\bullet(X)$. Let
\begin{itemize}
\item $u\colon\Omega_\ell\to\mathbb C^{N}\setminus\, \{{\rm diagonals}\}$ be a single-valued and holomorphic function, which defines a coherent ordering of canonical coordinates on $\Omega_\ell$, 
\item  $\Psi\colon \Omega_\ell\to GL(N,\mathbb C)$ be  a single-valued determination  of the $\Psi$-matrix on $\Omega_\ell$, 
\item $Y_0$ be  a solution of \eqref{16.07.17-1}-\eqref{16.07.17-2} in Levelt  form at $z=0$.
\end{itemize}
Then, the corresponding Stokes and central connection matrices $(S,S_-,C)$ are constant on $\Omega_\ell$.
\end{teorema}
\esh

\bsh
\begin{teorema}\cite{CDG}\label{iso2.1}
Let $p\in\mathcal B_X\setminus\mathcal K_X$ be a semisimple coalescing point of $QH^\bullet(X)$, and let $\ell$ be an oriented line admissible at $p$. Let $\Omega\subseteq QH^\bullet(X)\setminus\mathcal K_X$ be an open simply connected neighborhood of $p$, on which an ordering $u\colon \Omega\to\mathbb C^{N}$ of canonical coordinates and a holomorphic branch $\Psi\colon\Omega\to GL(N,\mathbb C)$ of the $\Psi$-matrix have been fixed. If $\Omega$ is sufficiently small, the monodromy data $(S,S_-, C)$ computed at $p$ are the same computed at any other point of $\Omega$ at which $\ell$ is admissible. In particular, they are the same data computed\footnotemark$\ $ in any $\ell$-chamber with non-empty intersection with $\Omega$.
\end{teorema}
\esh
\footnotetext[16]{Here, the ordering of the canonical coordinates and the holomorphic branch of $\Psi$ on any $\ell$-chamber are the ones prolonged from $\Omega$.}

\begin{defi}Let $p\in QH^\bullet(X)\setminus\mathcal K_X$ be a semisimple point of the Frobenius manifold $QH^\bullet(X)$.
We call \emph{monodromy data} (or \emph{monodromy local moduli}) of $QH^\bullet(X)$ at $p$ the tuple $(\mu, R,S,C)$, where
\begin{itemize}
\item the tensor $\mu$ is the grading operator defined by equation \eqref{06.07.17-3},
\item the matrix $R$ is the matrix associated with the operator $c_1(X)\cup(-)\colon H^\bullet(X,\mathbb C)\to H^\bullet(X,\mathbb C)$ with respect to the fixed basis $(T_\alpha)_{\alpha=1}^N$,
\item the Stokes and central connection matrices $(S,C)$ are defined as in Definition \ref{defistokcenconn}.
\end{itemize}
The monodromy data define local invariants of the Frobenius structure, as described in Theorems \ref{06.07.17-4}, \ref{iso2}, \ref{iso2.1}.
\end{defi}

As explained in the previous paragraphs, the definition of the monodromy data $(S,C)$ is subordinate to many non-canonical choices, namely:\begin{enumerate}
\item the choice of an oriented line $\ell(\phi)=\{z=\rho e^{i\phi},~\rho\in\mathbb{R}\}$ in the complex plane, with slope $\phi\in[0;2\pi[$;
\item the choice of a point in the fiber $\Pi^{-1}(1)$ of the universal cover $\Pi\colon\mathcal R\to \mathbb C^*$, in order to fix the sheet with principal value of the argument $\arg z$ in the interval $[0;2\pi[$; 
\item the choice of {an} ordering of canonical coordinates on each $\ell$-chamber $\Omega_\ell$;
\item the choice of the signs of the square roots (\ref{5luglio2017-1}) defining the normalized idempotent vielbein $(f_i)_{i=1}^{N}$, and hence the matrix $\Psi$ on each $\ell$-chamber $\Omega_\ell$;
\item the choice of a {solution} $Y_0$ in the Levelt form corresponding to the same exponent $R$.
\end{enumerate} 

A detailed analysis of the effects of different choices on the numerical values of the data $(S,C)$ has been developed in \cite{CDG}: in particular, the freedom  in the choices (1)-(5) above can be quantified through the action of suitable groups on the set of Stokes and Central Connection matrices. Here, we briefly summarize the main results, and we refer the reader to \cite{CDG} for more details.

\bsh
\begin{teorema}\cite{CDG} Let $p\in QH^\bullet(X)\setminus\mathcal K_X$ be a semisimple point, and let $(S,C)$ be the Stokes and Central Connection matrices computed at $p$ with respect to some choices of normalizations (2)-(5) above. All other possible values of the data $(S,C)$ corresponding to different choices can be obtained through the actions of the following groups.
\begin{itemize}
\item \emph{Action of the symmetric group} $\frak S_{N}$: the permutation $\tau$ corresponding to the re-ordering 
\[(u_1,\dots, u_{N})\mapsto (u_{\tau(1)},\dots, u_{\tau(N)})
\]
acts on the monodromy data as follows
\[S\mapsto PSP^{-1},\quad C\mapsto CP^{-1},\quad P_{ij}:=\delta_{j\tau(i)}.
\]
\item  \emph{Action of the group} $(\mathbb Z/2\mathbb Z)^{N}$: different choices of signs in equations \eqref{5luglio2017-1} correspond to transformations
\[S\mapsto ISI,\quad C\mapsto CI,\]
where $I$ is a diagonal matrix with entries $\pm1$.
\item \emph{Action of the group} $\mathcal C_0(X)$: a different choice $Y_0\mapsto Y_0G$, with $G\in\mathcal C_0(X)$, of the solution in Levelt  form at $z=0$ acts on the monodromy data as
\[S\mapsto S,\quad C\mapsto G^{-1}C.
\]
\item \emph{Action of the Galois group }$Deck(\Pi)\cong \mathbb Z$: different choices of a base point in the fiber $\Pi^{-1}(1)$, namely of the principal determination of the argument $\arg z$, correspond to the transformations
\[S\mapsto S,\quad C\mapsto M_0^{-k}C,\quad k\in\mathbb Z,
\] where $M_0=\exp(2\pi i \mu)\exp(2\pi i R)$.
\end{itemize}
\end{teorema}
\esh

\begin{defi}[Triangular order]
\label{30luglio2016-2}Let $p\in QH^\bullet(X)\setminus\mathcal K_X$ be a semisimple point, and let $S$ be the Stokes matrix  computed at $p$ with respect to some admissible oriented line $\ell$. We say that $(u_1(p),...,u_N(p))$ are in \emph{triangular} order with respect to the line $\ell$ whenever $S$ is upper triangular. 
\end{defi}

Notice that the triangularity of an ordering of canonical coordinates only depends on the choice of the oriented line $\ell$, according to point (3) of Theorem \ref{19.03.18-1}. In general, triangular orders of canonical coordinates at a semisimple point $p$ are not unique\footnote{This is the case for a semisimple coalescing point $p$, and for  the points of all $\ell$-chambers intersecting a sufficiently small neighborhood of $p$.}. Among all possible triangular {orderings} of the canonical coordinates,  a particularly convenient one is  the \emph{lexicographical order} w.r.t an admissible line $\ell(\phi)$, defined as follows. 

\begin{defi}[Lexicographical order]\label{30luglio2016-1}  Let $p\in QH^\bullet(X)\setminus\mathcal K_X$ be a semisimple point, and let $\ell$ be an admissible oriented line. Let us consider the rays starting from the points $u_1(p),\dots, u_N(p)$ in the complex plane
\[L_j:=\left\{u_j(p)+\rho e^{i\left(\frac{\pi}{2}-\phi\right)}\colon \rho\in\mathbb R_+\right\},\quad j=1,\dots, n,
\]and for any complex number $z_0$ let us define the oriented line
\[L_{z_0,\phi}:=\left\{z_0+\rho e^{-i\phi}\colon\rho\in\mathbb R\right\}
\]where the orientation is induced by $\mathbb R$. In this way we have a natural total order $\preceq$ on the points of $L_{z_0,\phi}$. We can choose $z_0$, with $|z_0|$ sufficiently large, so that the intersections
$L_{j}\cap L_{z_0,\phi}=:\left\{p_j\right\}
$ are non-empty. The canonical coordinates $u_j(p)$'s are in the \emph{$\ell$-lexicographical order} if
\[p_1\preceq p_2\preceq p_3\preceq\dots\preceq p_N.
\]
The definition does not depend on the choice of $z_0\in\mathbb C$, with $|z_0|$ sufficiently large.
\end{defi}

\begin{defi}
If $\ell$ has been chosen and $(u_1,...,u_n)$ are in lexicographical order as in Definition \ref{30luglio2016-1}, we say that the monodromy data $S$ and $C$ computed with respect to the above $\ell$ at the above point $u$ are {\it monodromy data in lexicographical order}.
\end{defi}

\begin{oss}Observe that if $u_1,\dots, u_N$ are in the lexicographical order with respect to the admissible line $\ell(\phi)$, then:
\begin{enumerate}
\item the Stokes matrix is in upper triangular form;
\item the nearest Stokes rays to the positive half-line $\operatorname{pr}(\ell_+(\phi))$ are of the form
\[R_{i,i+1}\subseteq \Pi_{\rm left}(\phi),
\quad\quad
R_{j,j-1}\subseteq \Pi_{\rm right}(\phi),\quad
\]where $1\leq i\leq N-1$ and $2\leq j\leq N$. 
\end{enumerate}

In general, condition $(1)$  alone   does not imply that the canonical coordinates are in the lexicographical order: it does if and only if the number of nonzero entries of the Stokes matrix $S$ is maximal (and equal to  $\frac{N(N+1)}{2}$) (see \cite{CDG} for further details).
\end{oss}

At this point, we can finally describe 
\begin{enumerate}
\item how  the monodromy data $(S,C)$, computed at a semisimple point $p\in QH^\bullet(X)\setminus\mathcal B_X$, change when the oriented line $\ell$, admissible at $p$, changes direction $\phi$;
\item or, dually, how the values of the pair $(S,C)$ computed in different $\ell$-chambers, for some fixed choice of the line $\ell$, are related to each other.
\end{enumerate}
In both cases, this is described by the action of the \emph{braid group} $\mathcal B_{N}$. Recall that this group is generated by $N$ elementary braids $\beta_{12},\beta_{23},\dots,\beta_{N-1,N}$ with the relations
\[\beta_{i,i+1}\beta_{j,j+1}=\beta_{j,j+1}\beta_{i,i+1}\quad\text{for }i+1\neq j,j+1\neq i,
\]
\[\beta_{i,i+1}\beta_{i+1,i+2}\beta_{i,i+1}=\beta_{i+1,i+2}\beta_{i,i+1}\beta_{i+1,i+2}.
\]

\begin{defi}
\label{10novembre2018-2}
Let $p\in QH^\bullet(X)\setminus \mathcal B_X$ be a semisimple point. Let $\ell$ be an oriented line admissible at $p$, and let $(S, C)$ be the monodromy data computed at $p$ with respect to the line $\ell$ and in the $\ell$-lexicographical order. For any elementary braid $\beta_{i,i+1}\in\mathcal B_{N}$ let us define the transformed set of matrices $(S^{\beta_{i,i+1}}, C^{\beta_{i,i+1}})$ through the equations 
\begin{equation}
\label{stokesbraid1}
S^{\beta_{i,i+1}}:= A^{\beta_{i,i+1}}(S)~ S~ A^{\beta_{i,i+1}}(S),
\quad C^{\beta_{i,i+1}}:= C~ (A^{\beta_{i,i+1}})^{-1},
\end{equation}
where
\[
\left(A^{\beta_{i,i+1}}(S)\right)_{hh}=1,\quad\quad h=1,\dots, N\quad h\neq i,i+1,
\]
\[\left(A^{\beta_{i,i+1}}(S)\right)_{i+1,i+1}=-s_{i,i+1},
\]
\[
\left(A^{\beta_{i,i+1}}(S)\right)_{i,i+1}=\left(A^{\beta_{i,i+1}}(S)\right)_{i+1,i}=1.
\]
For a generic braid $\beta$, which is a product of $m$ elementary braids $\beta=\beta_{i_1,i_{1}+1}\dots\beta_{i_m,i_{m}+1}$, the action is
\begin{equation}
\label{connbraid2-29luglio}
S\mapsto S^\beta:= A^\beta(S)\cdot S \cdot A^\beta(S)^T,\quad C\mapsto C^\beta:=C\cdot (A^\beta)^{-1},
\end{equation}
where
\[A^\beta(S)=A^{\beta_{i_m,i_{m}+1}}\left(S^{\beta_{i_{m-1},i_{m-1}+1}}\right)\cdot
... \cdot A^{\beta_{i_2,i_{2}+1}}\left(S^{\beta_{i_1,i_{1}+1}}\right)\cdot A^{\beta_{i_1,i_{1}+1}}(S).
\]This defines an action (on the right) of the braid group on the set of the data $(S,C)$. Notice that the triangularity of the Stokes matrices is preserved. In the paper, we will often call {\it mutations} of $S$ and $C$ any of the matrices $S^\beta$, $C^\beta$ obtained by the action above of the braid group.

\end{defi}

\bsh
\begin{teorema}\cite{dubro1,dubro2,CDG}
Let $p\in QH^\bullet(X)\setminus \mathcal B_X$ be a semisimple point. Let $\ell$ be an oriented line admissible at $p$, and let $(S_{\rm lex}, C_{\rm lex})$ be the monodromy data computed at $p$ with respect to the line $\ell$ and in the $\ell$-lexicographical order. Without loss of generality, we can suppose that the Stokes rays are all distinct\footnotemark. 

\begin{enumerate}

\item Let us consider a point $p'\in QH^\bullet(X)\setminus \mathcal B_X$ in another $\ell$-chamber, and let us consider a path $\gamma\colon [0,1]\to QH^\bullet(X)$ with $\gamma(0)=p$ and $\gamma(1)=p'$. If along $\gamma$ the only Stokes ray\footnotemark $\ R_{i,i+1}(t)$, with $t\in[0,1]$, crosses the line $\ell$ in the clockwise direction, or equivalently the eigenvalue $u_i(t)$ rotates in the counter-clockwise direction with respect to the eigenvalue $u_{i+1}(t)$, then the monodromy data at $p'$ in the $\ell$-lexicographical order are given by
\[(S_{\rm lex}^{\beta_{i,i+1}}, C_{\rm lex}^{\beta_{i,i+1}}).
\]
If along $\gamma$ more Stokes rays cross $\ell$ (or, equivalently, more eigenvalues $u_i$'s rotate with respect to each other) then the resulting monodromy data are obtained from the composition of elementary braids transformations as in \eqref{connbraid2-29luglio}.

\item Let us consider a counter-clockwise rotation of the line $\ell$ into another oriented line $\ell'$ admissible at $p$. If, during the rotation, the line $\ell$ crosses only the Stokes ray $R_{i,i+1}$, then the values of the monodromy data $(S_{\rm lex}', C_{\rm lex}')$ computed at $p$ with respect to the line $\ell'$ and in the $\ell'$-lexicographical order are given by 
\[S_{\rm lex}'=S_{\rm lex}^{\beta_{i,i+1}},\quad C_{\rm lex}'=C_{\rm lex}^{\beta_{i,i+1}}.
\]If during the rotation the line $\ell$ crosses more Stokes rays, then the resulting data are obtained from the composition of elementary braids transformations as in \eqref{connbraid2-29luglio}.
\end{enumerate}
\end{teorema}
\esh
\footnotetext[18]{If this is not the case, we can choose another point in the same $\ell$-chamber (so that the corresponding monodromy data $(S_{\rm lex}, C_{\rm lex})$ are the same) with this property.}
\footnotetext[19]{Here the labeling of the Stokes rays $R_{i,i+1}(t)$ is the one prolonged from the initial point at $t=0$.}

\bsh
\begin{cor}\cite{CDG}\label{centerbraidlemma}
The braid corresponding to a complete counter-clockwise $2\pi$-rotation of $\ell$ is 
\[(\beta_{12}\beta_{23}\dots\beta_{N-1,N})^{N},
\]and its acts on the monodromy data as follows:
\begin{itemize}
\item trivially on Stokes matrices,
\item the central connection matrix is transformed as $C\mapsto M_0^{-1}C$.
\end{itemize}
\end{cor}
\esh

\subsubsection{Inverse Problem for Frobenius manifolds}\label{invfrob}
The monodromy data $(\mu,R,S,C)$, defined as in the previous paragraphs, can be interpreted as \emph{local moduli} for the semisimple Frobenius manifold structures. Let us briefly recall how they define a sort of \emph{coordinate system} in the space of solutions of WDVV equations, and how the Frobenius structure can be locally reconstructed from their knowledge. 

The starting point of this observation is the following. Assume we are already given a semisimple Frobenius manifold, i.e. (for the interests of this paper) a quantum cohomology of a smooth projective variety $X$. Let $$Y_0(z,u)=\left(\sum_{p=0}^\infty\phi_{p}(u)z^p\right)z^{\mu}z^R,\quad \phi_p(u):=\left(\phi_{i\alpha,p}(u)\right)_{i,\alpha=1}^N,\quad \phi_0(u)=\Psi(u)$$ be a solution of the system \eqref{16.07.17-2} in Levelt form. Then the following parametric formulae hold true, provided that $\prod_{i=1}^N\phi_{i1,0}(u)\neq 0$:
\begin{align}
\label{12.11.18-1}\eta_{\alpha\beta}&=\sum_{i=1}^N\phi_{i\alpha,0}(u)\phi_{i\beta,0}(u),\\
\frac{\partial}{\partial t^1}&=\sum_{i=1}^N\frac{\partial}{\partial u_i},\\
E&=\sum_{i=1}^Nu_i\frac{\partial}{\partial u_i},\\
t_\alpha(u)&=\sum_{i=1}^N\phi_{i\alpha,0}(u)\phi_{i1,1}(u),\quad t_\alpha:=\eta_{\alpha\beta}t^\beta,\\
\label{12.11.18-5}F(t(u))=\frac{1}{2}\left\{t^\alpha t^\beta\sum_{i=1}^n\phi_{i\alpha,0}(u)\right.&\left.\phi_{i\beta,1}(u)-\sum_{i=1}^N\left(\phi_{i1,1}(u)\phi_{i1,2}(u)+\phi_{i1,3}(u)\phi_{i1,0}(u)\right)\right\}.
\end{align}
This means that from the knowledge of the only matrices $\phi_0(u),\phi_1(u),\phi_2(u),\phi_3(u)$ we are able to locally reconstruct the Frobenius manifold structure. The crucial point is that these matrices can be reconstructed from the only datum of $(\mu, R, S, C)$ attached to a chamber, through a Riemann-Hilbert boundary value problem (RH b.v.p.).

Assume indeed that, conversely, we are given the datum of $(\mu, R,S,C)$, of a fixed point $u^{(0)}=(u^{(0)}_1,\dots, u^{(0)}_N)\in\mathbb C^N$ and an admissible oriented line $\ell$ such that
\begin{itemize}
\item properties (3),(4),(5),(6) of Theorem \ref{19.03.18-1} are satisfied\footnote{Here we assume that $\eta$ is already known. If this is not the case, then only points (3) and (4) of Theorem \ref{19.03.18-1} must be satisfied. The metric $\eta$ can thus be reconstructued through equation \eqref{12.11.18-1}.}, 
\item $u^{(0)}_i\neq u^{(0)}_j$ for $i\neq j$.
\end{itemize}

Let $D$ be a sufficiently small  disc in the complex plane with center $z=0$, and denote by 
\begin{itemize}
\item $P_R$ and $P_L$ the external parts of $D$ on the right and on the left, respectively, of the oriented line $\ell$,
\item $\tilde{\ell}_+$ and $\tilde{\ell}_-$ the parts of $\ell_+$ and $\ell_-$, respectively, on the common border of $P_R$ and $P_L$,
\item $\partial D_R$ and $\partial D_L$ the parts of $\partial D$ which border $P_R$ and $P_L$ respectively.
\end{itemize}

\begin{figure}[ht!]
\centering
\def\svgscale{1}
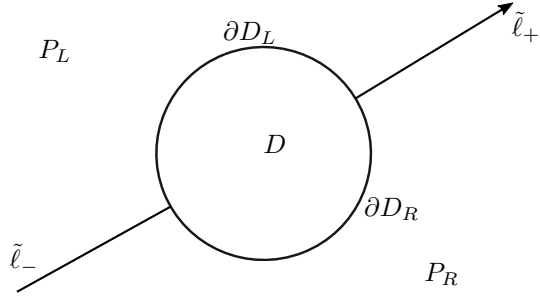
\caption{Boundary for the Riemann-Hilbert problem.}
\label{RH}
\end{figure}

Consider now the following\newline

\noindent
{\bf Problem. }Find a piecewise analytic function $\Phi(z)$ defined in $\mathbb C$ with
\[\Phi(z)=\begin{sistema}
\Phi_R(z),\quad z\in P_R,\\
\Phi_L(z),\quad z\in P_L,\\
\Phi_0(z),\quad z\in D,
\end{sistema}
\]such that the following conditions are satisfied:
\begin{enumerate}
\item $\Phi_R\in\mathcal O(P_R)\cap\mathcal C^0(\overline{P_R}),\quad \Phi_L\in\mathcal O(P_L)\cap\mathcal C^0(\overline{P_L}),\quad \Phi_0\in\mathcal O(D)\cap\mathcal C^0(\overline{D})$,
\item $\Phi_L(z)=\Phi_R(z)e^{zU}Se^{-zU},\quad z\in \tilde{\ell}_+$,
\item $\Phi_L(z)=\Phi_R(z)e^{zU}S^Te^{-zU},\quad z\in \tilde{\ell}_-$,
\item $\Phi_0(z)=\Phi_R(z)e^{zU}C^{-1}z^{-R}z^{-\mu},\quad z\in\partial D_R$,
\item $\Phi_0(z)=\Phi_L(z)e^{zU}S^{-1}C^{-1}z^{-R}z^{-\mu},\quad z\in\partial D_L$,
\item $\Phi_{L/R}(z)\to\mathbbm 1$ if $z\to \infty$ in $P_{L/R}$ (uniformly in closed sub-sectors).
\end{enumerate}

\bsh
\begin{teorema}[\cite{Malg3, malbook},\cite{Miwa}, \cite{dubro2}]\label{teomiwmal}
If the RH b.v.p. above has a solution for the point $u^{(0)}$, then 
\begin{enumerate}
\item the solution is unique,
\item it admits a unique solution for all $u$ in a sufficiently small neighborhood of $u^{(0)}$, analytically depending on $u$.
\end{enumerate}
The solution analytically continues to a meromorphic function on the universal cover of $\mathbb C^N\setminus \Delta$, where 
\[\Delta:=\left\{u\in\mathbb C^N\colon u_i=u_j,\quad \text{for some }i\neq j\right\}.
\]
\end{teorema}
\esh
The poles of the analytic continuation are along a hypersurface $\Theta$ in $\mathbb C^N\setminus \Delta$ (the \emph{Malgrange's divisor}), along which the RH b.v.p. is not solvable. 
Given a solution $\Phi$ of the above Riemann-Hilbert b.v.p. we can construct \emph{both} the system \eqref{16.07.17-1}-\eqref{16.07.17-2} \emph{and} three of its solutions, by setting
\[Y_{L/R}(z,u):=\Phi_{L/R}(z,u)e^{zU},\quad Y_0(z,u):=\Phi_0(z,u)z^\mu z^R.
\]
Indeed, from the asymptotic expansions 
$$\Phi_{L/R}(z)=\mathbbm 1+\frac{1}{z}F_1+O\left(\frac{1}{z^2}\right),\quad z\to \infty,\quad z\in P_{L/R},
$$
and the convergent Taylor series representation 
$$
\Phi_0(z)=\sum_{p=0}^\infty\phi_pz^p,\quad z\to 0,
$$
we deduce that
\begin{align}
\label{13.11.18-1}\frac{\partial Y_{L/R}}{\partial z}\ Y^{-1}_{L/R}&=U+\frac{1}{z}[F_1,U]+O\left(\frac{1}{z^2}\right),\\
\label{13.11.18-2}\frac{\partial Y_0}{\partial z}\ Y^{-1}_0&=\frac{1}{z}\left(\phi_0\mu\phi_0^{-1}+O\left(z\right)\right),\\
\label{13.11.18-3}\frac{\partial Y_{L/R}}{\partial u_i}\ Y^{-1}_{L/R}&=zE_i+[F_1, E_i]+O\left(\frac{1}{z}\right),\\
\label{13.11.18-4}\frac{\partial Y_0}{\partial u_i}\ Y^{-1}_0&=\frac{\partial \phi_0}{\partial u_i}\phi_0^{-1}+O(z),
\end{align}
and from the independence of the matrices $S,C$ w.r.t. both $z$ and $u$, we conclude that the r.h.s.'s of \eqref{13.11.18-1}-\eqref{13.11.18-2}, and of \eqref{13.11.18-3}-\eqref{13.11.18-4} respectively, are equal. Hence, the functions $Y_{L/R}, Y_0$ are solutions of the system of equations
\begin{align*}
\frac{\partial Y}{\partial u_i}&=(zE_i+V_i)Y,\quad V_i(u):=[F_1(u),E_i]\equiv\frac{\partial \phi_0}{\partial u_i}\phi_0^{-1},\\
\frac{\partial Y}{\partial z}&=\left(U+\frac{1}{z}V(u)\right)Y,\quad V(u):=[F_1(u), U]\equiv \phi_0\mu\phi_0^{-1}.
\end{align*}
Then using the formulae \eqref{12.11.18-1}-\eqref{12.11.18-5} the Frobenius structure can be reconstructed. More precisely, we have the following

\bsh
\begin{teorema}\cite{dubro1, dubro0, dubro2}\label{teodubinv}
If the RH b.v.p. above admits solution for $u^{(0)}$, and if 
\[\prod_{i=1}^N\phi_{i1,0}(u^{(0)})\neq 0,
\]then the formulae \eqref{12.11.18-1}-\eqref{12.11.18-5} define an analytic Frobenius structure on a sufficiently small neighborhood of the point $u^{(0)}$.
\end{teorema}
\esh

\begin{oss}\label{ossinv}
Notice that if we want to interpret the monodromy data $(\mu, R,S,C)$ as a system of coordinates on the space of solutions of WDVV equations, we have to keep track of the freedom and ambiguities up to which they are defined. For example, according to Theorem \ref{06.07.17-4}, the 4-tuples 
\[(\mu, R,S,C),\quad \text{and}\quad (\mu, G^{-1}RG, S, GC),\ G\in\mathcal G(X)
\]allow us to locally reconstruct the same $\ell$-chamber of the quantum cohomology $QH^\bullet(X)$ through the parametric formulae \eqref{12.11.18-1}-\eqref{12.11.18-5}.
\end{oss}

\begin{oss}
In the very interesting paper \cite{sabbah_insp}, C. Sabbah addressed the problem of extending the validity of both Theorem \ref{teomiwmal} and of Theorem \ref{teodubinv} in order to allow $u^{(0)}$ to be a \emph{semisimple coalescing point}, i.e. with
\[u^{(0)}_i=u^{(0)}_j,\quad \text{for some }i\neq j.
\] 
It is proved that if the RH b.v.p. is solvable at $u^{(0)}$, and if the corresponding matrix $$V_0:=\phi_0(u^{(0)})\cdot \mu\cdot \phi_0(u^{(0)})^{-1}$$ satisfies the following conditions:
\begin{enumerate}
\item it is of the form
\[V_0=[F_1,U_0],\quad U_0:={\rm diag}(u^{(0)}_1,\dots, u^{(0)}_N),
\]
\item and it admits an eigenvector with non-zero components w.r.t. the standard basis of $\mathbb C^N$, 
\end{enumerate}
then the RH b.v.p. is solvable in a sufficiently small neighborhood of $u^{(0)}$, analytically continues to a meromorphic function on the universal cover of the complement of a hypersurface $\Theta$ in $\mathbb C^N$. Consequently, in the spirit of \cite{CDG}, we have that the monodromy data define a sistem of local moduli for the Frobenius structure also at semisimple coalescing points.
\end{oss}

The analytic continuation of the Frobenius structure can be obtained by changing the given monodromy data $(\mu, R,S,C)$ to those associated with another chamber, by the action of the braid group, and then solving again the Riemann-Hilbert problem with the new data.

\vskip 0.2 cm 
In \cite{guzzetti2} the above formulae have been used to construct in closed form the solution $F(t)$ of the WDVV equations for $N = 3$ in some relevant cases:
\begin{itemize}
\item
Using  the five algebraic solutions of the Painlev\'e  VI equation, the three polynomial solutions of the WDVV equations are obtained, corresponding to the Frobenius structure on the orbit space of Coxeter groups,   
plus two algebraic solutions.   

\item
 Using the Painlev\'e VI transcendent associated with the monodromy data of the quantum cohomology or $\mathbb{P}_\mathbb{C}^2$, the 
Kontsevich's solution  of
the WDVV equations  is obtained  which generates 
 the numbers $N_k$  of rational curves $
\mathbb{P}_\mathbb{C}^1\longrightarrow \mathbb{P}_\mathbb{C}^2$ of degree $k$ passing through $3k-1$ generic
points. Namely 
$$
F(t^1,t^2,t^3)= 
{1\over 2}\left[ (t^1)^2 t^3+t^1(t^2)^2\right]+ \sum_{k=1}^{\infty} {N_k \over
(3k-1)!} (t^3)^{3k-1} e^{k t^2}. 
$$
This procedure shows that  the $N_k$'s can be computed as 
an application of the isomonodromic deformation approach to Frobenius manifolds. 
\end{itemize}
\newpage
\section{Helix Theory in Triangulated Categories}\label{Helixsec}

In this Section we recall basic facts about exceptional objects and collections, more general semiorthogonal decompositions in triangulated categories as well as about the operations of mutations naturally defined on these objects. Without claim to completeness and originality, our aim is to give a summary, self-contained as much as possible, of the so-called \emph{Helix theory} in triangulated categories as developed in algebraic geometry literature (\cite{goru}, \cite{ru}, \cite{BK}, \cite{helix} and references therein). Particular emphasis will be given to results stressing analogies and similarities with certain aspects of the analytic theory of Frobenius manifolds explained in the previous Sections. In order to reach the widest audience as possible (we think of the mathematicians working on integrable systems, asymptotic analysis, ordinary differential equations in complex domains, isomonodromic deformations theory, etc.), for the convenience of the reader we also have included proofs of the main results that will be used in the subsequent part of the paper. The experienced reader can easily skip several subsections.

\subsection{Prerequisites} We only assume the reader to be comfortable with the notions of derived and triangulated categories. The main reference, to which we will refer for these basics definitions and the labeling of the axioms (e.g. the axioms ${\rm TR}i$ with $i=1,\dots,4$), is \cite{gelman}. Other useful and complete references are \cite{BBR} and \cite{huy}. In what follows we use the notations $\mathscr D$ for a triangulated category, and we denote its shift (or translation) functor by $[1]\colon\mathscr D\to\mathscr D$. As usual, we denote by $[n]\colon\mathscr D\to\mathscr D$ the subsequent $n$ compositions of the shift functor $[1]$.

\subsection{Notations and preliminaries} 
Let $\mathbb K$ be a field\footnote{Here we work on a general ground field $\mathbb K$ but, starting from Section \ref{geocase} we will specialize to the case $\mathbb K=\mathbb C$.}.  We denote by $\text{GrVect}^{<\infty}_{\mathbb K}$ the category of finite dimensional $\mathbb Z$-graded vector spaces: in what follows we will denote the $p$-th degree of $V^\bullet$ by $\text{Gr}^p(V^\bullet)$ or $V^p$.  $\text{GrVect}^{<\infty}_{\mathbb K}$ is a triangulated category, the shift being defined by
\[\text{Gr}^p(V^\bullet[k]):=\text{Gr}^{p+k}(V^\bullet),\quad p,k\in\mathbb Z,
\]
 and we also have operations of \emph{tensor product} and \emph{dualization} with the usual gradations
\[\text{Gr}^p(V^\bullet\otimes W^\bullet):=\bigoplus_{i+j=p}\text{Gr}^i(V^\bullet)\otimes \text{Gr}^j(W^\bullet),\quad\quad \text{Gr}^p\left((V^\bullet)^*\right):=\left(\text{Gr}^{-p}(V^\bullet)\right)^*.
\]
The category $\text{GrVect}^{<\infty}_{\mathbb K}$ is equivalent to the bounded derived category of finite dimensional $\mathbb K$-vector spaces, denoted by $\mathcal D^b(\mathbb K)$: the equivalence is realized by the functors
\[\Phi\colon\text{GrVect}^{<\infty}_{\mathbb K}\to\mathcal D^b(\mathbb K)\colon  V^\bullet\mapsto \bigoplus_{i\in\mathbb Z}(\text{Gr}^iV^\bullet)[-i],\quad\text{with zero differentials,}\]\[ H^\bullet\colon\mathcal D^b(\mathbb K)\to \text{GrVect}^{<\infty}_{\mathbb K}\colon F^\bullet\mapsto H^\bullet(F^\bullet).
\]
Let $\mathscr D$ be a triangulated category. We will assume that $\mathscr D$ is a $\mathbb K$-\emph{linear category of finite type} (or $\Hom$-\emph{finite}), i.e. that
\[\Hom^{\bullet}(X,Y):=\bigoplus_{i\in\mathbb Z}\Hom^i(X,Y) 
\]is a finite dimensional graded $\mathbb K$-vector space for all $X,Y\in\Ob(\mathscr D)$, and where we posed $\Hom^i(X,Y):=\Hom(X,Y[i])\ \text{for any }i\in\mathbb Z$. Sometimes, it will be useful to consider the category $\mathscr D$ to be $\mathscr D^b(\mathbb K)$-enriched, by identifying the graded vector spaces $\Hom^{\bullet}(X,Y)$ with the associated complex through the equivalence $\Phi$ above.

\begin{defi}
Let $V^\bullet$ be a finite dimensional graded $\mathbb K$-vector space and $X$ be an object in a $\mathbb K$-linear triangulated category $\mathscr D$. We define the tensor product $V^\bullet \otimes X$, an object of $\mathscr D$, as a solution of a universal problem, by requiring
\[\Hom^\bullet (Y,V^\bullet\otimes X)=V^\bullet \otimes\Hom^\bullet (Y,X)\quad\forall \,Y\in\Ob(\mathscr D).
\]Such a universal problem admits a solution: the tensor product can be constructed as
\[V^\bullet \otimes X:=\bigoplus_{i}V^i\otimes X[-i],
\]where
\[V^i\otimes X[-i]:=\underbrace{X[-i]\oplus\dots\oplus X[-i]}_{\dim_{\mathbb K}V^i \text{ times}}.
\]
\end{defi}

\begin{oss}
We can define an analogous operation of tensor product $-\otimes-\colon\mathcal D^b(\mathbb K)\times\mathscr D\to\mathscr D$ by composition with the cohomology functor in the first entry:
\[\xymatrix{
\mathcal D^b(\mathbb K)\times\mathscr D\ar[rr]^{H^\bullet\times \mathbbm 1_{\mathscr D}\quad}&&\text{GrVect}^{<\infty}_{\mathbb K}\times\mathscr D \ar[rr]^{\quad\quad\quad -\otimes -}&&\mathscr D.
}
\]In this way, the object $F^\bullet\otimes X$ depends only on the quasi-isomorphism class of $F^\bullet$.
\end{oss}
\begin{lemma}\label{lemalge}If $V^\bullet\in\Ob(\operatorname{GrVect^{<\infty}_{\mathbb K}})$, $X\in\Ob(\mathscr D)$ and if $j,k\in\mathbb Z$, then
\[V^\bullet[j]\otimes X[k]=(V^\bullet\otimes X)[j+k], \quad (V^\bullet [j])^*=(V^\bullet)^*[-j].
\]
\end{lemma}
\proof For the first equality it is easy to see that the r.h.s. solves the universal problem which defines the l.h.s.. The second equality can be trivially deduced by a direct comparison of the gradings. \endproof

\begin{defi}
If $\mathscr D$ and $\mathscr E$ are two $\mathbb K$-linear triangulated categories, a covariant exact functor $F\colon \mathscr D\to\mathscr E$ is called \emph{linear} if
\[F(V^\bullet\otimes X)=V^\bullet\otimes F(X)
\]for any graded vector space $V^\bullet$ and any object $X$. Analogously, a contravariant functor $F\colon\mathscr D^{\text{op}}\to\mathscr E$ is \emph{linear} if it satisfies
\[F(V^\bullet\otimes X)=(V^\bullet)^*\otimes F(X)
\]for any graded vector space $V^\bullet$ and any object $X$.
\end{defi}
So, in particular, the bifunctor $\Hom^\bullet(-,-)\colon\mathscr D\times\mathscr D^{\text{op}}\to\operatorname{GrVect_{\mathbb K}^{<\infty}}$ is bilinear:
\[\Hom^\bullet(W^\bullet\otimes X,V^\bullet \otimes Y)=(W^\bullet)^* \otimes V^\bullet\otimes\Hom^\bullet(X,Y)
\]for any $X,Y\in\Ob(\mathscr D)$ and any graded vector spaces $V^\bullet$ and $W^\bullet$.
Thus, for any $X,Y\in\Ob(\mathscr D)$, we have the identifications
\[\End(\Hom^\bullet(X,Y))=\Hom^\bullet(\Hom^\bullet(X,Y)\otimes X,Y)=\Hom^\bullet(X,\Hom^\bullet(X,Y)^*\otimes Y).
\]Hence, the identity morphism $\operatorname{id}\colon \Hom^\bullet(X,Y)\to \Hom^\bullet(X,Y)$ induces two canonical morphisms
\[j^*(X,Y)\colon\Hom^\bullet(X,Y)\otimes X\to Y,
\]
\[j_*(X,Y)\colon X\to \Hom^\bullet(X,Y)^*\otimes Y.
\]
\begin{prop}\label{25.03.17.1}
Let $E\in\Ob(\mathscr D)$ be a generic object. Let us define the functors
\begin{align*}\Phi_E\colon\mathcal D^b(\mathbb K)\to\mathscr D\colon&\  V^\bullet\mapsto V^\bullet\otimes E,\\
\Phi_E^*\colon\mathscr D\to\mathcal D^b(\mathbb K)\colon&\ X\mapsto\Hom^\bullet(X,E)^*,\\
\Phi_E^!\colon\mathscr D\to\mathcal D^b(\mathbb K)\colon&\ X\mapsto\Hom^\bullet(E,X).
\end{align*}
We have the adjunctions $\Phi_E^*\dashv \Phi_E\dashv \Phi^!_E$.
\end{prop}
\proof
This is a simple check of the definition of adjoint functors. Notice that the unity of the adjunction $\Phi^*_E\dashv\Phi_E$ and counity of the adjunction $\Phi_E\dashv \Phi_E^!$ are given by the morphisms $j_*(-,E)$ and $j^*(E,-)$ respectively.
\endproof

\begin{defi}[Generated triangulated subcategory]If $\Omega\subseteq\Ob(\mathscr D)$, we denote by $\langle\Omega\rangle$ the smallest full triangulated subcategory of $\mathscr D$ containing all objects of $\Omega$.\end{defi}

\begin{defi}
If $\mathcal A,\mathcal B\subseteq\Ob(\mathscr D)$ we define the set
\[\mathcal A *\mathcal B:=\left\{X\in\Ob(\mathscr D)\colon A\to X\to B\to A[1],\ \text{for some }A\in\mathcal A,\ B\in\mathcal B\right\}.
\]
Notice by the octahedral axiom (TR4) that the operation $*$ is associative.
\end{defi}

The subcategory $\langle\Omega\rangle$ is obtained by taking the closure with respect to shifts and cones. More precisely, we have the following
\begin{prop}\label{24-03-17.1}Let $\Omega\subseteq\Ob(\mathscr D)$, and let us define
\[\Omega_1:=\left\{X[n]\colon X\in\Omega, n\in\mathbb Z\right\}, \quad \Omega_r:=\underbrace{\Omega_1*\dots*\Omega_1}_{r\ times}.
\]Then
\[\langle\Omega\rangle\equiv\bigcup_{r\in\mathbb N^*}\Omega_r.
\]
\end{prop}

\subsection{Exceptional Objects and Mutations}
Let $\mathscr D$ be a $\mathbb K$-linear triangulated category.
\begin{defi}[Exceptional Object, Pair and Collection]An object $E\in\Ob(\mathscr D)$ is called \emph{exceptional} if $\Hom^\bullet(E,E)$ is a 1-dimensional $\mathbb K$-algebra generated by the identity morphism.\\
An ordered pair $(E_1,E_2)$ of exceptional objects of $\mathscr D$ is called \emph{exceptional} or \emph{semiorthogonal} if
\[\Hom^\bullet(E_2,E_1)=0.
\]More in general, an ordered collection $(E_1,E_2,\dots, E_k)$ of exceptional objects of $\mathscr D$ is called \emph{exceptional} or \emph{semiorthogonal} if
\[\Hom^\bullet(E_j,E_i)=0\quad\text{whenever }i<j.
\] 
An exceptional collection is said to be \emph{full} if it generates $\mathscr D$, i.e. any full triangulated subcategory containing all objects $E_i$ is equivalent to $\mathscr D$ via the inclusion functor.
\end{defi}

\begin{prop}[\cite{Bondal}]\label{propbondal}
Let $E\in\Ob(\mathscr D)$ be a generic object. Then $E$ is exceptional if and only if the functor
\begin{align*}\Phi_E\colon\mathcal D^b(\mathbb K)\to\mathscr D\colon&\  V^\bullet\mapsto V^\bullet\otimes E\end{align*}
is fully faithful. In particular, the category $\langle E\rangle\equiv \operatorname{Im}\Phi_E$ is equivalent to the category $\mathcal D^b(\mathbb K)$.
\end{prop}
\proof
Using the notations of Proposition \ref{25.03.17.1}, $\Phi_E$ is fully faithful if and only if the natural transformation $\Phi_E^!\Phi_E\Leftarrow \mathbbm 1_{\mathcal D^b(\mathbb K)}$ is a natural isomorphism. This holds if and only if $\Hom^\bullet(E,E)\cong\mathbb K$.
\endproof

\begin{oss}Given an exceptional collection in $\mathscr D$, there are several operations generating other such collections. Indeed, the group $\operatorname{Aut}(\mathscr D)$ of isomorphism classes of auto-equivalences of the category $\mathscr D$ acts on the set of exceptional collections: the element $\Psi\in\operatorname{Aut}(\mathscr D)$ acts in the obvious way, by associating with the exceptional collection $\frak E:=(E_1,\dots, E_n)$ the collection $\Psi\frak E:=(\Psi E_1,\dots,\Psi E_n)$. Analogously, the additive group $\mathbb Z^n$ acts on the sets of exceptional collection of length $n$ by shifts: if $\frak E:=(E_1,\dots,E_n)$ is an exceptional collection, then also $\frak E[{\bold k}]:=(E_1[k_1],E_2[k_2],\dots, E_n[k_n])$ is exceptional for any $(k_1,\dots,k_n)\in\mathbb Z^n$. The actions of both $\operatorname{Aut}(\mathscr D)$ and $\mathbb Z^n$ preserve the fullness of an exceptional collection. 
\end{oss}

In what follows, we are going to define a nontrivial action of the \emph{braid group} $\mathcal B_{n}$ on the set of (full) exceptional collections of length $n$.

\begin{defi}[Orthogonal complements]Let $\mathscr A$ be a full triangulated subcategory of $\mathscr D$. We introduce two full triangulated subcategories 
$^\perp \mathscr A$ and $\mathscr A^\perp$ defined by
\[^\perp \mathscr A:=\left\{X\in\Ob(\mathscr D)\colon\Hom(X,A)=0\text{ for all }A\in\Ob(\mathscr A)\right\},
\]
\[ \mathscr A^\perp:=\left\{X\in\Ob(\mathscr D)\colon\Hom(A,X)=0\text{ for all }A\in\Ob(\mathscr A)\right\}.
\]These subcategories are called respectively \emph{left} and \emph{right orthogonals} to $\mathscr A$ in $\mathscr D$. \end{defi}

\begin{oss}It is easy to see that, if $\mathscr A=\langle E\rangle$ is the smallest triangulated subcategory containing an object $E\in \Ob(\mathscr D)$, the following characterization of the orthogonal complements $^\perp \langle E\rangle$ and $\langle E\rangle^\perp$ holds:
\[^\perp\langle E\rangle\equiv\ ^\perp E,\quad \text{where }^\perp E:=\left\{X\in\Ob(\mathscr D)\colon\Hom^\bullet(X,E)=0\right\},
\]
\[ \langle E\rangle^\perp\equiv E^\perp,\quad \text{where }E^\perp:=\left\{X\in\Ob(\mathscr D)\colon\Hom^\bullet(E,X)=0\right\}.
\]
\end{oss}
\begin{defi}[Mutations of objects]
\label{10novembre2018-1}
Let $E\in\Ob(\mathscr D)$ be an exceptional object. 
For any $X\in\Ob(\mathscr D)$ we can define two new objects 
\[
\mathbb L_EX\in\Ob(E^\perp),\quad \mathbb R_EX\in\Ob(^\perp E)
\]
called respectively \emph{left} and \emph{right mutations} of $X$ with respect to $E$. These two objects are defined as the cones
\[\mathbb L_E(X):=\operatorname{Cone}\left(\Phi_E\Phi_E^!(X)\to X\right),\quad \mathbb R_E(X):=\operatorname{Cone}\left(X\to\Phi_E\Phi_E^*(X)\right)[-1],
\]where the functors $\Phi_E,\Phi^*_E,\Phi^!_E$ are the ones introduced in Proposition \ref{25.03.17.1}. We thus have the distinguished triangles
\begin{equation}
\label{triangleleft}
\xymatrix{
\mathbb L_EX[-1]\ar[r]&\Hom^\bullet(E,X)\otimes E\ar[r]^{\quad\quad\quad\quad j^*}&X\ar[r]&\mathbb L_EX}
\end{equation}
\begin{equation}
\label{triangleright}
\xymatrix{
\mathbb R_EX\ar[r]&X\ar[r]^{ j_*\quad\quad\quad\quad}&\Hom^\bullet(X,E)^*\otimes E\ar[r]&\mathbb R_EX[1]}
\end{equation}
extending the canonical morphisms $j^*(E,X)$ and $j_*(X,E)$.
\end{defi} 
By applying the functor $\Hom^\bullet(E,-)$ to \eqref{triangleleft}, and the functor $\Hom^\bullet(-,E)$ to \eqref{triangleright}, and using the fact that $E$ is exceptional, we obtain the \emph{orthogonality relations}
\begin{equation}\label{orthleftright}
\Hom^\bullet(E,\mathbb L_EX)=0,\quad\Hom^\bullet(\mathbb R_EX,E)=0.
\end{equation}

\begin{oss}
Our definitions of $\mathbb L_EX$ and $\mathbb R_EX$ differ from the original ones given in \cite{helix} by a shift operator. In particular, 
\begin{itemize}
\item what here we denote $\mathbb L_EX$, in \cite{helix} is $\mathbb L_EX[1]$,
\item and our $\mathbb R_EX$ in \cite{helix} is $\mathbb R_EX[-1]$.
\end{itemize}
The reason of such a choice will be explained later. In what follows we will reformulate some results of \cite{helix}, adapted to our definition, and we leave to the reader the small modifications of some proofs.
\end{oss}

In general, the third term in a distinguished triangle is not \emph{canonically} defined by the other two terms. In this case, however, the objects $\mathbb L_EX$ and $\mathbb R_EX$ are \emph{unique up to unique} isomorphism because of the orthogonality relations. Indeed, let us assume that we have the two following distinguished triangles\footnote{Here we use the axiom TR2 of triangulated category.}:
\[\xymatrix{
\Hom^\bullet(E,X)\otimes E\ar[r]^{ \quad\quad\quad\quad j^*}\ar@{=}[d]&X\ar[r]\ar@{=}[d]&S\ar[r]\ar@{.>}[d]_{h}&\left(\Hom^\bullet(E,X)\otimes E\right)[1]\ar@{=}[d]\\
\Hom^\bullet(E,X)\otimes E\ar[r]^{ \quad\quad\quad\quad j^*}&X\ar[r]^{\alpha\quad}&\mathbb L_EX\ar[r]&\left(\Hom^\bullet(E,X)\otimes E\right)[1]
}
\]Then, by axiom TR3 of triangulated category there exists a morphism $h\colon S\to \mathbb L_EX$, which necessarily is an isomorphism (the other two vertical maps being the identities). We want to show that $h$ is unique. Let us apply the functor $\Hom^\bullet(-,\mathbb L_EX)$ to the first line of the diagram: recalling that the shift functor $T=(-)[1]$ is an auto-equivalence, and using the orthogonality relations \eqref{orthleftright} we obtain
\begin{align*}
\Hom^\bullet\left(\left(\Hom^\bullet(E,X)\otimes E\right)[1],\mathbb L_EX\right)&=\Hom^\bullet\left(\Hom^\bullet(E,X)\otimes E,\mathbb L_EX[-1]\right)\\&=\Hom^\bullet(E,X)^*\otimes\Hom^\bullet(E,\mathbb L_EX[-1])\\
&=0.
\end{align*}
So we get the long exact sequence
\[\xymatrix{
\dots\ar[r]&0\ar[r]&\Hom^\bullet(S,\mathbb L_EX)\ar[r]&\Hom^\bullet(X,\mathbb L_EX)\ar[r]&\dots
}
\]Since $\alpha\circ j^*=0$, there exists $h\in\Hom^\bullet(S,\mathbb L_EX)$ as above, which must be unique by injectivity. Analogously one shows that $\mathbb R_EX$ is unique up to a unique isomorphism. Notice, in particular, that $\mathbb L_EE=\mathbb R_EE=0$.

Moreover, as a consequence of axiom TR1, we necessarily have that
\begin{equation}\label{leftproj}\mathbb L_EX=X\quad\text{for all }X\in E^\perp,
\end{equation}
\begin{equation}\label{rightproj}\mathbb R_EX=X\quad\text{for all }X\in{^\perp E},
\end{equation}
which means that the operations $\mathbb L_E,\mathbb R_E$ are \emph{projections} onto the subcategories $E^\perp$ and ${^\perp E}$. Some other useful properties of these projections are summarized in the following

\bsh
\begin{prop}[\cite{helix}]\label{RLproperties}Let $\mathscr D$ be a $\mathbb K$-linear triangulated category, and $E$ an exceptional object.
\begin{enumerate}
\item For any object $X\in\Ob (\mathscr D)$ and any pair of integer $k,\ell \in\mathbb Z$ we have
\[\mathbb L_{E[k]}\left(X[\ell]\right)=\left(\mathbb L_EX\right)[\ell],\quad \mathbb R_{E[k]}\left(X[\ell]\right)=\left(\mathbb R_EX\right)[\ell].
\]
\item
If $E'\in\Ob(^\perp E), E''\in\Ob(E^\perp)$ and $X\in\Ob(\mathscr D)$, the following bifunctor isomorphisms hold
\begin{equation}
\label{functiso1}\Hom^\bullet(E',X)=\Hom^\bullet(E',\mathbb R_EX)=\Hom^\bullet(E',\mathbb L_EX)=\Hom(\mathbb L_EE',\mathbb L_EX),
\end{equation}
\begin{equation}
\label{functiso2}\Hom^\bullet(X,E'')=\Hom^\bullet(\mathbb L_EX,E'')=\Hom^\bullet(\mathbb R_EX,E'')=\Hom(\mathbb R_EX,\mathbb R_EE'').
\end{equation}
\item The functors
\[
\xymatrix{
\mathscr D\ar[rr]^{X\mapsto \mathbb R_EX}&&^\perp E\\
\mathscr D\ar[rr]^{X\mapsto \mathbb L_EX}&& E^\perp
}
\]are respectively the right adjoint functor to the inclusion $\xymatrix{^\perp E\ \ar@{^{(}->}[r]&\mathscr D}$, and the left adjoint functor to the inclusion $\xymatrix{ E^\perp\ \ar@{^{(}->}[r]&\mathscr D}$.
\item The following identities hold
\[\mathbb  L_E\circ \mathbb R_E=\mathbb L_E,\quad \mathbb R_E\circ \mathbb L_E=\mathbb R_E.
\]
\item The restrictions
\[\mathbb L_E|_{^\perp E}\colon^\perp E\to E^\perp\quad\text{and}\quad \mathbb R_E|_{E^\perp}\colon E^\perp\to {^\perp E}
\]are functors inverse to each other, establishing an isomorphism between these two subcategories.
\item The following isomorphism holds functorially on $X$ and $E$
\[\Hom^\bullet(\mathbb L_EX[-1],E)^*=\Hom^\bullet(E,\mathbb R_EX).
\]
\end{enumerate}
\end{prop}
\esh

\proof
Let us prove (1). Applying Lemma \ref{lemalge} we have that
\[\Hom^\bullet (E[k],X[\ell])\otimes E[k]=\left(\Hom^\bullet (E,X)\otimes E\right)[\ell],
\]and that
\[
\Hom^\bullet(X[\ell],E[k])^*\otimes E[k]=\left(\Hom^\bullet(X,E)^*\otimes E\right)[\ell].
\]
Moreover, it is easily seen that the following diagrams are commutative:
\[
\xymatrix{X[\ell]\ar[rr]^{j_*(X[\ell],E[k])\quad\quad\quad\quad\quad}\ar[rrd]_{j_*(X,E)[\ell]\quad}&&\Hom^\bullet(X[\ell],E[k])^*\otimes E[k]\ar@{=}[d]\\
&&\left(\Hom^\bullet(X,E)^*\otimes E\right)[\ell]}\quad\quad
\xymatrix{\Hom^\bullet(E[k],X[\ell])\otimes E[k]\ar[rr]^{\quad\quad\quad\quad j^*(E[k],X[\ell])}\ar@{=}[d]&& X[\ell]\\
(\Hom(E,X)\otimes X)[\ell]\ar[rru]_{\quad j^*(E,X)[\ell]}&&}
\]

Thus, applying the shift functor $(-)[\ell]$ to \eqref{triangleleft}, \eqref{triangleright} we obtain
\[\xymatrix{
\mathbb L_EX[\ell-1]\ar[r]&(\Hom^\bullet(E,X)\otimes E)[\ell]\ar[rr]^{\quad\quad\quad\quad (-1)^\ell j^*[\ell]}&&X[\ell]\ar[r]&\mathbb L_EX[\ell],}
\]
\[\xymatrix{
\mathbb R_EX[\ell]\ar[r]&X[\ell]\ar[rr]^{(-1)^\ell j_*[\ell]\quad\quad\quad\quad}&&(\Hom^\bullet(X,E)^*\otimes E)[\ell]\ar[r]&\mathbb R_EX[\ell+1].}
\]
Recalling now that, by the axiom TR1 of triangulated category, one can change the sign of any two morphisms in a distinguished triangle, we conclude. For points (2),(3),(4),(5),(6) we refer to \cite{helix}, where the reader can find and easily adapt the proofs by keeping track of the difference of shiftings in the definition of left and right mutations.
\endproof

\begin{defi}
If $(E_1,E_2)$ is an exceptional pair, we define its \emph{left} and \emph{right mutations} to be the pairs
\[\mathbb L(E_1,E_2):=(\mathbb L_{E_1}E_2,E_1)\quad\text{and}\quad \mathbb R(E_1,E_2):=(E_2,\mathbb R_{E_2}E_1)
\]respectively. \end{defi}

\bsh
\begin{prop}The pairs $\mathbb L(E_1,E_2), \mathbb R(E_1,E_2)$ are exceptional. Moreover, $\mathbb L,\mathbb R$ act on the set of exceptional pairs as inverse transformations:
\[\mathbb L\circ \mathbb R(E_1,E_2)=(E_1,E_2)\quad\text{and}\quad \mathbb R\circ \mathbb L(E_1,E_2)=(E_1,E_2).
\]
\end{prop}
\esh
\proof[Proof]
The first statement follows from the orthogonality relations \eqref{orthleftright}. From relations \eqref{leftproj},\eqref{rightproj} we have that
\[\mathbb R_{E_1}E_2=E_2\quad\text{and}\quad \mathbb L_{E_2}E_1=E_1,
\]so that, by Proposition \eqref{RLproperties}, we deduce
\[\Hom^\bullet(\mathbb L_{E_1}E_2[-1],E_1)^*=\Hom^\bullet(E_1,\mathbb R_{E_1}E_2)=\Hom^\bullet(E_1,E_2),
\]
\[\Hom^\bullet(E_2,\mathbb R_{E_2}E_1[1])=\Hom^\bullet(\mathbb L_{E_2}E_1,E_1)^*=\Hom^\bullet(E_1,E_2)^*.
\]It follows that the triangles
\[
\xymatrix{
\mathbb L_{E_1}E_2[-1]\ar[r]&\Hom^\bullet(E_1,E_2)\otimes E_1\ar[r]&E_2\ar[r]&\mathbb L_{E_1}E_2}\]
\[\xymatrix{
\mathbb L_{E_1}E_2[-1]\ar[r]&\Hom^\bullet(\mathbb L_{E_1}E_2[-1],E_1)^*\otimes E_1\ar[r]&\mathbb R_{E_1}\mathbb L_{E_1}E_2\ar[r]&\mathbb L_{E_1}E_2
}
\]can be canonically identified, i.e. $\mathbb R\circ \mathbb L(E_1,E_2)=(E_1,E_2).$ Similarly the other identity follows.
\endproof

For a more general exceptional sequence, we give the following

\begin{defi}
Let $(E_0,\dots,E_k)$ be an exceptional collection in $\mathscr D$. For $1\leq i\leq k$ we define
\[\mathbb L_i(E_0,\dots,E_k):=(E_0,\dots, \mathbb L_{E_{i-1}}E_i,E_{i-1},\dots,E_k),
\]
\[\mathbb R_i(E_0,\dots,E_k):=(E_0,\dots, E_i,\mathbb R_{E_i}E_{i-1},\dots,E_k).
\]
\end{defi}

\bsh
\begin{prop}[\cite{BK}]
The mutations preserve the exceptionality and satisfy
\begin{equation}\label{invbraid}
\mathbb L_i\mathbb R_i=\mathbb R_i\mathbb L_i=\operatorname{Id}
\end{equation}
\begin{equation}\label{braid1}
\mathbb L_i\mathbb L_j=\mathbb L_j\mathbb L_i\quad\text{for}\quad |i-j|>1
\end{equation}
\begin{equation}\label{braid2}
\mathbb L_{i+1}\mathbb L_i\mathbb L_{i+1}=\mathbb L_i\mathbb L_{i+1}\mathbb L_i\quad\text{for}\quad1<i<k.
\end{equation}
So the braid group $\mathcal B_{k+1}$ acts by the mutations on exceptional objects of length $(k+1)$.
\end{prop}
\esh

\proof[Proof]
The fact that exceptionality is preserved, and the first two relations \eqref{invbraid},\eqref{braid1} follow from the previous results. The only non-trivial relation is \eqref{braid2}. Let $(A,B,C)$ be an exceptional triple. We have to show the commutativity of the diagram
\[
\xymatrix{
&(A,B,C)\ar[dl]_{\mathbb L_1}\ar[dr]^{\mathbb L_2}&\\
(\mathbb L_AB,A,C)\ar[d]_{\mathbb L_2}&&(A,\mathbb L_BC,B)\ar[d]^{\mathbb L_1}\\
(\mathbb L_AB,\mathbb L_AC,A)\ar[dr]_{\mathbb L_1}&&(\mathbb L_A\mathbb L_BC,A,B)\ar[dl]^{\mathbb L_2}\\
&(\mathbb L_{\mathbb L_AB}\mathbb L_AC,\mathbb L_AB,A)&
}
\]So we have to prove that
\[\mathbb L_{\mathbb L_AB}\mathbb L_AC=\mathbb L_A\mathbb L_BC.
\]Applying the exact linear functor $\mathbb L_A|_{^\perp A}$ to the canonical triangle
\[\mathbb L_BC[-1]\to\Hom^\bullet(B,C)\otimes B\to C\to \mathbb L_BC,
\]and recalling that $\Hom^\bullet(B,C)=\Hom^\bullet(\mathbb L_AB,\mathbb L_AC)$ by Proposition \eqref{RLproperties} we find the triangle
\[\mathbb L_A\mathbb L_BC[-1]\to\Hom^\bullet(\mathbb L_AB,\mathbb L_AC)\otimes \mathbb L_AB\to \mathbb L_AC\to \mathbb L_A\mathbb L_BC.
\]
\endproof

\begin{oss}\label{23-04-17.1}
In the previous exposition, we have followed the main references on the subject and we have defined a \emph{left} action of the braid group on the set of eceptional collections in a $\mathbb K$-linear triangulated category $\mathscr D$. In what follows, in order to establish a perfect correspondence between Helix theory and the theory of local monodromy invariants for quantum cohomologies of Fano manifolds, it will be convenient to consider the braid group $\mathcal B_{n+1}$ as acting on the \emph{right} on the set of exceptional collections of length $n+1$: if we denote by $\beta_{i,i+1}$ with $1\leq i\leq n$ the generators of the braid group, satisfying the relations 
\[\beta_{i,i+1}\beta_{j,j+1}=\beta_{j,j+1}\beta_{i,i+1},\quad |i-j|>1,
\]
\[\beta_{i,i+1}\beta_{i+1,i+2}\beta_{i,i+1}=\beta_{i+1,i+2}\beta_{i,i+1}\beta_{i+1,i+2},
\]and if $\frak E=(E_0,\dots, E_{n})$ is an exceptional collection, we will define 
\[\frak E^{\beta_{i,i+1}}:=\mathbb L_i\frak E.
\] 
\end{oss}

\subsection{Semiorthogonal decompositions,  admissible subcategories, and mutations functors}
Let $\mathscr D$ be a $\mathbb K$-linear triangulated category. In this section we introduce some definitions generalizing the ones of (full) exceptional collection and of left/right mutations with respect to them.

\begin{defi}[\cite{BK, BO, BO2}] A sequence $\mathscr A_1,\dots,\mathscr A_n$ of full triangulated subcategories of $\mathscr D$ is said to be \emph{semiorthogonal} if
\[\mathscr A_i\subseteq\mathscr A_j^\perp\quad \text{for all }i<j.
\]A semiorthogonal sequence $\mathscr A_1,\dots,\mathscr A_n$ is said to define a \emph{semiorthogonal decomposition} of $\mathscr D$ if one of the following equivalent conditions holds:
\begin{enumerate}
\item $\mathscr D$ is generated by the $\mathscr A_i$,  i.e. $\mathscr D=\langle \mathscr A_i\rangle_{i=1}^n$;
\item for any $X\in \Ob(\mathscr D)$ there exists a chain of morphisms
\[
\xymatrix@C-20pt{0\ar@{=}[r]&X_n\ar[rr]&& X_{n-1}\ar[rr]\ar[dl]&& \dots\ar[rr]&&X_2\ar[rr]&&X_1\ar[rr]\ar[dl]&&X_0\ar@{=}[r]\ar[dl]&X\\
&&A_n\ar@{.>}[ul]& && && &A_2\ar@{.>}[ul]& & A_1\ar@{.>}[ul]&&}
\]with $A_i\in\Ob(\mathscr A_i)$.
\end{enumerate}
\end{defi}

The equivalence of $(1),(2)$ immediately follows from Proposition \ref{24-03-17.1}. The chain of morphisms of point (2) is usually called a \emph{filtration} of the object $X$.
\begin{defi}[Filtrations and Postnikov systems]\label{postsyst}Given an object $X\in\Ob(\mathscr D)$, we call \emph{filtration} of $X$, the datum of a set of objects $\left\{X_i\right\}_{i=1}^m$, and a chain of morphisms
\[0=X_m\to X_{m-1}\to\dots\to X_2\to X_1\to X_0=X.
\]These morphisms induce on the cones a family of arrows in the opposite direction, which fit into the diagram
\[
\xymatrix@C-20pt{0\ar@{=}[r]&X_n\ar[rr]&& X_{n-1}\ar[rr]\ar[dl]&& \dots\ar[rr]\ar[dl]&&X_2\ar[rr]&&X_1\ar[rr]\ar[dl]&&X_0\ar@{=}[r]\ar[dl]&X\\
&&L_n\ar@{.>}[ul]& &L_{n-1}\ar@{.>}[ul]\ar@{.>}[ll]^{\delta_{n-1}}& &\ar@{.>}[ll]\dots&& L_2\ar@{.>}[ul]\ar@{.>}[ll]^{\delta_{2}}& & L_1\ar@{.>}[ul]\ar@{.>}[ll]^{\delta_{1}}&&}
\]
where $\delta_{i+1}\circ \delta_i=0$, the bottom triangles are commutative, and the top triangles are distinguished (dashed arrows have degree 1). This diagram is called \emph{(right) Postnikov system} (or \emph{(right) Postnikov tower}), and the object $X$ is called the \emph{canonical convolution} of the Postnikov system.
\end{defi}

\bsh
\begin{prop}\label{30.03.17-1}
The Postnikov system induced by a semiorthogonal decomposition is functorial, i.e. given $X,X'\in\Ob(\mathscr D)$ and a morphism $f\colon X\to X'$ there exists a \emph{unique} prolongation to their Postnikov systems:
\[
\xymatrix@R-20pt{
&A_{n}\ar@{.>}[dl]\ar[dd] | \hole&&&&A_2\ar@{.>}[ld]\ar[dd] | \hole&&A_1\ar@{.>}[ld]\ar[dd] | \hole&\\
0\ar[rr]\ar[dd]&&X_{n-1}\ar[lu]\ar[dd]&\dots&X_2\ar[rr]\ar[dd]&&X_1\ar[rr]\ar[dd]\ar[ul]&&X\ar[dd]^{f}\ar[lu]\\
&A'_n\ar@{.>}[dl]&&&&A_2'\ar@{.>}[ld]&&A_1'\ar@{.>}[ld]&\\
0\ar[rr]&&X'_{n-1}\ar[lu]&\dots&X_2'\ar[rr]&&X_1'\ar[rr]\ar[ul]&&X'\ar[lu]
}
\]In particular, the Postnikov system of an object $X$ is unique up to a unique isomorphism.
\end{prop} 
\esh

\proof
Since $X_1\in\langle\mathscr A_2,\dots,\mathscr A_n\rangle$, we have that $\Hom^\bullet(X_1,A_1')=0$: thus we have the isomorphisms $\Hom(X_1,X_1')\cong\Hom(X_1,X')$ and $\Hom(A_1,A_1')\cong\Hom(X,A_1')$. Consequently there exists a unique morphism of distinguished triangles from $X_1\to X\to A_1\to X_1[1]$ to $X_1'\to X'\to A_1'\to X_1'[1]$ which fits into the diagram. By induction one concludes.\endproof

\begin{oss}\label{spectral}Let $\Phi\colon\mathscr D\to\mathscr A$ be a covariant cohomological functor with values in an abelian category $\mathscr A$, and set $\Phi^q(X):=\Phi(X[q])$ for any object $X\in\Ob(\mathscr D)$, $q\in\mathbb Z$. Given an object $X\in\Ob (\mathscr D)$, there exists a spectral sequence converging to $\Phi^\bullet(X)$. Let us realize $X$ as the canonical convolution of a Postnikov system, as in Definition \ref{postsyst}, and delete the $X_0=X$ term:
\[
\xymatrix@C-20pt{0\ar@{=}[r]&X_n\ar[rr]&& X_{n-1}\ar[rr]\ar[dl]&& \dots\ar[rr]\ar[dl]&&X_2\ar[rr]&&X_1\ar[dl]&\\
&&L_n\ar@{.>}[ul]& &L_{n-1}\ar@{.>}[ul]\ar@{.>}[ll]^{\delta_{n-1}}& &\ar@{.>}[ll]\dots&& L_2\ar@{.>}[ul]\ar@{.>}[ll]^{\delta_{2}}& & L_1\ar@{.>}[ul]\ar@{.>}[ll]^{\delta_{1}}&&}
\]By applying the functor $\Phi$ to this diagram we obtain a bigraded exact couple $(D,E,i,j,k)$ 
\[\xymatrix{D_1^{\bullet,\bullet}\ar[rr]^{i}&&D_1^{\bullet,\bullet}\ar[dl]^{j}\\
&E_1^{\bullet,\bullet}\ar[ul]^{k}&}
\]where $E_1^{p,q}:=\Phi^q(L_{p+1})$, $D_1^{p,q}:=\Phi^q(X_p)$, and the morphism $i,j,k$ have degree $(-1,1), (0,0), (1,0)$ respectively. We thus obtain a spectral sequence $(E^{p,q}_1,d_1:=kj)$ which can be shown to converge to $\Phi^{p+q}(X)$. For further details see \cite{gelman}, Ex. III.7.3c and Ex. IV.2.2a.
\end{oss}

Let us now introduce the strictly related notion of admissibility of a subcategory.
\begin{defi}[\cite{Bondal, BK}]A full triangulated subcategory $\mathscr A$ of $\mathscr D$ is called 
\begin{itemize}
\item \emph{left admissible} if the inclusion functor $i\colon\mathscr A\to\mathscr D$ admits a left adjoint functor \\$ i^*\colon\mathscr D\to\mathscr A$;
\item \emph{right admissible} if the inclusion functor $i\colon\mathscr A\to\mathscr D$ admits a right adjoint functor $i^!\colon\mathscr D\to\mathscr A$;
\item \emph{admissible} if it is both left and right admissible.
\end{itemize} 
\end{defi}

\begin{lemma}[\cite{Bondal}]\label{2-4-17.1}Let $\mathscr A,\mathscr B$ be two full triangulated subcategories of $\mathscr D$.
\begin{enumerate}
\item Let $\mathscr D=\langle \mathscr A,\mathscr B\rangle$ be a semiorthogonal decomposition of $\mathscr D$. Then $\mathscr A$ is left admissible and $\mathscr B$ is right admissible.
\item Conversely, if $\mathscr A$ is left admissible and $\mathscr B$ is right admissible, then $\langle\mathscr A,\ ^\perp\mathscr A\rangle$ and $\langle\mathscr B^\perp, \mathscr B\rangle$ are semiorthogonal decompositions of $\mathscr D$.
\end{enumerate}
\end{lemma}

\proof For any object $X\in\Ob(\mathscr D)$ there exists a distinguished triangle
\[B\to X\to A\to B[1],
\]with $A\in\Ob(\mathscr A)$ and $B\in\Ob(\mathscr B)$. By Proposition \ref{30.03.17-1}, such a distinguished triangle is unique up to unique isomorphism. So, for point (1), the associations
\[i_{\mathscr A}^*(X):=A,\quad i_{\mathscr B}^!(X)=B,
\]are well defined and are respectively left/right adjoint functors to the inlcusions $i_{\mathscr A}, i_{\mathscr B}$. For point (2), given any object $X\in\Ob(\mathscr D)$, by the properties of adjoint functors we have two morphisms
\[X\to i_{\mathscr A}i_{\mathscr A}^*(X),\quad i_{\mathscr B}i_{\mathscr B}^!(X)\to X.
\]By completing them to a distinguished triangle, and using the semiorthogonality condition, it is easily seen that the completing objects are respectively in $^\perp\mathscr A$ and $\mathscr B^\perp$.
\endproof

\begin{cor}[\cite{Bondal}]
If $\mathscr A_1,\dots,\mathscr A_n$ is a semiorthogonal sequence of full triangulated subcategories of $\mathscr D$ such that
\begin{itemize}
\item $\mathscr A_1,\dots,\mathscr A_k$ are left admissible,
\item $\mathscr A_{k+1},\dots,\mathscr A_n$ are right admissible,
\end{itemize}
then
\[\langle\mathscr A_1,\dots,\mathscr A_k,\ ^\perp\langle\mathscr A_1,\dots,\mathscr A_k\rangle\cap \langle\mathscr A_{k+1},\dots,\mathscr A_n\rangle^\perp,\mathscr A_{k+1},\dots,\mathscr A_n\rangle
\]is a semiorthogonal decompositon.
\end{cor}

\begin{cor} If $(\mathscr A_1,\dots,\mathscr A_n)$ is a semiorthogonal sequence of full admissible triangulated subcategories of $\mathscr D$, then the following are equivalent
\begin{enumerate}
\item $\mathscr D=\langle\mathscr A_1,\dots,\mathscr A_n\rangle$ is a semiorthogonal decomposition,
\item $\bigcap_{j=1}^n\mathscr A_j^\perp=0$,
\item $\bigcap_{j=1}^n\ ^\perp\mathscr A_j=0$.
\end{enumerate} 
\end{cor}

\begin{defi}[Mutations functors]
Let $\mathscr A$ be a full triangulated subcategory of $\mathscr D$. 
\begin{itemize}
\item Let us assume that $\mathscr A$ is \emph{left admissible}. Then, we define a functor $\mathbb R_{\mathscr A}\colon\mathscr D\to\mathscr D$, called \emph{right mutation functor with respect to }$\mathscr A$ as follows: for any $X\in\Ob(\mathscr D)$ we define
\[\mathbb R_{\mathscr A}(X):=\operatorname{Cone}\left(X\to ii^*(X)\right)[-1],
\]where $i\colon\mathscr A\to\mathscr D$ is the inclusion functor and $i^*$ is its left adjoint.
\item Let us assume that $\mathscr A$ is \emph{right admissible}. Then, we define a functor $\mathbb L_{\mathscr A}\colon\mathscr D\to\mathscr D$, called \emph{left mutation functor with respect to }$\mathscr A$ as follows: for any $X\in\Ob(\mathscr D)$ we define
\[\mathbb L_{\mathscr A}(X):=\operatorname{Cone}\left(ii^!(X)\to X\right),
\]where $i\colon\mathscr A\to\mathscr D$ is the inclusion functor and $i^!$ is its right adjoint.
\end{itemize}
\end{defi}

\bsh
\begin{prop}
Let $(E_1,\dots, E_k)$ be an exceptional collection in $\mathscr D$. Then the subcategory $\langle E_1,\dots, E_k\rangle$ is admissible, and moreover
\[\mathbb R_{\langle E_1,\dots, E_k\rangle}=\mathbb R_{E_k}\circ\mathbb R_{E_{k-1}}\circ\dots\circ\mathbb R_{E_1},
\]
\[\mathbb L_{\langle E_1,\dots, E_k\rangle}=\mathbb L_{E_1}\circ\mathbb L_{E_2}\circ\dots\circ\mathbb L_{E_k}.
\]In particular, the r.h.s. depend only on $\langle E_1,\dots, E_k\rangle$, and not on the exceptional collection $(E_1,\dots, E_k)$.
\end{prop}
\esh
\proof
Let us proceed by induction on the length $k$ of the exceptional collection. If $k=1$, then the statement is obvious for the results of the previous Section. Let us assume that it is true for all exceptional collections of length $k-1$. Then, we have two distinguished triangles
\[\mathbb R_{\langle E_1,\dots E_{k-1}\rangle}X\to X\to F,\quad \text{with }F\in\Ob \langle E_{1},\dots, E_{k-1}\rangle,
\]
\[\mathbb R_{E_k}\mathbb R_{\langle E_1,\dots E_{k-1}\rangle}X\to \mathbb R_{\langle E_1,\dots E_{k-1}\rangle}X\to F',\quad\text{with }F'\in\Ob \langle E_{k}\rangle.
\]We can fit these triangles into a bigger diagram which, by the octahedral axiom TR4, has exact column and rows:
\[\xymatrix{
\mathbb R_{E_k}\mathbb R_{\langle E_1,\dots E_{k-1}\rangle}X\ar[d]\ar[r] & X\ar@{=}[d]\ar@{-->}[r]& F'' \ar@{-->}[d]\\
\mathbb R_{\langle E_1,\dots E_{k-1}\rangle}X\ar[r]\ar[d] & X\ar[r]\ar@{-->}[d] & F\ar@{-->}[d]\\
F'\ar@{-->}[r] & 0\ar@{-->}[r] & F'[1]
}
\] Here the upper-left square is commutative. Focusing on the right column, we have that $F\in\Ob\langle E_1,\dots, E_{k-1}\rangle$, $F'\in\Ob\langle E_k\rangle$, and consequently $F''\in\Ob\langle E_1,\dots, E_k\rangle$, being the subcategory triangulated, and thus closed by taking cones. The association $X\mapsto F''$ define a left adjoint for the inclusion $\langle E_1,\dots, E_k\rangle\to\mathscr D$, and $\langle E_1,\dots, E_k\rangle$ is left admissible. A similar argument shows that $\langle E_1,\dots, E_k\rangle$ is right admissible.
\endproof

\bsh
\begin{prop}
Let $\mathscr A$ be an admissible full triangulated subcategory of $\mathscr D$.
\begin{enumerate}
\item Both functors $\mathbb L_\mathscr A, \mathbb R_\mathscr A$ are vanishing if restricted to $\mathscr A$.
\item For any $X\in\Ob(\mathscr D)$ we have that $\mathbb L_\mathscr A(X)\in\Ob(\mathscr A^\perp)$ and $\mathbb R_\mathscr A(X)\in\Ob(^\perp\mathscr A)$.
\item The restricted functors $\mathbb L_\mathscr A|_{^\perp\mathscr A},\mathbb R_\mathscr A|_{\mathscr A^\perp}$ induce mutually inverse equivalences $^\perp\mathscr A\to \mathscr A^\perp$ and $\mathscr A^\perp\to\ ^\perp\mathscr A$, respectively.
\item If $\Psi\in\operatorname{Aut}(\mathscr D)$ is an auto-equivalence of $\mathscr D$, then
\[\Psi\circ\mathbb L_{\mathscr A}=\mathbb L_{\Psi(\mathscr A)}\circ\Psi, \quad \Psi\circ\mathbb R_{\mathscr A}=\mathbb R_{\Psi(\mathscr A)}\circ\Psi.
\]
\end{enumerate}
\end{prop}
\esh

\subsection{Saturatedness and Serre Functors} 

\begin{defi}[\cite{BK}]
A triangulated $\mathbb K$-linear category $\mathscr D$ is \emph{saturated} if and only if any (covariant/contravariant) cohomological functor of finite type, i.e. any functor
\[F\colon\mathscr D\to\operatorname{Vect}_{\mathbb K}^{<\infty},\quad F\colon\mathscr D^{\operatorname{op}}\to\operatorname{Vect}_{\mathbb K}^{<\infty}
\]such that
\begin{itemize}
\item $F$ takes distinguished triangles into exact sequences,
\item $\sum_{i\in\mathbb Z}\dim_{\mathbb K}F(A[i])<\infty$ for any object $A\in\operatorname{Ob}(\mathscr D)$,
\end{itemize}
is representable. This is equivalent to the requirement that any exact functor $\Phi\colon \mathscr D\to\mathcal D^b(\mathbb K)$ is representable, the category $\mathscr D$ being $\mathcal D^b(\mathbb K)$-enriched by seeing $\Hom^\bullet(X,Y)$ as a complex with trivial differentials.
\end{defi}

The following results describe how the properties of admissibility and saturatedness interact with one other.

\bsh
\begin{prop}\label{2-4-17.2}Let $\mathscr D$ be a $\mathbb K$-linear triangulated category.
\begin{enumerate}

\item If $\mathscr D$ is saturated, and $\mathscr A\subseteq\mathscr D$ is left (or right) admissible, then $\mathscr A$ is saturated.
\item If $\mathscr A$ is a saturated category, imbedded in $\mathscr D$ as a full triangulated subcategory, then $\mathscr A$ is admissible.
\item If $\mathscr D=\langle\mathscr A_1,\dots,\mathscr A_n\rangle$ is a semiorthogonal decomposition, and $\mathscr D$ is saturated, then each $\mathscr A_i$ is admissible.
\end{enumerate}
\end{prop}
\esh

\proof For a proof of the points (1) and (2), see \cite{Bondal}, \cite{BK} and \cite{kuz1}. For the point (3) let us proceed by induction on the length of the semiorthogonal decomposition. Since $\mathscr D=\langle \mathscr A_n^\perp, \mathscr A_n\rangle$, by point Lemma \ref{2-4-17.1} it follows that $\mathscr A_n$ is right admissible, and $\mathscr A_n^\perp$ left admissible. Hence $\mathscr A_n$ and $\mathscr A_n^\perp$ are saturated (by (1)), and admissible (by (2)). A simple inductive argument completes the proof.
\endproof

\bsh
\begin{prop}\label{2-4-17.3}Let $\mathscr D$ be a $\mathbb K$-linear triangulated category.
\begin{enumerate}

\item Let $\mathscr A$ be an admissible subcategory of $\mathscr D$. Suppose that both $\mathscr A$ and $\mathscr A^\perp$ are saturated. Then also $\mathscr D$ is saturated.

\item Let $\mathscr D=\langle\mathscr A_1,\dots,\mathscr A_n\rangle$ be a given semiorthogonal decomposition. Then $\mathscr D$ is saturated if and only if each $\mathscr A_i$ is saturated.
\end{enumerate}
\end{prop}
\esh

\proof For a proof of point (1) see \cite{BK}. For point (2), let us suppose that $\mathscr D$ is saturated. Then by (3) and (1) of Proposition \ref{2-4-17.2}, we have that each $\mathscr A_i$ is saturated. Vice versa, an inductive argument completes the proof, using (1). 
\endproof

\begin{defi}[\cite{BK}]A \emph{Serre functor} in a $\mathbb K$-linear category of finite type $\mathscr D$ is a $\mathbb K$-linear auto-equivalence $\kappa \colon\mathscr D\to\mathscr D$ such that there exist bi-functorial isomorphisms of $\mathbb K$-vector spaces
\[\xymatrix{\eta_{A,B}\colon\Hom(A,B)\ar[r]^{\cong}&\Hom(B,\kappa(A))^*}
\]for any two objects $A,B\in\operatorname{Ob}(\mathscr D)$.
\end{defi}

If $\mathscr D$ is a $\mathbb K$-linear triangulated category for which there exists a Serre functor, then it is automatically compatible with the triangulated structure. 

\bsh
\begin{prop}[\cite{BK}] Any Serre functor on a triangulated $\mathbb K$-linear category is exact, i.e.
\begin{itemize}
\item it commutes with shift operators,
\item it takes distinguished triangles into distinguished triangles.
\end{itemize}
Moreover, we have that
\begin{enumerate}
\item the category $\mathscr D$ has a Serre functor if and only if all functors $\Hom(X,-)^*,\Hom(-,X)^*$, for any object $X\in\operatorname{Ob}(\mathscr D)$, are representable. 
\item Any two Serre functors $\kappa_1$ and $\kappa_2$ are connected by a canonical functor isomorphism.
\end{enumerate}
\end{prop}
\esh

Because of the previous Proposition, it is clear that any $\Hom$-finite saturated category admits a Serre functor, since the functors $\Hom(X,-)^*,\Hom(-,X)^*$ are cohomological of finite type, for any object $X\in\operatorname{Ob}(\mathscr D)$.

\bsh
 \begin{prop}[\cite{BK}]\label{serre}
 Let $\mathscr D$ be a triangulated $\Hom$-finite $\mathbb K$-linear category admitting a full exceptional collection $\frak E=(E_0,\dots,E_n)$. Then $\mathscr D$ is saturated, and hence it has a Serre functor.
 \end{prop}
 \esh
 \proof
 By Proposition \ref{2-4-17.3}, we already know that $\mathscr D$ is saturated if and only if each subcategory generated by an $E_i$ is saturated. It is easily seen that the category generated by an exceptional object $E$ is saturated: an exact functor $\Phi\colon\mathscr D\to \mathcal D^b(\mathbb K)$ is represented by $\Phi(E)^*$ if $\Phi$ is covariant, $\Phi(E)$ otherwise. \endproof

\subsection{Dual Exceptional Collections and Helices} In the subsequent subsections we will always assume that the $\mathbb K$-linear triangulated category $\mathscr D$ admits a full exceptional collection $\frak E:=(E_0,E_1,\dots,E_n)$.
\subsubsection{Left and Right Dual Exceptional Collections}
Starting from the full exceptional collection $\frak E:=(E_0,E_1,\dots,E_n)$, we can define other two collections
\[^\vee\frak E:=(^\vee E_0,{^\vee}E_1,\dots,{^\vee}E_n),\quad \frak E^\vee:=(E_0^\vee,E_1^\vee,\dots,E_n^\vee),
\]
called respectively \emph{left} and \emph{right dual collections}, defined by iterated mutations
\[^\vee E_k:=\mathbb R_{E_n}\mathbb R_{E_{n-1}}\dots \mathbb R_{E_{n-k+1}}E_{n-k},
\]
\[E_k^\vee:=\mathbb L_{E_0}\mathbb L_{E_1}\dots \mathbb L_{E_{n-k-1}}E_{n-k}
\]for $k=0,1,\dots,n$. Adopting the conventions of Remark \ref{23-04-17.1}, we can define the dual exceptional collections through the action of the braids
\[^\vee\frak E:=\frak E^{\beta'},\quad\beta':=(\beta^{-1}_{12}\beta^{-1}_{23}\dots\beta^{-1}_{n,n+1})(\beta^{-1}_{12}\beta^{-1}_{23}\dots\beta^{-1}_{n-1,n})\dots\beta^{-1}_{12},
\]
\[\frak E^\vee:= \frak E^\beta,\quad\beta:=(\beta_{n,n+1},\beta_{n-1,n}\dots,\beta_{12})(\beta_{n,n+1},\beta_{n-1,n}\dots,\beta_{23})\dots\beta_{n,n+1}.
\]

Notice that we have

\begin{equation}\label{dualorthright}
\Hom^\bullet(E_h, E_k^\vee)=\begin{sistema}
0\quad\text{if }h=0,\dots, n-k-1,\quad\text{by definition of left mutation}\\
\\
\mathbb K\quad\text{if }h=n-k,\\
\\
0\quad\text{if }h=n-k+1,\dots,n,\quad\text{by iteration of \eqref{functiso1}}
\end{sistema}
\end{equation}

\begin{equation}\label{dualorthleft}\Hom^\bullet(^\vee E_k,E_h)=\begin{sistema}
0\quad\text{if }h=0,\dots, n-k-1,\quad\text{by iteration of \eqref{functiso2}}\\
\\
\mathbb K\quad\text{if }h=n-k,\\
\\
0\quad\text{if }h=n-k+1,\dots,n,\quad\text{by definition of right mutation}
\end{sistema}
\end{equation}
$ $\\
where the graded vector space $\mathbb K$ is concentrated in degree $0$. In other words, we have that
\begin{itemize}
\item $\Hom^\alpha(E_h, E_k^\vee)$ vanishes except for $\alpha=0$ and $h=n-k$ (in which case it is $\mathbb K$),
\item $\Hom^\alpha(^\vee E_k,E_h)$ vanishes except for $\alpha=0$ and $h=n-k$ (in which case it is $\mathbb K$).
\end{itemize}
These orthogonality relations actually define the left and right dual collections uniquely up to \emph{unique} isomorphisms: this is a consequence of Yoneda Lemma, as the following results shows.

\bsh
\begin{prop}[\cite{helix}]\label{propdual}
If $\mathscr D$ is a $\mathbb K$-linear triangulated category generated by the exceptional collection $(E_0,\dots,E_n)$, then for any $k=0,1,\dots, n$ we have that:
\begin{itemize}
\item the object $^\vee E_k$ represents the covariant functor
\[X\mapsto \Hom^\bullet(E_{n-k},\mathbb L_{E_{n-k+1}}\dots \mathbb L_{E_{n}}X);
\]
\item the object $E^\vee_k$ represents the contravariant functor
\[X\mapsto\Hom^\bullet(E_{n-k},\mathbb L_{E_{n-k+1}}\dots \mathbb L_{E_{n}}X)^*.
\]
\end{itemize}
In particular, for any object $X\in\operatorname{Ob}(\mathscr D)$, we get the functorial isomorphisms
\[\Hom^\bullet(^\vee E_k,X)=\Hom^\bullet(X,E^\vee_k)^*.
\]
\end{prop}
\esh

\proof
Observing that $\mathbb L_{E_0}\dots \mathbb L_{E_n}X=0$ for any object $X$, since it is an object of the subcategory $\langle E_0,\dots,E_n\rangle^\perp=\mathscr D^\perp=0$, and applying the functors $\Hom^\bullet(^\vee E_k,-)$ and $\Hom^\bullet(-, E_k^\vee)$ to the triangle
\[\Hom^\bullet(E_h,\mathbb L_{E_{h+1}}\dots \mathbb L_{E_{n}}X)\otimes E_h\to \mathbb L_{E_{h+1}}\dots \mathbb L_{E_{n}}X\to \mathbb L_{E_{h}}\dots \mathbb L_{E_{n}}X
\]starting from $h=0$ up to $h=n-k-1$, we iteratively obtain 
\[\Hom^\bullet(^\vee E_k,\mathbb L_{E_{h}}\dots \mathbb L_{E_{n}}X)=\Hom^\bullet(\mathbb L_{E_{h}}\dots \mathbb L_{E_{n}}X, E_k^\vee)=0
\]for any $h\in\left\{0,\dots,n-k\right\}$. So, at the step $h=n-k$, applying $\Hom^\bullet(^\vee E_k,-)$, we get
\begin{align*}
\Hom^\bullet(^\vee E_k,X)&=\Hom^\bullet(^\vee E_k, \mathbb L_{E_{n-k+1}}\dots \mathbb L_{E_n}X)\quad\quad\quad (\text{by iteration of }\eqref{functiso1})\\
&=\Hom^\bullet( E_{n-k}, \mathbb L_{E_{n-k+1}}\dots \mathbb L_{E_n}X)\otimes \Hom^\bullet(^\vee E_k,E_{n-k})\\
&=\Hom^\bullet( E_{n-k}, \mathbb L_{E_{n-k+1}}\dots \mathbb L_{E_n}X)
\end{align*}
because of orthogonality relations \eqref{dualorthleft}. Analogously, applying $\Hom^\bullet(-, E_k^\vee)$ to the same triangle, we get
\begin{align*}
\Hom^\bullet(X, E_k^\vee)&=\Hom^\bullet(\mathbb L_{E_{n-k+1}}\dots \mathbb L_{E_n}X, E_k^\vee)\quad\quad\quad(\text{by iteration of }\eqref{functiso2})\\
&=\Hom^\bullet(E_{n-k}, \mathbb L_{E_{n-k+1}}\dots \mathbb L_{E_n}X)^*\otimes \Hom^\bullet(E_{n-k},E^\vee_k)\\
&=\Hom^\bullet(E_{n-k}, \mathbb L_{E_{n-k+1}}\dots \mathbb L_{E_n}X)^*
\end{align*}
because of orthogonality relations \eqref{dualorthright}.
\endproof

\subsubsection{Geometrical Dual Exceptional Collection}\label{geomdualcoll} In the geometric case, i.e. the case in which the triangulated category is the derived category of coherent sheaves $\mathcal D^b(X)$ of a smooth complex projective variety $X$, then we can introduce a third notion of \emph{geometrical dual} exceptional collection. Given an object $E^\bullet\in{\rm Ob}(\mathcal D^b(X))$, we define its dual to be the object
\[(E^\bullet)^*:={\bf R}\HOM_{\mathcal O_X}^{\bullet}(E^\bullet,\mathcal O_X),
\]
where ${\bf R}\HOM_{\mathcal O_X}^\bullet(E^\bullet,-)\colon D^+(\mathcal {QC}oh(X))\to D^+(\mathcal {QC}oh(X))$ denotes the right derived functor of the left exact functor 
\[\HOM_{\mathcal O_X}^{\bullet}(E^\bullet,-)\colon K^+(\mathcal{QC}oh(X))\to K^+(\mathcal{QC}oh(X)),
\]
\[\HOM_{\mathcal O_X}^{j}(E^\bullet,\mathcal L^\bullet):=\prod_p\HOM_{\mathcal O_X}(E^p,\mathcal L^{p+j}),\quad \text{with }df:=f\circ d_{E}+(-1)^{j+1}d_{\mathcal L}\circ f.
\]
Because of the smoothness assumption on the variety $X$, each object $E^\bullet\in{\rm Ob}(\mathcal D^b(X))$ is isomorphic to a bounded complex of locally free sheaves. Consequently, we have that
\begin{enumerate} 
\item $(E^\bullet)^*$ is a well defined object of $\mathcal D^b(X)$,
\item each object is \emph{reflexive}, i.e. $(E^\bullet)^{**}=E^\bullet$.
\end{enumerate}
Using the reflexivity property, the following definition is well posed.

\begin{defi}
Let $\frak E=(E_0,\dots, E_n)$ be a full exceptional collection in $\mathcal D^b(X)$. We define the \emph{geometrical dual exceptional collection} $\frak E^*$ as the collection
\[\frak E^*:=(E_n^*,\dots, E_0^*).
\]
\end{defi}

 \subsubsection{Helices}Following \cite{helix}, we introduce the 
\begin{defi}[Helix] If $(E_0,\dots, E_n)$ is a full exceptional collection, we call \emph{helix} the infinite collection $(E_i)_{i\in\mathbb Z}$ defined by the iterated mutations
\[E_{i+n+1}=\mathbb R_{E_{i+n}}\dots \mathbb R_{E_{i+1}}E_i,\quad E_{i-n-1}=\mathbb L_{E_{i-n}}\dots \mathbb L_{E_{i-1}}E_i,\quad i\in\mathbb Z.
\]Such a helix is said to be \emph{of period }$n+1$, and any family of $n+1$ consecutive objects $(E_i, E_{i+1},\dots, E_{i+n})$ is called \emph{helix foundation}. The braid group $\mathcal B_{n+1}$ acts on the set of helices of period $n+1$: the mutations functors $\mathbb L_i,\mathbb R_i$ act on the helix by replacing all the pairs 
\[(E_{i-1+k(n+1)}, E_{i+k(n+1)}),\quad\text{with }k\in\mathbb Z
\]with their left/right mutations. In this way, the mutation of a helix is still a helix.
\end{defi}

\bsh
\begin{prop}[\cite{helix}]Let $(E_{i-n-1},\dots, E_{i-1})$ and $(E_i,\dots, E_{i+n})$ be two consequent foundations of an helix in a $\mathbb K$-ilnear triangulated category $\mathscr D$. For any object $X\in\Ob(\mathscr D)$ the following functorial isomorphisms hold:
\begin{equation}\label{24.07.17-1}\Hom^\bullet(E_i,X)=\Hom^\bullet(X,E_{i-n-1})^*.
\end{equation}
In particular, we deduce the \emph{periodicity condition}
\begin{equation}\label{periodicitaserre}\Hom^\bullet(E_i,E_j)=\Hom^\bullet(E_{i-n-1},E_{j-n-1}).
\end{equation}
\end{prop}
\esh

\proof Notice that, if we consider two consecutive foundations of a helix of period $n+1$, i.e.
\[E_{i-n-1},\dots, E_{i-1},\quad E_i,\dots, E_{i+n},
\]the collection $(F_0,\dots, F_n)$, defined by the relations
\[F_k:=\mathbb L_{E_i}\dots \mathbb L_{E_{i+n-k-1}}E_{i+n-k}=\mathbb R_{E_{i-1}}\dots \mathbb R_{E_{i-k}}E_{i-k-1},
\]is at the same time right dual collection of $(E_i,\dots, E_{i+n})$ and left dual collection of $(E_{i-n-1},\dots, E_{i-1})$. Thus, by Proposition \ref{propdual}, we deduce that for any object $X\in\operatorname{Ob}(\mathscr D)$ it holds the functorial isomorphism
\[\Hom^\bullet(E_i,X)=\Hom^\bullet(X,E_{i-n-1})^*.
\]Applying it to $X=E_j$, we obtain
\[\Hom^\bullet(E_i,E_j)=\Hom^\bullet(E_{j},E_{i-n-1})^*;
\]analogously we have
\[\Hom^\bullet(E_j,X)=\Hom^\bullet(X,E_{j-n-1})^*,\]and dualizing we get
\[\quad\Hom^\bullet(E_j,X)^*=\Hom^\bullet(X,E_{j-n-1}).
\]If we take $X=E_{i-n-1}$ in the last isomorphism, we finally obtain the periodicity condition
\[\Hom^\bullet(E_i,E_j)=\Hom^\bullet(E_{i-n-1},E_{j-n-1}).
\]
\endproof

We already know, by Proposition \ref{serre}, that the category $\mathscr D$ admits a Serre functor $\kappa\colon\mathscr D\to\mathscr D$ (unique up to a canonical isomorphism): from the result above we deduce that if $\frak E:=(E_0,\dots, E_n)$ is an exceptional collection, then 
\[\kappa(E_i)=E_{i-n-1},\quad i\in \mathbb Z,
\]for any exceptional object of the helix generated by $\frak E$. Remarkably, the knowledge of the action of $\kappa$ on such a helix is enough to reconstruct its action on the whole category $\mathscr D$.

\bsh
\begin{cor}[\cite{helix}]\label{corhelixser}
The action on the set of full exceptional collections $\frak E$'s (of length $n+1$) of the central element of the braid group $\mathcal B_{n+1}$
\[\beta^2=(\beta_{n,n+1}\beta_{n-1,n}\dots\beta_{23}\beta_{12})^{n+1},\quad \beta:=(\beta_{n,n+1},\beta_{n-1,n}\dots,\beta_{12})(\beta_{n,n+1},\beta_{n-1,n}\dots,\beta_{23})\dots\beta_{n,n+1},
\]can be extended to a Serre functor $\kappa$ of the category $\mathscr D$.
\end{cor}
\esh

\proof
We have to show that given an object $X\in\Ob(\mathscr D)$, the image $\kappa(X)$ is uniquely determined by the images of the objects of an exceptional collection $(E_0,\dots, E_n)$. Let us consider $X$ as the canonical convolution of the Postnikov system whose associated complex is
\[0\to V_0^\bullet\otimes E_0\to V_1^\bullet\otimes E_1\to\dots\to V_n^\bullet\otimes E_n\to 0,\quad V_k^\bullet:=\Hom^\bullet(\mathbb R_{E_{k-1}}\dots\mathbb R_{E_0}X, E_k)^*.
\]The differentials are given by an element of
\[\bigoplus_{p}\Hom^1(V_p^\bullet\otimes E_p, V_{p+1}^\bullet\otimes E_{p+1})\cong \bigoplus_{p,\alpha}\Hom^{-\alpha}(V_p^\bullet,V_{p+1}^\bullet)\otimes \Hom^{\alpha+1}(E_p, E_{p+1}).
\]By \eqref{periodicitaserre}, the same element defines differentials of the complex
\[0\to V_0^\bullet\otimes E_{-n-1}\to V_1^\bullet\otimes E_{-n}\to\dots\to V_n^\bullet\otimes E_{-1}\to 0,\]
whose canonical convolution defines an object $\kappa(X)$. In order to show that $$\Hom^\bullet(Y, X)\cong \Hom^\bullet(\kappa(X), Y)^*,$$ we use the procedure described in Remark \ref{spectral} for both the linear covariant cohomological functors $\Hom^\bullet(Y, -)$ and $\Hom^\bullet(\kappa(-), Y)^*$. The first one is computed through the spectral sequence whose first sheet is
\[E_1^{p,q}=\Hom^q(Y, V^\bullet_p\otimes E_p)=\bigoplus_\alpha V^{-\alpha}_p\otimes \Hom^{q+\alpha}(Y, E_p), 
\]and for the second one we have
\[E_1^{p,q}=\bigoplus_\alpha V^{-\alpha}_p\otimes \left(\Hom^{-q-\alpha}(E_{p-n-1}, Y)\right)^*.
\]By \eqref{24.07.17-1} one concludes.
\endproof

\subsubsection{$m$-Blocks}

Following B.V. Karpov and D.Y. Nogin (\cite{karpov-nogin, karpluck}) we introduce the definition of \emph{m-blocks} as a particular class of exceptional collections.

\begin{defi}[\cite{karpov-nogin}]\label{blockcoll}If $\mathscr D$ is a triangulated $\mathbb K$-linear category of finite type, and if $\frak E=(E_1,\dots,E_k)$ is an exceptional collection, we will say that $\frak E$ is a \emph{block} if 
\[\Hom^\bullet(E_i,E_j)=0\quad\text{whenever }i\neq j.
\]
More in general, an $m$-block is an exceptional collection 
\[(\frak E_1,\dots,\frak E_m)=( E_{11},\dots,E_{1\alpha_1},E_{21},\dots,E_{2\alpha_2},\dots, E_{m1},\dots, E_{m\alpha_m})
\]such that all subcollections $\frak E_j=( E_{j1},\dots, E_{j\alpha_j})$ are blocks. We will call \emph{type} of the $m$-block $\frak E$ the $m$-tuple $(\alpha_1,\dots,\alpha_m)$.
\end{defi}

A close notion of \emph{levelled exceptional collection} has been introduced by L. Hille in \cite{hillealg}. The following result is an immediate consequence of the vanishing condition defining an $m$-block.

\begin{prop}\label{propmutblock}
Let $\frak E$ be an $m$-block exceptional collection of type $(\alpha_1,\dots, \alpha_m)$. The left/right mutations of two objects in a same block $\frak E_j$ act just as permutations. More concretely, if $E_{h}, E_{h+1}$ are objects in a same block of an $m$-block exceptional collection
\[\frak E:=(E_1,\dots,\underbrace{\dots, E_h,E_{h+1}, \dots\ }_{\text{same block}}, \dots, E_n),
\]
then 
\[\frak E^{\beta_{h,h+1}}=(E_1,\dots,\underbrace{\dots, E_{h+1},E_{h}, \dots\ }_{\text{same block}}, \dots, E_n).
\]

\end{prop}

\newpage
\section{Non-symmetric orthogonal geometry of Mukai lattices}\label{Mukaisec}

In this Section, companion of the previous one, we discuss properties of classes of the exceptional objects of a triangulated category $\mathscr D$ in its Grothendieck group $K_0(\mathscr D)$, the \emph{exceptional bases}. This can be interpreted as a sort of \emph{linearization} procedure, translating objects of an abstract category into elements of a more elementary algebraic structure (a lattice or vector space), more suitable for explicit computations. In particular, we summarize the main points of the geometry of lattices endowed with unimodular bilinear non-symmetric forms, the so-called \emph{Mukai lattices}, introduced and extensively studied by A. L. Gorodentsev \cite{Go3,Go2.2}. 
\subsection{Grothendieck Group and Mukai Lattices}
Let $\mathscr D$ be a (small) triangulated category. Let us denote by $[\mathscr D]$ the set of isomorphism classes of objects of $\mathscr D$.
\begin{defi}The \emph{Grothendieck group} $K_0(\mathscr D)$ is the group defined as the quotient of the free abelian group generated by $[\mathscr D]$ over the following \emph{Euler relations}:
\[[B]=[A]+[C],
\] whenever there is a triangle in $\mathscr D$
\[\xymatrix{
A\ar[r]&B\ar[r]&C\ar[r]&A[1].
}
\]This group is the solution of the following universal problem: 
\emph{to find an abelian group $X$ and a function $[-]\colon [\mathscr D]\to X$ such that, given a function $\varphi\colon [\mathscr D]\to G$, with values in an abelian group $G$, and preserving the Euler relations, there exists a unique group homomorphism $\overline{\varphi}\colon X\to G$ making the following diagram commutative}
\[\xymatrix{
[\mathscr D]\ar[rr]^{\varphi}\ar[rd]_{[-]}&&G\\
&X\ar[ru]_{\overline{\varphi}}&
}
\]
\end{defi}
A triangulated functor $F\colon\mathscr D\to\mathscr D'$ induces a group homomorphism between $K_0(\mathscr D)$ and $K_0(\mathscr D')$, by sending $[E]$ to $[F(E)]$.
If $\mathscr D$ is $\mathbb K$-linear, we can naturally define the so called \emph{Grothendieck--Euler--Poincaré pairing} 
\[\chi(E,F):=\sum_{i}(-1)^i\dim_{\mathbb K}\Hom ^i(E,F),
\]for any objects $E,F\in\Ob\mathscr D$.

\begin{oss}Let us write some useful identities valid in Grothendieck groups. First of all, note that $[0]=0$. Moreover, from the distinguished triangle $A\to A\oplus B\to B\to A[1]$ we have
\[[A\oplus B]=[A]+[B],
\]whereas from $A\to 0\to A[1]\to A[1]$ we deduce that $[A[1]]=-[A]$. 
\end{oss}

\bsh
\begin{lemma}\label{projofmut}
If $E\in\Ob\mathscr D$ is an exceptional object, then for any object $X\in\Ob\mathscr D$ the following identities hold in the Grothendieck group $K_0(\mathscr D)$:
\[[\mathbb L_EX]=[X]-\chi(E,X)\cdot [E],
\]
\[[\mathbb R_EX]=[X]-\chi(X,E)\cdot [E].
\]
\end{lemma}
\esh

\proof
From the distinguished triangle defining $\mathbb L_EX$ we have
\begin{align*}
[\mathbb L_EX]&=[X]-[\Hom^\bullet(E,X)\otimes E]\\
&=[X]-\left[\bigoplus_i E[-i]^{\oplus\dim_{\mathbb K}\Hom^i(E,X)}\right]\\
&=[X]-\left(\sum_i(-1)^i\dim_{\mathbb K}\Hom^i(E,X)\right)[E].
\end{align*}
Analogously, we have
\begin{align*}
[\mathbb R_EX]&=[X]-[\Hom^\bullet(X,E)^*\otimes E]\\
&=[X]-\left[\bigoplus_i E[-i]^{\oplus\dim_{\mathbb K}(\Hom^{-i}(X,E))^*}\right]\\
&=[X]-\left(\sum_i(-1)^i\dim_{\mathbb K}\Hom^{-i}(X,E)\right)[E].
\end{align*}
\endproof

Let us assume that $\mathscr D$ admits a full exceptional collection $(E_0,\dots, E_n)$: it then follows that $K_0(\mathscr D)$ is freely generated by $([E_0],\dots,[E_n])$. In this case $(K_0(\mathscr D),\chi(\cdot ,\cdot))$ admits a structure of \emph{exceptional unimodular Mukai lattice}:

\begin{defi}[Mukai Lattice]A \emph{unimodular Mukai lattice} is a finitely generated free $\mathbb Z$-module $V$ endowed with a unimodular bilinear (not necessarily symmetric) form $\langle\cdot, \cdot\rangle\colon V\times V\to\mathbb Z$. An element $e\in V$ will be said to be \emph{exceptional} if $\langle e,e\rangle=1$. A $\mathbb Z$-basis $\varepsilon:=(e_0,\dots,e_n)$ of the Mukai lattice is called \emph{exceptional} if
\[\langle e_i,e_i\rangle=1\quad\text{for all }i,\quad\langle e_j,e_i\rangle=0\quad\text{for }j>i.
\]In other words, the Gram matrix must be of the upper triangular form
\[\begin{pmatrix}
1&&&&\\
0&1&&*&\\
0&0&1&&\\
\vdots&\vdots&\vdots&\ddots&\\
0&0&0&\dots&1
\end{pmatrix}.
\] A Mukai lattice is called \emph{exceptional} if it admits an exceptional basis.
\end{defi}
It is thus clear that the projection on $K_0(\mathscr D)$ of a full exceptional collection in $\mathscr D$ is an exceptional basis.

\begin{defi}[Mutations of exceptional bases]Let $(V,\langle\cdot,\cdot\rangle)$ be an exceptional Mukai lattice of rank $n+1$. If $\varepsilon:=(e_0,\dots,e_n)$ is an exceptional basis we define for any $1\leq i\leq n$
\[\mathbb L_i\varepsilon:=(e_0,\dots, \mathbb L_{e_{i-1}}e_i,e_{i-1},\dots, e_n),\quad \mathbb L_{e_{i-1}}e_i:=e_i-\langle e_{i-1},e_i\rangle\cdot e_{i-1},
\]
\[\mathbb R_i\varepsilon:=(e_0,\dots,e_i,\mathbb R_{e_i}e_{i-1} ,\dots,e_n),\quad \mathbb R_{e_i}e_{i-1}:=e_{i-1}-\langle e_{i-1},e_i\rangle\cdot e_i.
\]In particular, we still get exceptional basis, called \emph{left} and \emph{right mutations} of $\varepsilon$. It is easy to see that this defines an action\footnote{In what follows we will use the same conventions and notations of Remark \ref{23-04-17.1} for the action of the braid group on the set of exceptional bases of a Mukai lattice.} of the braid group $\mathcal B_n$ on the set of exceptional bases of $V$. \end{defi}

Note that, accordingly to Lemma \ref{projofmut}, the projection on the Grothendieck group of the mutation of a full exceptional collection $(E_0,\dots, E_n)$ coincides with the corresponding mutation of the exceptional basis $([E_0],\dots,[E_n])$.

\begin{defi}[Left and right dual exceptional bases] Given an exceptional basis $\varepsilon=(e_0,\dots, e_n)$ of an exceptional Mukai lattice $(V,\langle\cdot,\cdot\rangle)$, we define two other exceptional bases $$^\vee\varepsilon=(^\vee e_0,\dots, ^\vee e_n)\quad\text{ and }\quad\varepsilon^\vee=(e_0^\vee ,\dots, e_n^\vee )$$ called respectively \emph{left} and \emph{right dual exceptional bases} defined through the action of the braids 
\[^\vee\varepsilon:=\varepsilon^{\beta'},\quad\beta':=(\beta^{-1}_{12}\beta^{-1}_{23}\dots\beta^{-1}_{n,n+1})(\beta^{-1}_{12}\beta^{-1}_{23}\dots\beta^{-1}_{n-1,n})\dots\beta^{-1}_{12},
\]
\[\varepsilon^\vee:= \varepsilon^\beta,\quad\beta:=(\beta_{n,n+1},\beta_{n-1,n}\dots,\beta_{12})(\beta_{n,n+1},\beta_{n-1,n}\dots,\beta_{23})\dots\beta_{n,n+1}.
\]
Notice, in particular that we have the following orthogonality relations
\begin{equation}\label{23-04-17.2}\langle e_h, e_k^\vee \rangle=\delta_{h,n-k},\quad \langle ^\vee e_k, e_h \rangle=\delta_{h,n-k}.
\end{equation}
\end{defi}

\bsh
\begin{prop}\label{propbraidmukai}
If $\varepsilon=(e_0,\dots, e_n)$ is an exceptional basis of $(V,\langle\cdot ,\cdot\rangle)$, and $G$ is the Gram matrix of $\langle\cdot ,\cdot\rangle$ with respect to $\varepsilon$, i.e.
\[G_{hk}:=\langle e_h ,e_k\rangle\quad 0\leq h,k\leq n,
\]then the Gram matrix  
\begin{itemize}
\item with respect to the exceptional basis $\mathbb L_i\varepsilon$ is given by $H^i\cdot G \cdot H^i$, where
\[H^i=\begin{pmatrix}
1&&&&&&&\\
&\ddots&&&&&&\\
&&1&&&&&\\
&&&-G_{i-1,i}&1&&&\\
&&&1&0&&&\\
&&&&&1&&\\
&&&&&&\ddots&\\
&&&&&&&1\\
\end{pmatrix},
\]the entry $-G_{i-1,i}$ being in the position $(i-1,i-1)$;
\item with respect to the exceptional basis $\mathbb R_i\varepsilon$ is given by $K^i\cdot G\cdot K^i$, where
\[K^i=\begin{pmatrix}
1&&&&&&&\\
&\ddots&&&&&&\\
&&1&&&&&\\
&&&0&1&&&\\
&&&1&-G_{i-1,i}&&&\\
&&&&&1&&\\
&&&&&&\ddots&\\
&&&&&&&1\\
\end{pmatrix},
\]the entry $-G_{i-1,i}$ being in the position $(i,i)$;
\item with respect to both the right and left dual exceptional basis $\varepsilon^\vee$ and $^\vee\varepsilon$ is given by
\[J\cdot G^{-T}\cdot J,
\]where $J$ is the anti-diagonal matrix $$J=\begin{pmatrix}
&&1\\
&\iddots&\\
1&&
\end{pmatrix}.$$
\end{itemize}
\end{prop}
\esh

\proof Lemma \ref{projofmut} implies that
\[(\mathbb L_i\varepsilon)_k=\sum_a (H^i)^a_ke_a,\quad (\mathbb R_i\varepsilon)_k=\sum_a (K^i)^a_ke_a,
\]from which the first two points immediately follow. If we define a matrix $X$ such that
\[e^\vee_k=\sum_aX^h_ae_a,
\]then the orthogonality relations \eqref{23-04-17.2} imply that
\[G\cdot X=J,
\]so that the Gram matrix with respect to the basis $\varepsilon^\vee$ is given by
\[(G^{-1}\cdot J)^T\cdot G\cdot(G^{-1}\cdot J)=J\cdot G^{-T}\cdot J.
\]The computations for the basis $^\vee\varepsilon$ are identical, and are left as an exercise for the reader. 
\endproof

\subsubsection{Geometrical Dual Exceptional Basis}\label{geomdualbas} In the case of the derived category $\mathcal D^b(X)$ of a smooth projective variety, we can also introduce the following notion.

\begin{defi}
Let $\frak E=(E_0,\dots, E_n)$ be a full exceptional collection in $\mathcal D^b(X)$. We define the \emph{geometrical dual} of the projected exceptional basis $[\frak E]:=([E_0],\dots, [E_n])$ as the basis
\[[\frak E^*]:=([E_n^*],\dots, [E_0^*]).
\]
\end{defi}

Notice that, because of the properness of $X$, the exceptionality of the basis $[\frak E^*]$ can be deduced also from Grothendieck--Hirzebruch--Riemann--Roch Theorem: for any $E^\bullet,F^\bullet\in{\rm Ob}(\mathcal D^b(X))$ one has 
\begin{equation}\label{GHRR}\chi(E^\bullet,F^\bullet)=\int_X{\rm ch}(E^{\bullet*}){\rm ch}(F^{\bullet}){\rm td}(X),
\end{equation}
where
\[\chi(E^\bullet,F^\bullet):=\sum_i(-1)^i\dim_\mathbb C\Hom^i_{\mathcal D^b(X)}(E^\bullet,F^\bullet),
\]
\[{\rm td}(X):=\prod_{j=1}^{\dim_\mathbb CX}\frac{\delta_j}{1-\exp(-\delta_j)},\quad \delta_j\text{ Chern roots of }TX,
\]
and for a complex of locally free sheaves $L^\bullet$
\[{\rm ch}(L^\bullet):=\sum_i(-1)^i{\rm ch}(L^i).
\]
Using \eqref{GHRR} and the reflexivity property, we deduce that $\chi([E_i^*], [E_j^*])=\chi([E_j], [E_i])=0$.

\bsh
\begin{prop}\label{propdualexccoll}
Let $\frak E=(E_0,\dots, E_n)$ be a full exceptional collection in $\mathcal D^b(X)$. If $G$ is the Gram matrix associated with  the exceptional basis $[\frak E]$ of $(K_0(X),\chi(\cdot,\cdot))$, then the Gram matrix associated with the geometrical dual basis $[\frak E^*]$ is given by
\[J\cdot G^T\cdot J,
\]where $J$ is the anti-diagonal matrix
$$J=\begin{pmatrix}
&&1\\
&\iddots&\\
1&&
\end{pmatrix}.$$
\end{prop}
\esh

\subsection{Isometries and canonical operator}
In the previous section, we have seen that in any triangulated category $\mathscr D$ also the group $\operatorname{Aut}(\mathscr D)$ of isomorphism classes of auto-equivalences acts on the set of full exceptional collections. This action projects onto the Grothendieck group $K_0(\mathscr D)$ through the actions of \emph{isometries} preserving the Grothendieck--Euler--Poincaré form, and hence acting on the set of exceptional bases.  

\begin{defi}[Isometries] Given two Mukai lattices $(V_1,\langle\cdot,\cdot\rangle_1),(V_2,\langle\cdot,\cdot\rangle_2)$, any $\mathbb Z$-linear map $\phi\colon V_1\to V_2$ such that
\[\langle x,y\rangle_1=\langle\phi(x),\phi(y)\rangle_2,\quad x,y\in V_1
\]is called \emph{isometry} for the Mukai structures. If $\phi$ is invertible, then we will say that the Mukai structures $(V_1,\langle\cdot,\cdot\rangle_1),(V_2,\langle\cdot,\cdot\rangle_2)$ are isometrically isomorphic. 

The set of all $\mathbb Z$-linear isometric automorphisms $\phi\colon V\to V$ of a Mukai lattice $(V,\langle\cdot,\cdot\rangle)$ is denoted by $\operatorname{Isom}_{\mathbb Z}(V,\langle\cdot,\cdot\rangle)$, or simply $\operatorname{Isom}_{\mathbb Z}$ if no confusion arises. 
\end{defi}

Since Serre functors are prototypical and important auto-equivalences in $\mathbb K$-linear triangulated categories, their projections on the Grothendieck group play a particularly important role.

\begin{defi}[Canonical operator] Given a Mukai lattice $(V,\langle\cdot,\cdot\rangle)$, we call \emph{canonical operator} the unique $\mathbb Z$-linear operator $\kappa\colon V\to V$ satisfying the property
\[\langle x,y\rangle=\langle y,\kappa(x)\rangle,\quad x,y\in V.
\]
\end{defi}

Although Serre functors do not always exist in $\mathbb K$-linear triangulated categories, at the level of Mukai structures this existence problem always admits a solution. 

\begin{defi}[Left and right correlations] Given a Mukai lattice $(V,\langle\cdot,\cdot\rangle)$ there are two well defined operators between $V$ and its dual $V^*:=\Hom_{\mathbb Z}(V,\mathbb Z)$, called respectively \emph{left} and \emph{right correlations}:
\[\lambda\colon V\to V^*\colon x\mapsto \langle x,\cdot\rangle,
\]
\[\rho\colon V\to V^*\colon x\mapsto \langle\cdot,x\rangle.
\]Because of the unimodularity of the pairing $\langle\cdot,\cdot\rangle$, both left and right correlations $\lambda, \rho$ define isomorphisms of abelian groups.
\end{defi}

\bsh
\begin{prop}[\cite{Go3, Go2.2}]\label{26.06.17-1}Let $(V,\langle\cdot,\cdot\rangle)$ be a Mukai lattice.
\begin{enumerate}
\item There exists a unique canonical operator $\kappa\colon V\to V$, and it is defined in terms of left and right canonical correlations as the composition
\[\rho^{-1}\circ\lambda\colon V\to V.
\]
\item Given any basis $(e_1,\dots, e_n)$ (not necessarily exceptional) of $V$, with respect to which the Gram matrix of the pairing $\langle\cdot,\cdot\rangle$ is $G$, then the matrix associated with the canonical operator $\kappa$ is given by
\[\kappa=G^{-1}\cdot G^T.
\]
\end{enumerate}
\end{prop}
\esh
\proof An exercise for the reader.\endproof

\subsection{Adjoint operators and canonical algebra}
Let us consider a Mukai lattice $(V,\langle\cdot,\cdot\rangle)$. 
\begin{defi}[Left and right adjoint operators]Let $\phi\in\End_{\mathbb Z}(V)$. We define two new operators $^\vee\phi$ and $\phi^\vee$ called respectively \emph{left} and \emph{right adjoint to }$\phi$ through the following identities:
\[\langle^\vee\phi (x),y\rangle=\langle x,\phi (y)\rangle,
\]
\[\langle\phi(x),y\rangle=\langle x,\phi^\vee(y)\rangle,
\]
for any $x,y\in V$. Fixed a (non-necessarily exceptional) basis $(e_0,\dots, e_n)$ of $V$, in terms of matrix representation we have
\[^\vee\phi=G^{-T}\cdot\phi^T\cdot G^T,\quad \phi^\vee=G^{-1}\cdot\phi^T\cdot G.
\]
\end{defi}
Because of the non-symmetry of the pairing $\langle\cdot,\cdot\rangle$, in general one has $^\vee\phi\neq\phi^\vee$.

\begin{defi}[Canonical algebra]An endomorphism $\phi\in\End_{\mathbb Z}(V)$ is called \emph{reflexive} if $^\vee\phi=\phi^\vee$. The subalgebra $\mathcal A\subseteq \End_{\mathbb Z}(V)$ of all reflexive operators of $V$ is called \emph{canonical algebra}. \end{defi}

The proofs of the following Proposition is straightforward, and is left as an exercise for the reader.
\begin{prop}[\cite{Go3, Go2.2}]\label{24-04-17.1}
Let $\phi\in\End_{\mathbb Z}(V)$. The following conditions are equivalent:
\begin{enumerate}
\item $^\vee\phi=\phi^\vee$,
\item $\phi=\phi^{\vee\vee}$,
\item $\phi=^{\vee\vee}\phi$,
\item $\phi\kappa=\kappa\phi$.
\end{enumerate}
Hence, the canonical algebra $\mathcal A$ coincides with the center $\mathcal Z(\kappa)$ of the canonical operator.
\end{prop}

\begin{prop}[\cite{Go3, Go2.2}]\label{24-04-17.2}
The following sets are contained in the canonical algebra $\mathcal A$:
\begin{enumerate}
\item $\mathcal A^+:=\left\{\phi\in\End_\mathbb Z(V)\colon ^\vee\phi=\phi^\vee=\phi\right\}$, whose elements are called \emph{self-adjoint operators};
\item $\mathcal A^-:=\left\{\phi\in\End_\mathbb Z(V)\colon ^\vee\phi=\phi^\vee=-\phi\right\}$, whose elements are called \emph{anti-self-adjoint operators};
\item $\operatorname{Isom}(V,\langle\cdot,\cdot\rangle)\equiv\left\{\phi\in\operatorname{Aut}_\mathbb Z(V)\colon ^\vee\phi=\phi^\vee=\phi^{-1}\right\}$.
\end{enumerate}
\end{prop}

Given any field $\mathbb K$, we can extend scalars from $\mathbb Z$ to $\mathbb K$, by considering the vector space $V\otimes_\mathbb Z\mathbb K$ endowed with the non-symmetric bilinear form $\langle\cdot,\cdot\rangle$, extended by $\mathbb K$-bilinearity. All previous definitions (and notations) can be trivially adapted to this extension of scalars. 

\begin{prop}[\cite{Go3, Go2.2}]
If $\mathbb K$ is a field of characteristic not equal to 2, then the following direct sum of $\mathbb K$-vector spaces holds
\[\mathcal A_\mathbb K=\mathcal A^+_{\mathbb K}\oplus\mathcal A^-_{\mathbb K}.
\]
\end{prop}

\begin{prop}[\cite{Go3, Go2.2}]
The Lie algebra of the complex Lie group $\operatorname{Isom}_{\mathbb C}$ is equal to $\mathcal A^-_{\mathbb C}$.
\end{prop}

\subsection{Isometric classification of Mukai structures}
The following Proposition underlines the importance of the canonical operator for the isometric classification of Mukai structures.

\bsh
\begin{prop}[\cite{Go3, Go2.2, GoAlg}]\label{classisomukai}Let $V$ be a free $\mathbb Z$-module of finite rank, and let $\langle\cdot,\cdot\rangle_1,\langle\cdot,\cdot\rangle_2$ two non-symmetric unimodular bilinear forms defining two Mukai lattice structures on $V$.
\begin{enumerate}
\item The two Mukai structures share the same canonical operator if and only if there exists an invertible operator $\psi\in\mathcal A^+_{\langle\cdot,\cdot\rangle_1}\cap\mathcal A^+_{\langle\cdot,\cdot\rangle_2}$ and such that
\[\langle x, y\rangle_1=\langle x, \psi(y)\rangle_2,\quad x,y\in V.
\]
\item If $\mathbb K$ is an algebraically closed field of characteristic zero, the two Mukai vector spaces $(V\otimes_{\mathbb Z}\mathbb K,\langle\cdot,\cdot\rangle_1)$ and $(V\otimes_{\mathbb Z}\mathbb K,\langle\cdot,\cdot\rangle_2)$ are isometrically isomorphic if and only if there exists an isomorphism $\phi\in\operatorname{Aut}_{\mathbb K}(V\otimes_{\mathbb Z}\mathbb K)$ such that
\[\phi\circ\kappa_1=\kappa_2\circ\phi.\]
\end{enumerate}
\end{prop}
\esh

\proof
If $\lambda_1,\lambda_2$ and $\rho_1,\rho_2$ denote respectively the left and right correlations for the two Mukai structures on $V$, then the operator $\psi:=\rho_2^{-1}\rho_1=\lambda_2^{-1}\lambda_1$ satisfies all properties of point (1). Indeed, since $\kappa=\rho_1^{-1}\lambda_1=\rho_2^{-1}\lambda_2$, we have that
\[\psi\kappa=(\rho_2^{-1}\rho_1)(\rho_1^{-1}\lambda_1)=\rho_2^{-1}\lambda_1=(\rho_2^{-1}\lambda_2)(\lambda_2^{-1}\lambda_1)=\kappa\psi,
\]
\[\langle x,\psi y\rangle_2=[\rho_1(y)](x)=\langle x,y\rangle_1,\quad x,y\in V.
\]
For point (2), if the two structures are isometric then existence of an isomorphism $\phi$ intertwining the canonical operators is clear. Hence, let us suppose to have two different Mukai structures on the same vector space $V\otimes_{\mathbb Z}\mathbb K$ sharing the same canonical operator. By point (1), we deduce existence of a self-dual isomorhism $\psi\in\operatorname{Aut}_{\mathbb K}(V\otimes_{\mathbb Z}\mathbb K)$ such that
\[\langle x, y\rangle_1=\langle x, \psi(y)\rangle_2,\quad x,y\in V\otimes_{\mathbb Z}\mathbb K.
\]
The field $\mathbb K$ being algebraically closed, a polynomial $p\in\mathbb K[X]$ can be constructed in such a way that the operator $\alpha:=p(\psi)$ satisfies $\alpha^2=\psi$ (see Lemma 16.2 of \cite{GoAlg}). Such an operator $\alpha$ is self-adjoint, since
\[\alpha^\vee=p(\psi)^\vee=p(\psi^\vee)=p(\psi)=\alpha,
\]and it clearly satisfies the condition $\langle\alpha (x),\alpha (y)\rangle_2=\langle x,y\rangle_1$.
\endproof

In particular, Proposition \ref{classisomukai} implies that a non-degenerate non-symmetric bilinear form over algebraically closed field of characteristic zero is uniquely determinated by Jordan normal form of its canonical operator.

\begin{defi}[]Given a Mukai lattice $(V,\langle\cdot,\cdot\rangle)$, and an algebraically closed field $\mathbb K$ of characteristic zero, the Mukai space $(V\otimes_\mathbb Z\mathbb K,\langle\cdot,\cdot\rangle)$ will be called \emph{decomposable}, if there exist two subspaces $U,V$ such that
\begin{enumerate}
\item $V\otimes_\mathbb Z\mathbb K= U\oplus V$,
\item the restrictions $\langle\cdot,\cdot\rangle$ to $U$ and $V$ are nondegenerate,
\item $U$ and $V$ are \emph{bi-orthogonal}, namely $\langle u,v\rangle=\langle v,u\rangle=0$ for all $u\in U$ and $v\in V$.
\end{enumerate}
The space will be called \emph{indecomposable} if it is not decomposable.
\end{defi}

The following result gives a complete classification of all indecomposable Mukai structures over an algebraically closed field $\mathbb K$ of characteristic zero.

\bsh
\begin{teorema}[\cite{Go3, Go2.2, GoAlg}]\label{13.07.17-5}
Let $(V,\langle\cdot,\cdot\rangle)$ be a Mukai lattice, and let $\mathbb K$ be an algebraically closed field of characteristic zero. If $(V\otimes_\mathbb Z\mathbb K,\langle\cdot,\cdot\rangle)$ is indecomposable, then it is isometrically isomorphic to one of the following Mukai spaces.
\begin{enumerate}
\item {\rm{Space of type $U_n$:}} consider the coordinate space $\mathbb K^n$ endowed with the non-degenerate bilinear form whose Gram matrix with respect to the standard basis is
\[G=\begin{pmatrix}
&&&&&1\\
&&&&-1&1\\
&&&1&-1&\\
&&\iddots&\iddots&&\\
&(-1)^{n-2}&(-1)^{n-3}&&&\\
(-1)^{n-1}&(-1)^{n-2}&&&&
\end{pmatrix}.
\]
In this case, by Proposition \ref{26.06.17-1}, we have that the canonical operator is of the form
\[\kappa=G^{-1}\cdot G^T=(-1)^{n-1}\mathbbm 1+M,\quad M^{n-1}\neq 0,\quad M^n=0,
\]and its Jordan form is 
\[
 J_n((-1)^{n-1}):=
 \begin{pmatrix}
(-1)^{n-1}&1&&&\\
&(-1)^{n-1}&1&&\\
&&\ddots&\ddots&\\
&&&(-1)^{n-1}&1\\
&&&&(-1)^{n-1}
\end{pmatrix}.
\]
\item {\rm Space of type $W_n(	\lambda)$ with $\lambda\neq (-1)^{n-1}$:} consider the coordinate space $\mathbb K^{2n}$ endowed with the non-degenerate bilinear form whose Gram matrix with respect to the standard basis is 
\[G=
\begin{pmatrix}
0&\mathbbm 1_n\\
J_n(\lambda)&0
\end{pmatrix},\quad J_n(\lambda):=\begin{pmatrix}
\lambda&1&&&\\
&\lambda&1&&\\
&&\ddots&\ddots&\\
&&&\lambda&1\\
&&&&\lambda
\end{pmatrix}.
\]In this case, by Proposition \ref{26.06.17-1} the canonical operator is of the form
\[\kappa=G^{-1}\cdot G^T=\begin{pmatrix}
J_n(\lambda)^{-1}&0\\
0&J_n(\lambda)^T
\end{pmatrix},
\]and its Jordan form consists of two $n\times n$ blocks with nonzero inverse eigenvalues.
\end{enumerate}

The two types of space $U_n, W_n(\lambda)$ with $\lambda\neq(-1)^{n-1}$ are not isometrically isomorphic.
\end{teorema}
\esh

Together with Theorem \ref{13.07.17-5}, the following result gives a complete classification of all Mukai spaces over algebraically closed field of characteristic zero.

\bsh
\begin{teorema}[\cite{malcev} Chapter VI-VII, \cite{Go3, Go2.2}]\label{classcompleta}
Let $(V,\langle\cdot,\cdot\rangle)$ is a vector space over an algebraically closed field $\mathbb K$ of characteristic zero endowed with a non-degenerate bilinear form. If $f\in\End(V)$ is an isometry, then $V$ splits as a bi-orthogonal direct sum of subspaces $V_\lambda$, where
\begin{enumerate}
\item for $\lambda=\pm 1$, $V_\lambda$ is the root space
\[V_\lambda:=\bigoplus_{n\in\mathbb N}\ker\left(f-\lambda\mathbbm 1\right)^n,
\]and the restriction of $\langle\cdot,\cdot\rangle$ on $V_\lambda$ is non-degenerate;
\item for $\lambda\neq\pm 1$, the space $V_\lambda$ is the sum of isotropic root subspaces
\[\left(\bigoplus_{n\in\mathbb N}\ker\left(f-\lambda\mathbbm 1\right)^n\right)\oplus\left(\bigoplus_{n\in\mathbb N}\ker\left(f-\lambda^{-1}\mathbbm 1\right)^n\right),
\]and the restriction of $\langle\cdot,\cdot\rangle$ on $V_\lambda$ defines a non-degenerate pairing between these two subspaces.
\end{enumerate}
\end{teorema}
\esh

\subsection{Geometric case: the derived category $\mathcal D^b(X)$}\label{geocase} In previous sections, we have treated the general case of a $\mathbb K$-linear triangulated category $\mathscr D$. Now we consider the case of the bounded derived category of coherent sheaves on a projective complex variety $X$. In order to work with a $\Hom$-finite derived category, we assume that $X$ is \emph{smooth}: in this way each object is a \emph{perfect complex}, i.e. locally quasi-isomorphic to a bounded complex of locally free sheaves of finite rank on $X$.

The condition of existence of a full exceptional collection in $\mathcal D^b(X)$ imposes strict conditions on the topology and the geometry of $X$. The key property is a result of \emph{motivic decomposition} for the rational Chow motive of $X$. We refer the reader to Appendix \ref{chowmot} for basic notions on Chow motives.

\begin{defi} Let $X$ be a smooth projective variety over $\mathbb C$. We will say that the rational Chow motive of $X$,  denoted by $\frak h(X)_{\mathbb Q}\in{\rm CHM}(\mathbb C)_{\mathbb Q}$, is \emph{discrete} (or of \emph{Lefschetz type}) if it is polynomial in the Lefschetz motive $\mathbb L$, i.e. if it admits a decomposition as a direct sum
\[\frak h(X)_{\mathbb Q}\cong \bigoplus_{i=1}^n\mathbb L^{\otimes a_i},\quad a_i\in\left\{0,\dots, \dim_{\mathbb C}X\right\},
\]where by convention $\mathbb L^0:=\frak h({\rm Spec}(\mathbb C))_{\mathbb Q}$. The integer $n$ will be called \emph{length} of the motive $\frak h(X)_{\mathbb Q}$.
\end{defi}

\bsh
\begin{teorema}[\cite{gorchorlov}, \cite{marcotabu1}, \cite{bernabolo1}]\label{13.07.17-1}
If $X$ is a smooth projective variety over $\mathbb C$ with an exceptional collection in $\mathcal D^b(X)$, then the rational Chow motive $\frak h(X)_{\mathbb Q}$ is discrete. 
\end{teorema}
\esh

There exist many proofs in the literature of this fact, all differing in techniques. In \cite{gorchorlov} the statement was proved using $K$-motives. In \cite{marcotabu1} (see also \cite{marcotabu2}) a more general statement was proved (assuming that $X$ is a smooth proper Deligne--Mumford stack over ${\rm Spec}(\mathbb K)$, for a perfect field $\mathbb K$) using the connection between Chow and non-commutative motives discovered by M. Kontsevich (see \cite{tabu}). In the case of smooth projective varieties, we can deduce Theorem \ref{13.07.17-1} easily follows from the following result, essentially due to S. Kimura.

\bsh
\begin{teorema}[\cite{kimura}]\label{12.07.17-1}Let $X$ be a smooth projective variety over $\mathbb C$. The following conditions are equivalent:
\begin{enumerate}
\item the Grothendieck group $K_0(X)_\mathbb Q$ is a finite dimensional $\mathbb Q$-vector space;
\item the rational Chow motive of $X$ is \emph{discrete}.
\end{enumerate}
Furthermore, if these conditions hold true, the length of $\frak h(X)_{\mathbb Q}$ coincides with $\dim_{\mathbb Q}K_0(X)_{\mathbb Q}$.
\end{teorema}
\esh

\proof The proof of the fact that (1) implies (2) follows from the main result of \cite{kimura}. Conversely, if 
\[\frak h(X)_{\mathbb Q}\cong \bigoplus_{i=1}^n\mathbb L^{\otimes a_i},\quad a_i\in\left\{0,\dots, \dim_{\mathbb C}X\right\},
\]using the properties of the Lefschetz motive $\mathbb L$ we deduce that
\begin{align*}CH^r(X)_{\mathbb Q}&\cong \Hom_{\rm CHM(\mathbb C)_{\mathbb Q}}(\mathbb L^{\otimes r}, \frak h(X)_{\mathbb Q})\\
&\cong\bigoplus_{i=1}^n\Hom_{\rm CHM(\mathbb C)_{\mathbb Q}}(\mathbb L^{\otimes r}, \mathbb L^{\otimes a_i})\cong \mathbb Q^{N(r)},
\end{align*}
where $N(r):={\rm card}\left\{i\colon a_i=r\right\}$. Since the Chern character ${\rm ch}\colon K_0(X)_{\mathbb Q}\to CH^\bullet(X)_{\mathbb Q}$ is an isomorphism (not preserving the gradation), we conclude that $K_0(X)_{\mathbb Q}$ is finite dimensional.
\endproof

\bsh
\begin{cor}[\cite{marcotabu1}, \cite{minifolds}]\label{corollifexce}
Let $X$ be a smooth projective variety over $\mathbb C$ such that $K_0(X)$ is free of finite rank. Then:
\begin{enumerate}
\item X is of Hodge--Tate type, i.e. $h^{p,q}(X)=0$ if $p\neq q$.
\item The cycle maps $c_X^r\colon CH^r(X)_{\mathbb Q}\to H^{2r}(X,\mathbb Q)$ are isomorphisms.
\item The forgetful morphism $K_0(X)_{\mathbb Q}\to K^{\rm top}_0(X^{\rm an})_{\mathbb Q}$ is an isomorphism.
\item ${\rm Pic}(X)$ is free of finite rank and $c_1\colon {\rm Pic}(X)\to H^2(X,\mathbb Z)$ is an isomorphism.
\item $H_1(X,\mathbb Z)=0$.
\end{enumerate}
\end{cor}
\esh

\proof Since the Grothendieck group $K_0(X)$ is free and of finite rank, the conditions (1), (2) of Theorem \ref{12.07.17-1} hold true. 
By the universal property of the Chow motives, any Weil cohomological functor $H^\bullet$ with values in $\rm GrVect_{\mathbb Q}^{<\infty}$ factorizes through $\rm CHM(\mathbb C)_{\mathbb Q}$: hence, using the same notation of Theorem \ref{12.07.17-1} and its proof, we have
\[H^\bullet(X)\cong \bigoplus_{i=1}^n H^\bullet(\mathbb L)^{\otimes a_i},
\]and consequently 
\[H^{2r}(X)\cong \mathbb Q^{N(r)},\quad\text{since }H^{2i}(\mathbb L)\cong \begin{sistema}
\mathbb Q,\ i=2\\
\\
0,\ \text{otherwise}.
\end{sistema}
\]By taking the Hodge realization of the rational Chow motive $\frak h(X)_{\mathbb Q}$, we deduce point (1). Point (2) follows from the fact that ${\rm im}(c_X^r)\subseteq H^{r,r}(X)$ and that $CH^r(X)_{\mathbb Q}$ and $H^{2r}(X;\mathbb Q)$ have the same dimension $N(r)$. Statement (3) follows from the commutative diagram 
\[\xymatrix{
K_0(X)_{\mathbb Q}\ar[rr]^{\rm ch}\ar[d]&& CH^\bullet(X)_{\mathbb Q}\ar[d]\\
K^{\rm top}_0(X^{\rm an})_{\mathbb Q}\ar[rr]^{\quad\rm ch^{top}}&& H_{dR}^\bullet(X^{\rm an};\mathbb Q)
}
\]and from (2). Here the vertical arrows denote the natural forgetful morphisms, $\rm ch$ denotes the Chern character as defined in \cite{fulton_int}, and $\rm ch^{top}$ denotes the topological version of the Chern character. As shown in \cite{minifolds}, the freeness of $K_0(X)$ implies that ${\rm Pic}(X)$ is free. From the exponential long exact sequence we have that
\[\xymatrix{H^1(X,\mathcal O_X)\ar[r]&{\rm Pic}(X)\ar[r]^{c_1}& H^1(X,\mathbb Z)\ar[r]& H^2(X,\mathcal O_X),}
\]
and by point (1) we have that $h^{0,1}(X)=h^{0,2}(X)=0$. Hence the first Chern class map is an isomoprhism. It follows that Pic$(X)$ is of finite rank. For the last statement, by the Universal Coefficient Theorem we have a (non-canonical) isomorphism
\[H^2(X,\mathbb Z)\cong \mathbb Z^{\beta_2(X)}\oplus H_1(X,\mathbb Z)^{\rm tors},
\]and by (4) we deduce that $H_1(X,\mathbb Z)$ is a finitely generated torsion-free abelian group, and hence free. By point (1) we have that $h^{1,0}(X)=0$.
\endproof

\bsh
\begin{cor}[\cite{marcotabu1}, \cite{kuz5}]\label{13.07.17-2}If $X$ is a smooth projective variety over $\mathbb C$ such that $\mathcal D^b(X)$ admits a full exceptional collection, then $X$ is of Hodge--Tate type, and the length of the collection is equal to $\sum_{p}h^{p,p}(X)$.
\end{cor}
\esh

\begin{oss}
The proof of Corollary \ref{13.07.17-2} which appears in \cite{kuz5} is not based on motivic decomposition techniques, rather on an \emph{additivity property of Hochschild homology} (\cite{kuz2.2}). If $X$ admits a semiorthogonal decomposition $\langle\mathscr A_i\rangle_{i=1}^n$ of $\mathcal D^b(X)$, then Hochschild homology admits the decomposition
\begin{equation}\label{13.07.17-3}HH_\bullet(X)\cong\bigoplus_{i=1}^nHH_\bullet(\mathscr A_i).
\end{equation}Using the following properties of Hochschild homology 
\begin{itemize}
\item $HH_\bullet ({\rm Spec}(\mathbb K))\cong \mathbb K$ (Example 1.17 of \cite{kuz5}), 
\item and Hochschild-Kostant-Rosenberg Theorem
\[HH_\bullet(X)\cong\bigoplus_{q-p=k}H^q(X,\Omega^p_X),
\]
\end{itemize}
by Proposition \ref{propbondal} and the decomposition \eqref{13.07.17-3} one easily concludes.
\end{oss}

The following result shows that, in the isometric classification of Mukai structures, for all varieties $X$ the vector spaces $(K_0( X)_\mathbb C,\chi(\cdot,\cdot))$ correspond just to one of the possible cases, namely bi-orthogonal sums of $U_n$-type spaces. It is based on the \emph{dévissage property} of coherent sheaves on $X$ (see e.g. \cite{shafa2}, Section II.6.3.3 and \cite{chrissginz}, Section 5.9).

\bsh
\begin{prop}[\cite{chrissginz}]\label{13.07.17-4}
Let $X$ be a smooth projective variety over $\mathbb C$, and let $E$ be a rank $d$ vector bundle on $X$. The endomorphism
\[\varphi_E\colon K_0(X)_{\mathbb C}\to K_0(X)_{\mathbb C}\colon [\mathscr F]\mapsto [\mathscr F\otimes (E-d\cdot \mathcal O_X)]
\]is nilpotent. In particular, 
\[\varphi_E^{\dim_\mathbb CX+1}=0.
\]
\end{prop}
\esh

\bsh
\begin{cor}\label{corcancappa}Let $X$ be a smooth projective variety over $\mathbb C$. The canonical morphism $\kappa\colon K_0(X)_{\mathbb C}\to K_0(X)_{\mathbb C}$ is of the form
\[\kappa=(-1)^{\dim_\mathbb CX}\mathbbm 1+ M,\quad \text{with }M\text{ nilpotent}.
\]Hence, the Mukai vector space $(K_0(X)_{\mathbb C},\chi)$ is isomorphic to a direct sums of irreducible Mukai spaces of type $U_n$. In particular, a necessary condition for the irreducibility of $(K_0(X)_{\mathbb C},\chi)$ is that 
\[\sum_{j=0}^{\dim_{\mathbb C}X}\beta_{j}(X)\equiv \dim_\mathbb CX+1\quad (\rm mod\ 2).
\]
\end{cor}
\esh

\proof
The canonical morphism $\kappa$ is defined by $\kappa([\mathscr F])=(-1)^{\dim_\mathbb CX}[\mathscr F\otimes \omega_X]$. By Proposition \ref{13.07.17-4}, the morphism $\kappa+(-1)^{\dim_{\mathbb C}X+1}\mathbbm 1$ is nilpotent. The last two statements follow from the classification of indecomposable Mukai spaces, Theorem \ref{13.07.17-5}, and from Theorem \ref{classcompleta}.
\endproof

Hence, in the isometric classification, the case of Projective Spaces is of particular importance as the following results show.

\bsh
\begin{teorema}[\cite{Go3, Go2.2}]\label{04.08.17-1}
The Mukai spaces $(K_0(\mathbb P^{k-1})_\mathbb C,\chi)$ are indecomposable of type $U_k$. Their isometry group $${\rm Isom}_\mathbb C(K_0(\mathbb P^{k-1}_\mathbb C)_\mathbb C,\chi)$$ has two connected components. The identity component is a unipotent abelian algebraic group of dimension\footnotemark $\ [\frac{k}{2}]$.
\end{teorema}
\esh
\footnotetext[24]{Here $[x]$ denotes the integer part of $x\in\mathbb R$.}

Notice indeed that the necessary condition of Corollary \ref{corcancappa} is satisfied in the case of complex Projective Spaces. From Corollary \ref{corcancappa} and Theorem \ref{04.08.17-1}, we deduce the following result.

\bsh
\begin{cor}\label{04.08.17-2}
Let $X$ be a smooth projective variety. The identity component of the isometry group
\[{\rm Isom}_\mathbb C(K_0(X)_\mathbb C,\chi)_0
\]is unipotent and abelian.
\end{cor}
\esh

\newpage
\section{The Main Conjecture}\label{chmainconj}

In the occasion of the 1998 ICM in Berlin, the second author formulated a conjecture relating two apparently different and unrelated aspects of the geometry of Fano varieties, namely their enumerative geometry (quantum cohomology)  and their derived category of coherent sheaves.

\subsection{Original version of the Conjecture and known results}

Let us recall the original statement of the Conjecture formulated by the second author.

\bsh
\begin{conj}[\cite{dubro0}]\label{originalconj}Let $X$ be a Fano variety.
\begin{enumerate}
\item The quantum cohomology $QH^\bullet(X)$ is semisimple if and only if the category $\mathcal D^b(X)$ admits a full exceptional collection $(E_1,\dots, E_n)$.
\item The Stokes matrix $S$, computed with respect to a fixed oriented line $\ell$ admissible for the system \eqref{16.07.17-2} and in the lexicographical order with respect to $\ell$, is equal to the Gram matrix of the Grothendieck--Euler--Poincaré product with respect to a full exceptional collection in $\mathcal D^b(X)$,
\[S_{ij}:=\chi(E_i,E_j).
\]
\item The central connection matrix, connecting the solution $Y^{(0)}_{\rm right}$ of Theorem \ref{16.07.17-3} with the topological-enumerative solution $Y_{\rm top}:=\Psi\cdot Z_{\rm top}$ of Proposition \ref{18.07.17-2}, is of the form
\[C=C'\cdot C'',
\]where the columns of $C''$ are the components of the Chern characters ${\rm ch}(E_i)$, and the matrix $C'$ represents an endomorphism of $H^\bullet(X,\mathbb C)$ commuting with $c_1(X)\cup(-)\colon H^\bullet(X,\mathbb C)\to H^\bullet(X,\mathbb C)$.
\end{enumerate}
\end{conj}
\esh

Over the years, the Conjecture \ref{originalconj} has been partially verified in specific cases by several authors, as follows.

\begin{enumerate}
\item In \cite{guzzetti1} the third author proved point (2) of Conjecture \ref{originalconj} for projective spaces. He performed a detailed analysis of the action of the braid group on the set of monodromy data. By this action, the Stokes matrix computed at a point of the small quantum cohomology was transformed into a Stokes matrix (associated with another chamber of the manifold) which coincides with the  {\it inverse} of the Gram matrix $\chi(E_i,E_j)$, where $\{E_j=\mathcal{O}(j-1)\}_{j=1}^{k}$ is the Beilinson exceptional collection. See also Remark \ref{04.08.17-4}.

\item The results of \cite{guzzetti1} were recovered by S. Tanabé in \cite{tanabe}, who showed how to calculate the Stokes matrices of quantum cohomologies of projective spaces in terms of a certain hypergeometric group. Furthermore, in \cite{morvdp}, J. A. Cruz Morales and M. Van der Put showed another method to obtain the same results of \cite{guzzetti1} for projective spaces, and also for the case of weighted projective spaces, using multisummation techniques.

\item In \cite{ueda1} K. Ueda extended the results of \cite{guzzetti1} to all complex Grassmannians. His proof relies on a conjecture of K. Hori and C. Vafa (\cite{horivafa}), rigorously proved by A. Bertram, I. Ciocan-Fontanine, and B. Kim (\cite{BCFK}), relating quantum cohomology of Grassmannians with quantum cohomology of projective spaces. The analysis of the action of the braid group is not treated. Note that in \cite{ueda1} the  \emph{phenomenon of coalescence} for the isomonodromic system \eqref{16.07.17-2} is not recognised: a priori, the monodromy data at points of small quantum cohomology of  Grassmannians with  coalescence could be not  well defined. A rigorous analysis of this point has been developed in \cite{CDG0}, and adapted to the geometry of Frobenius manifolds in \cite{CDG}.
\item In \cite{ueda2} K. Ueda proved point (2) of Conjecture \ref{originalconj} for cubic surfaces, using a toric degeneration of the surfaces and A. Givental's mirror results (\cite{givemirtor}).
\item In \cite{bm4} A. Bayer and Yu. I. Manin proved point (1) of Conjecture \ref{originalconj} for all del Pezzo surfaces.
\item Out of the 106 deformation classes of smooth Fano threefolds (see \cite{fano31,fano32}, \cite{shigemori1,shigemori2}), only 59 satisfy the condition of vanishing odd cohomology, necessary for the semisimplicity of the quantum cohomology.  In \cite{ciolli, cio5}, G. Ciolli proved validity of point (1) of Conjecture \ref{originalconj} for 36 out of these 59 families.
\item A. Bayer proved in \cite{bay4} that the family of varieties for which point (1) of Conjecture \ref{originalconj} holds true is closed under blow-ups at any number of points. Furthermore, Bayer also suggested to drop any reference to the condition of being Fano in the statement of  Conjecture \ref{originalconj}. No explicit result is available in the non-Fano case for points (2)-(3) of Conjecture \ref{originalconj}. We plan to address this problem in a future publication \cite{inprep}.
\item The results of Y. Kawamata \cite{kawa1, kawa2, kawa3} confirm validity of point (1) of Conjecture \ref{originalconj} for projective toric manifolds.
\item In \cite{golyshev} V. Golyshev proved validity of point (2) of Conjecture \ref{originalconj} for minimal Fano three-folds, i.e. with minimal cohomology 
\[H^{2k+1}(X,\mathbb Z)\cong 0,\quad H^{2k}(X,\mathbb Z)\cong\mathbb Z.
\]
\item In \cite{iwakitaka}, K. Iwaki and A. Takahashi proved validity of point (2) of Conjecture \ref{originalconj} for a class of orbifold projective lines $\mathbb P^1_A$  whose quantum cohomology is known to be semisimple, being isomorphic to the Frobenius manifold constructed from the theory of primitive forms for the polynomial $f_A(x_1,x_2,x_3):=x_1^{a_1}+x_2^{a_2}+x_3^{a_3}-x_1x_2x_3$. Their proof is based on Kontsevich's Homological Mirror Symmetry, an equivalence of triangulated categories between $\mathcal D^b(\mathbb P^1_A)$ and the directed Fukaya category $\mathcal {DF}uk(\mathbb C^3,f_A)$. See Section \ref{KHMS} for more details.
\item In \cite{pech} C. Pech showed that the (small) quantum cohomology of odd symplectic Grassmannians of lines $IG(2,2n+1)$ is semisimple. These results, together with the results of A. Kutznetsov \cite{kuz6} adpated, \emph{mutatis mutandis}, to the odd case, confirm the validity of the point (1) of Conjecture \ref{originalconj}.

\item In \cite{gms}, S. Galkin, A. Mellit and M. Smirnov proved validity of point (1) of Conjecture \ref{originalconj} for the symplectic isotropic Grassmannian $IG(2,6)$. The importance of this result is due to the fact that it underlines the need of considering the whole big quantum cohomology for the formulation of the conjecture, the small quantum locus being contained in the caustic (recall Definition \ref{17agosto2018-1}). This result has been generalized for all symplectic isotropic Grassmannian $IG(n,2n)$: on the one hand it is known that these Grassmannians admit full exceptional collections (\cite{kuz6}, \cite{samokhin}), on the other hand it has been proved by N. Perrin that their (big) quantum cohomology is generically semisimple (see \cite{perrin}). See also \cite{cmmps}.
\end{enumerate}

\subsection{Relationships with Zaslow's Conjecture for solutions of $tt^*$-equations} Speculations on some relationships between monodromy phenomena of differential equations and objects in derived categories already appeared in the literature some years before Conjecture \ref{originalconj}. In this regard, it is particularly worth mentioning the inspiring paper \cite{zas} by E. Zaslow. 

Together with the system of linear differential equations \eqref{06.07.17-2}, associated with any semisimple Frobenius manifold there is a second isomonodromic problem, formulated in \cite{dubro93}, and arising from the work of S. Cecotti and C. Vafa on \emph{topological-antitopological fusion} (\cite{cecvafa}). From an analytical point of view, the $tt^*$-equations can be interpreted as the isomonodromicity conditions of a linear differential equation on $\mathbb P^1_{\mathbb C}$ with \emph{two irregular singularities} (of Poincaré rank 1) at both $z=0$ and $z=\infty$ (see equation 2.18 of \cite{dubro93}). In order to suitably describe the monodromy of their solutions it is then necessary to introduce Stokes matrices and a central connection matrix. For an extended analytical asymptotical study of certain cases of $tt^*$-equations, reduced to the Toda lattice (with opposite sign), we refer the reader to the detailed papers \cite{guestits1,guestits2,guestits3}, \cite{guestlin1,guestlin2}. 
 
Remarkably, Cecotti and Vafa predicted existence of solutions of the $tt^*$-equations for which the corresponding Stokes matrices have integer entries. Basing on this fact, a  classification of two-dimensional $N=2$ supersymmetric field theories was subsequently developed (\cite{cecotti-vafa}). From a physical point of view, indeed, given a $N=2$ massive theory (i.e. whose ``chiral ring'' is a semisimple Frobenius algebra)  the entries of the Stokes matrices can be interpreted as the numbers of solitons of minimal energy (\emph{Bogomolnyi solitons}) between the vacua states in the infrared limit. The idea of the classification is based on the information about the number of vacua states and of solitons between them. See Section 4 and 6 of \cite{cecotti-vafa} for further details.

Basing on a previous observation of M. Kontsevich, in \cite{zas} E. Zaslow conjectured a similar relationship between the Stokes phenomenon arising from the 
$tt^*$-isomonodromic problem attached to the topological $\sigma$-model with target a Kähler manifold $X$ and the Gram matrices associated with objects of 
(full) exceptional collections in the derived categories of coherent sheaves $\mathcal D^b(X)$. This is  strictly related to Conjecture \ref{originalconj}, because 
 the Stokes matrices of both the system \eqref{06.07.17-2} and of the $tt^*$-problem are believed to be \emph{equal}. Such an equality was already conjectured 
 by C. Vafa\footnote{Private communication to the second author.}. For the $tt^*$-Toda case, in the  papers \cite{guestlin1}, \cite{guestits3}, the two 
 systems of differential equations (very different from an analytical point of view!) are shown to be related through an Iwasawa factorization of loop groups 
 (see Chapter 12 of \cite{guestbook}), and  the equality of Stokes data  for the $tt^*$-Toda equations is proved (in \cite{guestits3} details 
  are worked out for the $4\times 4$ system, while  in \cite{guestho} for the $n\times n$ case). We expect that the Iwasawa factorisation 
  method of \cite{guestits3}  may equally  work in order to prove the equality of Stokes matrices in the general case, not only for the  $tt^*$-Toda equations.  A similar approach was already used
   by A. Bobenko and A. Its in \cite{bobits}, where the reader can find explicit formulae relating the monodromy data of (the linear system associated with) 
   Painlevé III and its data associated by the approach of J. Dorfmeister, F. Pedit, and H. Wu (\cite{dpw}), based on the Iwasawa decomposition of a twisted 
   loop group of $SL(2,\mathbb C)$. To the best  of our knowledge, these are the only references in which such a relationships between the two Stokes phenomena is discussed with some analytical details.


\subsection{Gamma classes, graded Chern character, and morphisms \textnormal{$\textcyr{D}^\pm_X$}} Let $X$ be a smooth projective variety of (complex) dimension $d$ with odd-vanishing cohomology, $V$ be a complex vector bundle of rank $r$ on $X$, and let $\delta_1,\dots,\delta_r$ be the Chern roots of the bundle $V$, so that
\[c_k(V)=\sigma_k(\delta_1,\dots,\delta_r),\quad 1\leq k\leq r
\]where $\sigma_k$ is the $k$-th elementary symmetric polynomial.
Starting from the Taylor series expansion of the functions $\Gamma(1\pm t)$ near $t=0$, and applying the Hirzebruch's construction of characteristic classes, we can define two characteristic classes $\widehat\Gamma^\pm(V)\in H^\bullet(X,\mathbb C)$ by
\[\widehat\Gamma^\pm(V):=\prod_{j=1}^r\Gamma(1\pm\delta_j).
\]In particular we will denote by $\widehat\Gamma^\pm(X)$ the characteristic class $\widehat\Gamma^\pm(TX)$.\newline

For any object $E\in\Ob \mathcal D^b(X)$ we define a \emph{graded} version of (the Grothendieck's definition of) the Chern character: being $X$ smooth, the object $E$ is isomorphic in $\mathcal D^b(X)$ to a (bounded) complex of locally free sheaves $F^\bullet$. We thus define
\[\operatorname{Ch}(E):=\sum_j(-1)^j\operatorname{Ch}(F^j),
\]
where
\[\operatorname{Ch}(F^j):=\sum_{h}e^{2\pi i \alpha_h},\quad\text{the }\alpha_h\text{'s being the Chern roots of }F^j.
\]The definition is well posed, since it can be easily shown to be independent of the bounded complex $F^\bullet$ of locally free sheaves.

Let us now define two morphisms $\textcyr{D}^\pm_X\colon K_0(X)\to H^\bullet(X,\mathbb C)$ given by
\[ \textcyr{D}^\pm_X(E):=\frac{i^{\overline{d}}}{(2\pi)^\frac{d}{2}}\widehat\Gamma^\pm(X)\cup \exp(\pm\pi i c_1(X))\cup\operatorname{Ch}(E), \]
where $\overline{d}\in\left\{0,1\right\}$ is the residue class $d\text{ mod}(2)$.

\subsection{Refined statement of the Conjecture}$\quad $
\bsh
\begin{conj}\label{congettura}
Let $X$ be a smooth Fano variety of Hodge--Tate type, namely such that $h^{p,q}(X)=0$ for any $p\neq q$. 
\begin{enumerate}
\item The quantum cohomology $QH^\bullet(X)$ is semisimple if and only if there exists a full exceptional collection in the derived category of coherent sheaves $\mathcal D^b(X)$. 

\item If $QH^\bullet(X)$ is semisimple, then for any oriented line $\ell$ (of slope $\phi\in[0;2\pi[$) in the complex plane there is a correspondence between $\ell$-chambers and founded helices, i.e. helices with a marked foundation, in the derived category $\mathcal D^b(X)$. 

\item The monodromy data computed in a $\ell$-chamber $\Omega_\ell$, in the lexicographical order, are related to the following geometric data of the corresponding exceptional collection $\frak E_\ell=(E_1,\dots, E_n)$ (the marked foundation):
\begin{enumerate}
\item the Stokes matrix is equal to the \emph{inverse} of the Gram matrix of the Grothendieck-Poincaré-Euler product on $K(X)_{\mathbb C}$, computed with respect to the exceptional basis $([E_i])_{i=1}^n$
\[ S^{-1}_{ij}=\chi(E_i, E_j);
\]
\item the Central Connection matrix $C\equiv C^{(0)}$,  connecting the solution $Y_{\rm right}^{(0)}$ of Theorem \ref{16.07.17-3} with the topological-enumerative solution $Y_{\rm top}=\Psi\cdot Z_{\rm top}$ of Proposition \ref{18.07.17-2}, coincides with the matrix associated with the $\mathbb C$-linear morphism 
\[\textnormal{\textcyr{D}}^-_X\colon K_0(X)_{\mathbb C}\to H^\bullet(X,\mathbb C)\colon E\mapsto \frac{i^{\bar d}}{(2\pi)^\frac{d}{2}}\widehat\Gamma^-_X\cup\exp(-\pi i c_1(X))\cup{\rm Ch}(E),
\]where $d=\dim_\mathbb C X$, and $\bar d$ is the residue class $d $ {\rm (mod 2)}. The matrix is computed with respect to the exceptional basis $([E_i])_{i=1}^n$ and any pre-fixed basis $(T_\alpha)_{\alpha}$ in cohomology (see Section \ref{notations}).
\end{enumerate}
\end{enumerate}
\end{conj}
\esh

\begin{oss} Let $X$ be a smooth projective variety with semisimple quantum cohomology. From the original Conjecture \ref{originalconj}, in \cite{dubro0} it was conjectured that there exists  an atlas of $QH^\bullet(X)$ whose charts, denoted $Fr(S,C)$, are expected to be in one-to-one correspondence with exceptional collections in $\mathcal D^b(X)$. Point (2) of Conjecture \ref{congettura} clarifies this. In order to have such a correspondence, each of the charts discussed in \cite{dubro0} should cover a single $\ell$-chamber. The correspondence with exceptional collections is not one-to-one, since two foundations of the same helix, obtained one from another by iterated applications of the Serre functor (or its inverse), are associated with monodromy data computed with respect to other solutions $Y^{(k)}_{\rm left/right}$ of Theorem \ref{16.07.17-3}. In other words, the choice of the foundation of the helix corresponds to the choice of the branch of the logarithmic factors $z^\mu z^{c_1(X)\cup(-)}$ of $Y_0(z)$.
\end{oss}

\begin{oss}\label{teoblock}
The Main Theorem of \cite{CDG} implies a constraint on the kind of exceptional collections associated with the monodromy data in a neighborhood of a semisimple coalescing point of the quantum cohomology. If the eigenvalues $u_i$'s coalesce, at some semisimple point $t_0\in QH^\bullet(X)$, to $s< n$ values $\lambda_1,\dots,\lambda_s$ with multiplicities $p_1,\dots, p_s$ (with $p_1+\dots +p_s=n$, here $n$ is the sum of the Betti numbers of $X$, equivalently, it is the dimension of $QH^{\bullet}(X)$ as a Frobenius manifold), then the corresponding monodromy data can be expressed in terms of Gram matrices and characteristic classes of objects of a full $s$-block exceptional collection of type $(p_1,\dots,p_s)$ (see Definition \ref{blockcoll}). Let $\mathcal U_{\epsilon_1}(t_0)$ be a  polydisc centered at $t_0$, of poly-radius $\epsilon_1$ suffciently small as in the statement of Theorem 4.1 of \cite{CDG}. Once we fix 
\begin{itemize}
\item an oriented line $\ell$, 
\item a continuous labeling of canonical coordinates on $\mathcal U_{\epsilon_1}(t_0)$,
\item a continuous branch of the $\Psi$-matrix on $\mathcal U_{\epsilon_1}(t_0)$,
\end{itemize}
the monodromy invariants $(S,C)$ computed in each $\ell$-chamber intersecting $\mathcal U_{\epsilon_1}(t_0)$ are constant and the same, being equal to the ones computable at the semisimple coalescence point $t_0$. Notice that on $\mathcal U_{\epsilon_1}(t_0)$ we have $\prod_jp_j!$ possible $\ell$-triangular orders (Definition \ref{30luglio2016-2}), each of which is $\ell$-lexicographical in exactly one $\ell$-chamber intersecting $\mathcal U_{\epsilon_1}(t_0)$. By Proposition \ref{propmutblock} the exceptional collections attached to these $\ell$-chambers, as in Conjecture \ref{congettura}, differ just by a permutation of the objects of each single block.
\end{oss}

\begin{oss}
\label{04.08.17-4}
 One of the main differences between the statements of Conjecture \ref{originalconj} and Conjecture \ref{congettura} is the point concerning the Stokes matrix $S$. The identification of $S$ with the inverse of the Gram matrix is forced by the Grothendieck--Hirzebruch--Riemann--Roch Theorem and the constraints of monodromy data, as it will be evident from Proposition \ref{18.07.17-1} and Corollary \ref{corcongsempl}. 
However, in case of projective spaces, the necessity of considering the inverse of the Gram matrix is not evident. In \cite{guzzetti1}, an explicit braid was constructed relating  the Stokes matrix $S$ computed at a point of the small quantum cohomology, to another Stokes matrix $S^\beta$  (associated with another $\ell$-chamber). The latter was proved to  coincide with the {\it inverse} of the  Gram matrix of the  Beilinson exceptional collection $\frak B=(\mathcal O(i))_{i=0}^{k-1}$. 
Then, using the identity
\[[S^{-1}]^\beta=PS^TP,\quad \beta=\beta_{12}(\beta_{23}\beta_{12})(\beta_{34}\beta_{23}\beta_{12})\dots(\beta_{k-1,k}\dots\beta_{12}),\quad P^\alpha_\beta=\delta_{\alpha+\beta,1+k},
\]where $S$ is any $k\times k$ Stokes matrix (see \cite{zas}), together with the numerical coincidence
\[PG^TP=G,\quad \text{for }G_{ab}:=\chi(\mathcal O(a-1),\mathcal O(b-1))=\binom{k-1+b-a}{b-a},\quad 1\leq a,b\leq k,
\] it was proved  that $G$ (or equivalently $S$) and its inverse are in the same orbit under the action of the braid group.  In other words,   the result of  \cite{guzzetti1} says that there  is an $\ell$-chamber where the Stokes matrix coincides with the  Gram matrix of the  Beilinson exceptional collection $\frak B$.

For a generic smooth projective variety $X$ the following questions naturally arise.\newline

{\bf Problem 1. }If point (3.a) of Conjecture \ref{congettura} holds true, so that $S^{-1}=G$ with $G$ Gram matrix of $\chi$ computed with respect to an exceptional basis $\varepsilon$, does the Stokes matrix $S=G^{-1}$ coincide with the Gram matrix of another exceptional basis $\varepsilon'$?\newline 

{\bf Problem 1.bis. }Given a full exceptional collection $\frak E$ in $\mathcal D^b(X)$ with associated Gram matrix $G$, does it exist another exceptional collection $\frak E'$ with associated Gram matrix $G^{-1}$?\newline

Problems 1 and 1.bis admit affermative solutions. By Proposition \ref{propbraidmukai} and Proposition \ref{propdualexccoll} it immediately follows that the Gram matrix corresponding to both $(\frak E^*)^\vee$ and $^\vee(\frak E^*)$ is equal to $G^{-1}$. This guarantees, for example, that for all other cases of smooth projective varieties for which point (2) of Conjecture \ref{originalconj} was proved (e.g. \cite{ueda2,ueda1}, \cite{golyshev}), point (3.a) of Conjecture \ref{congettura} holds true. We still do not know the answer to the following  question.\newline

{\bf Problem 2.}  If Problems 1/1.bis admit affermative solutions, are $\varepsilon/\frak E$ and $\varepsilon'/\frak E'$ in the same orbit under the action of the braid group?\newline

Finally, let us notice that the solutions $\varepsilon'/\frak E'$ to Problem 1 and 1.bis have been explicitly constructed basing on the geometrical nature of the triangulated category under consideration, for which a further duality operation is available (Section \ref{geomdualcoll} and Section \ref{geomdualbas}). In a more general context, we do not know the answer to the following  question.\newline

{\bf Problem 3. }Let $\frak E$ be a full exceptional collection in a triangulated category $\mathscr D$ with associated Gram matrix $G$. Does there  exist another full exceptional collection $\frak E'$ in $\mathscr D$ with Gram matrix $G^{-1}$? If yes, are $\frak E$ and $\frak E'$ in the same orbit under the action of the braid group?
\end{oss}

\begin{oss}\label{sasha}
During the final revision of the present work, a very interesting paper \cite{tar-var} by V. Tarasov and A. Varchenko appeared\footnote{The first author thank A. Varchenko for very useful and instructive discussions.}. The main focus of \cite{tar-var} is the study of equivariant quantum cohomologies of the cotangent bundles $T^*\mathcal F_\lambda$ of partial flag varieties $\mathcal F_\lambda$, with $\lambda\in\mathbb Z^N_{\geq 0}$ such that $|\lambda|=n$, parametrizing chains of subspaces $$0=F_0\subset F_1\subset\dots\subset F_{N-1}\subset F_N\equiv\mathbb C^n,$$ with $\dim F_i/F_{i-1}=\lambda_i$ for $i=1,\dots, N$. The varieties $T^*\mathcal F_\lambda$ being not projective,  a good definition of their quantum cohomology is available only in the equivariant case, $QH^\bullet_T(T^*\mathcal F_\lambda)$. Here the torus $T$ is defined as $T:=A\times \mathbb C^*$,  where $A\subseteq GL_n(\mathbb C)$ denotes the subgroup of diagonal matrices which acts on $\mathbb C^n$, and hence on $T^*\mathcal F_\lambda$, whereas $\mathbb C^*$ acts by multiplication in each fiber. The quantum multiplication by divisors on $QH^\bullet_T(T^*\mathcal F_\lambda)$ is described in \cite{mo}, and an analogue family of the deformed connections \eqref{14.11.18-1} is defined (see also \cite{bmo}): in the notations of \cite{tar-var}
\[\nabla^{\rm quant}_{\lambda, q,\kappa, i}:=\kappa{q}_i\frac{\partial}{\partial{q}_i}-D_i*_{ q},\quad i=1,\dots, N,
\]where
\begin{itemize}
\item $(q_1,\dots, q_N)\in(\mathbb C^*)^N$ represent the parameters of deformation corresponding to the small equivariant quantum cohomology,
\item $D_i*_{ q}\colon QH^\bullet_T(T^*\mathcal F_\lambda)\to QH^\bullet_T(T^*\mathcal F_\lambda)$ denotes the operators of quantum multiplication by the divisors 
$$D_i:=\sum_{j=1}^{\lambda_i}\gamma_{i,j},\quad i=1,\dots, N,$$ where $\left\{\gamma_{i,j}\right\}_{j=1}^{\lambda_i}$ is the set of Chern roots of the vector bundle $E_i$ on $\mathcal F_\lambda$ with fiber $F_{i}/F_{i-1}$,
\item $\kappa\in\mathbb C^*$ denotes the spectral parameter of deformation of the flat connection (analogue to our parameter $z$).
\end{itemize}
The system of equations defining flat sections of $\nabla^{\rm quant}_{\lambda, q,\kappa, i}$ are called \emph{equivariant quantum differential equations} of $T^*\mathcal F_\lambda$. Among other results, Tarasov and Varchenko proved that the leading term of the asymptotic expansion of $q$-hypergeometric solutions of this system of equations can be written as the equivariant $\Gamma$-class of $T^*\mathcal F_\lambda$ multiplied by exponentials of the equivariant first Chern classes of the vector bundles $E_i$'s, namely the characteristic classes
\[\widehat{\Gamma}^+_{T^*\mathcal F_\lambda}\cup\prod_{j=1}^N\exp\left(\pi i (\lambda_j-n)c_1(E_j)\right)q_j^{c_1(E_j)}.
\]
See Appendix B of \cite{tar-var} for more details. The study of relationships between Conjecture \ref{congettura} and the work of \cite{tar-var} in the equivariant setting will be the object of a future research project. 
\end{oss}

\bsh
\begin{prop}\label{18.07.17-1}Let $X$ be a smooth projective variety of complex dimension $d$.
\begin{enumerate}
\item Let $E,F\in\Ob\mathcal D^b(X)$. Then 
\[\textnormal{\textcyr{D}}^\pm_X(E)\cup\textnormal{\textcyr{D}}^\mp_X(F)=\frac{(-1)^{\overline{d}}}{(2\pi)^d}\operatorname{Td}(X)\cup\operatorname{Ch}(E)\cup\operatorname{Ch}(F)\cup e^{-i\pi c_1(X)},
\]where $\operatorname{Td}(X)\in H^\bullet(X,\mathbb C)$ is the graded Todd characteristic class
\[\operatorname{Td}(X):=\prod_{j=1}^d\frac{2\pi i\delta_j}{1-e^{-2\pi i\delta_j}},
\]where $\delta_1,\dots,\delta_d$ are the Chern roots of the tangent bundle $TX$.

\item Let us naturally identify the tangent bundle $TQH^\bullet(X)$ with $H^\bullet(X,\mathbb C)$. Then 
for any $E\in\Ob\mathcal D^b(X)$ we have
\[e^{-i\pi\mu}(\textnormal{\textcyr{D}}^\pm_X(E))=i^d\textnormal{\textcyr{D}}^\mp_X(E^*),
\]where $\mu\in\End(H^\bullet(X,\mathbb C))$ is the grading operator defined in \eqref{06.07.17-3}.
\item If $\mu$ is the grading operator defined as before, $E,F\in\Ob\mathcal D^b(X)$ then
\[\int_Xe^{-i\pi\mu}(\textnormal{\textcyr{D}}^\pm_X(E))\cup e^{i\pi c_1(X)}\cup\textnormal{\textcyr{D}}^\pm_X(F)=\chi(E,F).
\]
\end{enumerate}
\end{prop}
\esh

\proof
From the well known relation
\[\Gamma(z)\Gamma(1-z)=\frac{\pi}{\sin (\pi z)}
\]we get
\[\Gamma(1+z)\Gamma(1-z)=\frac{2\pi i z}{1-e^{-2\pi i z}}e^{-i\pi z}.
\]Thus, 
\[\textnormal{\textcyr{D}}^+_X(E)\cup\textnormal{\textcyr{D}}^-_X(F)=\frac{(-1)^{\overline{d}}}{(2\pi)^d}\left(\prod_{j=1}^d\frac{2\pi i \delta_j}{1-e^{-2\pi i \delta_j}}e^{-i\pi \delta_j}\right)\cup\operatorname{Ch}(E)\cup\operatorname{Ch}(F),
\]and we conclude the proof of (1) since $c_1(X)=\sum_j\delta_j$. For (2) notice that if $\phi\in H^\bullet(X,\mathbb C)$, $\phi=\sum_p\phi_p$ with $\phi_p\in H^{2p}(X,\mathbb C)$ then
\begin{equation}\label{expmu}e^{-i\pi\mu}(\phi)=i^d\sum_p(-1)^p\phi_p,
\end{equation}
and one easily concludes. For the last point (3), we can apply (1),(2) and the Grothendieck--Hirzebruch--Riemann--Roch Theorem as follows
\begin{align*}
\int_Xe^{-i\pi\mu}(\textnormal{\textcyr{D}}^\pm_X(E))\cup e^{i\pi c_1(X)}\cup\textnormal{\textcyr{D}}^\pm_X(F)&=\\
i^d\int_X\textnormal{\textcyr{D}}^\mp_X(E^*)\cup e^{i\pi c_1(X)}\cup\textnormal{\textcyr{D}}^\pm_X(F)&=\quad\quad \text{(by (2))}\\
\frac{(-1)^{\overline{d}}}{(2\pi)^d}i^d\int_X\operatorname{Td}(X)\cup\operatorname{Ch}(E^*)\cup\operatorname{Ch}(F)&=\quad\quad\text{(by (1))}\\
\frac{(-1)^{\overline{d}}}{(2\pi)^d}i^d(2\pi i)^d\chi(E,F)&\quad\quad\text{(by GHRR).}
\end{align*}
\endproof

\bsh
\begin{cor}\label{corcongsempl}
Let $X$ be a Fano smooth projective variety for which points {\rm (3b)} of the Conjecture \ref{congettura} holds true. Then also point {\rm (3a)} holds true.
\end{cor}
\esh

\proof The Stokes and central connection matrices must satisfy the constraint, which is equivalent to
\[(e^{-i\pi \mu}C)^T\eta e^{i\pi R}C=S^{-1},
\] with $R=c_1(X)\cup(-)$. By point (3) of Proposition \ref{18.07.17-1} we conclude.
\endproof

\bsh
\begin{teorema}\label{03.08.17-1} Let $X$ be a smooth Fano variety of Hodge--Tate type for which Conjecture \ref{congettura} holds true. Then, all admissible operations on the monodromy data have a geometrical counterpart in the derived category $\mathcal D^b(X)$, as summarized in Table \ref{enrichconjintro} at the end of the Introduction. In particular, we have the following:
\begin{enumerate}
\item Mutations of the monodromy data $(S,C)$ correspond to mutations of the exceptional basis.

\item The monodromy data computed with respect to other solutions $Y^{(k)}_{\rm left/right}$, i.e. $(S, C^{(k)})$, are associated, as in the points {\rm (3a)-(3b)} of Conjecture \ref{congettura}, with different foundations of the helix, related to the marked one by an iterated application of the Serre functor $(\omega_X\otimes - )[\dim_{\mathbb C}X]\colon \mathcal D^b(X)\to  \mathcal D^b(X)$. 
\item The group ${\mathcal C}_0(X)$ is isomorphic to a subgroup of the identity component of the isometry group ${\rm Isom}_{\mathbb C}(\chi)$: more precisely, the morphism
\begin{equation}\label{isomgruppi}{\mathcal C}_0(X)\to {\rm Isom}_{\mathbb C}(\chi)_0\colon A\mapsto \left(\textnormal{\textcyr{D}}^-_X\right)^{-1}\circ A\circ \textnormal{\textcyr{D}}^-_X
\end{equation}defines a monomorphism. In particular, ${\mathcal C}_0(X)$ is abelian.
\end{enumerate}  
\end{teorema}
\esh

\proof
Claim (1) immediately follows from the definition of the action of the braid group on the monodromy data, and on the exceptional bases (Proposition \ref{propbraidmukai}). For claim (2), recall that the pairs of monodromy data $(S, C^{(k)})$ are related one to another by a power of the generator of the center of the braid group, as in Corollary \ref{centerbraidlemma}. Hence one concludes by Corollary \ref{corhelixser}. Point (3) follows from the identification of $S$ with the inverse of the Gram matrix, from the constraint (5) of Theorem \ref{19.03.18-1}, and from the definition of ${\mathcal C}_0(X)$: indeed, if $A\in{\mathcal C}_0(X)$, then
\begin{align*}(C^{-1}AC)^TS^{-1}(C^{-1}AC)&=C^TA^TC^{-T}S^{-1}C^{-1}AC\\
&=C^TA^T\eta e^{\pi i\mu} e^{\pi i R}AC\\
&=C^TA^T\eta e^{\pi i\mu} Ae^{\pi i R}C\\
&= C^T\eta e^{\pi i\mu} e^{\pi i R}C\\
&=S^{-1}.
\end{align*}
Moreover, since ${\mathcal C}_0(X)$ is unipotent (and since we are working in characteristic zero) it is connected (\cite{demagabr}, Prop. IV.2.4.1). From Corollary \ref{04.08.17-2}, we deduce that ${\mathcal C}_0(X)$ is abelian. 
\endproof

\begin{oss}
From Corollary \ref{centerbraidlemma}, Corollary \ref{corcancappa} and point (2) of Theorem \ref{03.08.17-1}, it follows that all eigenvalues of the monodromy matrix $M_0$ are equal to $(-1)^{\dim_{\mathbb C}X}$: indeed the inverse matrix $M_0^{-1}\in {\mathcal C}_0(X)$ corresponds to the canonical operator $\kappa\in {\rm Isom}_{\mathbb C}(\chi)_0$. Notice that for a smooth projective variety with odd-vanishing cohomology, this can also be obtained directly from the identity $M_0=\exp(2\pi i \mu)\exp(2 \pi i R)$, with $R=c_1(X)\cup (-)\colon H^\bullet(X,\mathbb C)\to H^\bullet(X,\mathbb C)$.  

From this simple fact, one can also deduce validity of several \emph{Diophantine constraints} on the entries of the Stokes matrices $S$, computed in any $\ell$-chambers of $QH^\bullet(X)$ with respect to any line $\ell$. Indeed, if we set
\[p_S(t):=\det(t\mathbbm 1_N-S^TS^{-1})=\sum_{j=0}^Np_j(S)t^j,\quad\text{ where $N:=\sum h^{p,q}(X)$,}
\]from the identity $p(t)=(t-(-1)^{\dim_{\mathbb C}X})^N$, we deduce the constraints (polynomial in the entries of $S$)
\[p_j(S)=(-1)^{(\dim_{\mathbb C}X+1)(N-j)}\binom{N}{j},\quad j=0,\dots, N.
\]Since $\det S=1$, we have that $$p_S\left(\frac{1}{t}\right)=\frac{(-1)^N}{t^N}p_S(t),$$so that $p_{N-j}(S)=(-1)^Np_j(S)$. Furthermore, since $p_0(S)=(-1)^{(\dim_\mathbb CX+1)N}$, we deduce that
\begin{itemize}
\item if $N$ is odd, then necessarily $\dim_\mathbb CX$ must be even,
\item the number of independent Diophantine constraints is $[\frac{N}{2}]$. 
\end{itemize}
As discussed in Section 5 of \cite{bondalsympl}, $p_1(S),\dots, p_{[\frac{N}{2}]}(S)$ actually define algebraically independent polynomials in the entries of $S$, and they freely generate the algebra of polynomial functions defined on the space of matrices $S$ invariant under the action of both $\mathcal B_N$ and $(\mathbb Z/2\mathbb Z)^{N}$ .

In the case $N=3$, necessarily $\dim_\mathbb CX$ must be even, and all  the Diophantine constraints reduce to
\begin{equation}\label{pseudomarkov}\quad a^2+b^2+c^2=abc,
\end{equation}where \[S=\begin{pmatrix}
1&a&b\\
0&1&c\\
0&0&1
\end{pmatrix}.
\]
By reducing equation \eqref{pseudomarkov} modulo 3, all integer solutions of equation \eqref{pseudomarkov} are of the form
\[a=3x,\quad b=3y,\quad c=3z,
\]where $(x,y,z)$ are integer solutions of the Markov equation
\[x^2+y^2+z^2=3xyz.
\]It is well known that all integer solutions of the Markov equation are equivalent (up to the action of  braid group $\mathcal B_3$) to the triple $(1,1,1)$. We refer the reader to the monograph \cite{aigner} for a modern survey of properties of Markov numbers and related open problems. 

In the case $N=4$, for the Stokes matrix
\begin{equation}\label{06.05.18-3}S=\begin{pmatrix}
1&a&b&c\\
0&1&d&e\\
0&0&1&f\\
0&0&0&1
\end{pmatrix},
\end{equation}the Diophantine constraints reduce to
\begin{align}
\label{06.05.18-1}a^2+b^2+c^2+d^2+e^2+f^2-abd-ace-bcf-def+acdf&=4(1-(-1)^{\dim_{\mathbb C}X}),\\
\label{06.05.18-2}(af-be+cd)^2=8(1-(-1)^{\dim_{\mathbb C}X}).
\end{align}
The polynomials in the lhs of \eqref{06.05.18-1}-\eqref{06.05.18-2} are generators of the algebra of polynomial functions on the space of matrices \eqref{06.05.18-3} which are invariant under the action of both the braid group $\mathcal B_4$ and of $(\mathbb Z/2\mathbb Z)^4$. In the case of surfaces ($\dim_\mathbb CX=2$), in the preprint \cite{dtdvvdb}, it is shown that a solution \eqref{06.05.18-3} of equations \eqref{06.05.18-1}-\eqref{06.05.18-2} is equivalent (modulo the action of $\mathcal B_4$ and $(\mathbb Z/2\mathbb Z)^4$) to exactly one of the following solutions
\[\begin{pmatrix}
1&2&2&4\\
0&1&0&2\\
0&0&1&2\\
0&0&0&1
\end{pmatrix},\quad
\begin{pmatrix}
1&n&2n&n\\
0&1&3&3\\
0&0&1&3\\
0&0&0&1
\end{pmatrix},\quad n\in\mathbb N.
\]The first solution corresponds to $\mathbb P^1\times\mathbb P^1$ (for which the Conjecture \eqref{congettura} holds true, as shown in Proposition \ref{congprod}): it is the Gram matrix of the full exceptional collection
\[(\mathcal O, \mathcal O(1,0),\mathcal O(0,1),\mathcal O(1,1)).
\]For the second family of solutions, for $n=0$ and $n=1$ we obtain the Gram matrices of full exceptional collections respectively in $\mathbb P^2\coprod {\rm pt}$, and the blow-up of $\mathbb P^2$ at one point. Remarkably, for all other values of $n$ it is believed that these matrices cannot be realized as Gram matrices of full exceptional collections on a smooth projective variety, although they arise in Grothendieck groups of \emph{non-commutative surfaces} (see e.g. \cite{dtdvp} for $n=2$, \cite{belpres} for the case $n\geq 2$, and \cite{belpresvan} for a comparison of their constructions).
\end{oss}

\bsh
\begin{prop}\label{congprod}The class of Fano variety for which Conjecture \ref{congettura} holds true is closed under finite products.
\end{prop}
\esh
\proof
Let $X,Y$ be Fano varieties for which Conjecture \ref{congettura} holds true, an let us define the canonical projections
\[\xymatrix{
&X\times Y\ar[dr]^{\pi_2}\ar[dl]_{\pi_1}&\\
X&&Y
}
\]If $X, Y$ have semisimple quantum cohomology, then also the tensor product of Frobenius manifolds $QH^\bullet(X\times Y)=QH^\bullet(X)\otimes QH^\bullet(Y)$ is semisimple. Furthermore, if $(E_0,\dots, E_n)$ and $(F_0,\dots, F_m)$ are full exceptional collections in $\mathcal D^b(X)$ and $\mathcal D^b(Y)$, respectively, then the collection $(E_{i}\boxtimes F_{j})_{(i,j)}$, indexed by all pairs $(i,j)$, is a full exceptional collection for $\mathcal D^b(X\times Y)$ (see e.g. \cite{kuz4}). Here we set
\[E\boxtimes F:=\pi_1^*E\otimes\pi_2^*F.
\]The order of the objects is intended to be the lexicographical one on the pairs $(i,j)$. Using the identities
\begin{align*}\widehat\Gamma^\pm_{X\times Y}&=\pi_1^*\widehat\Gamma^\pm_{X}\cup \pi_2^*\widehat\Gamma^\pm_{Y},\\
{\rm Ch}\left(E_{i}\boxtimes F_{j}\right)&=\pi_1^*{\rm Ch}(E_{i})\cup\pi_2^*{\rm Ch}(F_{j}),\\
c_1(X\times Y)&=\pi_1^*c_1(X)+\pi_2^*c_1(Y),\\
\dim(X\times Y)&=\dim X+\dim Y,
\end{align*}
and recalling that if $M,M'$ are two semisimple Frobenius manifolds we have that
\[S_{M\otimes M'}=S_M\otimes S_{M'},\quad C_{M\otimes M'}=C_M\otimes C_{M'}
\](see \cite{dubro2}, Lemma 4.10), we easily conclude.
\endproof

\subsection{Relations with Kontsevich's Homological Mirror Symmetry}\label{KHMS} The validity of the Conjecture \ref{congettura}, at least of its points (1) and (3a), can be heuristically deduced from M. Kontsevich's idea of \emph{Homological Mirror Symmetry} (\cite{hms, konens}). More precisely, Conjecture \ref{congettura} establishes an explicit relationship between two different geometrical aspects of the same Fano manifold $X$, the \emph{symplectic} one (the $A$-side) and the \emph{complex} one (the $B$-side), which can be connected through the study of an object \emph{mirror dual} to $X$. 

Although Mirror Symmetry phenomena were originally studied in the case of Calabi--Yau varieties,  several mirror conjectural correspondences have been generalized also to the Fano setting by the works of A. Givental (\cite{givehomgeo, give2, givemirtor}), M. Kontsevich (\cite{konens}), K. Hori and C. Vafa (\cite{horivafa}). If $X$ is a Fano manifold satisfying the semisimplicity condition of Conjecture \ref{congettura} then its mirror dual is conjectured to be a pair $(V,f)$ (the \emph{Landau--Ginzburg model}), where 
\begin{itemize}
\item $V$ is a non-compact Kähler manifold (whose symplectic form will be denoted by $\omega$), 
\item and $f\colon V\to\mathbb C$ is a holomorphic function which defines a Lefschetz fibration, i.e. $f$ admits only isolated non-degenerate critical points $\left\{p_1,\dots, p_n\right\}$, with only $A_1$-type singularities (i.e. Morse-type), and whose fibers are symplectic submanifolds of $V$.
\end{itemize}
With such an object one can associate two different categories, codifying respectively symplectic and complex geometrical properties of the the pair $(V, f)$. Let us briefly recall their constructions.
The symplectic geometry, also called $A$-\emph{side} or \emph{Landau--Ginzburg A-model}, is described by a Fukaya-type $\mathcal A_{\infty}$-category, originally introduced by M. Kontsevich and later by K. Hori, and whose explicit and rigorous construction has been formalized by P. Seidel (\cite{seidel1, seidel2, seidel3}).     
On the fibration $f\colon V\to\mathbb C$ one can consider a symplectic transport, by considering as horizontal spaces the symplectic orthogonal complement of vertical subspaces, i.e.
\[\mathcal H_p:=(\ker df_p)^{\perp\omega},\quad p\in V.
\]
For a fixed regular value $z_0\in\mathbb C$, by choosing $n$ paths $\gamma_i$ connecting $z_0$ with the critical values\footnote{We assume that the critical values, and the paths are numbered in clockwise order around the regular value $z_0$.}  $z_i:=f(p_i)$ for $i=1,\dots, n$, so that one can symplectically transport along the arc $\gamma_i$ the vanishing cycles at $p_i$. In this way one obtains a Lagrangian disc $D_i\subseteq V$ fibered above $\gamma_i$ (such a disc is called the \emph{Lefschetz thimble} over $\gamma_i$), and whose boundary is a Lagrangian sphere $L_i$ in the fiber $f^{-1}(z_0)$. Assuming genericity conditions, in particular that all the paths intersect each other only at $z_0$ and that all Lagrangian spheres intersect transversally in $f^{-1}(z_0)$, one can introduce the so called \emph{directed Fukaya category of} $(f,\left\{\gamma_i\right\})$.

\begin{defi}[\cite{seidel1, seidel2}]The directed Fukaya category $\mathcal {F}uk(V,f,\left\{\gamma_i\right\})$ is defined as the $\mathcal A_\infty$-category whose objects are the Lagrangian spheres $L_1,\dots, L_n$ and whose morphisms are given by
\[
\Hom(L_i,L_j):=\begin{sistema}
CF^\bullet(L_i,L_j;\mathbb C)\cong \mathbb C^{|L_i\cap L_j|},\quad\text{ if }i<j,\\
\\
\mathbb C\cdot {\rm Id},\quad\quad\quad\quad\quad\quad\quad\quad\quad\quad\text{ if }i=j,\\
\\
0,\quad\quad\quad\quad\quad\quad\quad\quad\quad\quad\quad\quad\text{ if }i>j,
\end{sistema}
\]where the Floer cochain complex $CF^\bullet(L_i,L_j;\mathbb C)$ with complex coefficients, the differential $m_1$, the composition $m_2$ and all other higher degree products $m_k$'s are defined in terms of Floer Lagrangian (co)homology in the fiber $f^{-1}(z_0)$. 
\end{defi}
The directed Fukaya category is unique up to quasi-isomorphism, and the derived category $\mathcal {DF}uk(V,f)$ only depends on $f\colon V\to\mathbb C$ (\cite{seidel1}, Corollary 6.5). Furthermore, the objects $(L_1,\dots, L_n)$ define a full exceptional collection of $\mathcal {DF}uk(V,f)$, and to different choices of paths $\left\{\gamma_i\right\}$ (actually inside the same \emph{Hurwtiz equivalence} class) there correspond different choices of full exceptional collections, related one to another by operations called \emph{mutations} (not totally coinciding with the ones discussed in Section \ref{Helixsec}). For more details the reader can consult the cited references. 

The second category associated with the pair $(V,f)$, encoding its complex geometrical aspects, is the so called \emph{triangulated category of singularities} defined by D. Orlov (\cite{orlov1, orlov3}). If $Y$ is an algebraic variety over $\mathbb C$, in what follows we denote by $\frak{Perf}(Y)$ the full triangulated subcategory of $\mathcal D^b(X)$ formed by  perfect complexes, i.e. objects locally isomorphic to a bounded complex of coherent sheaves of finite type: in particular, if $Y$ is smooth, then $\frak{Perf}(Y)\equiv \mathcal D^b(Y)$.
\begin{defi}[\cite{orlov1, orlov3}] We define the \emph{triangulated category of singularities of }$(V,f)$ as the disjoint union
\[\mathcal D_{\rm sing}(V,f):=\coprod_{z\in\mathbb C}\mathcal D_{\rm sing}(V_z),\quad V_z:=f^{-1}(z),
\]where we introduced the quotient category,
\[\mathcal D_{\rm sing}(V_z):=\mathcal D^b(V_z)/\frak{Perf}(V_z).
\]Such a quotient is defined by localizing the category $\mathcal D^b(V_z)$ with respect to the class of morphisms $s$ embedding into an exact triangle
\[\xymatrix{
X\ar[r]^{s}&Y\ar[r]&Z\ar[r]&X[1],}\quad \text{with }Z\in\Ob(\frak{Perf}(V_z)).
\]In particular,  note that $\mathcal D_{\rm sing}(V_z)$ is non-trivial only at the critical values $z_1,\dots, z_n$ of $f$.
\end{defi}

The crucial point in our discussion is the following homological formulation of Mirror Symmetry in the Fano case. 

\begin{conj}[Homological Mirror Symmetry, \cite{konens}]\label{hms}
Let $X$ be a Fano variety. There exist equivalences of triangulated categories as follows:

\[\xymatrix{A\text{\scshape-Model}&&B\text{\scshape-Model}\\\ 
\mathcal{F}uk(X)\ar[rrd]&&\mathcal D^b(X)\ar[lld]\\
\mathcal{DF}uk(V,f)\ar[rru]&& \mathcal D_{sing}(V,f)\ar[llu]}\]
\end{conj}

It is believed that the net of equivalences described above could be recast in terms of isomorphy of Frobenius manifolds structures, associated with $X$ and $(V,f)$, respectively. More precisely, we have that
\begin{itemize}
\item the Frobenius manifold related to the symplectic geometry of $X$ (the $A$-side) is the quantum cohomology $QH^\bullet(X)$;
\item the Frobenius manifold associated with $(V,f)$, and encoding information about its complex geometrical aspects (the $B$-side), is the Frobenius manifold structure defined on the space of miniversal unfoldings of $f$. The general construction is well-defined thanks to the works of A. Douai and C. Sabbah \cite{douaisab1,douaisab2,sabbah}, C. Hertling \cite{hertling, hertt*} and also of S. Barannikov's construction of Frobenius structures arising from semi-infinite variations of Hodge structures (\cite{barann1}). These efforts can be seen as a generalization of the construction of K. Saito \cite{Saito2}, who considered the case of germs of functions defined on $\mathbb C^n$.
\end{itemize}
The $A$-model and $B$-model of the Landau--Ginzburg mirror $(V,f)$ are conjectured to \emph{numerically} related in the following way:

\begin{conj}\label{conjABLG}
The Stokes matrix of the $B$-model Frobenius manifold associated with $(V,f)$ equals the Gram matrix of the Grothendieck--Euler--Poincaré product $\chi(\cdot,\cdot)$ product on $\mathcal{DF}uk(V,f)$.
\end{conj}

Putting together Conjecture \ref{conjABLG} and Conjectures \ref{hms}, it is clear that points (1) and (3.a) of Conjecture \ref{congettura} should (heuristically) follow.

\begin{oss}
In a recent preprint \cite{kuzsmi}, A. Kuznetsov and M. Smirnov formulated an intriguing conjecture, probably very close to Conjecture \ref{congettura}. Their work focuses on the case of a complex Fano variety $X$ of Picard rank one and index $m$ (i.e. $-K_X=mH$ for the ample generator $H$ of ${\rm Pic}(X)$). It is conjectured that a necessary condition for the semisimplicity of the small quantum cohomology in $QH^\bullet(X)$ is existence of a full \emph{Lefschetz collection}, i.e. an exceptional collection of the form
\[(\underbrace{E_1,E_2,\dots, E_{\tau_0}}_{\text{starting block}};E_1(H),\dots, E_{\tau_1}(H);\dots; E_{1}((m-1)H),\dots, E_{\tau_{m-1}}((m-1)H)),
\]where $\tau_0\geq \tau_1\geq\dots\geq\tau_{m-1}\geq 0$, which is \emph{minimal} with respect to the partial order on the set of Lefschetz collections defined by inclusion of the starting blocks, and whose \emph{residual category}
\[\langle E_1,E_2,\dots, E_{\tau_{m-1}};\dots; E_{1}((m-1)H),\dots, E_{\tau_{m-1}}((m-1)H)\rangle^\perp
\]is generated by a completely orthogonal exceptional collection.

Kuznetsov and Smirnov also suggest that the structure of the residual category can be deduced from HMS Conjecture \ref{hms}, as exposed in the previous paragraphs. Namely, it is expected the existence of a \emph{$\mu_m$-equivariant}\footnote{Here $\mu_m$ denotes the group of $m$-th roots of unity.} Landau--Ginzburg model $(V,f)$ for $X$, i.e. with a $\mu_m$-action on $V$ and $f$ equivariant with respect to the standard $\mu_m$-action on $\mathbb C$. Consequently, the nonzero critical values of $f$ can be partitioned into a number of free $\mu_m$-orbits. Under the equivalence $\mathcal{DF}uk(V,f)\cong \mathcal D^b(X)$, it is expected that
\begin{itemize}
\item the objects of the rectangular subcollection
\[(E_1,E_2,\dots, E_{\tau_{m-1}};\dots; E_{1}((m-1)H),\dots, E_{\tau_{m-1}}((m-1)H))
\]are generated by the thimbles associated with critical points of $f$ with \emph{nonzero} critical values;
\item the residual category is generated by thimbles associated with critical points over $0\in\mathbb C$. By semisimplicity, the critical points of $f$ are isolated and simple: the corresponding vanishing cycles and their symplectic transport do not intersect neither over a neighborhood of $0$ nor over a neighborhood of the chosen regular point in $\mathbb C$. From this, the complete orthogonality of the corresponding exceptional objects in $\mathcal{DF}uk(V,f)$ follows.
\end{itemize}

 In \cite{kuzsmi} the validity of this conjecture for complex Grassmannians $\mathbb G(r,k)$ is also shown, with $r$ a prime number, under the assumption of completeness of Fonarev's collection (a Lefschetz collection introduced in \cite{fonarev}). The validity of the completeness assumption is proved for $r=3$.
 
 We believe that further investigations are needed, in order to understand any contingent relationship between the Conjecture of \cite{kuzsmi} and Conjecture \ref{congettura}. In particular, we believe that, after identifying the set of critical values of $f$ with the spectrum of the operator $\mathcal U$, it has to be clarified whether the conjecture of Kuznetsov and Smirnov can be justified through the study of the Stokes phenomenon of the equations \eqref{06.07.17-1}-\eqref{06.07.17-2} at semisimple points of the small quantum cohomology in either the discriminant locus (defined by $\prod_iu_i=0$) or the coalescence locus (i.e. $u_i=u_j$ for some $i\neq j$), or their intersection. Furthermore, we also wonder about relationships (if any) between the minimal Lefschetz collection of \cite{kuzsmi} and the $m$-block exceptional collections of Remark \ref{teoblock}.
\end{oss}

\subsection{Galkin-Golyshev-Iritani  $\Gamma$-Conjecture II and its relationship with Conjecture \ref{congettura}}\label{relggi}Few months after the beginning of our research project, two very interesting papers by S. Galkin, V. Golyshev and H. Iritani appeared (\cite{gamma1,gamma2}). In \emph{loc. cit.}, the authors proposed two conjectures, called $\Gamma$-conjectures, describing the exponential asymptotics of flat sections for an extended deformed connection, that for clarity we will denote $\widehat\nabla^{\rm GGI}$. This connection is defined on the pull-back of the tangent bundle $TQH^\bullet(X)$ along the natural projection $\pi\colon \mathbb C^*\times QH^\bullet(X)\to QH^\bullet(X)$, in an analogous way to our connection $\widehat\nabla$ (Section \ref{extdefconn}), although with a difference along the tangential direction of the spectral parameter $\lambda\in\mathbb C^*$, namely
\begin{align*}
\widehat\nabla^{\rm GGI}_{\alpha}&=\frac{\partial}{\partial t^\alpha}+\frac{1}{\lambda}(T_\alpha\circ),\\
\widehat\nabla^{\rm GGI}_{\frac{\partial}{\partial\lambda}}&=\frac{\partial}{\partial\lambda}-\frac{1}{\lambda^2}\mathcal U+\frac{1}{\lambda}\mu.
\end{align*}
The two connections $\widehat\nabla, \widehat\nabla^{\rm GGI}$ can be identified by setting $\lambda=z^{-1}$. Despite this difference, since Galkin, Golyshev and Iritani focused their attention on flat \emph{vector fields}, rather than flat \emph{differentials} defining deformed flat coordinates, the isomonodromic linear differential system considered in \cite{gamma1,gamma2}, namely 
\begin{equation}\label{ggi}\widehat\nabla^{\rm GGI}_{\frac{\partial}{\partial\lambda}}Y=0,\quad Y\in \Gamma(\pi^*TQH^\bullet(X)),
\end{equation} is exactly the same as the one considered in the present paper, i.e. 
\begin{equation}\label{noi}\widehat\nabla_{\frac{\partial}{\partial z}}\xi=0,\quad \xi\in\Gamma(\pi^*T^*QH^\bullet(X)),
\end{equation} provided we set $\lambda=-z^{-1}$. Clearly, in order to  obtain solutions of \eqref{noi}, or equivalently of its gauge equivalent form \eqref{06.07.17-1}-\eqref{06.07.17-2}, starting from solutions of \eqref{ggi} one has to specify a sign $(\pm)$ in the identification formula
\begin{equation}\label{identvariab}\lambda=e^{\pm i\pi}z^{-1},\quad\lambda,z\in\mathcal R.
\end{equation}The two choices lead to two solutions of \eqref{noi}, differing from each other by a right multiplication by the monodromy matrix $M_0^{\pm 1}$. In this way, the monodromy and Stokes phenomenon of solutions of \eqref{ggi} can be described analogously to our equation \eqref{noi}, by introducing a Stokes matrix and a central connection matrix that, for clarity we will denote by $(S^{\rm GGI},C^{\rm GGI})$. More precisely, let us fix an admissible direction $\phi$ for the system \eqref{ggi}, and let us denote by
\[Y_{{\rm left}, \phi}^{\rm GGI},\quad Y_{{\rm right}, \phi}^{\rm GGI},
\] the solutions of the system \eqref{ggi} of asymptotic expansion
\[\Psi^{-1}(\mathbbm 1+O(\lambda))\exp\left(-\frac{1}{\lambda}U\right),\quad U={\rm diag}(u_1,\dots, u_n),
\]respectively in the sectors 
\[\phi<\arg\lambda<\phi+\pi,\quad \phi-\pi<\arg\lambda<\phi.\] 
The corresponding Stokes and central connection matrices, as defined in \cite{gamma1}, are defined by the equations
\[Y_{{\rm right}, \phi}^{\rm GGI}=Y_{{\rm left}, \phi}^{\rm GGI}\ S^{\rm GGI}_\phi,
\]
\[Y_{{\rm right}, \phi}^{\rm GGI}=Y^{\rm GGI}_{0}\ C^{\rm GGI}_\phi.
\]
The central connection matrix $C^{\rm GGI}$ is apparently defined in an analogous way to our matrix $C$: it is computed with respect to to the solution of \eqref{ggi} given at the point $0\in QH^\bullet(X)$ by
\begin{equation}\label{topggi}Y^{\rm GGI}_0(\lambda)=\Theta_{\rm top}(-\lambda^{-1})\lambda^{-\mu}\lambda^{c_1(X)\cup},
\end{equation}where $\Theta_{\rm top}(z)\equiv \Theta_{\rm top}(z,0)$ is the series given in Proposition \ref{18.07.17-2} (see Section 2.3 of \cite{gamma1}).

\begin{conj}[$\Gamma$-Conjecture II, \cite{gamma1}]\label{gammaggi}
Let $X$ be a Fano variety with semisimple quantum cohomology. The entries of the columns of the central connection matrix $C^{\rm GGI}$ (computed at $0\in QH^\bullet(X)$, with respect to an admissible line, and a branch of the $\Psi$-matrix) are the components (with respect to a prefixed basis of $H^{\bullet}(X,\mathbb C)$) of the characteristic classes 
\[\frac{1}{(2\pi)^{\frac{d}{2}}}\widehat\Gamma^+_X\cup{\rm Ch}(E_i),\quad d=\dim_{\mathbb C}X,
\]for an exceptional collection $(E_1,\dots, E_n)$.
\end{conj}

We want now to show the equivalence of point (3b) of Conjecture \ref{congettura} with Conjecture \ref{gammaggi}, by pointing out a few subtleties. Let us focus on the power-logarithmic term of the solution \eqref{topggi}, namely $\lambda^{-\mu}\lambda^{c_1(X)\cup}$. Using the identification \eqref{identvariab}, and point (4) of Lemma 2.1 of \cite{CDG}, we have the following identity: 
$$
(e^{\pm i\pi}z^{-1})^{-\mu}(e^{\pm i\pi}z^{-1})^{R}=z^{\mu}e^{\mp i\pi\mu}z^{-R}e^{\pm i\pi R}$$
\begin{equation}
\label{17agosto2018-6}
=z^\mu z^R\underbrace{e^{\mp i\pi\mu}e^{\pm i\pi R}}_{K_\pm},
\end{equation}
$R$ denoting the operator $c_1(X)\cup(-)$. Thus, using the identification of variables \eqref{identvariab}, the solution \eqref{topggi} is equal to 
\[\Theta_{\rm top}(z)z^\mu z^RK_\pm.
\]The matrix $K_\pm$ responsible for the difference between the predicted forms of the central connection matrix of both Conjecture \ref{congettura} and Conjecture \ref{gammaggi} is an element of the $(\eta,\mu)$ parabolic orthogonal group $\mathcal G(X)$ (it follows from point (4) of Lemma 2.1 of \cite{CDG}), but it \emph{is not an element of }$\mathcal C_0(X)$. This means that the solution \eqref{topggi}, under the identification of the systems \eqref{ggi} and \eqref{noi}, is not a solution in the \emph{natural} Levelt  form at $z=0$, i.e. the one dictated by the limit of system \eqref{noi} at the classical limit point \eqref{classical1}-\eqref{classical2}. In order to understand to which adjoint orbit in $\frak g(X)$ it corresponds the choice done in \cite{gamma1}, let us factor the power-logarithmic term $\lambda^{-\mu}\lambda^{c_1(X)\cup}$ as follows
$$
(e^{\pm i\pi}z^{-1})^{-\mu}(e^{\pm i\pi}z^{-1})^{R}=e^{\mp i\pi\mu}z^{\mu}e^{\pm i\pi R}z^{-R}
$$
$$=e^{\mp i\pi\mu}\underbrace{z^{\mu}e^{\pm i\pi R}z^{-\mu}}_{\text{ polynomial in }z}z^\mu z^{-R}
$$
\begin{equation}
\label{17agosto2018-7}
=e^{\mp i\pi\mu}P_{e^{\pm i\pi R}}(z)z^\mu z^{-R},
\end{equation}
where we used for $P$ the notations of Theorem \ref{06.07.17-4}, point (2).  By combining \eqref{17agosto2018-6} and  \eqref{17agosto2018-7}, we have
$$
z^\mu z^RK_\pm=e^{\mp i\pi \mu}P_{e^{\pm i\pi R}}(z)~z^\mu z^{-R}.
$$
Thus,  the Levelt form chosen in \cite{gamma1} has exponent $-R$,  the opposite to the natural one usually used for ``good'' Frobenius manifolds.

Now there is a further subtlety, which apparently could create an incompatibility of Conjecture \ref{congettura} and Conjecture \ref{gammaggi}. 

\bsh
\begin{prop}\label{propgammapiu}
Let $X$ be a smooth projective variety for which Conjecture \ref{congettura} holds true. In particular let the central connection matrix $C$ (computed with respect to some choice of $\Psi$ and of an oriented line) be equal to the matrix associated with the morphism $\textnormal{\textcyr{D}}^-_X$ and with some exceptional collection $\frak E=(E_1,\dots, E_n)$. Then there exists another choice of the fundamental solution of equation \eqref{noi} in natural Levelt form at $z=0$, obtained from the topological-enumerative solution by multiplication to the right by an element of $\mathcal C_0(X)$, with respect to which the central connection matrix has  entries  given by the components of  
\[\frac{i^{\bar{d}}}{(2\pi)^\frac{d}{2}}\widehat\Gamma^\pm_X\cup{\rm Ch}(E_i),\quad d=\dim_\mathbb CX.
\]\end{prop}
\esh

\proof
Let $f(t)\in \mathbb C[\![t]\!]$ be a formal power series of the form 
\[
f(t)=1+O(t),\quad f(t)\tilde f(t)=1,\quad \tilde f(t):=f(-t),\]
and let us introduce the corresponding characteristic class by applying the Hirzebruch construction (\cite{hirz})
\[\lambda_f:=\prod_{j=1}^d f(\delta_j),\quad \delta_j\text{'s Chern roots of }TX.
\]We claim that $\lambda_f\cup\colon H^\bullet(X,\mathbb C)\to H^\bullet(X,\mathbb C)$ is an element of $\mathcal C_0(X)$. Indeed it is clearly of the form 
\[\lambda_f\cup(-)=\mathbbm 1_{H^\bullet(X,\mathbb C)}+\Delta,\quad \Delta\text{ is $\mu$-nilpotent},
\]and it is a $\left\{\cdot,\cdot\right\}$-isometry as the following computation shows:
\begin{align*}
\left\{\lambda_f\cup a,\lambda_f\cup b\right\}&=\int_Xe^{i\pi\mu}\left[\lambda_f\cup a\right]\cup\lambda_f\cup b\\
&=\int_X\lambda_{\tilde f}\cup\lambda_f\cup e^{i\pi\mu}(a)\cup b\\
&=\int_Xe^{i\pi\mu}(a)\cup b\\
&=\left\{a,b\right\}.
\end{align*}
Finally, because of the commutativity of $H^\bullet(X,\mathbb C)$ the operator $\lambda_f\cup (-)$ commutes with the multiplication by the first Chern class $c_1(X)\cup(-)$, and consequently it is an element of $\mathcal C_0(X)$. The statement follows from the choices 
\[f(t)=e^{\pi i t}\frac{\Gamma(1+t)}{\Gamma(1-t)},\quad f(t)=e^{\pi i t},
\]
respectively. 
\endproof

So, if we consider an even dimensional Fano variety $X$, we would have to handle  two different fundamental solutions (only one of which is in natural Levelt form at $z=0$) such that the corresponding central connection matrices are both of the type described by Conjecture \ref{gammaggi}. Despite the apparent incompatibility, this is due to a \emph{difference in the exceptional collections} considered.

\bsh
\begin{teorema}\label{teoequivgamma2}
Point $\rm (3b)$ of  Conjecture \ref{congettura} and Conjecture \ref{gammaggi} are equivalent. More precisely, if the central connection matrix $C_{-\phi+2\pi}$ of the system \eqref{noi} is of the form prescribed by point \emph{(3.b)} of Conjecture \ref{congettura} for some exceptional collection $\frak E$, then the central connection matrix $C_\phi^{\rm GGI}$ of the system \eqref{ggi} is of the form prescribed by $\Gamma$-conjecture {\rm II} for the \emph{dual} exceptional collection $\frak E^*$.
\end{teorema}
\esh

\proof
If we identify the systems \eqref{ggi} and \eqref{noi} by by $\lambda=e^{\pi i}z^{-1}$, then we have the identifications
\[Y_{{\rm left}, \phi}^{\rm GGI}(\lambda)\equiv Y_{{\rm left},-\phi}(z),
\]
\[Y_{{\rm right}, \phi}^{\rm GGI}(\lambda)\equiv Y_{{\rm right},-\phi+2\pi}(z),
\]and we thus obtain the identity
\begin{equation}\label{identconn0}C_{-\phi+2\pi}=K_+C^{\rm GGI}_\phi.
\end{equation}
If the entries of the columns of $C_\phi^{\rm GGI}$ are the components of the characteristic classes
\[\frac{1}{(2\pi)^{\frac{d}{2}}}\widehat\Gamma^+_X\cup {\rm Ch}(E_i),
\]as predicted by Conjecture \ref{gammaggi} for an exceptional collection $(E_i)_i$, then from identities \eqref{expmu}, \eqref{identconn0} we deduce that the entries of the columns of $C_{-\phi+2\pi}$ are the components of 
\[\frac{i^d}{(2\pi)^\frac{d}{2}}\widehat\Gamma^-_X\cup{\rm Ch}(E_i')\cup\exp(-\pi i c_1(X)),\quad E_i'=E_i^*.
\]
Analogously, if we identify the systems \eqref{ggi} and \eqref{noi} by $\lambda=e^{-\pi i}z^{-1}$, then we have the following identifications
\[Y_{{\rm left}, \phi}^{\rm GGI}(\lambda)\equiv Y_{{\rm left},-\phi-2\pi}(z),
\]
\[Y_{{\rm right}, \phi}^{\rm GGI}(\lambda)\equiv Y_{{\rm right},-\phi}(z).
\]
Consequently, we have that
\begin{equation}\label{identconn1}C_{-\phi}=K_-C_\phi^{\rm GGI}.
\end{equation}
Notice that equations \eqref{identconn0}- \eqref{identconn1} are coherent with Corollary \ref{centerbraidlemma}, because of the identity $K_+K_-^{-1}=M_0^{-1}$ (easily seen using again point (4) of Lemma 2.1 of \cite{CDG}). It follows that the entries of the columns of $C_{-\phi}$ are the components of 
\[\frac{i^d}{(2\pi)^\frac{d}{2}}\widehat\Gamma^-_X\cup{\rm Ch}(E_i')\cup\exp(-\pi i c_1(X)),\quad \kappa(E_i')=E_i^*,
\]where $\kappa=(\omega_X\otimes-)[\dim_\mathbb CX]$ denotes the Serre functor.
\endproof

In conclusion, we would like to alert the reader to some delicate points in the reconstruction procedure of the Frobenius structure on $QH^\bullet(X)$, through an inverse problem starting from its monodromy data as described in Section \ref{monlocalmod}. Indeed, despite the equivalence described in Theorem \ref{teoequivgamma2} above, the slogan ``$\Gamma$-conjecture II refines the conjecture of \cite{dubro0}'' could potentially lead to wrong results, if misunderstood or interpreted \emph{ad litteram}, because of different choices of Levelt forms w.r.t. the natural ones usually chosen in the theory of Frobenius manifolds.

So, for example, if we consider a smooth projective variety $X$, and we want to locally reconstruct the Frobenius structure of some $\ell$-chambers of $QH^\bullet(X)$, we could set at least six different RH b.v.p.'s as in Section \ref{invfrob} (at a fixed point $u^{(0)}=(u^{(0)}_1,\dots, u^{(0)}_N)\in\mathbb C^N$) associated with data $(\mu, R,S,C)$, where
\begin{enumerate}[(I)]
\item either $R=c_1(X)\cup$, and the columns of $C$ being the components of the characteristic classes
\[\frac{i^{\overline{d}}}{(2\pi)^\frac{d}{2}}\widehat\Gamma^\pm_X\cup{\rm Ch}(E_i)\cup e^{-\pi i c_1(X)},
\]
\item or $R=c_1(X)\cup$, and the columns of $C$ being the components of the characteristic classes
\[\frac{i^{\overline{d}}}{(2\pi)^\frac{d}{2}}\widehat\Gamma^\pm_X\cup{\rm Ch}(E_i),
\]
\item or $R=-c_1(X)\cup$, and the columns of $C$ being the components of the characteristic classes
\[\frac{1}{(2\pi)^\frac{d}{2}}\widehat\Gamma^\pm_X\cup{\rm Ch}(E_i),
\]
\end{enumerate}
where $(E_1,\dots, E_N)$ is an exceptional collection, and the matrix $S$ can be computed in each case through the constraint (5) of Theorem \ref{19.03.18-1}. Proposition \ref{propgammapiu} and Theorem \ref{teoequivgamma2} guarantee that if one of the above RH b.v.p.'s is solvable at $u^{(0)}$, then \emph{all of them are solvable}, and their solutions can be used for the reconstruction of the Frobenius structure (see also Remark \ref{ossinv}). Nevertheless, the solutions $\Phi=(\Phi_L,\Phi_R,\Phi_0)$ have \emph{different enumerative meaning}: the coefficient of the solution $\Phi_0$, for example, are Gromov-Witten invariants with descendants (namely the series $\Psi\cdot \Theta_{\rm top}$, using the notation of Proposition \ref{18.07.17-2}) only in case [I], whereas in case [II] and [III] they are related to them only up to the action of the group $\mathcal C_0(X)$.

\begin{oss}
Basing on $\Gamma$-conjecture II of S. Galkin, V. Golyshev and H. iritani, in the recent paper \cite{dubtyp}, F. Sanda and Y. Shamoto formulated an analogue of Conjecture \ref{originalconj} (1)-(2), called by the authors \emph{Dubrovin type conjecture}, concerning the case of Fano manifolds with non-necessarily semisimple quantum cohomology. In such a case, it is conjectured that the role played by exceptional collections should be played by more general semiorthogonal decompositions of the derived categories. We plan to address this case, with some explicit examples, in a future publication \cite{cs}. 
\end{oss}

\newpage
\section{Proof of the Main Conjecture for Projective Spaces}\label{computationsproj}
\subsection{Notations and preliminaries}\label{notationsproj}
In what follows
\begin{itemize}
\item the symbol $\mathbb P$ will stand for $\mathbb P^{k-1}_{\mathbb C}$, $k\geq 2$;
\item we denote by  $\sigma$ the generator of the 2-nd cohomology group $H^2(\mathbb P,\mathbb C)$, so that
\[H^\bullet(\mathbb P,\mathbb C)\cong \frac{\mathbb C[\sigma]}{(\sigma^k)}.
\]We also assume that $\sigma$ is normalized so that
\[\int_{\mathbb P}\sigma^{k-1}=1.
\]In this way $\sigma$ coincides with the hyperplane class, so that $c_1(\mathbb P)=k\sigma$.
\end{itemize}

The flat coordinates $t^1,\dots,t^k$ for the quantum cohomology of $\mathbb P$ are the coordinates with respect to the homogeneous base $$(1,\sigma,\sigma^2,\dots,\sigma^{k-1}),$$ the matrix of Poincar\'e metric being constant
\[\eta_{\alpha\beta}=\eta\left(\frac{\partial}{\partial t^\alpha},\frac{\partial}{\partial t^\beta}\right)=\delta_{\alpha+\beta,k+1}.
\]
Observe that the unity vector field is $e=\frac{\partial}{\partial t^1}$, and the Euler vector field is
\[E=\sum_{\alpha\neq 2}(1-q_\alpha)t^\alpha\frac{\partial}{\partial t^\alpha}+k\frac{\partial}{\partial t^2},\quad q_h=h-1\quad\text{for }h=1,\dots,k.
\]
If $\zeta$ is a column vector whose components are the components of the gradient of a deformed flat coordinate, w.r.t the frame $\left(\frac{\partial}{\partial t^\alpha}\right)_{\alpha=1}^k=(\sigma^{\alpha-1})_{\alpha=1}^k$, then it must satisfy the system
 \begin{empheq}[left=\empheqlbrace]{align*} 
        \partial_\alpha\zeta&=z\mathcal C_\alpha\zeta,\quad \alpha=1,\dots,k,\\
        \partial_z\zeta&=\left(\mathcal U+\frac{1}{z}\mu\right)\zeta.
    \end{empheq}
If we restrict to the locus of \emph{small quantum cohomology}, i.e. to the points $(0,t^2,0,\dots,0)$, the system above reduces to the system
\begin{align}\label{28.07.17-1}\partial_2\zeta&=z\mathcal C_2\zeta,\\ \label{28.07.17-2}\partial_z\zeta&=\left(\mathcal U+\frac{1}{z}\mu\right)\zeta,\end{align}
where at the point $(0,t^2,0,\dots,0)$
\[\mathcal U:=\begin{pmatrix}
0&&&&kq\\
k&0&&&\\
&k&0&&\\
&& \ddots&\ddots&\\
&&&k&0
\end{pmatrix},\quad q:=e^{t^2},\quad\mathcal C_2=\frac{1}{k}\mathcal U,
\quad\mu=\text{diag}\left(-\frac{k-1}{2},-\frac{k-3}{2},\dots,\frac{k-3}{2},\frac{k-1}{2}\right).
\]

The Stokes data of the system
\begin{equation}\label{eq0}
\frac{d\zeta}{dz}=\left(\mathcal U+\frac{1}{z}\mu\right)\zeta.
\end{equation}
 have been computed in \cite{guzzetti1}. Below, we review the main steps of \cite{guzzetti1}, and we complement them with the computation of the topological solution, of $\mathcal{C}_0(\mathbb{P})$ and of the central connection matrix.
 
The eigenvalues of the matrix $\mathcal U(0,t^2,0,\dots,0)$ are
\[u_h=ke^{\frac{2\pi i (h-1)}{k}}q^{\frac{1}{k}}\quad h=1,\dots,k,
\]and let us compute the corresponding eigenvectors $x_1,\dots, x_k$. The equations for $x_h=(x_h^1,\dots,x_h^k)$ read
\[kx_h^\ell=u_hx^{\ell+1}_h,\quad \ell=1,\dots,k-1,
\] 
\[kqx^k_h=u_hx^1_h.
\]By choosing $x_h^k=e^{\frac{i\pi (h-1)}{k}}$, we get all the entries
\[x_h^\ell=\left(\frac{u_h}{k}\right)^{k-\ell}x^k_h=q^{\frac{k-\ell}{k}}e^{(1-2\ell)i \pi\frac{(h-1)}{k}}\quad h,\ell=1,\dots,k.
\]Since the norm of the eigenvector $x_h$ is
\[\eta(x_h,x_h)=kq^{\frac{k-1}{k}},
\]we find (choosing signs of square roots) the orthogonal vectors $f_1,\dots,f_k$
\[f_h^\ell=k^{-\frac{1}{2}}q^{\frac{k+1-2\ell}{2k}}e^{(1-2\ell)i \pi\frac{(h-1)}{k}}\quad h,\ell=1,\dots,k.
\]
Thus the matrix $\Psi$ is given by
\[\Psi=\left(
\begin{array}{c|c|c|c}
&&&\\
f_1&f_2&\dots&f_k\\
&&&\\
\end{array}
\right)^{-1}.
\]
Instead of working with the differential equation \eqref{eq0}, in \cite{guzzetti1} the following  gauge equivalent system of differential equations for $\xi(z,t^2):=\eta\cdot \zeta(z,t^2)$ is considered
\begin{align}\label{eq2}
\partial_2\xi&=z\mathcal C_2^T\xi,\\
\label{eq2bis}
\partial_z\xi&=\left(\mathcal U^T-\frac{1}{z}\mu\right)\xi.
\end{align}
Now, the entries of the column vector $\xi$  are the components of the differential of a deformed flat coordinate.


A simple computation shows that with the following substitution
\[\xi_\alpha(z,t^2)=\frac{1}{k^{\alpha-1}}z^{\frac{k-1}{2}-\alpha+1}\vartheta^{\alpha-1}\Phi(z,t^2),
\]for any $\alpha=1,2,\dots,k$ and where $\vartheta:=z\frac{d}{dz}$, the system \eqref{eq2}-\eqref{eq2bis} is equivalent to the equations
\begin{align*}
\vartheta^k\Phi-(kz)^kq\Phi&=0,\\
\partial_2^k\Phi-z^kq\Phi&=0.
\end{align*}
The compatibility of these equations implies the following functional dependence of $\Phi$ on $(z,t^2)$: \[\Phi(t^2,z)=\Phi(q^{\frac{1}{k}}z).\]
Thus, the study of the system \eqref{eq2bis}, restricted to the point $t^2=0$, is equivalent to the study of the generalized hypergeometric equation
\begin{equation}
\label{eq3}
\vartheta^k\Phi(z)-(kz)^k\Phi(z)=0,
\end{equation}where $\vartheta:=z\frac{d}{dz}$. Given a solution $\Phi$ of \eqref{eq3}, the corresponding solution of equation \eqref{eq2bis} is given by
\begin{equation}
\label{17agosto2018-8}
\xi=
\begin{pmatrix}
z^{\frac{k-1}{2}}\Phi(z)\\
\vdots\\
\frac{1}{k^{\alpha-1}}z^{\frac{k-1}{2}-\alpha+1}\vartheta^{\alpha-1}\Phi(z)\\
\vdots\\
\frac{1}{k^{k-1}}z^{\frac{1-k}{2}}\vartheta^{k-1}\Phi(z)
\end{pmatrix}.
\end{equation}

\subsection{Computation of the Topological-Enumerative Solution} In this section, we use the characterization of the topologcal-enumerative solution described in Proposition \ref{18.07.17-2}. 
\begin{lemma}The formal\footnote{The components of the series $\Phi(z)$ w.r.t. the basis $(\sigma^0,\dots, \sigma^k)$ are actually convergent, according to Theorem \ref{06.07.17-4}.} series $\Phi(z)\in H^\bullet(\mathbb P,\mathbb C)[\log z]\llbracket z\rrbracket$
\[\Phi(z)=e^{k\sigma\log(z)}\sum_{n=0}^\infty f(n)z^{kn},\quad f(n)\in\mathbb C[\sigma]/(\sigma^k)\]
satisfies equation \[\vartheta^k\Phi(z)-(kz)^k\Phi(z)=0,
\]if and only if the coefficients $f(n)$ satisfy the following difference equation
\[(\sigma+n)^kf(n)=f(n-1),\quad n\geq 1.
\]
\end{lemma}

\proof[Proof]Observe that
\begin{align*}
\vartheta\Phi(z)&=ze^{k\sigma\log(z)}\left(\frac{k\sigma}{z}\sum_{n=0}^\infty f(n)z^{kn}+\sum_{n=0}^\infty f(n)knz^{kn-1}\right)\\
&=k e^{k\sigma\log(z)}\sum_{n=0}^\infty (\sigma+n)f(n)z^{kn}.
\end{align*}
By an inductive argument one can easily show that
\[\vartheta^\alpha\Phi(z)=k^\alpha e^{k\sigma\log(z)}\sum_{n=0}^\infty (\sigma+n)^\alpha f(n)z^{kn}.
\]So, using the fact that $\sigma^k=0$, we have that
\[\vartheta^k\Phi(z)=(kz)^k\Phi(z),
\]if and only if
\[(\sigma+n)^kf(n)=f(n-1)\quad \text{ for all }n\geq 1.
\]
\endproof

\bsh
\begin{prop} \label{19.07.17-2}
\begin{enumerate}
\item For any fixed value $f(0)\in H^\bullet (\mathbb P,\mathbb C)$, the corresponding formal solution of \eqref{eq3} is given by
\[\Phi(z)=f(0)\cdot\sum_{p=0}^{k-1}\left(\sum_{l=0}^{p}\frac{(k\log z)^{p-l}}{(p-l)!}a_{l}(z)\right)\sigma^p,\]
where, for $0\leq l\leq k-1$, we have introduced the notation
\begin{equation}\label{18.07.17-3}a_{l}(z):=\sum_{n=0}^\infty \alpha_{n,l}z^{kn},\quad \alpha_{0,l}:=\delta_{0,l},\quad \alpha_{n,l}:=\sum_{\substack{h_1+\dots+h_n=l\\ 0\leq h_i\leq k-1}}\left(\prod_{j=1}^n\frac{(-1)^{h_j}}{j^{k+h_j}}\binom{k-1+h_j}{h_j}\right).
\end{equation}
Representing $\Phi(z)=\sum_{i=1}^k\Phi_i(z)\sigma^{k-i}$, we deduce that each component 
\[\Phi_i(z):=f(0)\cdot\sum_{l=0}^{k-i}\frac{(k\log z)^{k-i-l}}{(k-i-l)!}a_{l}(z)
\] is a convergent solution of \eqref{eq3}. 
\item Another representation of the solution is given by the formula
\[\Phi(z)=f(0)e^{k\sigma\log z}\sum_{n=0}^\infty\frac{\Gamma(-\sigma-n)^k}{\Gamma(-\sigma)^k}e^{\pm k \pi i n}z^{kn}
\]for any choice of the sign $(\pm)$.
\item Moreover, if $f(0)=1$, the fundamental solution $\Xi$ of \eqref{eq2bis} given by
\begin{equation}\label{fundsol}\Xi_0(z)=
\begin{pmatrix}
z^{\frac{k-1}{2}}\Phi_1(z)&\dots&z^{\frac{k-1}{2}}\Phi_k(z)\\
\vdots&&\vdots\\
\frac{1}{k^{\alpha-1}}z^{\frac{k-1}{2}-\alpha+1}\vartheta^{\alpha-1}\Phi_1(z)&\dots&\frac{1}{k^{\alpha-1}}z^{\frac{k-1}{2}-\alpha+1}\vartheta^{\alpha-1}\Phi_k(z)\\
\vdots&&\vdots\\
\frac{1}{k^{k-1}}z^{\frac{1-k}{2}}\vartheta^{k-1}\Phi_1(z)&\dots&\frac{1}{k^{k-1}}z^{\frac{1-k}{2}}\vartheta^{k-1}\Phi_k(z)
\end{pmatrix}
\end{equation}is of the form
\[\Xi(z)=\eta \Theta_{\rm top}(z)z^\mu z^{c_1(\mathbb P)\cup(-)},\quad \Theta_{\rm top}(z)^\alpha_\gamma=\delta^\alpha_\gamma+\sum_{n=0}^\infty\sum_{\lambda}\sum_{\beta\in{\rm Eff}(\mathbb P)\setminus\left\{0\right\}}\langle\tau_n \sigma^\gamma, \sigma^\lambda\rangle^{\mathbb P}_{0,2,\beta}\eta^{\lambda\alpha}z^{n+1},
\]
\[\text{with }\quad\langle\tau_n \sigma^\gamma, \sigma^\lambda\rangle^{\mathbb P}_{0,2,\beta}:=\int_{[\overline{\mathcal M}_{0,2}(X,\beta)]^{\rm vir}}\psi_1^{n}\cup {\rm ev}_1^*(\sigma^\gamma)\cup {\rm ev}_2^*(\sigma^\lambda).
\] 
\end{enumerate}
\end{prop}
\esh

\proof[Proof]  From the identity
\[(1+\sigma)^{-1}=1-\sigma+\sigma^2-\dots+(-1)^{k-1}\sigma^{k-1},
\]one easily shows that if $n\geq 1$, then 
\[(n+\sigma)^{-k}=n^{-k}\left(1+\frac{\sigma}{n}\right)^{-k}=\sum_{h=0}^{k-1}\frac{(-1)^h}{n^{k+h}}\binom{k-1+h}{h}\sigma^{h}.
\]
As a consequence, we have that
\[f(n)=\sum_{l=0}^{k-1}f(0)\sigma^l\alpha_{n,l},
\]
where the numbers $\alpha_{n,l}\in\mathbb Q$ are defined as in \eqref{18.07.17-3}.
It follows that
\begin{align*}\Phi(z)&=f(0)e^{k\sigma\log z}\sum_{n=0}^\infty\sum_{l=0}^{k-1}f(0)\sigma^l\alpha_{n,l} z^{kn}\\
&=f(0)\left(\sum_{m=0}^{k-1}\frac{(k\log z)^m}{m!}\sigma^m\right)\cdot\left(\sum_{n=0}^\infty\sum_{l=0}^{k-1}\sigma^l\alpha_{n,l} z^{kn}\right)\\
&=f(0)\left(\sum_{m=0}^{k-1}\frac{(k\log z)^m}{m!}\sigma^m\right)\left(\sum_{l=0}^{k-1}a_{l}(z)\sigma^l\right)\\
&=f(0)\cdot\sum_{p=0}^{k-1}\left(\sum_{l=0}^{p}\frac{(k\log z)^{p-l}}{(p-l)!}a_{l}(z)\right)\sigma^p.
\end{align*}
This proves point $(1)$. For the second point, observe that also the functions
\[f_{\pm}(n):=\frac{\Gamma(-\sigma-n)^k}{\Gamma(-\sigma)^k}e^{\pm k\pi i n}
\]satisfy the relation
\[(\sigma+n)^kf_{\pm}(n)=f_{\pm}(n-1).
\]
For the last claim, if we write the solution $\Xi_0$ in the form
\[\Xi_0(z)=z^{-\mu}A(z)\eta z^{R},\quad R\equiv c_1(\mathbb P)\cup (-)\colon H^\bullet(\mathbb P,\mathbb C)\to H^\bullet(\mathbb P,\mathbb C),
\]by Proposition \ref{18.07.17-2} it is sufficient to prove that $A(z)$ is holomorphic in $z=0$ and $A(0)=\mathbbm 1$.
From the identity 
\[\Phi(z)=z^{k\sigma}\sum_{n=0}^\infty f(n)z^{kn},
\]we obtain for $1\leq\alpha\leq k$ the relation
\[\vartheta^{\alpha-1}\Phi(z)=z^{k\sigma}\left\{(k\sigma)^{\alpha-1}+\sum_{p=0}^{\alpha-2}\binom{\alpha-1}{p}k^{\alpha-1}\sigma^p\sum_{n=0}^\infty f(n)n^{\alpha-1-p}z^{kn}\right\},
\]and by definition of $A(z)$ we have the identity
\[\sum_{j=1}^k A(z)^\alpha_j\sigma^{j-1}=\frac{1}{k^{\alpha-1}}\left\{(k\sigma)^{\alpha-1}+\sum_{p=0}^{\alpha-2}\binom{\alpha-1}{p}k^{\alpha-1}\sigma^p\sum_{n=0}^\infty f(n)n^{\alpha-1-p}z^{kn}\right\}.
\]This shows that $A(z)$ is holomorphic in $z=0$, and furthermore that $A(0)=\mathbbm 1$.
\endproof

\subsection{Computation of the group $\mathcal{{C}}_0(\mathbb P)$}

Let us introduce the $k\times k$ matrices $J_i$, $i\geq 0$, defined by
\[(J_i)_{ab}:=\delta_{i,a-b}.
\]

\bsh
\begin{teorema}\label{teoctildecpn}
The group ${\mathcal C}_0(\mathbb P)$ is an abelian unipotent algebraic group of dimension $[\frac{k}{2}]$. In particular, the exponential map defines an isomorphism
\[{\mathcal C}_0(\mathbb P)\cong \underbrace{\mathbb C\oplus\dots\oplus\mathbb C}_{[\frac{k}{2}]\text{ copies}}.
\]With respect to the basis $(1,\sigma,\dots, \sigma^{k-1})$ of $H^\bullet(\mathbb P,\mathbb C)$, the group ${\mathcal C}_0(\mathbb P)$ is described as follows
\[{\mathcal C}_0(\mathbb P)=\left\{C\in GL(k,\mathbb C)\colon C=\sum_{i=0}^{k-1}\alpha_iJ_i,\quad\alpha_0=1,\quad 2\alpha_{2n}+\sum_{\substack{i+j=2n\\1\leq i,j}}(-1)^i\alpha_i\alpha_{j}=0, \quad 2\leq 2n\leq k-1\right\}.
\]
\end{teorema}
\esh

\proof If $C\in{\mathcal C}_0(\mathbb P)$, in order to have that $P(z):=z^{\mu}z^{R}Cz^{-R}z^{-\mu}$ is a polynomial in $z$, where $R$ is the operator of classical multiplication by the first Chern class $c_1(\mathbb P)$, the matrix $C$ must be of the form
\[C=\sum_{i=0}^{k-1}\alpha_iJ_i,\quad\alpha_0=1.
\]
We have that $C\in{\mathcal C}_0(\mathbb P)$ if and only if
\[\left(\sum_{i=0}^{k-1}(-1)^i\alpha_iz^iJ_i^T\right)\eta\left(\sum_{i=0}^{k-1}\alpha_iz^iJ_i\right)=\eta.
\]
The l.h.s is equal to
\[\eta+\sum_{i=1}^{k-1}\alpha_iz^i\eta J_i+\sum_{i=1}^{k-1}(-1)^i\alpha_iz^iJ_i^T\eta+\sum_{h=2}^{2k-2}\left(\sum_{\substack{i+j=h\\ 1\leq i,j\leq k-1}}(-1)^i\alpha_i\alpha_jJ^T_i\eta J_j\right)z^h,
\]
and using the relations
\begin{align}
\eta J_i&=J_i^T\eta,\\
\label{19.07.17-1}
 \left(J_iJ_j\right)_{ab}&=\delta_{i+j,a-b}=(J_{i+j})_{ab},\\
  J_h&=0\quad\text{ if }h\geq k,
\end{align}
we obtain the equation
\[\sum_{\substack{1\leq i\leq k-1\\ i \text{ even}}}2\alpha_iz^i\eta J_i+\sum_{h=2}^{k-1}\left(\sum_{\substack{i+j=h\\ 1\leq i,j\leq k-1}}(-1)^i\alpha_i\alpha_j\right)z^h\eta J_h=0.
\]
So, we have the following constraints on the constants $\alpha_i$'s: 
\begin{align*}
2\alpha_2-\alpha_1^2&=0,\\
2\alpha_{4}-2\alpha_1\alpha_3+\alpha_2^2&=0,\\
2\alpha_6-2\alpha_1\alpha_5+2\alpha_2\alpha_4-\alpha_3^2&=0,\\
\dots&\\
2\alpha_{2n}+\sum_{\substack{i+j=2n\\1\leq i,j}}(-1)^i\alpha_i\alpha_{j}&=0, \quad 2\leq 2n\leq k-1.
\end{align*}
The Lie algebra of the group is 
\[{\frak g}_0(\mathbb P)=\left\{C\in \frak{gl}(k,\mathbb C)\colon C=\sum_{i=0}^{k-1}\alpha_iJ_i,\quad \alpha_{\rm even}=0\right\},
\]which is abelian by \eqref{19.07.17-1}, coherently with Theorem \ref{03.08.17-1}. In characteristic zero the structure of unipotent abelian group is well-known: in particular, the exponential map defines an isomorphism of groups (see \cite{demagabr}, Ch. IV.2.4 Proposition 4.1).
\endproof

The following result immediately follows from Theorem \ref{03.08.17-1} and Theorem \ref{04.08.17-1}.

\bsh
\begin{cor}\label{corctildecpn}
The groups ${\mathcal C}_0(\mathbb P)$ and the identity component ${\rm Isom}_\mathbb C(K_0(\mathbb P)_\mathbb C,\chi)_0$ are isomorphic.
\end{cor}
\esh

Remarkably, notice that the equations obtained above for the group ${\mathcal C}_0(\mathbb P)$ essentially coincide with those obtained by A. Gorodentsev for ${\rm Isom}_\mathbb C(K_0(\mathbb P)_\mathbb C,\chi)_0$ in \cite{Go3, Go2.2}.

\subsection{Computation of the Central Connection Matrix}
Using the labeling of the canonical coordinates $u_1,\dots, u_n$ introduced in the section \ref{notationsproj}, we introduce the corresponding Stokes' rays:
\[ R_{rs}:=\left\{z=-i\rho(\overline{u_r}-\overline{u_s}),\quad\rho>0\right\}.
\]
At a generic point of the small quantum cohomology $(0,t^2,0,\dots,0)$, we have
\begin{align*}-i(\overline{u_r}-\overline{u_s})&=-ikq^{-\frac{1}{k}}\left(e^{-\frac{2\pi i (r-1)}{k}}-e^{-\frac{2\pi i (s-1)}{k}}\right)\\
&=2kq^{-\frac{1}{k}}\sin\left(\frac{\pi}{k}(s-r)\right)\exp\left(i\left[\frac{2\pi}{k}-\frac{\pi}{k}(r+s)\right]\right).
\end{align*}
So if $1\leq r<s\leq k$ the Stokes' rays at a generic point $(0,t^2,0,\dots,0)$ are

\begin{equation}
\label{stokesrayscpn}
 R_{rs}=\left\{z\colon z=\rho\exp\left(i\left[\frac{2\pi}{k}-\frac{\pi}{k}(r+s)-\frac{\Im (t^2)}{k}\right]\right)\right\},
\end{equation}
\[ R_{sr}=- R_{rs}.
\]
\begin{figure}[ht!]
\centering
\def\svgscale{0.8}
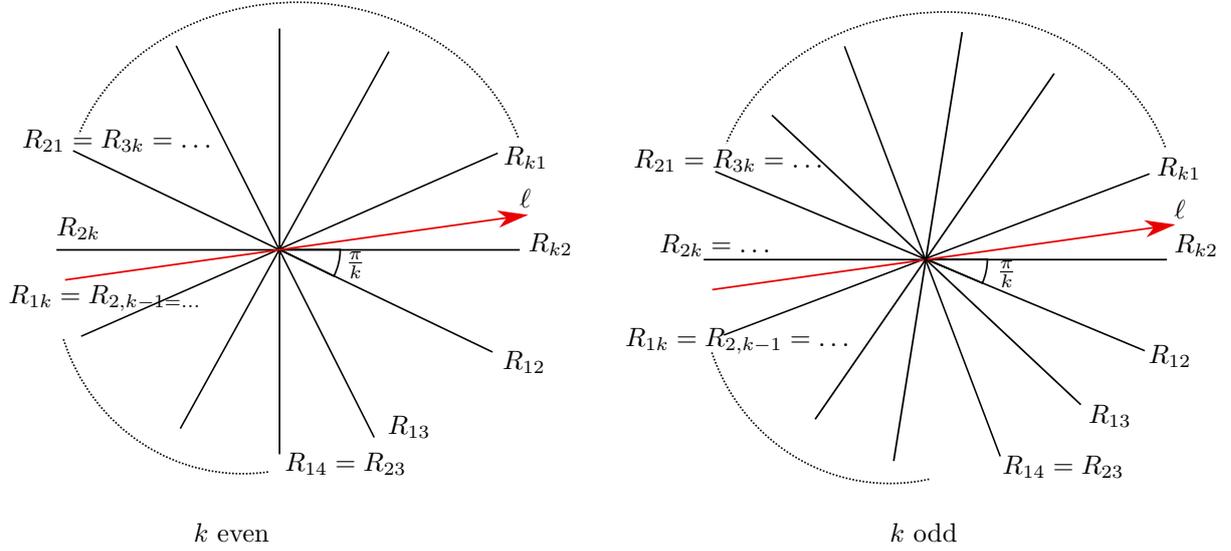
\caption{Configuration of Stokes rays for $k$ odd and $k$ even.}
\end{figure}

Since we want compute the central connection matrix at $t^2=0$ we have to fix an admissible line: following \cite{guzzetti1} we choose a line $\ell$ with slope $0<\phi<\frac{\pi}{k}$.

\bsh
\begin{prop}[\cite{guzzetti1}]\label{guz1} Let 
\[g(z):=
\begin{sistema}
\\
\frac{1}{(2\pi)^{\frac{k+1}{2}}}\int_\Lambda \Gamma(-s)^k z^{ks}ds,\quad\text{k even}\\
\\
\\
\frac{1}{(2\pi)^{\frac{k+1}{2}}i}\int_\Lambda \Gamma(-s)^k e^{-i\pi s}z^{ks}ds,\quad\text{k odd}\\
\\
\end{sistema}
\]where $\Lambda$ is a straight line going from $-c-i\infty$ to $-c+i\infty$, $c>0$. Fix a line $\ell$ with slope $0<\epsilon<\frac{\pi}{k}$. Then, for $k$ even, the fundamental solution $\Xi_R$, having asymptotic expansion
\[
\Xi=\eta \Psi^{-1}\left(\mathbbm 1+O\left(\frac{1}{z}\right)\right)e^{zU}\quad \text{on }\Pi_R,
\] 
is constructed by means of \eqref{17agosto2018-8} from the following basis of solutions of \eqref{eq3} 
\[
\Phi_R^{\alpha}(z)=
\begin{sistema}
(-1)^{\frac{k}{2}+1-\alpha}\sum_{h=0}^{k+1-2\alpha}(-1)^h\binom{k}{h}g\left(ze^{-i\pi+(\alpha+h-1)\frac{2\pi i}{k}}\right)\quad\text{if }1\leq \alpha\leq\frac{k}{2},\\
\\
(-1)^{\alpha-\frac{k}{2}-1}g\left(ze^{\frac{2\pi i}{k}(\alpha-\frac{k}{2}-1)}\right)\quad\text{if }\frac{k}{2}+1\leq\alpha\leq k.
\end{sistema}
\]For $k$ odd, he fundamental solution $\Xi_R$ is constructed with the basis
 \[\Phi_R^{\alpha}(z)=
\begin{sistema}
(-1)^{\frac{k+1}{2}-\alpha}\sum_{h=0}^{k+1-2\alpha}(-1)^h\binom{k}{h}g\left(ze^{\frac{2\pi i}{k}\left(\alpha-\frac{k+1}{2}+h\right)}\right)\quad\text{if }1\leq\alpha\leq\frac{k+1}{2},\\
\\
(-1)^{\alpha-\frac{k+1}{2}}g\left(ze^{\frac{2\pi i}{k}\left(\alpha-\frac{k+1}{2}\right)}\right)\quad\text{if }\frac{k+1}{2}+1\leq\alpha\leq k.
\end{sistema}
\]
\end{prop}
\esh
In \cite{guzzetti1}, the above basis of solutions is collected into a  row vector, which in transposed form looks as follows. For $k$ even,   
\[
\boldsymbol{\Phi}_R(z)^T:=
\begin{pmatrix}
(-1)^{\frac{k}{2}}\left(g(ze^{-i\pi})-\binom{k}{1}g\left(ze^{-i\pi+i\frac{2\pi}{k}}\right)+\dots-\binom{k}{k-1}g\left(ze^{i(\pi-\frac{2\pi}{k})}\right)\right)\\
\vdots\\
g\left(ze^{-\frac{4\pi i}{k}}\right)-\binom{k}{1}g\left(ze^{-\frac{2\pi i}{k}}\right)+\binom{k}{2}g(z)-\binom{k}{3}g\left(ze^{\frac{2\pi i}{k}}\right)\\
-g\left(ze^{-\frac{2\pi i}{k}}\right)+\binom{k}{1}g(z)\\
g(z)\\
-g\left(ze^{\frac{2\pi i}{k}}\right)\\
g\left(ze^{\frac{4\pi i}{k}}\right)\\
\vdots\\
(-1)^{\frac{k}{2}-1}g\left(ze^{\frac{2\pi i}{k}(\frac{k}{2}-1)}\right)
\end{pmatrix}
\]
where  the entry corresponding to $g(z)$ is the $n(k):=\left(\frac{k}{2}+1\right)$-th one.

For $k$ odd, 
\[\boldsymbol{\Phi}_R(z)^T:=
\begin{pmatrix}
(-1)^{\frac{k-1}{2}}\left(g\left(ze^{-\frac{2\pi i}{k}\left(\frac{k-1}{2}\right)}\right)-\binom{k}{1}g\left(ze^{-\frac{2\pi i}{k}\left(\frac{k-3}{2}\right)}\right)+\dots+\binom{k}{k-1}g\left(ze^{\frac{2\pi i}{k}\left(\frac{k-1}{2}\right)}\right)\right)\\
\vdots\\
g\left(ze^{-\frac{4\pi i}{k}}\right)-\binom{k}{1}g\left(ze^{-\frac{2\pi i}{k}}\right)+\binom{k}{2}g(z)-\binom{k}{3}g\left(ze^{\frac{2\pi i}{k}}\right)+\binom{k}{4}g\left(ze^{\frac{4\pi i}{k}}\right)\\
-g\left(ze^{-\frac{2\pi i}{k}}\right)+\binom{k}{1}g(z)-\binom{k}{2}g\left(ze^{\frac{2\pi i}{k}}\right)\\
g(z)\\
-g\left(ze^{\frac{2\pi i}{k}}\right)\\
g\left(ze^{\frac{4\pi i}{k}}\right)\\
\vdots\\
(-1)^{\frac{k-1}{2}}g\left(ze^{\frac{2\pi i}{k}\left(\frac{k-3}{2}\right)}\right)
\end{pmatrix}
\]and the entry corresponding to $g(z)$ is the $n(k):=\frac{k+1}{2}$-th one.

%
Now we compute the entries of the central connection matrix. We will denote by $\Phi_{\rm top}(z)$ the solution of Proposition \ref{19.07.17-2} corresponding to the choice $f(0)=1$. The computations will be done in cases, depending on the parity of $k$. \newline

{\bf CASE $k$ EVEN:} If $1\leq\alpha\leq \frac{k}{2}$, we have that

\begin{align*}
&\Phi^\alpha_R(z)=-\frac{2\pi i(-1)^{\frac{k}{2}+1-\alpha}}{(2\pi)^{\frac{k+1}{2}}}\sum_{h=0}^{k+1-2\alpha}\sum_{n=0}^\infty(-1)^h\binom{k}{h}\underset{s=n}{\operatorname{res}}\left(\Gamma(-s)^{k}z^{ks}e^{(\alpha+h-\frac{k}{2}-1)2\pi is}\right)ds\\
&=\frac{ i(-1)^{\frac{k}{2}-\alpha}}{(2\pi)^{\frac{k-1}{2}}}\sum_{h=0}^{k+1-2\alpha}\sum_{n=0}^\infty(-1)^h\binom{k}{h}\underset{w=0}{\operatorname{res}}\left(\Gamma(-w-n)^{k}z^{k(w+n)}e^{(\alpha+h-\frac{k}{2}-1)2\pi i(w+n)}\right)dw\\
&=\frac{i(-1)^{\frac{k}{2}-\alpha}}{(2\pi)^{\frac{k-1}{2}}}\sum_{h=0}^{k+1-2\alpha}\sum_{n=0}^\infty(-1)^h\binom{k}{h}\int_{\mathbb P}\left(\frac{\Gamma(-\sigma-n)^{k}z^{k(\sigma+n)}}{\Gamma(-\sigma)^k}\Gamma(1-\sigma)^ke^{(\alpha+h-\frac{k}{2}-1)2\pi i\sigma}\right)\\
&=\frac{i(-1)^{\frac{k}{2}-\alpha}}{(2\pi)^{\frac{k-1}{2}}}\sum_{h=0}^{k+1-2\alpha}(-1)^h\binom{k}{h}\int_{\mathbb P}\left\{\left(\Phi_{\rm top}(z)\cup e^{-k\pi i\sigma}\right)\cup\widehat{\Gamma}^{-}(\mathbb P)\cup {\rm Ch}(\mathcal O(\alpha+h-1))\right\}.\\
\end{align*}

If $\frac{k}{2}+1\leq\alpha\leq k$

\begin{align*}
&\Phi_R^\alpha(z)=-\frac{2\pi i (-1)^{\alpha-\frac{k}{2}-1}}{(2\pi)^{\frac{k+1}{2}}}\sum_{n=0}^\infty\underset{s=n}{\operatorname{res}}\left(\Gamma(-s)^kz^{ks}e^{(\alpha-\frac{k}{2}-1)2\pi is}\right)ds\\
&=\frac{i(-1)^{\alpha-\frac{k}{2}}}{(2\pi)^\frac{k-1}{2}}\sum_{n=0}^\infty\underset{w=0}{\operatorname{res}}\left(\Gamma(-w-n)^kz^{k(w+n)}e^{(\alpha-\frac{k}{2}-1)2\pi i(w+n)}\right)dw\\
&=\frac{i(-1)^{\alpha-\frac{k}{2}}}{(2\pi)^\frac{k-1}{2}}\sum_{n=0}^\infty\int_{\mathbb P}\left(\frac{\Gamma(-\sigma-n)^kz^{k(\sigma+n)}}{\Gamma(-\sigma)^k}\Gamma(1-\sigma)^{k}e^{(\alpha-\frac{k}{2}-1)2\pi i\sigma}\right)\\
&=\frac{i(-1)^{\alpha-\frac{k}{2}}}{(2\pi)^\frac{k-1}{2}}\int_{\mathbb P}\left\{\left(\Phi_{\rm top}(z)\cup e^{-k\pi i\sigma}\right)\cup\widehat{\Gamma}^{-}(\mathbb P)\cup {\rm Ch}(\mathcal O(\alpha-1))\right\}.\\
\end{align*}


{\bf CASE $k$ ODD:} If  $1\leq\alpha\leq\frac{k+1}{2}$ we have
\begin{align*}
&\Phi_R^\alpha(z)=-\frac{2\pi i(-1)^{\frac{k+1}{2}-\alpha}}{(2\pi)^{\frac{k+1}{2}}i}\sum_{n=0}^\infty\sum_{h=0}^{k+1-2\alpha}(-1)^h\binom{k}{h}\underset{s=n}{\operatorname{res}}\left(\Gamma(-s)^ke^{-i\pi s}z^{ks}e^{2\pi i s\left(\alpha-\frac{k+1}{2}+h\right)}\right)ds\\
&=\frac{(-1)^{\frac{k-1}{2}-\alpha}}{(2\pi)^{\frac{k-1}{2}}}\sum_{n=0}^\infty\sum_{h=0}^{k+1-2\alpha}(-1)^h\binom{k}{h}\underset{w=0}{\operatorname{res}}\left(\Gamma(-w-n)^ke^{-i\pi (w+n)}z^{k(w+n)}e^{2\pi i (w+n)\left(\alpha-\frac{k+1}{2}+h\right)}\right)dw\\
&=\frac{(-1)^{\frac{k-1}{2}-\alpha}}{(2\pi)^{\frac{k-1}{2}}}\sum_{n=0}^\infty\sum_{h=0}^{k+1-2\alpha}(-1)^h\binom{k}{h}\int_{\mathbb P}\left(-\frac{\Gamma(-\sigma-n)^k}{\Gamma(-\sigma)^k}\Gamma(1-\sigma)^ke^{-i\pi (\sigma+n)}z^{k(\sigma+n)}e^{2\pi i (\sigma+n)\left(\alpha-\frac{k+1}{2}+h\right)}\right)\\
&=\frac{(-1)^{\frac{k+1}{2}-\alpha}}{(2\pi)^{\frac{k-1}{2}}}\sum_{n=0}^\infty\sum_{h=0}^{k+1-2\alpha}(-1)^h\binom{k}{h}\int_{\mathbb P}\left(\frac{\Gamma(-\sigma-n)^k}{\Gamma(-\sigma)^k}\Gamma(1-\sigma)^ke^{-ki\pi (\sigma+n)}z^{k(\sigma+n)}e^{2\pi i \sigma\left(\alpha+h-1\right)}\right)\\
&=\frac{(-1)^{\frac{k+1}{2}-\alpha}}{(2\pi)^{\frac{k-1}{2}}}\sum_{h=0}^{k+1-2\alpha}(-1)^h\binom{k}{h}\int_{\mathbb P}\left\{\left(\Phi_{\rm top}(z)\cup e^{-ki\pi \sigma}\right)\cup\widehat{\Gamma}^{-}(\mathbb P)\cup{\rm Ch}(\mathcal O(\alpha+h-1))\right\}.
\end{align*}$ $\newline

If $\frac{k+1}{2}+1\leq\alpha\leq k$ we have
\begin{align*}
&\Phi_R^\alpha(z)=-\frac{2\pi i(-1)^{\alpha-\frac{k+1}{2}}}{(2\pi)^{\frac{k+1}{2}}i}\sum_{n=0}^\infty\underset{s=n}{\operatorname{res}}\left(\Gamma(-s)^ke^{-i\pi s}z^{ks}e^{2\pi i s \left(\alpha-\frac{k+1}{2}\right)}\right)ds\\
&=\frac{(-1)^{\alpha-\frac{k-1}{2}}}{(2\pi)^{\frac{k-1}{2}}}\sum_{n=0}^\infty\underset{w=0}{\operatorname{res}}\left(\Gamma(-w-n)^ke^{-i\pi (w+n)}z^{k(w+n)}e^{2\pi i (w+n) \left(\alpha-\frac{k+1}{2}\right)}\right)dw\\
&=\frac{(-1)^{\alpha-\frac{k-1}{2}}}{(2\pi)^{\frac{k-1}{2}}}\sum_{n=0}^\infty\int_{\mathbb P}\left(-\frac{\Gamma(-\sigma-n)^k}{\Gamma(-\sigma)^k}e^{-k(\sigma+n)\pi i}z^{k(\sigma+n)}\Gamma(1-\sigma)^ke^{2\pi i(\sigma+n)(\alpha-1)}\right)\\
&=\frac{(-1)^{\alpha-\frac{k+1}{2}}}{(2\pi)^{\frac{k-1}{2}}}\sum_{n=0}^\infty\int_{\mathbb P}\left(\frac{\Gamma(-\sigma-n)^k}{\Gamma(-\sigma)^k}e^{-kn\pi i}z^{k(\sigma+n)}e^{-ki\pi\sigma}\Gamma(1-\sigma)^ke^{2\pi i\sigma(\alpha-1)}\right)\\
&=\frac{(-1)^{\alpha-\frac{k+1}{2}}}{(2\pi)^{\frac{k-1}{2}}}\int_{\mathbb P}\left\{\left(\Phi_{\rm top}(z)\cup e^{-ki\pi\sigma}\right)\cup\widehat{\Gamma}^{-}(\mathbb P)\cup{\rm Ch}(\mathcal O(\alpha-1))\right\}.\\
\end{align*}

The form $\Phi_{\rm top}(z)\cup e^{-k\pi i\sigma}$ corresponds to the choice of another fundamental basis  $\widetilde{\Xi}_0$ in Levelt form at $z=0$, related to \eqref{fundsol} by the right multiplication by a matrix:

\[
\widetilde{\Xi}_0(z)=\Xi_0(z)\begin{pmatrix}
1&&&&&&\\
-k\pi i&1&&&&&\\
-\frac{k^2\pi^2}{2}&-k\pi i&1&&&&\\
\vdots&&&\ddots&&\\
\frac{(-k\pi i)^m}{m!}&\dots&\dots&\dots&1&&\\
\vdots&&&&&\ddots&\\
\frac{(-k\pi i)^{k-1}}{(k-1)!}&\dots&\dots&\dots&\dots&\dots&1
\end{pmatrix}.
\]

We claim that such a matrix is an element of the group $\mathcal{{C}}_0(\mathbb P)$ (see the previous section). Indeed if
\[\alpha_m:=\frac{(-k\pi i)^m}{m!}
\]
then, for $2\leq 2n\leq k-1$, we have that
\[2\alpha_{2n}+\sum_{\substack{i+j=2n\\1\leq i,j}}(-1)^i\alpha_i\alpha_{j}
=\frac{(-k\pi i)^{2n}}{(2n)!}\left(2+\sum_{j=1}^{2n-1}(-1)^j\binom{2n}{j}\right)=0.
\]
\newline\newline

\subsection{Reduction to Beilinson Form}\label{secbeilinson}

\begin{figure}[ht!]
\centering
\def\svgscale{0.5}
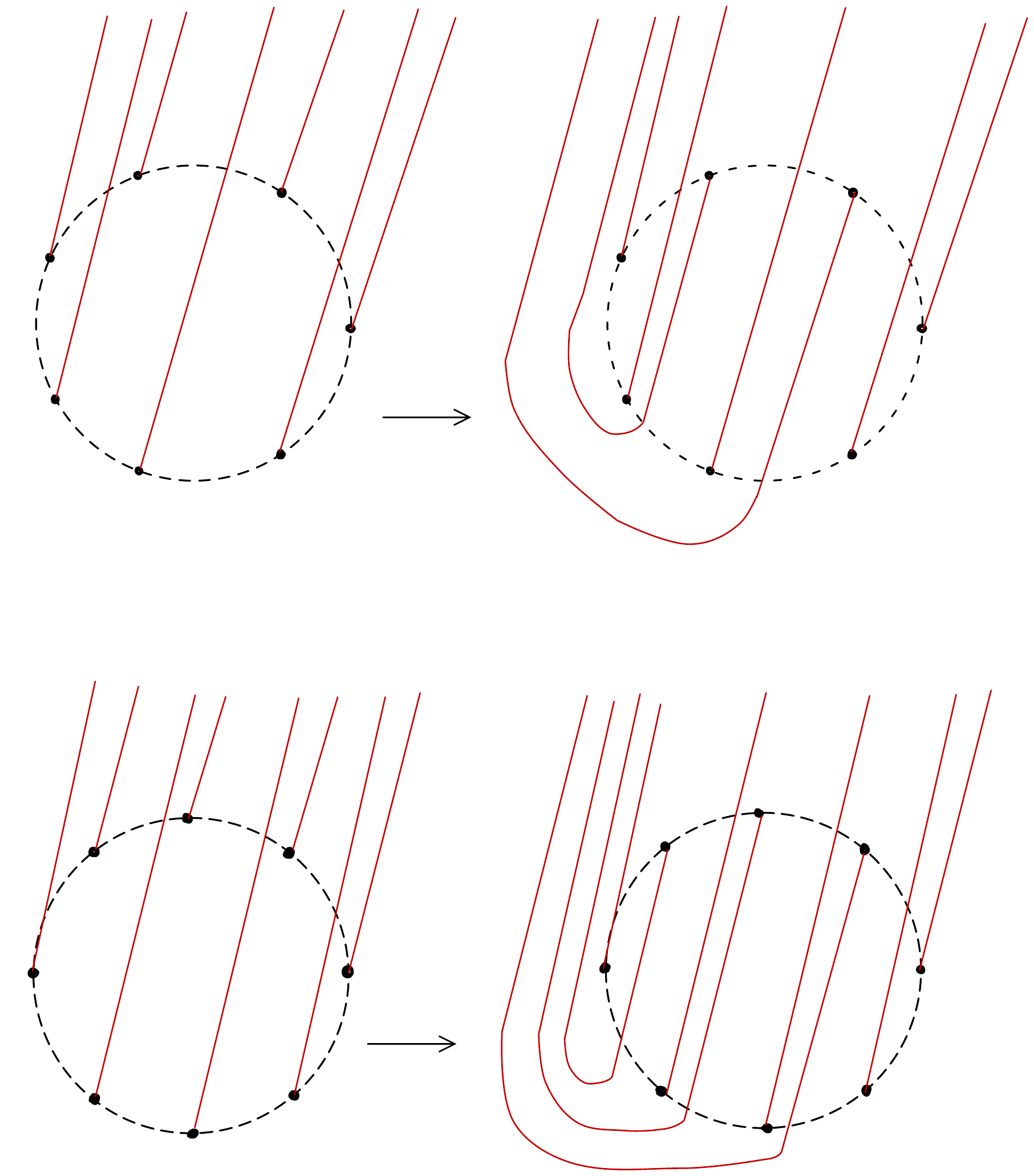
\caption{Action of the braid $\beta$ found by the third author in \cite{guzzetti1}: in the figure above we draw the case $k=7$, below the case $k=8$. Notice that the braid $\beta$ puts the canonical coordinates in counterclockwise order starting from the point $k\exp\left(\frac{2\pi i}{k}\right)$ in the complex plane.}
\label{treccia_davide}
\end{figure}

Let us recall that the canonical coordinates can always be reordered so that the corresponding Stokes matrix is upper triangular (\emph{the lexicographical order }w.r.t the line $\ell$). For the case of quantum cohomology of projective spaces, and for the choice of an admissible line $\ell$ with slope $0<\epsilon<\frac{\pi}{k}$, such an order is the one described in the left part of Figure \ref{treccia_davide}.
The matrices $P$ associated with this permutations are \cite{guzzetti1}

\begin{itemize}
\item for $k$ even
\[P=\left(\begin{array}{ccccccccccccccc}
&&&&&&&1&0&&&&&&\\
&&&&&&1&0&0&&&&&&\\
&&&&&0&0&0&1&0&&&&&\\
&&&&&1&0&0&0&0&&&&&\\
&&&&0&0&0&0&0&1&0&&&&\\
\vdots&\vdots&&\vdots&\vdots&\vdots&\vdots&\vdots&\vdots&\vdots&\vdots&\vdots&&\vdots&\vdots\\
&1&&&&&&&&&&&&0&\\
0&0&&&&&&&&&&&&0&1\\
1&0&&&&&&&&&&&&0&0\\
\end{array}\right)
\]where the 1 on the first row in on the $\frac{k}{2}+1$-th column;
\item for $k$ odd
\[P=\left(\begin{array}{ccccccccccccccc}
&&&&&&0&1&&&&&&&\\
&&&&&&0&0&1&&&&&&\\
&&&&&0&1&0&0&0&&&&&\\
&&&&&0&0&0&0&1&&&&&\\
&&&&0&1&0&0&0&0&0&&&&\\
\vdots&\vdots&&\vdots&\vdots&\vdots&\vdots&\vdots&\vdots&\vdots&\vdots&\vdots&&\vdots&\vdots\\
&1&&&&&&&&&&&&0&\\
0&0&&&&&&&&&&&&0&1\\
1&0&&&&&&&&&&&&0&0\\
\end{array}\right)
\]where the 1 on the first row in os the $\frac{k+1}{2}$-th column.
\end{itemize}
After such a renumeration of $u_1,\dots, u_n$, as a consequence of the computations of the preceding section, the central connection matrix (computed wrt $\Xi_0$ as in \eqref{fundsol}) is, for $k$ even
\[C_{\text{lex}}=\frac{i}{(2\pi)^{\frac{k-1}{2}}}
\begin{pmatrix}
\vdots& \vdots&\vdots&&\vdots\\
\pm \textcyr{D}^0&\mp\textcyr{D}^1&\pm\textcyr{D}^2&\dots&\mp\textcyr{D}^{k-1}\\
\vdots& \vdots&\vdots&&\vdots
\end{pmatrix}\cdot A_k,
\]where:
\begin{itemize} 
\item $\textcyr{D}^j$ is a column vector whose components are the components of the characteristic classes
\[\widehat{\Gamma}^-(\mathbb P)\cup{\rm Ch}(\mathcal O(j))\cup\exp(-\pi i c_1(\mathbb P));
\]
\item the sign $(+)$ is chosen if $\frac{k}{2}-1$ is even, $(-)$ if $\frac{k}{2}-1$ is odd;
\item the matrix $A_k$ is the $k\times k$ matrix
\[A_k:=
\left(\begin{array}{ccccccccccc}
0&0&0&0&0&&0&0&0&0&1\\
\\
&&&&&&0&0&1&0&\binom{k}{1}\\
\\
&&&&&&1&0&\binom{k}{1}&0&\binom{k}{2}\\
\\
&&&&&&\binom{k}{1}&0&\binom{k}{2}&0&\binom{k}{3}\\
&&&&&\dots&*&\vdots&*&\vdots&*\\
0&0&0&1&0&&*&\vdots&*&\vdots&*\\
0&1&0&\binom{k}{1}&0&\dots&*&\vdots&*&\vdots&*\\
1&\binom{k}{1}&0&\binom{k}{2}&0&\dots&*&\vdots&*&\vdots&*\\
0&0&1&\binom{k}{3}&0&&*&\vdots&*&\vdots&*\\
\vdots&\vdots&\vdots&\vdots&\vdots&&*&\vdots&*&\vdots&*\\
\vdots&\vdots&\vdots&\vdots&\vdots&\dots&*&\vdots&*&\vdots&*\\
\vdots&\vdots&\vdots&\vdots&\vdots&&*&\vdots&*&\vdots&*\\
&&&&&&\binom{k}{k-5}&0&\binom{k}{k-4}&0&\binom{k}{k-3}\\
\\
&&&&&&0&1&\binom{k}{k-3}&0&\binom{k}{k-2}\\
\\
&&&&&&0&0&0&1&\binom{k}{k-1}\\
\end{array}\right),
\]where the 1 of the first column is on the $(\frac{k}{2}+1)$-th row.
\end{itemize}$ $\newline

Analogously, if $k$ is odd then the central connection matrix in the lexicographical order is
\[C_{\text{lex}}=\frac{1}{(2\pi)^{\frac{k-1}{2}}}
\begin{pmatrix}
\vdots& \vdots&\vdots&&\vdots\\
\pm \textcyr{D}^0&\mp\textcyr{D}^1&\pm\textcyr{D}^2&\dots&\pm\textcyr{D}^{k-1}\\
\vdots& \vdots&\vdots&&\vdots
\end{pmatrix}\cdot A_k,
\]where:
\begin{itemize}
\item $\textcyr{D}^j$ is as before;
\item the sign $(+)$ is chosen if $\frac{k-1}{2}$ is even, $(-)$ if $\frac{k-1}{2}$ is odd;
\item the matrix $A_k$ is the $k\times k$ matrix
\[A_k:=\left(\begin{array}{cccccccccccc}
0&0&0&0&0&&&0&0&0&0&1\\
\\
&&&&&&&0&0&1&0&\binom{k}{1}\\
\\
&&&&&&&1&0&\binom{k}{1}&0&\binom{k}{2}\\
\\
&&&&&&&\binom{k}{1}&0&\binom{k}{2}&0&\binom{k}{3}\\
&&&&&&\dots&*&\vdots&*&\vdots&*\\
0&0&0&0&1&0&&*&\vdots&*&\vdots&*\\
0&0&1&0&\binom{k}{1}&0&\dots&*&\vdots&*&\vdots&*\\
1&0&\binom{k}{1}&0&\binom{k}{2}&0&\dots&*&\vdots&*&\vdots&*\\
0&1&\binom{k}{2}&0&\binom{k}{3}&0&&*&\vdots&*&\vdots&*\\
0&0&0&1&\binom{k}{4}&0&&*&\vdots&*&\vdots&*\\
0&0&0&0&0&1&&*&\vdots&*&\vdots&*\\
\vdots&\vdots&\vdots&\vdots&\vdots&&&*&\vdots&*&\vdots&*\\
\vdots&\vdots&\vdots&\vdots&\vdots&&\dots&*&\vdots&*&\vdots&*\\
\vdots&\vdots&\vdots&\vdots&\vdots&&&*&\vdots&*&\vdots&*\\
&&&&&&&\binom{k}{k-5}&0&\binom{k}{k-4}&0&\binom{k}{k-3}\\
\\
&&&&&&&0&1&\binom{k}{k-3}&0&\binom{k}{k-2}\\
\\
&&&&&&&0&0&0&1&\binom{k}{k-1}\\
\end{array}\right),
\]where the 1 of the first column is in the $\frac{k+1}{2}$-th row.
\end{itemize}

\begin{prop}[\cite{guzzetti1}]The action of the braid
\[\beta:=(\beta_{k-5,k-4}\beta_{k-6,k-5}\dots\beta_{12})(\beta_{k-6,k-5}\beta_{k-7,k-6}\dots\beta_{23})(\beta_{k-7,k-6}\dots\beta_{34})\dots
\]
\[\dots\beta_{\frac{k}{2}-2,\frac{k}{2}-1}(\beta_{k-3,k-2}\beta_{k-4,k-3}\dots\beta_{12})
\]for $k$ even, and
\[\beta:=(\beta_{k-5,k-4}\beta_{k-6,k-5}\dots\beta_{12})(\beta_{k-6,k-5}\beta_{k-7,k-6}\dots\beta_{23})(\beta_{k-7,k-6}\dots\beta_{34})\dots
\]
\[\dots(\beta_{\frac{k-3}{2},\frac{k-1}{2}}\beta_{\frac{k-5}{2},\frac{k-3}{2}})(\beta_{k-3,k-2}\beta_{k-4,k-3}\dots\beta_{12})
\]for $k$ odd, is represented by the multiplication of the matrix
\[ A^\beta(S):=\left(\begin{array}{ccccccccccccc}
0&0&0&0&0&&0&0&0&0&1&0&0\\
\\
&&&&&&0&0&1&0&\binom{k}{1}&0&0\\
\\
&&&&&&1&0&\binom{k}{1}&0&\binom{k}{2}&0&0\\
\\
&&&&&&\binom{k}{1}&0&\binom{k}{2}&0&\binom{k}{3}&0&0\\
&&&&&\dots&*&\vdots&*&\vdots&*&\vdots&\vdots\\
0&0&0&1&0&&*&\vdots&*&\vdots&*&\vdots&\vdots\\
0&1&0&\binom{k}{1}&0&\dots&*&\vdots&*&\vdots&*&\vdots&\vdots\\
1&\binom{k}{1}&0&\binom{k}{2}&0&\dots&*&\vdots&*&\vdots&*&\vdots&\vdots\\
0&0&1&\binom{k}{3}&0&&*&\vdots&*&\vdots&*&\vdots&\vdots\\
\vdots&\vdots&\vdots&\vdots&\vdots&&*&\vdots&*&\vdots&*&\vdots&\vdots\\
\vdots&\vdots&\vdots&\vdots&\vdots&\dots&*&\vdots&*&\vdots&*&\vdots&\vdots\\
\vdots&\vdots&\vdots&\vdots&\vdots&&*&\vdots&*&\vdots&*&\vdots&\vdots\\
&&&&&&\binom{k}{k-7}&0&\binom{k}{k-6}&0&\binom{k}{k-5}&0&0\\
\\
&&&&&&0&1&\binom{k}{k-5}&0&\binom{k}{k-4}&0&0\\
\\
&&&&&&0&0&0&1&\binom{k}{k-3}&0&0\\
\\
&&&&&&0&0&0&0&0&1&0\\
\\
0&0&0&0&0&&0&0&0&0&0&0&1\\
\end{array}\right)
\]
for $k$ even, and

\[A^\beta(S):=\left(\begin{array}{cccccccccccccc}
0&0&0&0&0&&&0&0&0&0&1&0&0\\
\\
&&&&&&&0&0&1&0&\binom{k}{1}&0&0\\
\\
&&&&&&&1&0&\binom{k}{1}&0&\binom{k}{2}&0&0\\
\\
&&&&&&&\binom{k}{1}&0&\binom{k}{2}&0&\binom{k}{3}&0&0\\
&&&&&&\dots&*&\vdots&*&\vdots&*&\vdots&\vdots\\
0&0&0&0&1&0&&*&\vdots&*&\vdots&*&\vdots&\vdots\\
0&0&1&0&\binom{k}{1}&0&\dots&*&\vdots&*&\vdots&*&\vdots&\vdots\\
1&0&\binom{k}{1}&0&\binom{k}{2}&0&\dots&*&\vdots&*&\vdots&*&\vdots&\vdots\\
0&1&\binom{k}{2}&0&\binom{k}{3}&0&&*&\vdots&*&\vdots&*&\vdots&\vdots\\
0&0&0&1&\binom{k}{4}&0&&*&\vdots&*&\vdots&*&\vdots&\vdots\\
0&0&0&0&0&1&&*&\vdots&*&\vdots&*&\vdots&\vdots\\
\vdots&\vdots&\vdots&\vdots&\vdots&&&*&\vdots&*&\vdots&*&\vdots&\vdots\\
\vdots&\vdots&\vdots&\vdots&\vdots&&\dots&*&\vdots&*&\vdots&*&\vdots&\vdots\\
\vdots&\vdots&\vdots&\vdots&\vdots&&&*&\vdots&*&\vdots&*&\vdots&\vdots\\
&&&&&&&\binom{k}{k-7}&0&\binom{k}{k-6}&0&\binom{k}{k-5}&0&0\\
\\
&&&&&&&0&1&\binom{k}{k-5}&0&\binom{k}{k-4}&0&0\\
\\
&&&&&&&0&0&0&1&\binom{k}{k-3}&0&0\\
\\
&&&&&&&0&0&0&0&0&1&0\\
\\
0&0&0&0&0&&&0&0&0&0&0&0&1\\
\end{array}\right)
\]for $k$ odd. Under the action of this braid, the Stokes matrix becomes
\[S^\beta=
\begin{pmatrix}
1&\binom{k}{1}&\binom{k}{2}&\binom{k}{3}&\binom{k}{4}&\dots&-\binom{k}{k-1}\\
& 1&\binom{k}{1}&\binom{k}{2}&\binom{k}{3}&\dots&-\binom{k}{k-2}\\
&&1&\binom{k}{1}&\binom{k}{2}&\dots&-\binom{k}{k-3}\\
&&&1&\binom{k}{1}&\dots&-\binom{k}{k-4}\\
&&&&\ddots&&\vdots\\
&&&&&&1
\end{pmatrix}.
\]
\end{prop}

Observe that, in both cases $k$ even/odd, we obtain that
\begin{equation}\label{treccia} A_k(A^\beta(S))^{-1}=
\begin{pmatrix}
0&\dots&&&&&&1\\
1&0&\dots&&&&&\binom{k}{k-1}\\
&1&0&\dots&&&&\binom{k}{k-2}\\
&&1&0&\dots&&&\binom{k}{k-3}\\
&&&\ddots&&&&\vdots\\
&&&&\ddots&&&\vdots\\
&&&&&&1&\binom{k}{1}
\end{pmatrix}.
\end{equation}
Indeed, observe that
\[A_k=\left(
\begin{array}{ccc|cc}
0&\dots&0&0&1\\ \hline
&&&0&*\\
& X &&\vdots&\vdots\\
&&&0&*\\ \hline
0&\dots&0&1&*
\end{array}
\right)\quad\text{and}\quad A^\beta(S)=
\left(
\begin{array}{ccc|cc}
&&&&\\
& X &&&\\
&&&&\\ \hline
&&&1&\\
&&&&1
\end{array}
\right).
\]
The matrix \eqref{treccia} is the matrix corresponding to the braid 
\[\beta':=\beta_{k-1,k}\beta_{k-2,k-1}\dots\beta_{12},
\]that is
\[A^{\beta'}(S^\beta).
\]
This is easily seen from the fact that 
\[\begin{pmatrix}
1&&&&&&&&\\
&\ddots&&&&&&&\\
&&1&&&&&&\\
&&&0&1&&&&\\
&&&1&x_1&&&&\\
&&&&&1&&&\\
&&&&&&1&&\\
&&&&&&&\ddots&\\
&&&&&&&&1
\end{pmatrix}
\begin{pmatrix}
1&&&&&&&&\\
&\ddots&&&&&&&\\
&&1&&&&&&\\
&&&1&&&\\
&&&&0&1&&&\\
&&&&1&x_2&&&\\
&&&&&&1&&\\
&&&&&&&\ddots&\\
&&&&&&&&1
\end{pmatrix}=\]\[
\begin{pmatrix}
1&&&&&&&&\\
&\ddots&&&&&&&\\
&&1&&&&&&\\
&&&0&0&1&&&\\
&&&1&0&x_1&&&\\
&&&0&1&x_2&&&\\
&&&&&&1&&\\
&&&&&&&\ddots&\\
&&&&&&&&1
\end{pmatrix}
\]and that $A^{\beta_1\beta_2}(S)=A^{\beta_2}(S^{\beta_1})A^{\beta_1}(S)$. As a consequence, we have that
\[A^{\beta'}(S^\beta)=
\begin{pmatrix}
0&\dots&&&&&&1\\
1&0&\dots&&&&&*\\
&1&0&\dots&&&&*\\
&&1&0&\dots&&&*\\
&&&\ddots&&&&\vdots\\
&&&&\ddots&&&*\\
&&&&&&1&*
\end{pmatrix}
\]and the entries $*$ are exactly those of the $k$-th column of $A^{\beta_{i,i+1}}(S^\beta)$, from the top to the bottom, namely 
\begin{align*}
-S^\beta_{1,k}&=\binom{k}{k-1}\\
-S^\beta_{2,k}&=\binom{k}{k-2}\\
\dots\\
-S^\beta_{k-1,k}&=\binom{k}{1}.
\end{align*}

We have thus obtained the following

\bsh
\begin{teorema}\label{risultato1}Consider the central connection matrix for the quantum cohomology of $\mathbb{P}^{k-1}_{\mathbb C}$,  connecting the fundamental matrix solution $\Xi_R$ from Proposition \ref{guz1} with the solution $\Xi_0$ in  \eqref{fundsol}  and set it in the lexicographical form $C_{\operatorname{lex}}$. After the action of the braid
\[\beta\beta':=(\beta_{k-5,k-4}\beta_{k-6,k-5}\dots\beta_{12})(\beta_{k-6,k-5}\beta_{k-7,k-6}\dots\beta_{23})(\beta_{k-7,k-6}\dots\beta_{34})\dots\]
\[\dots\beta_{\frac{k}{2}-2,\frac{k}{2}-1}(\beta_{k-3,k-2}\beta_{k-4,k-3}\dots\beta_{12})(\beta_{k-1,k}\beta_{k-2,k-1}\dots\beta_{12})
\]for $k$ even, and
\[\beta\beta':=(\beta_{k-5,k-4}\beta_{k-6,k-5}\dots\beta_{12})(\beta_{k-6,k-5}\beta_{k-7,k-6}\dots\beta_{23})(\beta_{k-7,k-6}\dots\beta_{34})\dots
\]
\[\dots(\beta_{\frac{k-3}{2},\frac{k-1}{2}}\beta_{\frac{k-5}{2},\frac{k-3}{2}})(\beta_{k-3,k-2}\beta_{k-4,k-3}\dots\beta_{12})(\beta_{k-1,k}\beta_{k-2,k-1}\dots\beta_{12})
\]for $k$ odd, the central connection matrix is
\[C_{\text{lex}}=\frac{i}{(2\pi)^{\frac{k-1}{2}}}
\begin{pmatrix}
\vdots& \vdots&\vdots&&\vdots\\
\pm \textnormal{\textcyr{D}}^0&\mp\textnormal{\textcyr{D}}^1&\pm\textnormal{\textcyr{D}}^2&\dots&\mp\textnormal{\textcyr{D}}^{k-1}\\
\vdots& \vdots&\vdots&&\vdots
\end{pmatrix},
\]for $k$ even, and
\[C_{\text{lex}}=\frac{1}{(2\pi)^{\frac{k-1}{2}}}
\begin{pmatrix}
\vdots& \vdots&\vdots&&\vdots\\
\pm \textnormal{\textcyr{D}}^0&\mp\textnormal{\textcyr{D}}^1&\pm\textnormal{\textcyr{D}}^2&\dots&\mp\textnormal{\textcyr{D}}^{k-1}\\
\vdots& \vdots&\vdots&&\vdots
\end{pmatrix},
\]for $k$ odd. Here
\begin{itemize}
\item$\textnormal{\textcyr{D}}^j$ is a column vector whose components are the components of the characteristic classes
\[\widehat{\Gamma}^-(\mathbb P)\cup{\rm Ch}(\mathcal O(j))\cup\exp(-\pi i c_1(\mathbb P));
\]
\item if $k$ is even, the sign $(+)$ is chosen if $\frac{k}{2}-1$ is even, $(-)$ if $\frac{k}{2}-1$ is odd;
\item if $k$ is odd, the sign $(+)$ is chosen if $\frac{k-1}{2}$ is even, $(-)$ if $\frac{k-1}{2}$ is odd.
\end{itemize}
The corresponding Stokes matrix (using the identity (5) of Theorem \ref{19.03.18-1}) is in the canonical form
\[s_{ij}=\binom{k}{i-j},\quad i<j.
\]After the conjugation by
\[(-1)^{\frac{k}{2}-1}\operatorname{diag}(1,-1,1,-1,\dots,1,-1)
\]if $k$ is even, or by
\[(-1)^{\frac{k-1}{2}}\operatorname{diag}(1,-1,1,-1,\dots,1,-1,1)
\] if $k$ is odd, the central connection matrix is in the canonical form
 \[C_{\text{lex}}=\frac{i}{(2\pi)^{\frac{k-1}{2}}}
\begin{pmatrix}
\vdots& \vdots&\vdots&&\vdots\\
\textnormal{\textcyr{D}}^0&\textnormal{\textcyr{D}}^1&\textnormal{\textcyr{D}}^2&\dots&\textnormal{\textcyr{D}}^{k-1}\\
\vdots& \vdots&\vdots&&\vdots
\end{pmatrix},
\]for $k$ even, and
\[C_{\text{lex}}=\frac{1}{(2\pi)^{\frac{k-1}{2}}}
\begin{pmatrix}
\vdots& \vdots&\vdots&&\vdots\\
\textnormal{\textcyr{D}}^0&\textnormal{\textcyr{D}}^1&\textnormal{\textcyr{D}}^2&\dots&\textnormal{\textcyr{D}}^{k-1}\\
\vdots& \vdots&\vdots&&\vdots
\end{pmatrix},
\]for $k$ odd. The corresponding Stokes matrix is in the form
\[s_{ij}=(-1)^{j-i}\binom{k}{i-j},\quad i<j.
\]
Namely, it is the inverse of the Gram matrix $\chi(\mathcal{O}(i-1),\mathcal{O}(j-1))$ with $i,j=1,\dots,k$.
\end{teorema}
\esh

%

\begin{figure}[ht!]
\centering
\def\svgscale{0.5}
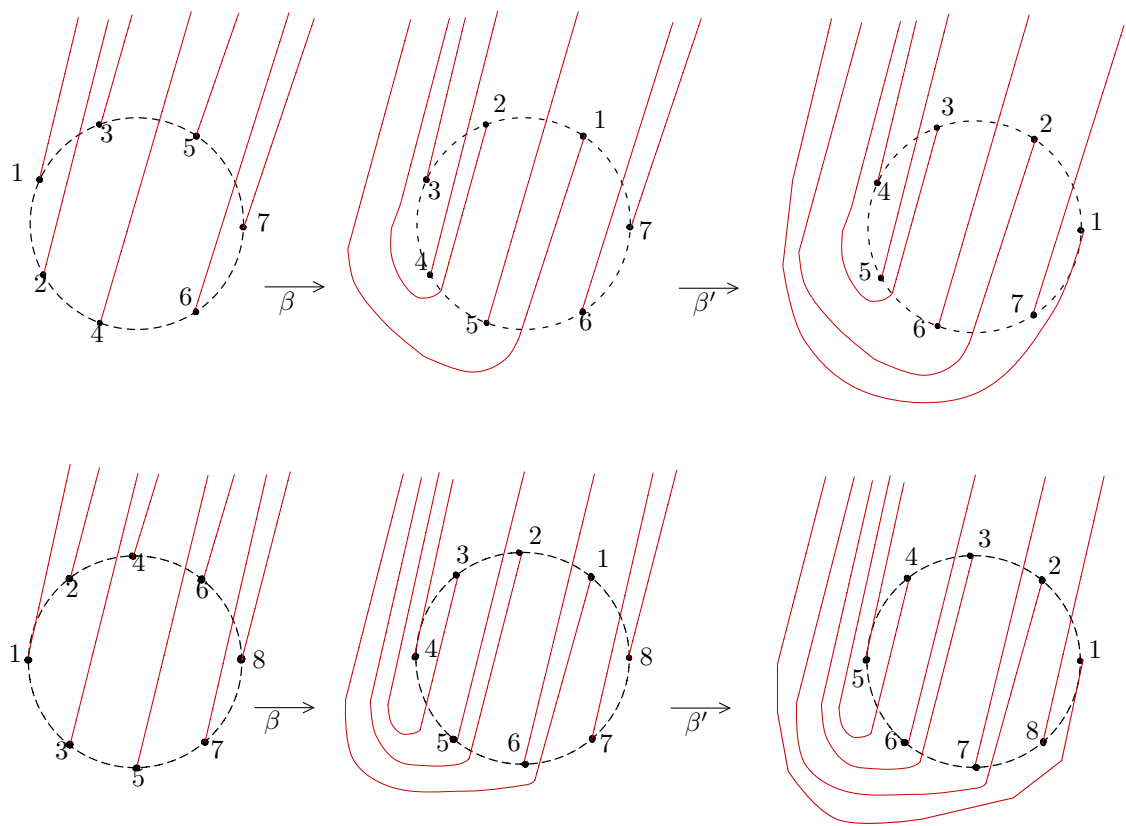
\caption{Action of the braid $\beta$ found by the third author in \cite{guzzetti1}, followed by the action of the remaining braid $\beta'$ determined above: in the figure above we draw the case $k=7$, below the case $k=8$. Notice that the braid $\beta$ puts the canonical coordinates in counterclockwise order starting from the point $k\exp\left(\frac{2\pi i}{k}\right)$ in the complex plane, and the braid $\beta'$ rearrange them in counterclockwise order starting from the point $k$ of the complex plane.}
\label{treccia_davide_completata}
\end{figure}

\begin{oss}\label{rktrecciadavidecompleta}
Notice that, the braid $\beta$ found by the third author in \cite{guzzetti1} puts the canonical coordinates in cyclic counterclockwise order (see Figure \ref{treccia_davide}). If we further act with the above braid $\beta'$, then the canonical coordinates dispose in cyclic counterclockwise order starting from the point $k$ in the complex plane (see Figure \ref{treccia_davide_completata}).
\end{oss}

\subsection{Mutations of the Exceptional Collections}$\quad $

\bsh
\begin{lemma}\label{bellatreccia}
The computed braid can be rewritten as the product
\[\beta\beta'=\beta_{12}(\beta_{34}\beta_{23}\beta_{12})(\beta_{56}\beta_{45}\beta_{34}\beta_{23}\beta_{12})\dots(\beta_{k-1,k}\beta_{k-2,k-1}\dots\beta_{12})
\]for $k$ even,
\[\beta\beta'=(\beta_{23}\beta_{12})(\beta_{45}\beta_{34}\beta_{23}\beta_{12})\dots(\beta_{k-1,k}\beta_{k-2,k-1}\dots\beta_{12})
\]for $k$ odd.
\end{lemma}
\esh

\proof[Proof]Consider the case $k$ even.
The only thing that we have to prove is that the braid
\begin{equation}\label{treccia1}(\beta_{k-5,k-4}\beta_{k-6,k-5}\dots\beta_{12})(\beta_{k-6,k-5}\beta_{k-7,k-6}\dots\beta_{23})(\beta_{k-7,k-6}\dots\beta_{34})\dots\beta_{\frac{k}{2}-2,\frac{k}{2}-1}\end{equation}
is equal to
\[\beta_{12}(\beta_{34}\beta_{23}\beta_{12})\dots(\beta_{k-5,k-4}\dots\beta_{12}).
\]
Note that the braid above ends with the product
\[\dots\beta_{\frac{k}{2}-2,\frac{k}{2}-1}\beta_{\frac{k}{2}-3,\frac{k}{2}-2}\beta_{\frac{k}{2}-2,\frac{k}{2}-1}.
\]
By the ``Yang--Baxter''  braid equations 
\begin{equation}
\label{10novembre2018-3}
\beta_{i,i+1}\beta_{i+1,i+2}\beta_{i,i+1}=\beta_{i+1,i+2}\beta_{i,i+1}\beta_{i+1,i+2},
\end{equation}
this product is equal to
\[\dots\beta_{\frac{k}{2}-3,\frac{k}{2}-2}\beta_{\frac{k}{2}-2,\frac{k}{2}-1}\beta_{\frac{k}{2}-3,\frac{k}{2}-2}.
\]Because of commutation relations, we can shift the first term on the left till we find 
\[\dots\beta_{\frac{k}{2}-3,\frac{k}{2}-2}\beta_{\frac{k}{2}-4,\frac{k}{2}-3}\beta_{\frac{k}{2}-3,\frac{k}{2}-2}\dots,
\]which is equal to
\[\dots\beta_{\frac{k}{2}-4,\frac{k}{2}-3}\beta_{\frac{k}{2}-3,\frac{k}{2}-2}\beta_{\frac{k}{2}-4,\frac{k}{2}-3}\dots.
\]
Again, starting from the first term, we can shift it on the left (until the commutation law allows), then use Yang--Baxter relations. Continuing this procedure, at the end we have eliminated the last term of \eqref{treccia1}, and we obtain a new first term:
\[\beta_{12}(\beta_{k-5,k-4}\beta_{k-6,k-5}\dots\beta_{12})(\beta_{k-6,k-5}\beta_{k-7,k-6}\dots\beta_{23})(\beta_{k-7,k-6}\dots\beta_{34})\dots\]\[\dots(\beta_{\frac{k}{2}-1,\frac{k}{2}}\beta_{\frac{k}{2}-2,\frac{k}{2}-1}\beta_{\frac{k}{2}-3,\frac{k}{2}-2}).
\]Now we continue the procedure of elimination of the last braid: we start from its first term, i.e. $\beta_{\frac{k}{2}-1,\frac{k}{2}}$, we shift it on the left, use Yang- Baxter relations, and so on, till we find
\[\beta_{12}\beta_{34}(\beta_{k-5,k-4}\beta_{k-6,k-5}\dots\beta_{12})(\beta_{k-6,k-5}\beta_{k-7,k-6}\dots\beta_{23})(\beta_{k-7,k-6}\dots\beta_{34})\dots\]\[\dots(\beta_{\frac{k}{2}-2,\frac{k}{2}-1}\beta_{\frac{k}{2}-3,\frac{k}{2}-2}).
\]Applying again the same procedure, before for $\beta_{\frac{k}{2}-2,\frac{k}{2}-1}$, and after for $\beta_{\frac{k}{2}-3,\frac{k}{2}-2}$, we have eliminated the last braid and we obtain
\[\beta_{12}(\beta_{34}\beta_{23}\beta_{12})(\beta_{k-5,k-4}\beta_{k-6,k-5}\dots\beta_{12})(\beta_{k-6,k-5}\beta_{k-7,k-6}\dots\beta_{23})(\beta_{k-7,k-6}\dots\beta_{34})\dots\]\[\dots(\beta_{\frac{k}{2},\frac{k}{2}+1}\beta_{\frac{k}{2}-1,\frac{k}{2}}\beta_{\frac{k}{2}-2,\frac{k}{2}-1}\beta_{\frac{k}{2}-3,\frac{k}{2}-2}\beta_{\frac{k}{2}-4,\frac{k}{2}-3}).
\]Iterating the same procedure, one obtains the braid
\[\beta_{12}(\beta_{34}\beta_{23}\beta_{12})(\beta_{56}\beta_{45}\beta_{34}\beta_{23}\beta_{12})\dots(\beta_{k-1,k}\beta_{k-2,k-1}\dots\beta_{12}).
\]The case $k$ odd is analogous, and we left the details to the reader.
\endproof

\begin{es}Consider for example $k=12$. We have that 
\[\beta\beta'=(\beta_{78}\beta_{67}\beta_{56}\beta_{45}\beta_{34}\beta_{23}\beta_{12})(\beta_{67}\beta_{56}\beta_{45}\beta_{34}\beta_{23})(\beta_{56}\beta_{45}\beta_{34})\beta_{45}\cdot\]\[\cdot(\beta_{9,10},\dots\beta_{12})(\beta_{11,12}\dots\beta_{12}).
\]We have to rearrange the first 4 braids. Let us apply the procedure described above:
\begin{align*}
(\beta_{78}\beta_{67}\beta_{56}\beta_{45}\beta_{34}\beta_{23}\beta_{12})(\beta_{67}\beta_{56}\beta_{45}\beta_{34}\beta_{23})(\beta_{56}\textcolor{red}{\beta_{45}}\beta_{34})\textcolor{red}{\beta_{45}}&=\\
(\beta_{78}\beta_{67}\beta_{56}\beta_{45}\beta_{34}\beta_{23}\beta_{12})(\beta_{67}\beta_{56}\beta_{45}\beta_{34}\beta_{23})(\beta_{56}\textcolor{green}{\beta_{34}}\beta_{45})\textcolor{green}{\beta_{34}}&=\\
(\beta_{78}\beta_{67}\beta_{56}\beta_{45}\beta_{34}\beta_{23}\beta_{12})(\beta_{67}\beta_{56}\beta_{45}\textcolor{green}{\beta_{34}}\beta_{23})\textcolor{green}{\beta_{34}}(\beta_{56}\beta_{45}\beta_{34})&=\\
(\beta_{78}\beta_{67}\beta_{56}\beta_{45}\beta_{34}\beta_{23}\beta_{12})(\beta_{67}\beta_{56}\beta_{45}\textcolor{blue}{\beta_{23}}\beta_{34})\textcolor{blue}{\beta_{23}}(\beta_{56}\beta_{45}\beta_{34})&=\\
(\beta_{78}\beta_{67}\beta_{56}\beta_{45}\beta_{34}\textcolor{blue}{\beta_{23}}\beta_{12})\textcolor{blue}{\beta_{23}}(\beta_{67}\beta_{56}\beta_{45}\beta_{34}\beta_{23})(\beta_{56}\beta_{45}\beta_{34})&=\\
(\beta_{78}\beta_{67}\beta_{56}\beta_{45}\beta_{34}\textcolor{red}{\beta_{12}}\beta_{23})\textcolor{red}{\beta_{12}}(\beta_{67}\beta_{56}\beta_{45}\beta_{34}\beta_{23})(\beta_{56}\beta_{45}\beta_{34})&=\\
\textcolor{red}{\beta_{12}}(\beta_{78}\beta_{67}\beta_{56}\beta_{45}\beta_{34}\beta_{23}\beta_{12})(\beta_{67}\beta_{56}\beta_{45}\beta_{34}\beta_{23})(\beta_{56}\beta_{45}\beta_{34}).
\end{align*}
Now we continue by eliminating the last braid, starting from its first term:
\begin{align*}
\beta_{12}(\beta_{78}\beta_{67}\beta_{56}\beta_{45}\beta_{34}\beta_{23}\beta_{12})(\beta_{67}\beta_{56}\beta_{45}\beta_{34}\beta_{23})(\textcolor{red}{\beta_{56}}\beta_{45}\beta_{34})&=\\
\beta_{12}(\beta_{78}\beta_{67}\beta_{56}\beta_{45}\beta_{34}\beta_{23}\beta_{12})(\beta_{67}\textcolor{red}{\beta_{56}}\beta_{45}\textcolor{red}{\beta_{56}}\beta_{34}\beta_{23})(\beta_{45}\beta_{34})&=\\
\beta_{12}(\beta_{78}\beta_{67}\beta_{56}\beta_{45}\beta_{34}\beta_{23}\beta_{12})(\beta_{67}\textcolor{blue}{\beta_{45}}\beta_{56}\textcolor{blue}{\beta_{45}}\beta_{34}\beta_{23})(\beta_{45}\beta_{34})&=\\
\beta_{12}(\beta_{78}\beta_{67}\beta_{56}\textcolor{blue}{\beta_{45}}\beta_{34}\textcolor{blue}{\beta_{45}}\beta_{23}\beta_{12})(\beta_{67}\beta_{56}\beta_{45}\beta_{34}\beta_{23})(\beta_{45}\beta_{34})&=\\
\beta_{12}(\beta_{78}\beta_{67}\beta_{56}\textcolor{green}{\beta_{34}}\beta_{45}\textcolor{green}{\beta_{34}}\beta_{23}\beta_{12})(\beta_{67}\beta_{56}\beta_{45}\beta_{34}\beta_{23})(\beta_{45}\beta_{34})&=\\
\beta_{12}\textcolor{green}{\beta_{34}}(\beta_{78}\beta_{67}\beta_{56}\beta_{45}\beta_{34}\beta_{23}\beta_{12})(\beta_{67}\beta_{56}\beta_{45}\beta_{34}\beta_{23})(\beta_{45}\beta_{34})&=\\
\beta_{12}\beta_{34}(\beta_{78}\beta_{67}\beta_{56}\beta_{45}\beta_{34}\beta_{23}\beta_{12})(\beta_{67}\beta_{56}\beta_{45}\beta_{34}\beta_{23})(\textcolor{red}{\beta_{45}}\beta_{34})&=\\
\beta_{12}\beta_{34}(\beta_{78}\beta_{67}\beta_{56}\beta_{45}\beta_{34}\beta_{23}\beta_{12})(\beta_{67}\beta_{56}\textcolor{red}{\beta_{45}}\beta_{34}\textcolor{red}{\beta_{45}}\beta_{23})\beta_{34}&=\\
\beta_{12}\beta_{34}(\beta_{78}\beta_{67}\beta_{56}\beta_{45}\beta_{34}\beta_{23}\beta_{12})(\beta_{67}\beta_{56}\textcolor{blue}{\beta_{34}}\beta_{45}\textcolor{blue}{\beta_{34}}\beta_{23})\beta_{34}&=\\
\beta_{12}\beta_{34}(\beta_{78}\beta_{67}\beta_{56}\beta_{45}\textcolor{blue}{\beta_{34}}\beta_{23}\textcolor{blue}{\beta_{34}}\beta_{12})(\beta_{67}\beta_{56}\beta_{45}\beta_{34}\beta_{23})\beta_{34}&=\\
\beta_{12}\beta_{34}(\beta_{78}\beta_{67}\beta_{56}\beta_{45}\textcolor{green}{\beta_{23}}\beta_{34}\textcolor{green}{\beta_{23}}\beta_{12})(\beta_{67}\beta_{56}\beta_{45}\beta_{34}\beta_{23})\beta_{34}&=\\
\beta_{12}\beta_{34}\textcolor{green}{\beta_{23}}(\beta_{78}\beta_{67}\beta_{56}\beta_{45}\beta_{34}\beta_{23}\beta_{12})(\beta_{67}\beta_{56}\beta_{45}\beta_{34}\beta_{23})\beta_{34}&=\\
\beta_{12}\beta_{34}\beta_{23}(\beta_{78}\beta_{67}\beta_{56}\beta_{45}\beta_{34}\beta_{23}\beta_{12})(\beta_{67}\beta_{56}\beta_{45}\textcolor{red}{\beta_{34}}\beta_{23})\textcolor{red}{\beta_{34}}&=\\
\beta_{12}\beta_{34}\beta_{23}(\beta_{78}\beta_{67}\beta_{56}\beta_{45}\beta_{34}\beta_{23}\beta_{12})(\beta_{67}\beta_{56}\beta_{45}\textcolor{blue}{\beta_{23}}\beta_{34})\textcolor{blue}{\beta_{23}}&=\\
\beta_{12}\beta_{34}\beta_{23}(\beta_{78}\beta_{67}\beta_{56}\beta_{45}\beta_{34}\textcolor{blue}{\beta_{23}}\beta_{12})\textcolor{blue}{\beta_{23}}(\beta_{67}\beta_{56}\beta_{45}\beta_{34}\beta_{23})&=\\
\beta_{12}\beta_{34}\beta_{23}(\beta_{78}\beta_{67}\beta_{56}\beta_{45}\beta_{34}\textcolor{green}{\beta_{12}}\beta_{23})\textcolor{green}{\beta_{12}}(\beta_{67}\beta_{56}\beta_{45}\beta_{34}\beta_{23})&=\\
\beta_{12}(\beta_{34}\beta_{23}\textcolor{green}{\beta_{12}})(\beta_{78}\beta_{67}\beta_{56}\beta_{45}\beta_{34}\beta_{23}\beta_{12})(\beta_{67}\beta_{56}\beta_{45}\beta_{34}\beta_{23}).
\end{align*}
At the final step, we have to eliminate the last braid, always starting from its first term:
\begin{align*}
\beta_{12}(\beta_{34}\beta_{23}\beta_{12})(\beta_{78}\beta_{67}\beta_{56}\beta_{45}\beta_{34}\beta_{23}\beta_{12})(\textcolor{red}{\beta_{67}}\beta_{56}\beta_{45}\beta_{34}\beta_{23})&=\\
\beta_{12}(\beta_{34}\beta_{23}\beta_{12})(\beta_{78}\textcolor{red}{\beta_{67}}\beta_{56}\textcolor{red}{\beta_{67}}\beta_{45}\beta_{34}\beta_{23}\beta_{12})(\beta_{56}\beta_{45}\beta_{34}\beta_{23})&=\\
\beta_{12}(\beta_{34}\beta_{23}\beta_{12})(\beta_{78}\textcolor{blue}{\beta_{56}}\beta_{67}\textcolor{blue}{\beta_{56}}\beta_{45}\beta_{34}\beta_{23}\beta_{12})(\beta_{56}\beta_{45}\beta_{34}\beta_{23})&=\\
\beta_{12}(\beta_{34}\beta_{23}\beta_{12})\textcolor{blue}{\beta_{56}}(\beta_{78}\beta_{67}\beta_{56}\beta_{45}\beta_{34}\beta_{23}\beta_{12})(\beta_{56}\beta_{45}\beta_{34}\beta_{23})&=\\
\beta_{12}(\beta_{34}\beta_{23}\beta_{12})\beta_{56}(\beta_{78}\beta_{67}\beta_{56}\beta_{45}\beta_{34}\beta_{23}\beta_{12})(\textcolor{red}{\beta_{56}}\beta_{45}\beta_{34}\beta_{23})&=\\
\beta_{12}(\beta_{34}\beta_{23}\beta_{12})\beta_{56}(\beta_{78}\beta_{67}\textcolor{red}{\beta_{56}}\beta_{45}\textcolor{red}{\beta_{56}}\beta_{34}\beta_{23}\beta_{12})(\beta_{45}\beta_{34}\beta_{23})&=\\
\beta_{12}(\beta_{34}\beta_{23}\beta_{12})\beta_{56}(\beta_{78}\beta_{67}\textcolor{green}{\beta_{45}}\beta_{56}\textcolor{green}{\beta_{45}}\beta_{34}\beta_{23}\beta_{12})(\beta_{45}\beta_{34}\beta_{23})&=\\
\beta_{12}(\beta_{34}\beta_{23}\beta_{12})\beta_{56}\textcolor{green}{\beta_{45}}(\beta_{78}\beta_{67}\beta_{56}\beta_{45}\beta_{34}\beta_{23}\beta_{12})(\beta_{45}\beta_{34}\beta_{23})&=\\
\beta_{12}(\beta_{34}\beta_{23}\beta_{12})\beta_{56}\beta_{45}(\beta_{78}\beta_{67}\beta_{56}\beta_{45}\beta_{34}\beta_{23}\beta_{12})(\textcolor{blue}{\beta_{45}}\beta_{34}\beta_{23})&=\\
\beta_{12}(\beta_{34}\beta_{23}\beta_{12})\beta_{56}\beta_{45}(\beta_{78}\beta_{67}\beta_{56}\textcolor{blue}{\beta_{45}}\beta_{34}\textcolor{blue}{\beta_{45}}\beta_{23}\beta_{12})(\beta_{34}\beta_{23})&=\\
\beta_{12}(\beta_{34}\beta_{23}\beta_{12})\beta_{56}\beta_{45}\textcolor{red}{\beta_{34}}(\beta_{78}\beta_{67}\beta_{56}\beta_{45}\beta_{34}\beta_{23}\beta_{12})(\beta_{34}\beta_{23})&=\\
\beta_{12}(\beta_{34}\beta_{23}\beta_{12})\beta_{56}\beta_{45}\beta_{34}(\beta_{78}\beta_{67}\beta_{56}\beta_{45}\beta_{34}\beta_{23}\beta_{12})(\textcolor{green}{\beta_{34}}\beta_{23})&=\\
\beta_{12}(\beta_{34}\beta_{23}\beta_{12})\beta_{56}\beta_{45}\beta_{34}(\beta_{78}\beta_{67}\beta_{56}\beta_{45}\textcolor{green}{\beta_{34}}\beta_{23}\textcolor{green}{\beta_{34}}\beta_{12})\beta_{23}&=\\
\beta_{12}(\beta_{34}\beta_{23}\beta_{12})\beta_{56}\beta_{45}\beta_{34}(\beta_{78}\beta_{67}\beta_{56}\beta_{45}\textcolor{blue}{\beta_{23}}\beta_{34}\textcolor{blue}{\beta_{23}}\beta_{12})\beta_{23}&=\\
\beta_{12}(\beta_{34}\beta_{23}\beta_{12})\beta_{56}\beta_{45}\beta_{34}\textcolor{blue}{\beta_{23}}(\beta_{78}\beta_{67}\beta_{56}\beta_{45}\beta_{34}\beta_{23}\beta_{12})\beta_{23}&=\\
\beta_{12}(\beta_{34}\beta_{23}\beta_{12})\beta_{56}\beta_{45}\beta_{34}\beta_{23}(\beta_{78}\beta_{67}\beta_{56}\beta_{45}\beta_{34}\textcolor{red}{\beta_{23}}\beta_{12})\textcolor{red}{\beta_{23}}&=\\
\beta_{12}(\beta_{34}\beta_{23}\beta_{12})\beta_{56}\beta_{45}\beta_{34}\beta_{23}(\beta_{78}\beta_{67}\beta_{56}\beta_{45}\beta_{34}\textcolor{green}{\beta_{12}}\beta_{23})\textcolor{green}{\beta_{12}}&=\\
\beta_{12}(\beta_{34}\beta_{23}\beta_{12})(\beta_{56}\beta_{45}\beta_{34}\beta_{23}\textcolor{green}{\beta_{12}})(\beta_{78}\beta_{67}\beta_{56}\beta_{45}\beta_{34}\beta_{23}\beta_{12}).
\end{align*}
\end{es}

In what follows we will denote by $\mathcal T$ the tangent sheaf of $\mathbb P$, by $\Omega$ the cotangent sheaf, and we will use the shorthands 
\[\bigwedge\nolimits^p\mathcal T(k):=\left(\bigwedge\nolimits^p\mathcal T\right)\otimes \mathcal O(k),\quad \bigwedge\nolimits^p\Omega(k):=\left(\bigwedge\nolimits^p\Omega\right)\otimes \mathcal O(k).
\]
The following formulae, due to R. Bott (\cite{bott}, \cite{okonek}, \cite{distler}), will be useful:
\begin{equation}\label{bott1}\dim_{\mathbb C} H^q\left(\mathbb P^{n}_{\mathbb C}, \bigwedge\nolimits^p\mathcal T(k)\right)=
\begin{sistema}
\binom{k+n+p+1}{p}\binom{k+n}{n-p},\quad\quad q=0,\quad k>-p-1,\\
\\
1,\quad\quad\quad\quad\quad\quad\quad\quad q=n-p,\quad k=-n-1,\\
\\
\binom{-k-p-1}{-k-n-1}\binom{-k-n-2}{p},\quad q=n,\quad k<-n-p-1,\\
\\
0,\quad\quad\quad\quad\quad\quad\quad\quad\quad\quad\quad\text{otherwise,}
\end{sistema}
\end{equation}

\begin{equation}\label{bott2}
\dim_{\mathbb C} H^q\left(\mathbb P^n_{\mathbb C}, \bigwedge\nolimits^p\Omega(k)\right)=
\begin{sistema}
\binom{k+n-p}{k}\binom{k-1}{p}, \quad\quad q=0,\quad 0\leq p\leq n,\quad k>p,\\
\\
1,\quad\quad\quad\quad\quad\quad\quad\quad k=0,\quad 0\leq q=p\leq n,\\
\\
\binom{-k+p}{-k}\binom{-k-1}{n-p},\quad q=n,\quad 0\leq p\leq n,\quad k<p-n,\\
\\
0,\quad\quad\quad\quad\quad\quad\quad\quad\quad\quad\quad\text{otherwise.}
\end{sistema}
\end{equation}

Consider Beilinson's exceptional collection $\frak B:=(\mathcal O,\mathcal O(1),\mathcal O(2),\dots,\mathcal O(k-1))$ in $\mathcal D^b\left(\mathbb P\right)$, with $\mathbb P=\mathbb P(V)$ ($\dim_{\mathbb C}V=k$), and the well known Euler exact sequence, together with its exterior powers
\begin{equation}\label{euler}\xymatrix@R=1pt{
0\ar[r]&\mathcal O\ar[r]&V\otimes \mathcal O(1)\ar[r]&\mathcal T\ar[r]&0,\\
0\ar[r]&\mathcal T\ar[r]&\bigwedge\nolimits^2V\otimes \mathcal O(2)\ar[r]&\bigwedge\nolimits^2\mathcal T\ar[r]&0,\\
&\vdots&\vdots&\vdots\\
0\ar[r]&\bigwedge\nolimits^{h-1}\mathcal T\ar[r]&\bigwedge\nolimits^hV\otimes \mathcal O(h)\ar[r]&\bigwedge\nolimits^h\mathcal T\ar[r]&0,\\
&\vdots&\vdots&\vdots\\
0\ar[r]&\bigwedge\nolimits^{k-2}\mathcal T\ar[r]&\bigwedge\nolimits^{k-1}V\otimes \mathcal O(k-1)\ar[r]&\mathcal O(k)\ar[r]&0.\\
}
\end{equation}
By Bott formulae \eqref{bott1}-\eqref{bott2}, we deduce that both $\Hom^\bullet(\mathcal O(h),\bigwedge\nolimits^h\mathcal T)$ and $\Hom^\bullet\left(\bigwedge\nolimits^{h-1}\mathcal T,\mathcal O(h)\right)$ are concentrated in degree 0 and they have the same dimension $\binom{k}{h}$. Hence, the short exact sequences \eqref{euler}, together with the identifications
\[\bigwedge\nolimits^hV=\Hom^\bullet\left(\mathcal O(h),\bigwedge\nolimits^h\mathcal T\right)=(\Hom^\bullet)^\vee\left(\bigwedge\nolimits^{h-1}\mathcal T,\mathcal O(h)\right),
\]allow us to explicitly compute successive right mutations of the sheaf $\mathcal O$: namely, denoting by $\sigma_{ij}$ the inverse braid $\beta_{ij}^{-1}$, for $0<h\leq k-1$ we have
\[\mathbb R_{[\mathcal O(1)\dots\mathcal O(h)]}\mathcal O=\left(\bigwedge\nolimits^h\mathcal T\right)[-h].
\]
Being the sheaf $\mathcal O(j)$ locally free, the functor $\mathcal O(j)\otimes (-)$ preserves the short exact sequences \eqref{euler}; moreover, observing that 
\[\Hom(\mathcal O(l),\mathcal O(m))\cong \Hom(\mathcal O(l+n),\mathcal O(m+n))
\]for all $l,m,n\in\mathbb Z$, we deduce that for $j<h\leq k-1$
\[\mathbb R_{[\mathcal O(j+1),\dots,\mathcal O(h)]}\mathcal O(j)=\left(\bigwedge\nolimits^{h-j}\mathcal T(j)\right)[j-h].
\]

\bsh
\begin{cor}\label{corexceptional1}
The central connection and the Stokes matrices of the quantum cohomology of $\mathbb P^{k-1}_{\mathbb C}$, computed at $0\in QH^{\bullet}(\mathbb P)$ and with respect to a line $\ell$ with slope $0<\epsilon<\frac{\pi}{k}$, corresponds (modulo action of $(\mathbb Z/2\mathbb Z)^k$) to the exceptional collections
\[\left(\mathcal O\left(\frac{k}{2}\right),\bigwedge\nolimits^1\mathcal T\left(\frac{k}{2}-1\right),\mathcal O\left(\frac{k}{2}+1\right),\bigwedge\nolimits^3\mathcal T\left(\frac{k}{2}-2\right),\dots,\mathcal O(k-1),\bigwedge\nolimits^{k-1}\mathcal T\right)
\]
for $k$ even, and
\[\left(\mathcal O\left(\frac{k-1}{2}\right),\mathcal O\left(\frac{k+1}{2}\right),\bigwedge\nolimits^2\mathcal T\left(\frac{k-3}{2}\right),\mathcal O\left(\frac{k+3}{2}\right),\bigwedge\nolimits^4\mathcal T\left(\frac{k-5}{2}\right),\dots,\mathcal O\left(k-1\right),\bigwedge\nolimits^{k-1}\mathcal T\right)
\]for $k$ odd.
\end{cor}
\esh

\proof[Proof] From Theorem \ref{risultato1} and from Lemma \ref{bellatreccia}, we have that the monodromy data computed at $0$ with respect to the line $\ell$ correspond to the exceptional collection
\[\frak B^{\sigma},\quad \sigma:=(\beta^{-1}_{12}\beta^{-1}_{23}\dots\beta^{-1}_{k-1,k})\dots(\beta^{-1}_{12}\beta^{-1}_{23}\beta^{-1}_{34}\beta^{-1}_{45}\beta^{-1}_{56})(\beta^{-1}_{12}\beta^{-1}_{23}\beta^{-1}_{34})\beta^{-1}_{12}
\]for $k$ even, and to
\[\frak B^{\sigma},\quad \sigma:=(\beta^{-1}_{12}\beta^{-1}_{23}\dots\beta^{-1}_{k-1,k})\dots(\beta^{-1}_{12}\beta^{-1}_{23}\beta^{-1}_{34}\beta^{-1}_{45})(\beta^{-1}_{12}\beta^{-1}_{23})
\]for $k$ odd. Using the previous observations, one obtains the collections above.\endproof

\subsection
{Monodromy data along the small quantum cohomology, and some results on the big quantum cohomology}\label{secstoksqcohcpn}$\quad $\newline
From Corollary \ref{corexceptional1}, we are able  to determine the monodromy data at \emph{any} point of the small quantum cohomology with respect to \emph{any} line $\ell$, together with the corresponding full exceptional collections.

The small quantum cohomology is identified with the set of points $(0,t^2,0,...,0)$, which can be represented in the $t^2$-plane. Fixed $\ell$, by formula \eqref{stokesrayscpn} we see that when $t^2$ varies, then some Stokes rays  cross $\ell$  whenever 
\begin{equation}
\label{lineforbidden}
\Im t^2=-k\phi +m\pi, \quad m\in\mathbb{Z}.
\end{equation}
Thus, the $t^2$-plane is divided into horizontal strips, whose boundary lines are  \eqref{lineforbidden}. These strips are the intersection of the small quantum cohomology with $\ell$-chambers. 

By the Isomonodromy Theorem \ref{iso2}, the monodromy data are constant in each horizontal strip. The data in different strips are related by a braid, as follows.  Passing from one strip  $\mathcal C_1(\ell)$ to an adjacent one   $\mathcal C_2(\ell)$, one Stokes ray crosses $\ell$, so the monodromy data change by a braid,  determined by $u_i,u_j$ associated with the ray crossing $\ell$, as in Section \ref{monlocalmod}.

  \begin{table}
\begin{tabular}{|c||c|c|c|}
\hline
& $S_\text{lex}$ & Exceptional Collection&Braid\\
\hline
\hline
$0<3\phi+\Im(t^1)<\pi$&$\left(
\begin{array}{ccc}
 1 & 3 & -3 \\
 0 & 1 & -3 \\
 0 & 0 & 1 \\
\end{array}
\right) $& $(\mathcal O(1),\mathcal O(2),\bigwedge^2\mathcal T)$&$id$\\
\hline
$\pi<3\phi+\Im(t^1)<2\pi$&$\left(
\begin{array}{ccc}
 1 & -3 & -6 \\
 0 & 1 & 3 \\
 0 & 0 & 1 \\
\end{array}
\right) $&$\left(\mathcal O(1),\bigwedge^1\mathcal T,\mathcal O(2)\right)$&$\omega_{1,3}$\\
\hline
$2\pi<3\phi+\Im(t^1)<3\pi$&$\left(
\begin{array}{ccc}
 1 & 3 & 3 \\
 0 & 1 & 3 \\
 0 & 0 & 1 \\
\end{array}
\right)$&{\color{red}$\left(\mathcal O,\mathcal O(1),\mathcal O(2)\right)$}&$\omega_{1,3}\omega_{2,3}$\\
\hline
$3\pi<3\phi+\Im(t^1)<4\pi$&$\left(
\begin{array}{ccc}
 1 & 3 & -6 \\
 0 & 1 & -3 \\
 0 & 0 & 1 \\
\end{array}
\right)$&$\left(\mathcal O,\Omega(2),\mathcal O(1)\right)$&$\omega_{1,3}\omega_{2,3}\omega_{1,3}$\\
\hline
$4\pi<3\phi+\Im(t^1)<5\pi$&$\left(
\begin{array}{ccc}
 1 & -3 & -3 \\
 0 & 1 & 3 \\
 0 & 0 & 1 \\
\end{array}
\right)$&$\left(\bigwedge^2\Omega(2),\mathcal O,\mathcal O(1)\right)$&$(\omega_{1,3}\omega_{2,3})^2$\\
\hline
$5\pi<3\phi+\Im(t^1)<6\pi$&$\left(
\begin{array}{ccc}
 1 & -3 & 6 \\
 0 & 1 & -3 \\
 0 & 0 & 1 \\
\end{array}
\right)$&$\left(\bigwedge^2\Omega(2),\Omega(1),\mathcal O\right)$&$(\omega_{1,3}\omega_{2,3})^2\omega_{1,3}$\\
\hline
$6\pi<3\phi+\Im(t^1)<7\pi$&$\left(
\begin{array}{ccc}
 1 & 3 & -3 \\
 0 & 1 & -3 \\
 0 & 0 & 1 \\
\end{array}
\right)$&$\left(\bigwedge^2\Omega(1),\bigwedge^2\Omega(2),\mathcal O\right)$&$(\omega_{1,3}\omega_{2,3})^3$\\
\hline
\end{tabular}
\caption{In this table we represent all possible Stokes matrices along the small quantum cohomology of $\mathbb P^2_\mathbb C$, in the $\ell$-lexicographical order for a line $\ell$ of slope $\phi$. We also write the corresponding (modulo shifts) exceptional collections associated with the monodromy data. Notice that the Beilinson exceptional collection $\frak B$ appears along the small quantum locus: it is obtained from the one of Corollary \ref{corexceptional1} by applying the braids $\omega_{1,3}\omega_{2,3}$.}
\label{collcp2}
\end{table}
 
Up to now we have fixed a line $\ell$ and considered the data in these ``static'' strips (which, remember, corresponds to  ``static'' $\ell$-chambers). If instead  we let $\ell$ rotate, say by increasing its slope $\phi$, then the  $\ell$-strips   \emph{glide} over the $t^2$-plane, according to equation \eqref{lineforbidden}. Consider a point $t^2\in\mathcal C_1(\ell)$, and let  the line $\ell$ vary by increasing its slope by $\Delta\phi$. Let $\ell^\prime$ be the admissible line after the rotation. Then, the strip  $\mathcal C_1(\ell)$ moves towards $\Im(t^2)\to-\infty$, so that, at the end of the  rotation, $t^2$ belongs to another strip, say $\mathcal C_2(\ell')=\mathcal C_2(\ell)-c$ for some positive constant $c=k\Delta\phi$. This process leads to the same braid actions obtained from the point of view described in the previous paragraph.

 In conclusion, if we know the data at a point of the small quantum cohomology with respect to some line $\ell$, then we can reconstruct the data at any  other point with respect to any other line.\newline


Starting from $0\in QH^{\bullet}(\mathbb P)$ with a line $\ell$ of slope $0<\phi<\frac{\pi}{k}$, we let increase $\phi$, so that the line $\ell$ rotates counter-clockwise. From the configuration of the Stokes rays,  it is easily seen that the \emph{first} crossing between $\ell$ and Stokes rays is described as follows:
\begin{itemize}
\item if $k\geq 4$ is \emph{even}, the line $\ell$ firstly crosses $\frac{k}{2}$ Stokes rays (which coincide) and the corresponding braid is
\[\omega_{1,k}:=\prod_{\substack{i=2\\ i\ \text{even}}}^k\beta_{i-1,i};
\]
\item if $k\geq 3$ is \emph{odd}, then the line $\ell$ firstly crosses $\frac{k-1}{2}$ Stokes rays (which coincide) and the corresponding braid is
\[\omega_{1,k}:=\prod_{\substack{i=3\\ i\ \text{odd}}}^k\beta_{i-1,i}.
\]
\end{itemize}

The \emph{second} crossing is:
\begin{itemize}
\item if $k\geq 4$ is \emph{even}, the line $\ell$ secondly cross $\frac{k}{2}-1$ Stokes rays (which coincide) and the corresponding braid is
\[\omega_{2,k}:=\prod_{\substack{i=3\\ i\ \text{odd}}}^{k-1}\beta_{i-1,i};
\]
\item if $k\geq 3$ is \emph{odd}, then the line $\ell$ firstly cross $\frac{k-1}{2}$ Stokes rays (which coincide) and the corresponding braid is
\[\omega_{2,k}:=\prod_{\substack{i=2\\ i\ \text{even}}}^{k-1}\beta_{i-1,i}.
\]
\end{itemize}
Furthermore, using symmetries of regular polygons, it is easy to see that the braids corresponding to subsequent crossings are \emph{alternatively} $\omega_{1,k}$ and $\omega_{2,k}$: in this ways, if we let $\ell$ rotate counterclockwise, and we have $N$ crossings in total, the resulting acting braid is the composition
\[\omega_{1,k}\omega_{2,k}\omega_{1,k}\omega_{2,k}\dots\]
with $N$ braids $\omega$ in total. Notice that after a complete rotation of $\ell$, the resulting braid is
\[(\omega_{1,k}\omega_{2,k})^k,
\]which, accordingly to Corollary \ref{centerbraidlemma} and using the ``Yang--Baxter'' braid relations \eqref{10novembre2018-3}, is easily seen to be  the central element $(\beta_{12},\dots,\beta_{k-1,k})^k$, .

\begin{figure}[ht!]
\centering
\def\svgscale{0.6}
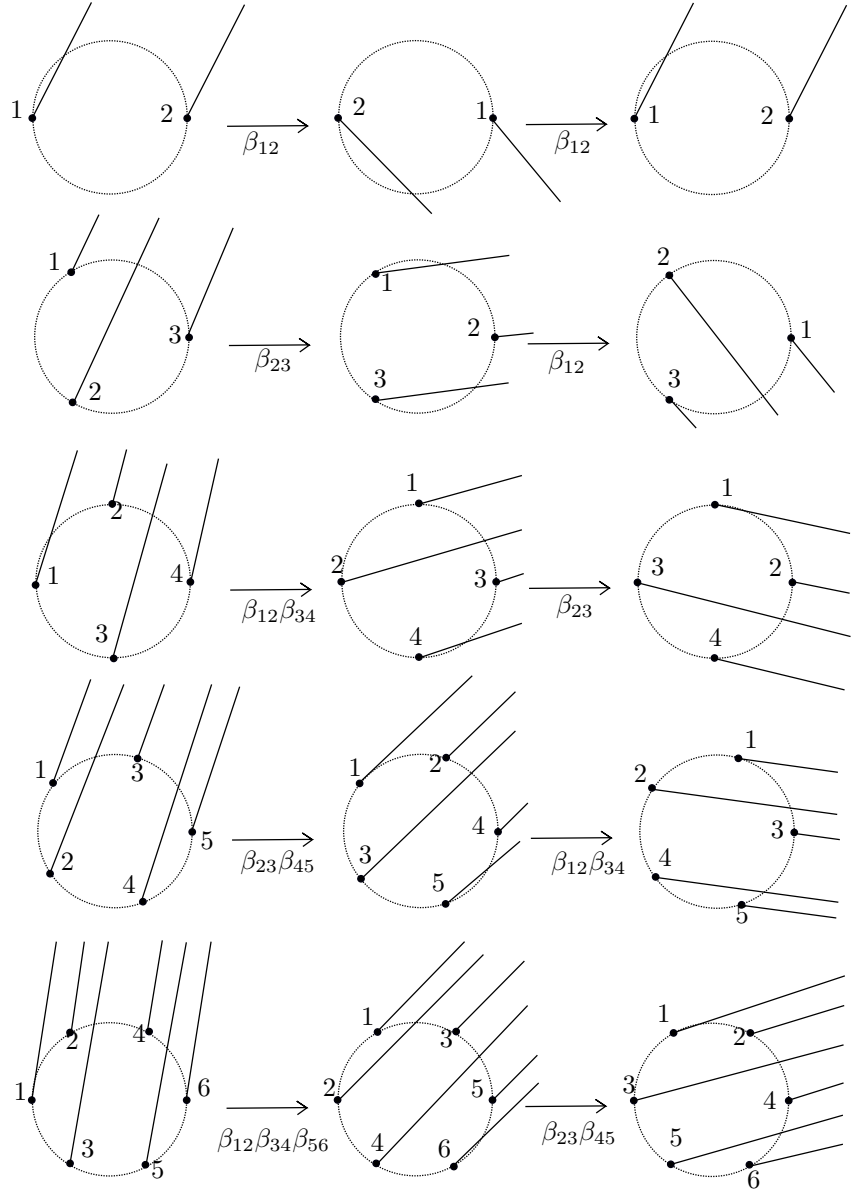
\caption{Here we represent the action of the braids $\omega_{1,k}, \omega_{2,k}$ for $2\leq k\leq 6$. In the left column the reader can find the canonical coordinates in the $\ell$-lexicographical order for $\ell$ of slope $\phi\in]0;\frac{\pi}{k}[$. In the central column we represent the action of the braid $\omega_{1,k}$, whereas in the right column the consecutive action of the braid $\omega_{2,k}$.}
\label{disegnosympolygons}
\end{figure}

\bsh
\begin{teorema}\label{teotreccebeilinson}
The braids of Lemma \ref{bellatreccia}, i.e.\[\beta\beta'=\beta_{12}(\beta_{34}\beta_{23}\beta_{12})(\beta_{56}\beta_{45}\beta_{34}\beta_{23}\beta_{12})\dots(\beta_{k-1,k}\beta_{k-2,k-1}\dots\beta_{12})
\]for $k$ even,
\[\beta\beta'=(\beta_{23}\beta_{12})(\beta_{45}\beta_{34}\beta_{23}\beta_{12})\dots(\beta_{k-1,k}\beta_{k-2,k-1}\dots\beta_{12})
\]for $k$ odd, which take the monodromy data computed at $0\in QH^\bullet(\mathbb P^{k-1}_{\mathbb C})$ (with respect to a line of slope $0<\phi<\frac{\pi}{k}$) to the data corresponding to the Beilinson's exceptional collection, are of the form
\[\omega_{1,k}\omega_{2,k}\omega_{1,k}\omega_{2,k}\dots\]
if and only if $k=2$ or $k=3$.
Thus, they do not correspond to analytic continuation along paths in the small quantum cohomology for $k\geq 4$.
\end{teorema}
\esh

\proof
For $k=2,3$ we have already shown that the braids $\beta\beta'$ are
\[\omega_{1,2}=\beta_{12}\quad\text{and}\quad \omega_{1,3}\omega_{2,3}=\beta_{23}\beta_{12}
\]respectively. So, let us suppose that $k\geq 4$ and that $\beta\beta'$ can be expressed as a product
\begin{equation}\label{ansatzproduct}\omega_{1,k}\omega_{2,k}\omega_{1,k}\omega_{2,k}\dots.\end{equation}
Let us start from the following observation: if a generic braid can be represented as a product of positive powers of elementary braids $\beta_{i,i+1}$, then any other of its factorizations in positive powers of elementary braids must consist of the same numbers of factors (this follows immediately from the relations defining the braid group $\mathcal B_n$). Thus, the product \eqref{ansatzproduct} should be a product of
\[\left(\frac{k}{2}\right)^2\\ \text{factors for $k$ even,}\quad
 \frac{k^2-1}{4}\\ \text{ factors for $k$ odd.}\]
We firstly consider the case $k$ \emph{even}: we are supposing existence of a number $n\in\mathbb N^*$ such that the product \eqref{ansatzproduct} contains $n$ times the braid $\omega_{1,k}$ and $n$ or $n-1$ times the  braid $\omega_{2,k}$. So, we must have
\[n\frac{k}{2}+m\left(\frac{k}{2}-1\right)=\left(\frac{k}{2}\right)^2
\]for some $n\in\mathbb N^*$ and $m\in\left\{n-1,n\right\}$, so that
\begin{equation}\label{kfromnandm}k=(n+m)\pm\frac{1}{2} (4(n+m)^2-16m)^{\frac{1}{2}}.
\end{equation}As a necessary condition we have that
\[4(n+m)^2-16m,\quad\text{with }m\in\left\{n-1,n\right\}
\]must be the square of some integer. Since 
\begin{itemize}
\item for $m=n$ the number $16(n^2-n)$ is a perfect square only for $n=1$,
\item for $m=n-1$ the number $16(n-1)^2+4$ is a perfect square only for $n=1$,
\end{itemize}
according to \eqref{kfromnandm} the only possible value of $k$ is $k=2$. Analogously, for the case $k\geq 3$ and \emph{odd}, if we suppose that it exists a number $n\in\mathbb N^*$ such that the product \eqref{ansatzproduct} contains $n$ times the braid $\omega_{1,k}$ and $n$ or $n-1$ times the  braid $\omega_{2,k}$, we necessarily must have
\[n\frac{k-1}{2}+m\frac{k-1}{2}=\frac{k^2-1}{4},\quad\text{with }m\in\left\{n-1,n\right\}\]\[\Longrightarrow n+m=\frac{k+1}{2},\quad\text{with }m\in\left\{n-1,n\right\}.
\]Thus, for any odd number $k\geq 3$, we have found a composition of $n$ times $\omega_{1,k}$ and $n$ or $n-1$ times $\omega_{2,k}$ whose length equals the length of $\beta\beta'$. In particular, we have that
\begin{itemize}
\item if $k=4n-1$ then $\omega_{1,k}$ and  $\omega_{2,k}$ appear the same number $n=m$ of times;
\item if $k=4n-3$ then $\omega_{1,k}$ appears $n$ times and $\omega_{2,k}$ appears $m=n-1$ times.
\end{itemize}Notice in particular that for $k=3$ we are in the first case, according to  what said at the beginning of the proof. We want now to show that $k=3$ is the only case in which the braid we have found is actually $\beta\beta'$. 

\begin{figure}[h]
\centering
\def\svgscale{0.4}
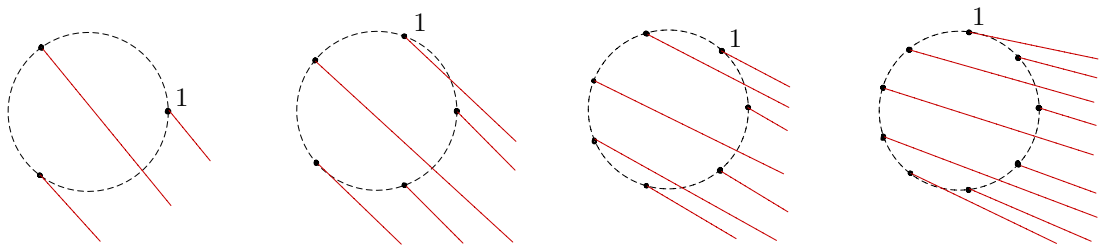
\caption{Configuration of the canonical coordinates, for $k$ odd ($k=3,5,7,9$) after the action of the candidate braid $\omega_{1,k}\omega_{2,k}\dots$. Notice that the final arrangement of the canonical coordinates is $(\dots,\underset{\text{$n$-th}}{1},\dots)$.}
\label{trecceomega}
\end{figure}

For this, notice that the braid $\beta\beta'$ takes the canonical coordinates in an ordered cyclic disposition $(u_1,u_2,\dots,u_k)$ starting from $1$ and going counter-clockwise along the regular $k$-gon formed by the canonical coordinates: we will denote this arrangement by the $k$-tuple $(1,2,3,\dots,k)$. Instead, the product of $\omega$'s we have found takes the canonical coordinates in another configuration: for example, the canonical coordinate $u_1$ is not taken in the first position but in the $n$-th in both cases $k=4n-1$ or $k=4n-3$: the corresponding $k$-tuple is of the form
\[(\dots,\underset{\text{$n$-th}}{1},\dots).\] Again we find that the only admissible case is $n=1$, and so $k=3$. This completes the proof. \endproof

\subsection{Symmetries and Quasi-Periodicity of Stokes matrices along the small quantum locus}\label{symcpn}
In this section we describe an interesting property of \emph{quasi-periodicity} of the Stokes matrices  computed at a point of the small quantum cohomology of $\mathbb P$ with respect to all possible admissible lines $\ell$. Because of the discussion at the beginning of the previous section, we can do the computation at any point, say, to fix ideas, at $0\in QH^\bullet(\mathbb P)$.  

For this let us introduce, following \cite{guzzetti1}, a new labelling of Stokes rays which is useful in order to describe  the \emph{Stokes factors} in which the matrix $S$ factorizes. Let us fix an admissible line $\ell$ in $\mathbb C$ of direction $\phi$ and choose an admissible direction $\arg z=\tau$  in the universal cover $\mathcal R$ (namely,  $\phi-\tau=0$ mod $2\pi$), which projects onto $\ell_+$. We label the Stokes rays in $\mathcal R$ as follows: the rays are labelled in counter-clockwise order (i.e. increasing the value of the argument) starting from the first one in $\Pi_{\text{right}}$ which will be $ R_0$. In this way
\[ R_0,\dots, R_{k-1}\subseteq \Pi_{\text{right}},
\]
\[ R_k,\dots, R_{2k-1}\subseteq \Pi_{\text{left}}.
\]The labeling is then extended to all integers, increasing the index in counter-clockwise direction, so to obtain a whole family $\left\{ R_i\right\}_{i\in\mathbb Z}$. For the choice of $\ell$ with slope $0<\phi<\frac{\pi}{k}$ we have that
\[\text{the ray } R_0\text{ projects onto } R_{1k},
\]
\[\text{the ray } R_1\text{ projects onto } R_{1,k-1},
\]
\[\dots
\]where we do not use  the lexicographical labelling but the original one (see equation \eqref{stokesrayscpn}). If we denote by $\mathscr S_j$ the sector in  $\mathcal R$ bounded by $ R_{j-1}$ and $ R_{j+k}$, then $\mathscr S_j$ has angular width of $\pi+\frac{\pi}{k}$ and consequently there exists a unique  solution $\Xi_j$ of the system \eqref{eq2bis} with the asymptotic expansion $\eta \Psi^{-1}Y_{\rm formal}(z,u)$ on $\mathscr S_j$. We define the Stokes factors to be the connection matrices $K_j$ such that
\[\Xi_{j+1}=\Xi_jK_j,\quad j\in\mathbb Z.
\]In this way we have that
\begin{equation}
\label{17agosto2018-9}
S=K_0K_1\dots K_{k-2}K_{k-1}.
\end{equation}
Moreover, notice that the first row of $\Xi_j(z)$ is equal to $z^{\frac{k-1}{2}}\boldsymbol{\Phi}_j(z)$, where  $\boldsymbol{\Phi}_j(z)$ is a {\it row} vector whose entries form a basis of independent solutions of equation \eqref{eq3}. Therefore, we have that  
$$
\boldsymbol{\Phi}_{j+1}=\boldsymbol{\Phi}_jK_j.$$
 Notice that if
\[
\boldsymbol{F}(z):=\left(\frac{1}{\sqrt{k}}\frac{1}{z^{\frac{k-1}{2}}}\exp(kz), \frac{1}{\sqrt{k}}\frac{e^\frac{i\pi}{k}}{z^{\frac{k-1}{2}}}\exp(ke^\frac{2\pi i}{k}z),\dots, \frac{1}{\sqrt{k}}\frac{e^{\frac{i\pi}{k}(k-1)}}{z^{\frac{k-1}{2}}}\exp(ke^{\frac{2\pi i}{k}(k-1)}z)\right)
\]
is the row vector whose entries are the first term of the asymptotic expansions of an actual basis of solution $\boldsymbol{\Phi}(z)$ of the  equation \eqref{eq3}, it is easily seen that
\[
\boldsymbol{F}(ze^{\frac{2\pi i}{k}})=
\boldsymbol{F}(z)T_F,\quad T_F=\begin{pmatrix}
0&&&&\dots&1\\
-1&0&&&&\\
&-1&0&&&\\
&&-1&&&\\
&&&\ddots&\ddots&\vdots\\
&&&&-1&0
\end{pmatrix}.
\]
As a consequence, if $\boldsymbol{\Phi}_m(z)$ is the unique genuine solution of the hypergeometric equation such that
\[
\boldsymbol{\Phi}_m(z)\sim \boldsymbol{F}(z)\quad z\to\infty\quad z\in\mathscr S_m,
\]
then 
\[
\boldsymbol{\Phi}_{m+2}(ze^{\frac{2\pi i}{k}})\sim \boldsymbol{F}(ze^{\frac{2\pi i}{k}})=F(z)T_F\quad z\in\mathscr S_m,
\]
so that
\[\boldsymbol{\Phi}_{m+2}(ze^{\frac{2\pi i}{k}})T_F^{-1}\sim \boldsymbol{F}(z)\quad z\in\mathscr S_m.
\] 
By uniqueness, this implies that
\[
\boldsymbol{\Phi}_{m+2}(ze^{\frac{2\pi i}{k}})T_F^{-1}=\boldsymbol{\Phi}_m(z)\quad z\in\mathcal R.
\]
We deduce from this identity the following properties of the Stokes factors.
\begin{lemma}[\cite{guzzetti1}]
\label{17agosto2018-10} 
For any $m,q\in\mathbb Z$ the following identity holds
\[K_{m+2q}=T_F^{-q}K_mT_F^q.
\]
\end{lemma}
\proof We have from the definitions of the $K_i$'s that
\begin{align*}\boldsymbol{\Phi}_{m+1}(z)&=\boldsymbol{\Phi}_m(z)K_{m}=\boldsymbol{\Phi}_{m+2}(ze^{\frac{2\pi i}{k}})T_F^{-1}K_m\\
&=\boldsymbol{\Phi}_{m+3}(ze^{\frac{2\pi i}{k}})K_{m+2}^{-1}T_F^{-1}K_m\\
&=\boldsymbol{\Phi}_{m+1}(z)T_FK_{m+2}^{-1}T_F^{-1}K_m.
\end{align*}Hence, $K_{m+2}=T_F^{-1}K_mT_F^1$. A simple inductive argument completes the proof.\endproof

From this one can deduce the following

\begin{teorema}[\cite{guzzetti1}] Let $\ell$ be an admissible line, and let us enumerate the rays as described above, and introduce the corresponding Stokes factors $K_i$'s. The Stokes matrix of the system \eqref{eq2bis}, and equivalently of the hypergeometric equation \eqref{eq3}, for $k>3$, is given by
\[S=\begin{sistema}
(K_0K_1T_F^{-1})^\frac{k}{2}T_F^{\frac{k}{2}}\equiv T_F^{\frac{k}{2}}(T_F^{-1}K_{k-2}K_{k-1})^\frac{k}{2},\quad k\text{ even}\\
\\
(K_0K_1T_F^{-1})^\frac{k-1}{2}K_0T_F^{\frac{k-1}{2}}\equiv T_F^{\frac{k-1}{2}}K_{k-1}(T_F^{-1}K_{k-2}K_{k-1})^\frac{k-1}{2},\quad k\text{ odd.}
\\
\end{sistema}\]Moreover, the two Stokes factors $K_{k-2}$ and $K_{k-1}$ are given by:
\begin{itemize}
\item for $k\ \bold{ even}$ we have
\[(K_{k-2})_{2,1}=-\binom{k}{1},\quad (K_{k-2})_{j,j}=1\ \text{for }j=1,\dots,k
\]
\[(K_{k-2})_{j,k-j+3}=\binom{k}{2j-3}\ \text{for }j=3,\dots,\frac{k}{2}+1
\]
\[(K_{k-1})_{j,j}=1\ \text{for }j=1,\dots,k,\quad (K_{k-1})_{j,k-j+2}=\binom{k}{2(j-1)}\ \text{for }j=2,\dots,\frac{k}{2}
\]and all other entries of $K_{k-2},K_{k-1}$ are zero.
\item for $k\ \bold{ odd}$ we have
\[(K_{k-2})_{2,1}=-\binom{k}{1},\quad (K_{k-2})_{j,j}=1\ \text{for }j=1,\dots,k
\]
\[(K_{k-2})_{j,k-j+3}=\binom{k}{2j-3}\ \text{for }j=3,\dots,\frac{k+1}{2}
\]
\[(K_{k-1})_{j,j}=1\ \text{for }j=1,\dots,k,\quad (K_{k-1})_{j,k-j+2}=\binom{k}{2(j-1)}\ \text{for }j=2,\dots,\frac{k+1}{2}
\]and all other entries of $K_{k-2},K_{k-1}$ are zero.
\end{itemize}
\end{teorema}

The above theorem was a crucial step introduced in \cite{guzzetti1} in order to explicitly compute the Stokes matrices and  prove point (3a) of Conjecture \ref{congettura}. 
With this results, we can now summarize the symmetries and quasi-periodicity relations of the Stokes matrices

\bsh
\begin{teorema}\label{quasi-period.th1}
Let $p$ be a point of the small quantum cohomology of $\mathbb P_{\mathbb C}^{k-1}$, let $\ell$ be an admissible line at $p$, and denote by $S_p(\ell)^{\text{lex}}$ the Stokes matrix computed at $p$, with respect to the line $\ell$ and in the $\ell$-lexicographical order. Then the entries of the matrices
\[S_p(\ell)^{\text{lex}}\quad\text{and}\quad S_p(e^\frac{2\pi i}{k}\ell)^{\text{lex}}
\]differ just by some signs. Moreover, also the entries
\[\left(S_p(\ell)\right)^{\text{lex}}_{j,j+1}\quad\text{and}\quad\left(S_p(e^\frac{\pi i}{k}\ell)\right)^{\text{lex}}_{j,j+1}
\]differ by some signs for any $i=1,\dots,k-1$. Thus, in particular the $(k-1)$-tuple
\[\left(\left|\left(S_p(\ell)^{\text{lex}}\right)_{1,2}\right|,\left|\left(S_p(\ell)^{\text{lex}}\right)_{2,3}\right|,\dots,\left|\left(S_p(\ell)^{\text{lex}}\right)_{k-1,k}\right|\right)
\]does not depend on $p$ and $\ell$. In particular, it is equal to
\[\left(\binom{k}{1},\dots,\binom{k}{k-1}\right).
\] 
\end{teorema}
\esh

\proof 
We  restrict to the case $k$ \emph{even}, the case $k$ \emph{odd} being analogous. 
In order to prove the theorem, it sufficies to compute monodromy data at $p=0\in QH^\bullet(\mathbb P)$. 
 For brevity, we will omit the index $p=0$ from the Stokes matrix $S_p(\ell)$.  So, if we  fix an admissible line $\ell$, and if $S(\ell)^{\text{lex}}$ is the Stokes matrix in the lexicographical form, then 
  $$
  S(e^\frac{2\pi i}{k}\ell)^{\text{lex}}= \left(S(\ell)^{\text{lex}}\right)^{\omega_1\omega_2}\quad\text{or}\quad\left(S(\ell)^{\text{lex}}\right)^{\omega_2\omega_1}.$$
In order to put 
\[S(\ell)=K_0K_1\dots K_{k-2}K_{k-1},
\]
of formula \eqref{17agosto2018-9} in lexicographical order, we act by conjugation $S^{\rm lex}(\ell)=PS(\ell)P^{-1}$, with  a \emph{unique} permutation matrix $P$ corresponding to the $\ell$-lexicographical order. We have
\[S(e^\frac{2\pi i}{k}\ell)=K_2K_3\dots K_kK_{k+1},
\]
and according to Lemma  \ref{17agosto2018-10} 
\[K_k=T_F^{-\frac{k}{2}}K_0T_F^{\frac{k}{2}},\quad K_{k+1}=T_F^{-\frac{k}{2}}K_1T_F^{\frac{k}{2}}.
\]
Thus,
\begin{align*}
S(e^\frac{2\pi i}{k}\ell)&=T_F^{\frac{k}{2}}\left(T_F^{-1}K_kK_{k+1}\right)^\frac{k}{2}\\
&=T_F^{\frac{k}{2}}\left(T_F^{-1}T_F^{-\frac{k}{2}}K_0K_1T_F^{\frac{k}{2}}\right)^\frac{k}{2}\\
&=T_F^{-1}(K_0K_1T_F^{-1})^\frac{k}{2}T_F^{\frac{k}{2}+1}\\
&=T_F^{-1}S(\ell)T_F.
\end{align*}
If we want to put the matrix $S(e^\frac{2\pi i}{k}\ell)$ in the lexicographical form, we have to conjugate it by a suitable permutation matrix $Q$ (corresponding to the lexicographical order with respect to the rotated line $e^\frac{2\pi i}{k}\ell$):
\begin{align*}S(e^\frac{2\pi i}{k}\ell)^{\text{lex}}&=Q\cdot S(e^\frac{2\pi i}{k}\ell)\cdot Q^{-1}\\
&=(QT_F^{-1})\cdot S(\ell)\cdot (QT_F^{-1})^{-1}.
\end{align*}

By definition of $T_F$, we clearly have that 
\[QT_F^{-1}=JP_1,
\]where $P_1$ is a permutation matrix and $J$ is a matrix of the form ${\rm diag}(\pm 1,\dots,\pm 1)$ (in particular, there will be $k-1$ times entries $(-1)$'s and just one entry $(+1)$, as in the matrix $T_F^{-1}$). Consequently, we have that
\[J^{-1}S(e^\frac{2\pi i}{k}\ell)^{\text{lex}}J=P_1\cdot S(\ell)\cdot P_1,
\]and since the lhs is upper triangular we conclude that $P=P_1$, by uniqueness of the lexicographical order. This proves the first statement.

For the second statement, it is sufficient to prove it just for the choice of $\ell$ with slope $0<\phi<\frac{\pi}{k}$. From the explicit expressions for the Stokes factors $K_{k-2}$ and $K_{k-1}$ of the previous Theorem, after some computations, one finds that the entries in the first upper-diagonals of the matrix $S(\ell)^{\text{lex}}$ are 

\[\begin{pmatrix}
1&-\binom{k}{1}&-\binom{k}{1}\binom{k}{k-2}+\binom{k}{k-1}&\binom{k}{1}^2-\binom{k}{2}&\dots&&\dots&&\\
&&&&&&&&\\
&1&\binom{k}{k-2}&-\binom{k}{1}&\dots&&\dots&&\\
&&&&&&&&\\
&&1&-\binom{k}{3}&\dots&&\dots&&\\
&&&&&&&&\\
&&&\ddots&\ddots&&&&\\
&&&&&&&&\\
&&&&1&-\binom{k}{k-3}&-\binom{k}{k-3}\binom{k}{2}+\binom{k}{k-1}&\binom{k}{k-3}\binom{k}{1}-\binom{k}{k-2}&\\
&&&&&&&&\\
&&&&&1&\binom{k}{2}&-\binom{k}{1}&\\
&&&&&&&&\\
&&&&&&1&-\binom{k}{k-1}&\\
&&&&&&&1&\\
\end{pmatrix},
\]
i.e. along the diagonals we have the general form
\[
\left(\begin{array}{cccccccccc}
&&&&&&&&&\\
&&\ddots&\ddots&&&&&&\\
&&&1&\color{red}{-\binom{k}{2n-1}}&-\binom{k}{2n-1}\binom{k}{k-2n}+\binom{k}{k-1}&\color{blue}{\binom{k}{2n-1}\binom{k}{1}-\binom{k}{2n}}&\dots&&\\
&&&&&&&&&\\
&&&&1&\binom{k}{k-2n}&\color{red}{-\binom{k}{1}}&\dots&&\\
&&&&&&&&&\\
&&&&&1&-\binom{k}{2n+1}&\dots&&\\
&&&&&&1&\ddots&&\\
\end{array}\right)
\]for $n=1,\dots,\frac{k}{2}-1$. Since the Stokes matrix $S(e^{\frac{\pi i}{k}}\ell)^{\text{lex}}$ is equal to $(S(\ell)^{\text{lex}})^{\omega_1}=A^{\omega_1}\cdot S(\ell)^{\text{lex}}\cdot A^{\omega_1}$, where
\[A^{\omega_1}=\begin{pmatrix}
0&1&&&&&\\
1&\binom{k}{1}&&&&&\\
&&0&1&&&\\
&&1&\binom{k}{3}&&&\\
&&&&\ddots&&\\
&&&&&0&1\\
&&&&&1&\binom{k}{k-1}
\end{pmatrix},
\]we find that
\[\left(S(e^{\frac{\pi i}{k}}\ell)^{\text{lex}}\right)_{2i+1,2i+2}=\binom{k}{2i+1}\quad\quad i=0,\dots,\frac{k}{2}-1
\]
and 
\begin{align*}\left(S(e^{\frac{\pi i}{k}}\ell)^{\text{lex}}\right)_{2i,2i+1}&=\left(S(\ell)^{\text{lex}}\right)_{2i-1,2i+2}-\left(S(\ell)^{\text{lex}}\right)_{2i-1,2i-1}\cdot \left(S(\ell)^{\text{lex}}\right)_{2i,2i+2}\\
&=\color{blue}{\binom{k}{2n-1}\binom{k}{1}-\binom{k}{2n}}\color{black}{+}\color{red}{\binom{k}{2n-1}\binom{k}{1}}\\
&=-\binom{k}{2n}
\end{align*}for $i=1,\dots,\frac{k}{2}-1$. This completes the proof.\endproof

\bsh
\begin{cor}\label{corbeilip12}
The Beilinson exceptional collection $\frak B$ corresponds to the monodromy data computed at some point of the small quantum cohomology of $\mathbb P^{k-1}_{\mathbb C}$ if and only if $k=2,3$.
\end{cor}
\esh

\proof Note that the inverse of the Gram matrix of the Grothendieck--Euler--Poincaré product, which would coincide with the Stokes matrix, has the following entries on the upper diagonal:
\[(-k,-k,\dots, -k,-k).
\]
\endproof
\begin{oss}
Note that the Corollary above cannot be deduced from Theorems \ref{risultato1} and \ref{teotreccebeilinson}. The reason is that \emph{a priori} the subgroup of $\mathcal B_k$ of braids fixing up to shifts the Beilinson exceptional collection $\frak B$ 
\[\left\{\beta\in \mathcal B_k\colon \frak B^\beta\equiv \frak B[{\bf m}] \right\},\quad {\bf m}:=(m_1,\dots, m_{k})\in\mathbb Z^k,
\]could be non-trivial. In general, it is still an open problem to study transitiveness and freeness of the braid group action on the set of exceptional collections. See \cite{helix} fur further details.
\end{oss}

\newpage
\section{Proof of the Main Conjecture for Grassmannians}\label{chcoalgrass}
In what follows
\begin{itemize}
\item $r,k$ will be natural numbers such that $0< r< k$.
\item We will denote by $\mathbb P$ the complex projective space $\mathbb P^{k-1}_{\mathbb C}$;
\item $\mathbb G$ will be the complex Grassmannian $\mathbb{G}(r,k)$ of $r$-planes in $\mathbb C^k$;
\item $\Pi$ will denote the cartesian product
\[\underbrace{\mathbb P\times\dots\times\mathbb P}_{r\text{ times}}.
\]
\item $\sigma\in H^2(\mathbb P,\mathbb C)$ will be the generator of the cohomology of $\mathbb P$, normalized so that
\[\int_{\mathbb P}\sigma^{k-1}=1.
\]We will denote the power $\sigma^h$, with $h\in\mathbb N$, by $\sigma_h$.
\end{itemize}

As an immediate generalization of the case of complex projective spaces, the Frobenius algebra structure defined by the (small) quantum cohomology of complex Grassmannians has been one of the first examples extensively studied both in physics \cite{vafa, wit1} and mathematical literature \cite{siebti, bert, buc1}. Here, we show how the validity of Conjecture \ref{congettura} for all complex Grassmannians can be directly deduced from the explicit results for projective spaces obtained in the previous Section. We also show validity of a property of \emph{quasi-periodicity} of the Stokes matrices along points of the small quantum cohomology, analogous to the one described in Theorem \ref{quasi-period.th1}. The main tool is an identification of the classical/quantum cohomology of the Grassmannian $\mathbb G$ with an \emph{exterior power} of the classical/quantum cohomology of $\mathbb P$. Such an identification has been described in the literature from many perspectives, and we briefly summarize it both in the classical and in the quantum case.

\subsection{Results on the Classical Cohomology of $\mathbb G$}\label{secclasscohgrass}
The complex Grassmannian $\mathbb G$ can be seen as a symplectic quotient. Let us consider the complex vector space $\Hom(\mathbb C^r,\mathbb C^k)$ endowed with its standard symplectic structure: if we introduce on $\Hom(\mathbb C^r,\mathbb C^k)$ coordinates $a_{ij}=x_{ij}+\sqrt{-1}y_{ij}$, for $1\leq i\leq k$ and $1\leq j\leq r$, then the standard symplectic structure is
\[\omega:=\sum_{i,j}dx_{ij}\wedge dy_{ij}.
\]Let us consider the action of $U(r)$ on $\Hom(\mathbb C^r,\mathbb C^k)$ defined by $g\cdot A:=A\circ g^{-1}$: this action is Hamiltonian and a moment map $\mu_{U(r)}\colon \Hom(\mathbb C^r,\mathbb C^k)\to \frak{u}(r)$ is given by
\[\mu_{U(r)}(A):=A^\dagger A-\mathbbm 1.
\]Since the subset $\mu^{-1}_{U(r)}(0)$ is the set of unitary $r$-frames in $\mathbb C^k$, we have clearly the identification
\[\mathbb G\cong\Hom(\mathbb C^r,\mathbb C^k)\sslash U(r):=\mu_{U(r)}^{-1}(0)/U(r).
\]
If $\mathbb T\subseteq U(r)$ is the subgroup of diagonal matrices, then $\mathbb T\cong U(1)^{\times r}$ is a maximal torus. Denoting by $\mu_\mathbb T\colon \Hom(\mathbb C^r,\mathbb C^k)\to \frak{u}(1)^{\times r}$ the composition of $\mu_{U(r)}$ and the canonical projection $\frak u(r)\to\frak u(1)^{\times r}$, we have that
$\mu_{\mathbb T}^{-1}(0)$ is the set of matrices $A\in M_{k,r}(\mathbb C)$ whose columns have unit length. Hence, we have
\[\Pi\cong\Hom(\mathbb C^r,\mathbb C^k)\sslash \mathbb T:=\mu_\mathbb T^{-1}(0)/\mathbb T.
\]Moreover, the quotient
\[\mu_{U(r)}^{-1}(0)/\mathbb T\]
can be identified with the flag manifold $\mathbb F:=\operatorname{Fl}(1,2,\dots,r,k)$ (for the identification we have to choose a Hermitian metric on $\mathbb C^k$, e.g. the standard one, compatible with the standard symplectic structure). Because of the inclusion $\mu_{U(r)}^{-1}(0)\subseteq \mu_\mathbb T^{-1}(0)$,
we have the following quotient diagram:
\[\xymatrix{
\mathbb F\ \ar[d]_{p}\ar@{^{(}->}[r]^{\iota}&\Pi\\
\mathbb G&
}
\]where $p$ is the canonical projecton, and $\iota$ the inclusion. Note that in this way there is also a natural rational map ``\emph{taking the span}''
\[\xymatrix{\Pi\ar@{-->}[r]&\mathbb G\colon\ (\ell_1,\dots,\ell_r)\mapsto \operatorname{span}\langle\ell_1,\dots,\ell_r\rangle},
\]whose domain is the image of $\iota$. On the manifold $\Pi$ we have $r$ canonical line bundles, denoted $\frak L_j$ for $j=1,\dots,r$, defined as the pull-back of the bundle $\mathcal O(1)$ on the $j$-th factor $\mathbb P$. If we denote $\frak V_1\subseteq\frak V_2\subseteq\dots\subseteq \frak V_r$ the tautological bundles over $\mathbb F$, we have that
\[\iota^*\frak L_j\cong (\frak V_j/\frak V_{j-1})^\vee.
\]
Denoting with the same symbol $x_i$ the Chern class $c_1(\frak L_i)$ on $\Pi$ and its pull-back $c_1(\iota^*\frak L_i)=\iota^*c_1(\frak L_i)$ on $\mathbb F$, we have
\[H^\bullet(\Pi,\mathbb C)\cong H^\bullet(\mathbb P,\mathbb C)^{\otimes r}\cong \frac{\mathbb C[x_1,\dots,x_r]}{\langle x_1^k,\dots x_r^k\rangle}\quad\text{(by Künneth Theorem),}
\]
\[H^\bullet(\mathbb F,\mathbb C)\cong\frac{\mathbb C[x_1,\dots,x_r]}{\langle h_{k-r+1},\dots, h_k\rangle},
\]where $h_j$ stands for the $j$-th complete symmetric polynomial in $x_1,\dots, x_r$. Since the classes $x_1,\dots, x_r$ are the Chern roots of the dual of the tautological bundle $\frak V_r$, we also have
\[H^{\bullet}\left(\mathbb G,\mathbb C\right)\cong\frac{\mathbb C[e_1,\dots, e_r]}{\langle h_{n-k+1},\dots,h_n\rangle}\cong \frac{\mathbb C[x_1,\dots,x_k]^{\frak S_k}}{\langle h_{n-k+1},\dots,h_n\rangle},
\]where the $e_j$'s are the elementary symmetric polynomials in $x_1,\dots, x_r$. This is the classical representation of the cohomology ring of the Grassmannian $\mathbb G$ with generators the Chern classes of the dual of the tautological vector bundle $\mathcal S$, and relations generated by the Segre classes of $\mathcal S$.\\
From this presentation of algebras, it is clear that any cohomology class of $\mathbb G$ can be \emph{lifted} to a cohomology class of $\Pi$: we will say that $\tilde \gamma\in H^\bullet(\Pi,\mathbb C)$ is the lift of $\gamma\in H^\bullet(\mathbb G,\mathbb C)$ if $p^*\gamma=\iota^*\tilde\gamma$. The following integration formula allow us to express the cohomology pairings on $H^\bullet(\mathbb G,\mathbb C)$ in terms of the cohomology pairings on $H^\bullet(\Pi,\mathbb C)$.

\bsh
\begin{teorema}[\cite{martin}] If $\gamma\in H^\bullet(\mathbb G,\mathbb C)$ admits the lift $\tilde\gamma\in H^\bullet(\Pi,\mathbb C)$, then
\begin{equation}\label{martin}\int_\mathbb G\gamma=\frac{(-1)^{\binom{r}{2}}}{r!}\int_\Pi\tilde\gamma\cup_{\Pi}\Delta^2,
\end{equation}where \[\Delta:=\prod_{1\leq i< j\leq r}(x_i-x_j).\]
\end{teorema}
\esh

\bsh
\begin{cor}[\cite{ell}] \label{ell-str}The linear morphism
\[\vartheta\colon H^\bullet(\mathbb G,\mathbb C)\to H^\bullet(\Pi,\mathbb C)\colon\gamma\mapsto\tilde\gamma\cup_{\Pi}\Delta
\]is injective, and its image is the subspace of antisymmetric part of $H^\bullet(\Pi,\mathbb C)$ with respect to the $\frak S_r$-action. Moreover
\[\vartheta(\alpha\cup_{\mathbb G}\beta)=\vartheta(\alpha)\cup_{\Pi}\tilde\beta=\tilde\alpha\cup_{\Pi}\vartheta(\beta).
\]
\end{cor}
\esh
\proof If $\vartheta(\gamma)=0$, then 
\[\int_\mathbb G\gamma\cup\gamma'=\frac{(-1)^{\binom{r}{2}}}{r!}\int_\Pi\left(\tilde\gamma\cup\Delta\right)\cup\left(\tilde\gamma'\cup\Delta\right)=0
\]for all $\gamma'\in H^\bullet(\mathbb G,\mathbb C)$. Then $\gamma=0$. Being clear that $\vartheta(\gamma)$ is antisymmetric, observe that any antisymmetric class is of the form $\tilde\gamma\cup\Delta$ with $\tilde\gamma$ symmetric in $x_1,\dots,x_r$. The last statement follows from the fact that the lift of a cup product is the cup product of the lifts.
\endproof

We can identify the antisymmetric part of $H^\bullet(\Pi,\mathbb C)\cong H^\bullet(\mathbb P,\mathbb C)^{\otimes r}$ with $\bigwedge\nolimits^rH^\bullet(\mathbb P,\mathbb C)$, using the identifications $i,j$ illustrated in the following diagram
\[\xymatrix{
H^\bullet(\mathbb P,\mathbb C)^{\otimes r}\ar[r]^{\pi}&\bigwedge\nolimits^rH^\bullet(\mathbb P,\mathbb C)\ar[dl]_{i}\\
[H^\bullet(\mathbb P,\mathbb C)^{\otimes r}]^{\text{ant}}\ar@{^{(}->}[u]\ar@/_1pc/[ur]_{j}&
}\]
where
\[\pi\colon H^\bullet(\mathbb P,\mathbb C)^{\otimes r}\to\bigwedge\nolimits^rH^\bullet(\mathbb P,\mathbb C)\colon\alpha_1\otimes\dots\otimes\alpha_r\mapsto\alpha_1\wedge\dots\wedge\alpha_r,
\]

\[i\colon\bigwedge\nolimits^rH^\bullet(\mathbb P,\mathbb C)\to [H^\bullet(\mathbb P,\mathbb C)^{\otimes r}]^{\text{ant}}\colon\alpha_1\wedge\dots\wedge\alpha_r\mapsto\sum_{\rho\in\frak S_r}\varepsilon(\rho)\alpha_{\rho(1)}\otimes\dots\otimes\alpha_{\rho(r)},
\]together with its inverse
\[j\colon[H^\bullet(\mathbb P,\mathbb C)^{\otimes r}]^{\text{ant}}\to\bigwedge\nolimits^rH^\bullet(\mathbb P,\mathbb C)\colon \alpha_1\otimes \dots\otimes \alpha_r\mapsto \frac{1}{r!}\alpha_1\wedge\dots\wedge\alpha_r.
\]
The Poincaré pairing $g^{\mathbb P}$ on $H^\bullet(\mathbb P,\mathbb C)$ induces a metric $g^{\otimes\mathbb P}$ on $H^\bullet(\mathbb P,\mathbb C)^{\otimes r}$ and a metric $g^{\wedge\mathbb P}$ on $\bigwedge\nolimits^rH^\bullet(\mathbb P,\mathbb C)$ given by
\[g^{\otimes\mathbb P}(\alpha_1\otimes\dots\otimes\alpha_r,\beta_1\otimes\dots\otimes\beta_r):=\prod_{i=1}^rg^{\mathbb P}(\alpha_i,\beta_i),
\]
\[g^{\wedge\mathbb P}(\alpha_1\wedge\dots\wedge\alpha_r,\beta_1\wedge\dots\wedge\beta_r):=\det\left(g^{\mathbb P}(\alpha_i,\beta_j)\right)_{1\leq i,j\leq r}.
\]
Using the identifications above, when $g^{\otimes\mathbb P}$ is restricted on the subspace $[H^\bullet(\mathbb P,\mathbb C)^{\otimes r}]^{\text{ant}}$ it coincides with $r!g^{\wedge\mathbb P}$ on $\bigwedge\nolimits^rH^\bullet(\mathbb P,\mathbb C)$. From the integration formula \eqref{martin}, we deduce the following result.

\bsh
\begin{cor}\label{corisometry}
The isomorphism
\[j\circ\vartheta\colon\left( H^\bullet(\mathbb G,\mathbb C), g^{\mathbb G}\right)\to\left(\bigwedge\nolimits^rH^\bullet(\mathbb P,\mathbb C), (-1)^{\binom{r}{2}}g^{\wedge\mathbb P}\right)
\]is an isometry.
\end{cor}
\esh
An additive basis of $H^\bullet(\mathbb G,\mathbb C)$ is given by the Schubert classes (Poincaré-dual to the Schubert cycles), given in terms of $x_1,\dots, x_r$ by the Schur polynomials
\[\sigma_\lambda:=\frac{\det\begin{pmatrix}
x_1^{\lambda_1+r-1}&x_1^{\lambda_2+r-2} &\dots& x_1^{\lambda_r}\\
x_2^{\lambda_1+r-1}&x_2^{\lambda_2+r-2} &\dots& x_2^{\lambda_r}\\
&&\vdots\\
x_r^{\lambda_1+r-1}&x_r^{\lambda_2+r-2} &\dots& x_r^{\lambda_r}
\end{pmatrix}
}{\det
\begin{pmatrix}
x_1^{r-1}&x_1^{r-2}&\dots&1\\
x_2^{r-1}&x_2^{r-2}&\dots&1\\
&&\vdots\\
x_r^{r-1}&x_r^{r-2}&\dots&1
\end{pmatrix}
}
\]where $\lambda$ is a partition whose corresponding Young diagram is contained in in a $r\times (k-r)$ rectangle. The lift of each Schubert class to $H^\bullet(\Pi,\mathbb C)$ is the Schur polynomial in $x_1,\dots, x_r$ (indeed each $x_i$ in the Schur polynomial has exponent at most $k-r<k$). Thus, under the identification above, the class $j\circ\vartheta(\sigma_\lambda)$ is $\sigma_{\lambda_1+r-1}\wedge\dots\wedge\sigma_{\lambda_r}\in\bigwedge\nolimits^rH^\bullet(\mathbb P,\mathbb C)$, $\sigma$ being the generator of $H^2(\mathbb P,\mathbb C)$.

Using the Künneth isomorphism $H^\bullet(\Pi,\mathbb C)\cong H^\bullet(\mathbb P,\mathbb C)^{\otimes r}$, the cup product $\cup_{\Pi}$ is expressed in terms of $\cup_{\mathbb P}$ as follows: 
\[\left(\sum_i\alpha_1^i\otimes\dots\otimes \alpha_r^i\right)\cup_{\Pi}\left(\sum_j\beta_1^j\otimes\dots\otimes \beta_r^j\right)=\sum_{i,j}(\alpha_1^i\cup_{\mathbb P}\beta_1^j)\otimes\dots\otimes(\alpha_r^i\cup_{\mathbb P}\beta_r^j).
\]

If $\gamma\in H^\bullet(\Pi,\mathbb C)^{\frak S_r}$, then $\gamma\cup_{\Pi}(-)\colon H^\bullet(\Pi,\mathbb C)\to H^\bullet(\Pi,\mathbb C)$ leaves invariant the subspace of anty-symmetric classes. Thus, $\gamma\cup_{\Pi}(-)$ induces an endomorphism $A_{\gamma}\in\End\left(\bigwedge\nolimits^rH^\bullet(\mathbb P,\mathbb C)\right)$ that acts on decomposable elements $\alpha=\alpha_1\wedge\dots\wedge\alpha_r$ as follows 
\begin{equation}\label{22-Feb-16-1}A_\gamma(\alpha)=j(\gamma\cup_{\Pi}i(\alpha))=\frac{1}{r!}\sum_{i,\rho}\varepsilon(\rho)(\gamma_1^i\cup_{\mathbb P}\alpha_{\rho(1)})\wedge\dots\wedge(\gamma_r^i\cup_{\mathbb P}\alpha_{\rho(r)}),
\end{equation}
where $\gamma_j^i\in H^\bullet(\mathbb P,\mathbb C)$ are such that
\[\gamma=\sum_{i}\gamma_1^i\otimes\dots\otimes\gamma_r^i.
\]

As an example, in the following Proposition we reformulate  in $\bigwedge\nolimits^rH^\bullet(\mathbb P,\mathbb C)$ the classical Pieri formula, expressing the multiplication by a special Schubert class $\sigma_\ell$ in $H^\bullet(\mathbb G,\mathbb C)$
\[\sigma_\ell\cup_{\mathbb G}\sigma_\mu=\sum_{\nu}\sigma_\nu,
\] where the sum is on all partitions $\nu$ which belong to the set $\mu\otimes\ell$ (the set of partitions obtained by adding $\ell$ boxes to $\mu$, at most one per column) and which are contained in the rectangle $r\times(k-r)$,  in terms of the multiplication by $\sigma_\ell=(\sigma)^\ell\in H^\bullet(\mathbb P,\mathbb C)$. We also make explicit the operation of multiplication by the classes $p_\ell\in H^\bullet(\mathbb G,\mathbb C)$ defined in terms of the special Schubert classes by
\[p_\ell:=-\left(\sum_{\substack{n_1+2n_2+\dots+rn_r=\ell \\ n_1,\dots, n_r\geq 0}}\frac{\ell(n_1+\dots+n_r-1)!}{n_1!\dots n_r!}\prod_{i=1}^r(-\sigma_i)^{n_i}\right),\quad \ell=0,\dots, k-1,
\]because of the nice form of their lifts $\tilde p_\ell\in H^\bullet(\Pi,\mathbb C)$.

\bsh
\begin{prop}\label{23-Feb-16-1}
If $\sigma_\mu\in H^\bullet(\mathbb G,\mathbb C)$ is a Schubert class then 
\begin{itemize}
\item the product $\sigma_\ell\cup_{\mathbb G}\sigma_\mu$ with a special Schubert class $\sigma_\ell$ is given by
\[j\circ\vartheta(\sigma_\ell\cup_{\mathbb G}\sigma_\mu)=\frac{1}{r!}\left(\sum_{\substack{i_1+\dots+i_r=\ell \\ i_1,\dots, i_r\geq 0}}\sum_{\rho\in\frak S_r}\bigwedge_{j=1}^r\sigma_{i_{\rho(j)}}\cup_{\mathbb P}\sigma_{\mu_j+r-j}\right);
\]
\item the product $p_\ell\cup_{\mathbb G}\sigma_\mu$ is given by
\[j\circ\vartheta(p_\ell\cup_{\mathbb G}\sigma_\mu)=\sum_{i=1}^r\sigma_{\mu_1+r-1}\wedge\dots\wedge(\sigma_{\mu_i+r-i}\cup_{\mathbb P}\sigma_\ell)\wedge\dots\wedge\sigma_{\mu_r}.
\]
\end{itemize}
\end{prop}
\esh

\proof From Corollary \eqref{ell-str} we have 
\[\vartheta(\sigma_\ell\cup_{\mathbb G}\sigma_\mu)=\tilde\sigma_\ell\cup_{\Pi}\vartheta(\sigma_\mu)
\]
If $\gamma=\tilde{\sigma}_\ell$ is the lift of the special Schubert class $\sigma_\ell\in H^\bullet(\mathbb G,\mathbb C)$, then 
\[\tilde\sigma_\ell=h_\ell(x_1,\dots, x_r)=\sum_{\substack{i_1+\dots+i_r=\ell \\ i_1,\dots, i_r\geq 0}}\sigma_{i_1}\otimes\dots\otimes\sigma_{i_r},
\]
and using \eqref{22-Feb-16-1} we easily conclude. Analogously, we have that

\[\tilde p_\ell=\sum_{i=1}^rx_i^\ell=\sum_{i=1}^r1\otimes\dots\otimes\underset{\text{$i$-th}}{\sigma_\ell}\otimes\dots\otimes 1,
\]and 
\[A_{\tilde p_\ell}(\alpha)=\sum_{i=1}^r\alpha_1\wedge\dots\wedge(\sigma_\ell\cup_{\mathbb P}\alpha_i)\wedge\dots\wedge\alpha_r.
\]\endproof

\bsh
\begin{cor}\label{corc1}For any $z\in\mathbb C^*$, any $t^2\sigma_1\in H^2(\mathbb G,\mathbb C)$, and any Schubert class $\sigma_\lambda\in H^\bullet(\mathbb G,\mathbb C)$, the following identity holds:
\[j\circ \vartheta(z^{t^2\sigma_1}\cup\sigma_\lambda)=\bigwedge_{j=1}^rz^{t^2\sigma}\cup\sigma_{\lambda_j+r-j}.
\]
\end{cor}
\esh
\proof
We have that
\begin{align*}
\bigwedge_{j=1}^rz^{t^2\sigma}\cup\sigma_{\lambda_j+r-j}&=\bigwedge_{j=1}^r\sum_{k_j=0}^\infty\frac{(\log z)^{k_j}}{k_j!}(t^2\sigma)^{k_j}\cup\sigma_{\lambda_j+r-j}\\
=&\sum_{k_1=0}^\infty\dots\sum_{k_r=0}^\infty\frac{(\log z)^{k_1+\dots+ k_r}}{k_1!\dots k_r!}\bigwedge_{j=1}^r(t^2\sigma)^{k_j}\cup\sigma_{\lambda_j+r-j}\\
&=\sum_{k=0}^\infty\frac{(\log z)^k}{k!}\sum_{k_1+\dots+k_r=k}\binom{k}{k_1\dots k_r}\bigwedge_{j=1}^r(t^2\sigma)^{k_j}\cup \sigma_{\lambda_j+r-j}\\
&=j\circ\vartheta\left(\sum_{k=0}^\infty\frac{(\log z)^k}{k!} ((t^2\sigma_1)^k\cup\sigma_\lambda)\right)\\
& =j\circ\vartheta(z^{t^2\sigma_1}\cup\sigma_\lambda).
\end{align*}
\endproof

\bsh
\begin{prop}\label{propmu} If $\mu^{\mathbb P}\in\End(H^\bullet(\mathbb P,\mathbb C))$ and $\mu^{\mathbb G}$ denotes the grading operator for the Projective Space and the Grassmannian respectively, defined as in \eqref{06.07.17-3}, then for all Schubert classes $\sigma_\lambda\in H^\bullet(\mathbb G,\mathbb C)$ the following identities hold:
\[j\circ\vartheta(\mu^\mathbb G\sigma_\lambda)=\sum_{j=1}^r\sigma_{\lambda_1+r-1}\wedge\dots\wedge\mu^{\mathbb P}\sigma_{\lambda_j+r-j}\wedge\dots\wedge\sigma_{\lambda_r},
\]
\[j\circ\vartheta (z^{\mu^\mathbb G}\sigma_\lambda)=\bigwedge_{j=1}^rz^{\mu^\mathbb P}\sigma_{\lambda_j+r-j},\quad z\in\mathbb C^*.
\]
\end{prop}
\esh
\proof
For the first identity notice that
\begin{align*}
(j\circ\vartheta)^{-1}&\left(\sum_{j=1}^r\sigma_{\lambda_1+r-1}\wedge\dots\wedge \left(\lambda_j+r-j-\frac{k-1}{2}\right)\sigma_{\lambda_j+r-j}\wedge\dots\wedge\sigma_{\lambda_r}\right)\\ &=\ \sigma_\lambda\cdot\left(\left(\sum_{j=1}^r\lambda_j\right)+r^2-\frac{r(r+1)}{2}-\frac{(k-1)r}{2}\right)\\
&=\ \sigma_\lambda\cdot\left(\left(\sum_{j=1}^r\lambda_j\right)-\frac{r(k-r)}{2}\right)\\
&=\ \mu^{\mathbb G}(\sigma_\lambda).
\end{align*}
For the second identity, we have that
\begin{align*}\bigwedge_{j=1}^rz^{\mu^\mathbb P}\sigma_{\lambda_j+r-j}&=\bigwedge_{j=1}^rz^{\mu^{\mathbb P}_{\lambda_j+r-j}}\cdot \sigma_{\lambda_j+r-j}\\
&=\exp\left(\log(z)\cdot\sum_{j=1}^r\mu^{\mathbb P}_{\lambda_j+r-j}\right)\cdot \bigwedge_{j=1}^r\sigma_{\lambda_j+r-j}\\
&=j\circ\vartheta\left(z^{\mu^\mathbb G}\sigma_\lambda\right).
\end{align*}
\endproof

\subsection{Quantum cohomology of $\mathbb G$ as exterior power of the quantum cohomology of $\mathbb P$}\label{secquantcohgrass} The identification in the classical cohomology setting of $H^\bullet(\mathbb G,\mathbb C)$ with the exterior power $\bigwedge^r H^\bullet(\mathbb P,\mathbb C)$ explained in the previous section, has been extended also to the quantum case from many different perspectives. 

The validity of this identification was already clear to physicists: it can be found already in Appendix A of \cite{cecotti-vafa}, in the description of the $\sigma$-model of $\mathbb G$ as ``\emph{equivalent up to $D$-terms}'' to a $\frak S_r$-quotient of a tensor product of $\sigma$-models of $\mathbb P$. Moreover, a similar relationship is outlined also in Appendix A of \cite{horivafa}. 

It was only with the paper  \cite{BCFK} that two different proofs of this identification for \emph{small} quantum cohomologies were given, the first one using localization and Grothendieck Quot schemes techinques, the second one being based on some explicit identities relating 3-points genus 0 Gromov--Witten invariants of $\mathbb G$ and $\mathbb P$ deduced from Vafa--Intriligator residue formula. 

Subsequently, these results were extended to a more general situation in \cite{BCFK2}: generalizing the conjecture by K. Hori and C. Vafa, the authors conjectured similar relationships between the Gromov--Witten theory in genus 0 of a non-abelian symplectic (or GIT, if you prefer) quotient $X\sslash G$ with that of the corresponding abelian quotient\footnote{Here, $X$ is assumed to be a smooth projective variety wth a linearized  action of a complex reductive group $G$, while $T$ denotes a maximal torus.} $X\sslash T$. These conjectures were shown to hold true for partial flag varieties. Subsequently, in \cite{CFKS}, the relationships described in \cite{BCFK, BCFK2} are summarized and extended in order to re-interpret them as identifications of Frobenius structures. In \cite{KS} two new tensorial operations of (germs of) Frobenius manifolds are introduced (symmetric and exterior powers): the quantum cohomology of $\mathbb G$ is described as \emph{$r$-th exterior power} of the quantum cohomology\footnote{More precisely, the germ of $QH^\bullet(\mathbb G)$ at a point of small quantum cohomology is identified with the germ of $QH^\bullet(\mathbb P)$ at a \emph{shifted} point of the small quantum cohomology. See Theorem \ref{teoBCFK} below. } of $\mathbb P$.

Contemporarily to the works mentioned above, in \cite{gatto1, gattobook1} L. Gatto described a new point of view for studying, in a greatly unified way, both classical and small quantum Schubert calculus.  The surprisingly simple Gatto's realization of Schubert calculus is based on the properties of Hasse--Schmidt derivations on the exterior algebra of a free module. Let $A$ be a commutative ring, and let $M$ be a free $A$-module of countable infinte rank. The $A$-modules $\bigwedge M$ and $\left(\bigwedge M\right) [\![ t]\!]$ are endowed with two natural $A$-algebras structures, the second one being defined through
\[\left(\sum_i\alpha_it^i\right)\wedge\left(\sum_j\beta_jt^j\right):=\sum_h\left(\sum_{i+j=h}\alpha_i\wedge\beta_j\right)t^{h}.
\]

\begin{defi}An $A$-module morphism $D\colon\bigwedge M\to\left(\bigwedge M\right) [\![ t]\!]$ is called a \emph{Hasse--Schmidt derivation} if it is an $A$-algebra homomorphism. i.e. such that
\begin{equation}\label{HSder}D(\alpha\wedge\beta)=D(\alpha)\wedge D(\beta),\quad \text{for all }\alpha,\beta\in \bigwedge M.
\end{equation}
A Hasse--Schmidt derivation can be expanded in series with respect to the indeterminate $t$, $D=\sum_{i\geq 0}D_it^i$ with $D_i\in{\rm End}_A(\bigwedge M)$ called \emph{components} of $D$.
If $\varepsilon:=(\varepsilon_i)_{i\geq 1}$ is a basis for the $A$-free module $M$, the unique Hasse--Schmidt derivation $\mathcal S$ such that $\mathcal S(\varepsilon_j)=\sum_{i\geq 0}\varepsilon_{i+j}t^i$ is called the ($\varepsilon$)-\emph{Schubert derivation}.
\end{defi} 

From the Schubert derivation $\mathcal S$ on $\bigwedge M$, L. Gatto was able to reconstruct both the classical ($A=\mathbb C$) and small quantum ($A=\mathbb C[q]$) Schubert calculus for $\mathbb G$ for all $(r,k)$, with $0<r<k$, at once. For a fixed $k$, if we denote by $M_k$ the $A$-submodule of $M$ generated by $(\varepsilon_1,\dots,\varepsilon_k)$, then the classical (resp. small quantum) cohomology of $\mathbb G$ can be realized (for all $r$ with $0<r<k$) as a quotient of the same commutative ring of endomorphisms of the exterior algebra of $M_k$. 
By varying $k$, one can realizes the totality of these rings as a quotient of the same ring of derivations on $\bigwedge M$. Notice that we can recover the abelian/non-abelian correspondence for Grassmannians  by indentifying the module $M_k$ with the classical/small quantum cohomology ring of $\mathbb P^{k-1}$.

Remarkably, in these realizations of the Schubert calculus, the content of the classical/small quantum Pieri rule is encoded in the single equation \eqref{HSder}, while classical/small quantum Giambelli can be deduced from it through simple algebraic manipulations evoking a formal \emph{integration by parts} procedure.

In subsequent works of D. Laksov and A. Thorup (\cite{lakthor1, lakthor2}) the results of Gatto were further extended in order to deal with more general Grassmann bundles and equivariant cohomology (see also \cite{laksov}). The point of view of Laksov and Thorup is quite different, being based on the fact that $\bigwedge^r M_k$ can be endowed with a structure of module over the ring of symmetric polynomials with coefficients in $A$. Their results were re-interpreted by L. Gatto and T. Santiago (\cite{gattosantiago1, gattosantiago2}) in terms of Hasse--Schmidt derivations, in a more unified framework. For further details, and many more applications of this formalism, we refer the  reader to the monograph \cite{gattobook2}.

The following isomorphism of the (small) quantum cohomology algebra of Grassmannians at a point $t\sigma_1=\log q\in H^2(\mathbb G,\mathbb C)$ is well-known (see \cite{wit1}, \cite{siebti}, \cite{bert}, \cite{buc1} for example)
\[QH_q^\bullet(\mathbb G)\cong\frac{\mathbb C[x_1,\dots,x_r]^{\frak S_r}[q]}{\langle h_{k-r+1},\dots,h_k-(-1)^{r-1}q\rangle},
\]while for the (small) quantum cohomology algebra of $\Pi$, being equal to the $r$-fold tensor product of the quantum cohomology algebra of $\mathbb P$ (see \cite{kauf}, \cite{dubro2}, \cite{manin}), we have 
\[QH_{q_1,\dots,q_r}^\bullet(\Pi)\cong \frac{\mathbb C[x_1,\dots,x_r][q_1,\dots, q_r]}{\langle x_1^k-q_1,\dots, x_r^k-q_r\rangle}.
\]
Following \cite{BCFK}, and interpreting now the parameters $q$'s just as formal parameters, if we denote by $\overline{QH}^\bullet_q(\Pi)$ the quotient of $QH_{q_1,\dots,q_r}^\bullet(\Pi)$ obtained by substituing $q_i=(-1)^{r-1}q$, and denoting the canonical projection by
\[[-]_q\colon QH_{q_1,\dots,q_r}^\bullet(\Pi)\to\overline{QH}^\bullet_q(\Pi),
\]
we can extend by linearity the morphisms $\vartheta,j$ of the previous section to morphisms
\[\overline{\vartheta}\colon QH_q^\bullet(\mathbb G)\to \overline{QH}^\bullet_q(\Pi),
\]
\[\overline{j}\colon \left[\overline{QH}^\bullet_q(\Pi)\right]^{\text{ant}}\to \left(\bigwedge^r H^\bullet (\mathbb P,\mathbb C)\right)\otimes_{\mathbb C} \mathbb C[q].
\]
Notice that the image under $\overline\vartheta$ of any Schubert class $\sigma_\lambda$ is equal to the classical product $\tilde\sigma_\lambda\cup_{\Pi}\Delta$, the exponents of $x_i$'s in the product $\sigma_\lambda(x)\prod_{i<j}(x_i-x_j)$ being less than $k$; as a consequence, the image of $\overline\vartheta$ is equal to the antisymmetric part with respect to the natural $\frak S_r$ action (permuting the $x_i$'s) $$\left[\overline{QH}^\bullet_q(\Pi)\right]^{\text{ant}}\cong [H^\bullet(\Pi,\mathbb C)]^{\text{ant}}\otimes_{\mathbb C}\mathbb C[q].$$ 

The following result, is a quantum generalization of Corollary \ref{ell-str}.

\bsh
\begin{teorema}[\cite{BCFK}] \label{teoBCFK}For any Schubert classes $\sigma_\lambda,\sigma_\mu\in H^\bullet (\mathbb G,\mathbb C)$ we have
\[\overline\vartheta(\sigma_\lambda*_{\mathbb G,q}\sigma_\mu)=\left[\vartheta(\sigma_\mu)*_{\Pi,q_1,\dots,q_r}\tilde \sigma_\lambda\right]_q.
\]   
\end{teorema}
\esh

Using the identification $\overline j$, we can deduce from the previous result the following generalization of Proposition \ref{23-Feb-16-1}.

\bsh
\begin{cor}\label{quantumprod}
If $\sigma_\mu\in H^\bullet (\mathbb G,\mathbb C)$ is a Schubert class then 
\begin{equation}\label{satakeprod}\overline{j}\circ\overline{\vartheta}(\sigma_\mu *_{\mathbb G,q} p_\ell)=\sum_{i=1}^r\sigma_{\mu_1+r-1}\wedge\dots\wedge \sigma_{\mu_i+r-i}*_{\mathbb P,(-1)^{r-1}q} \sigma_\ell\wedge\dots\wedge \sigma_{\mu_r}.\end{equation}
From this identity, it immediately follows that:
\begin{enumerate}
\item At the point $p=t^2\sigma_1\in H^2(\mathbb G,\mathbb C)$ of the small quantum cohomology of $\mathbb G$, the eigenvalues of the operator $$\mathcal U^\mathbb G_p:=c_1(\mathbb G)*_{q}(-)\colon H^\bullet(\mathbb G,\mathbb C)\to H^\bullet(\mathbb G,\mathbb C)$$ are given by the sums
\[u_{i_1}+\dots+u_{i_r},\quad 1\leq i_1<\dots<i_r\leq k,
\]where $u_1,\dots, u_k$ are the eigenvalues of the corresponding operator $\mathcal U^{\mathbb P}$ for projective spaces at the point $\hat p:= t^2\sigma_1+(r-1)\pi i\sigma_1\in H^2(\mathbb P,\mathbb C)$, i.e.
\[\mathcal U^{\mathbb P}_{\hat p}:=c_1(\mathbb P)*_{(-1)^{r-1}q}(-)\colon H^\bullet(\mathbb P,\mathbb C)\to H^\bullet(\mathbb P,\mathbb C).
\]
\item The spectrum $\operatorname{spec}(\mathcal U^{\mathbb G}_{p})$ is not simple if and only if the pair $(r,k)$ is such that
\[P_1(k)\leq r\leq k-P_1(k),
\]where $P_1(k)$ denotes the smallest prime divisor of $k$.
\item If $\pi_1,\dots,\pi_{n}$ denote the idempotents of the small quantum cohomology of the projective space $\mathbb P$ at the point $\hat p:= t^2\sigma_1+(r-1)\pi i\sigma_1\in H^2(\mathbb P,\mathbb C)$, then 
\begin{itemize}
\item the idempotents of the small quantum cohomology of $\mathbb G$ at $p=t^2\sigma_1\in H^2(\mathbb G;\colon\mathbb C)$ are 
\[(\overline{j}\circ\overline\vartheta)^{-1}\left(\varkappa_I\cdot \pi_I\right),\quad \pi_I:=\pi_{i_1}\wedge\dots\wedge \pi_{i_r},\]
with $1\leq i_1<\dots<i_r\leq k$, and where
\[ \varkappa_I:=\frac{1}{g^{\wedge\mathbb P}(\pi_I,\pi_I)}\det\begin{pmatrix}
g^{\mathbb P}(\pi_{i_1},\sigma_{r-1})&\dots&g^{\mathbb P}(\pi_{i_1},\sigma_{0})\\
\vdots&\ddots&\vdots\\
g^{\mathbb P}(\pi_{i_r},\sigma_{r-1})&\dots&g^{\mathbb P}(\pi_{i_r},\sigma_{0})
\end{pmatrix};
\]
\item the normalized idempotents are given by
\[(\overline{j}\circ\overline\vartheta)^{-1}\left(i^{\binom{r}{2}}\cdot f_I\right),\quad f_I:=\frac{\pi_I}{g^{\wedge\mathbb P}(\pi_I,\pi_I)^{\frac{1}{2}}}.
\]
\end{itemize}
\end{enumerate}
\end{cor}
\esh

\proof
Equation \eqref{satakeprod} is an immediate consequence of Theorem \ref{teoBCFK}; from this equality and from the value of the first Chern class $c_1(\mathbb G)=k\sigma_1$, one obtains point (1). For a proof of point (2), see \cite{cotti0}. The semisimplicity of the small quantum cohomology of the Grassmannian $\mathbb G$ is well known (see \cite{abrams}, \cite{cmp10}), so that existence of the idempotent vectors is guaranteed.
By Theorem \ref{teoBCFK} we deduce that the image of the idempotents $\alpha_1,\dots, \alpha_{\binom{k}{r}}$ of the small quantum cohomology of $\mathbb G$ under the map $\overline{j}\circ\overline\vartheta$ are scalar multiples of $\pi_I:=\pi_{i_1}\wedge\dots\wedge \pi_{i_r}$, with $1\leq i_1<\dots <i_r\leq k$, i.e. are of the form $\varkappa_I\cdot\pi_I$, for some constants $\varkappa_I\in\mathbb C^*$. Using Corollary \ref{corisometry}, from the equality $g^{\mathbb G}(\alpha_i,\alpha_i)=g^{\mathbb G}(\alpha_i,1)$, we find that necessarily 
\[\varkappa_I^2\cdot g^{\wedge\mathbb P}(\pi_I,\pi_I)=\varkappa_I\cdot g^{\wedge\mathbb P}(\pi_I,\sigma_{r-1}\wedge\dots\wedge\sigma_0),
\]and one concludes. By normalization, one obtains the expression for normalized idempotents.
\endproof

\subsection{Computation of the fundamental systems of solutions and monodromy data }\label{chconjgrass}
In what follows, if $V$ denotes a complex vector space and $\phi\in\End_\mathbb C(V)$, we denote by $\wedge^r\phi\in\End_\mathbb C(\bigwedge\nolimits ^r V)$ its $r$-exterior power: if a basis $(v_1,\dots, v_n)$ of $V$ is fixed, and  if $A$ denotes the matrix associated with $\phi$, then the matrix $\wedge^rA$ associated with $\wedge^r\phi$ is the one obtained by taking the $r\times r$ minors of $A$. The entries are disposed according to a pre-fixed ordering of the induced basis $(v_{i_1}\wedge\dots\wedge v_{i_r})_{1\leq i_1 <\dots < i_r\leq n}$ of $\bigwedge\nolimits^r V$.

\bsh
\begin{prop}\label{fundsol1} Let $Z^{\mathbb P}(z, t^2)$ be a solution of the system of differential equations \eqref{28.07.17-1}-\eqref{28.07.17-2}, i.e.
\begin{align}
\label{29.07.17-1}
\frac{\partial}{\partial t^2}Z^{\mathbb P}(z,t^2)&=z\mathcal C^{\mathbb P}_2(t^2)Z^{\mathbb P}(z,t^2),\quad\mathcal C^{\mathbb P}_2(t^2):=(\sigma)\circ^{\mathbb P}_{t^2\sigma},\\
\label{29.07.17-2}
\frac{\partial}{\partial z}Z^{\mathbb P}(z,t^2)&= \left(\mathcal U^\mathbb P(t^2)+\frac{1}{z}\mu^\mathbb P\right)Z^{\mathbb P}(z,t^2).
\end{align}
Then, the $r$-exterior power
\begin{equation}
\label{28.07.17-3}Z^{\mathbb G}(z, t^2):=\bigwedge\nolimits ^r\left( Z^\mathbb P(z,t^2+(r-1)\pi i)\right)
\end{equation}defines a solution of the system corresponding to the Grassmannian $\mathbb G$, namely
\begin{align}
\label{29.07.17-3}
\frac{\partial}{\partial t^2}Z^{\mathbb G}(z,t^2)&=z\mathcal C^{\mathbb G}_2(t^2)Z^{\mathbb G}(z,t^2),\quad\mathcal C^{\mathbb G}_2(t^2):=(\sigma_1)\circ^{\mathbb G}_{t^2\sigma_1},\\
\label{29.07.17-4}
\frac{\partial}{\partial z}Z^{\mathbb G}(z,t^2)&= \left(\mathcal U^\mathbb G(t^2)+\frac{1}{z}\mu^\mathbb G\right)Z^{\mathbb G}(z,t^2).
\end{align} Furthermore, if $Z^{\mathbb P}(z, t^2)$ is in Levelt form at $z=0$, then also \eqref{28.07.17-3} is in Levelt form at $z=0$.
\end{prop}
\esh

\proof Let us notice that
\[\left.\frac{\partial}{\partial t^2}(Z^{\mathbb G})^{A}_{B}\right|_{(z,t^2)}=\left.\sum_{a=1}^r\sum_{\ell=1}^k\det\begin{pmatrix}
(Z^{\mathbb P})^{\alpha_1}_{\beta_1}&\dots&(Z^{\mathbb P})^{\alpha_1}_{\beta_r}\\
\vdots&&\vdots\\
X^{\alpha_a}_\ell (Z^{\mathbb P})^{\ell}_{\beta_1}&\dots&X^{\alpha_a}_\ell (Z^{\mathbb P})^{\ell}_{\beta_r}\\
\vdots&&\vdots\\
(Z^{\mathbb P})^{\alpha_r}_{\beta_1}&\dots&(Z^{\mathbb P})^{\alpha_r}_{\beta_r}\\
\end{pmatrix}\right|_{(z,t^2+\pi i (r-1))},\quad X(z,t^2)=z\mathcal C_2^{\mathbb P}.
\]
By Corollary \ref{quantumprod}, the r.h.s. is easily seen to be equal to
\[\left(z\mathcal C_2^\mathbb G(t^2)Z^\mathbb G(z,t^2)\right)^A_B.
\]Analogously, we have that
\[\left.\frac{\partial}{\partial z}(Z^{\mathbb G})^{A}_{B}\right|_{(z,t^2)}=\left.\sum_{a=1}^r\sum_{\ell=1}^k\det\begin{pmatrix}
(Z^{\mathbb P})^{\alpha_1}_{\beta_1}&\dots&(Z^{\mathbb P})^{\alpha_1}_{\beta_r}\\
\vdots&&\vdots\\
W^{\alpha_a}_\ell (Z^{\mathbb P})^{\ell}_{\beta_1}&\dots&W^{\alpha_a}_\ell (Z^{\mathbb P})^{\ell}_{\beta_r}\\
\vdots&&\vdots\\
(Z^{\mathbb P})^{\alpha_r}_{\beta_1}&\dots&(Z^{\mathbb P})^{\alpha_r}_{\beta_r}\\
\end{pmatrix}\right|_{(z,t^2+\pi i (r-1))},
\]where we set $W(z,t^2)=\left(\mathcal U^\mathbb P(t^2)+\frac{1}{z}\mu^\mathbb P\right)$. Using Proposition \ref{propmu} and Corollary \ref{quantumprod}, one identifies the r.h.s. with 
\[\left[\left(\mathcal U^\mathbb G(t^2)+\frac{1}{z}\mu^\mathbb G\right)\cdot Z^\mathbb G(z, t^2)\right]^A_B.
\]For the last statement, notice that if
\[Z^\mathbb P(z,t^2)=\Phi(z,t^2)z^{\mu^\mathbb P}z^{c_1(\mathbb P)\cup(-)},\quad \Phi(-z,t^2)^T\eta^{\mathbb P}\Phi(z,t^2)=\eta^\mathbb P,
\]then using the generalized Cauchy-Binet identity for the minors of a product, and invoking Corollary \ref{corisometry}, Corollary \ref{corc1} and Proposition \ref{propmu}, one obtains that
\[Z^\mathbb G(z,t^2)=\widetilde\Phi(z,t)z^{\mu^\mathbb G}z^{c_1(\mathbb G)\cup (-)},\quad \widetilde\Phi(-z,t^2)^T\eta^\mathbb G\widetilde\Phi(z,t^2)=\eta^\mathbb G,\quad \widetilde\Phi(z,t^2)=\bigwedge\nolimits^r\Phi(z, t^2+\pi i(r-1)).
\]This concludes the proof.
\endproof

\bsh
\begin{cor} \label{corsoltop}
Let $Z^\mathbb P_{\rm top}(z,t^2)$ be the restriction to the small quantum locus of the topological-enumerative solution of $\mathbb P$. Then, the topological-enumerative solution of $\mathbb G$, restricted to the small quantum cohomology is given by
\[Z^\mathbb G_{\rm top}(z, t^2)=\left(\bigwedge\nolimits^rZ^\mathbb P_{\rm top}(z,t^2+\pi i(r-1))\right)\cdot e^{-\pi i(r-1)\sigma_1\cup (-)}.
\]
\end{cor}
\esh

\proof
According to Proposition \ref{18.07.17-2}, we have 
\[Z^\mathbb P_{\rm top}(z,t^2)=\Theta^\mathbb P_{\rm top}(z, t^2)z^{\mu^\mathbb P}z^{c_1(\mathbb P)\cup (-)},
\]and $\Theta_{\rm top}^\mathbb P$ is characterized  by the fact that
\[z^{-\mu^\mathbb P}\Theta_{\rm top}^\mathbb P(z, t^2)z^{\mu^\mathbb P}=\exp(t^2\sigma\cup(-))+ \sum_{i=1}^\infty A_iz^i,\quad A_i\in \frak{gl}(k,\mathbb C).
\]
Hence, from Proposition \ref{fundsol1}, we deduce that
\[\left(\bigwedge\nolimits^rZ^\mathbb P_{\rm top}(z,t^2+\pi i(r-1))\right)=H(z,t^2)z^{\mu^\mathbb G}z^{c_1(\mathbb G)\cup (-)},
\]where
\[z^{-\mu^\mathbb G}H(z,t^2)z^{\mu^\mathbb G}=\exp((t^2+\pi i (r-1))\sigma_1\cup(-))+ \sum_{i=1}^\infty A'_iz^i,\quad A'_i\in \frak{gl}\left(\binom{k}{r},\mathbb C\right),
\]by Proposition \ref{propmu} and Corollary \ref{corc1}. Using Proposition \ref{18.07.17-2}, we conclude.
\endproof

Let us consider a fixed choice $\Psi^\mathbb P(t^2)$ of the $\Psi$-matrix for $\mathbb P$ along points $t^2\sigma\in H^2(\mathbb P,\mathbb C)$ of the small quantum cohomology. By point (3) of Corollary \ref{quantumprod}, a choice of the $\Psi$-matrix for the Grassmannian $\mathbb G$ is given by the $r$-exterior power
\begin{equation}\label{29.07.17-5}\Psi^\mathbb G(t^2):=i^{\binom{r}{2}}\bigwedge\nolimits^r\Psi^\mathbb P(t^2+\pi i (r-1)),\quad t^2\sigma_1\in H^2(\mathbb G,\mathbb C).
\end{equation}If we set $Y^{\mathbb P/\mathbb G}:=\Psi^{\mathbb P/\mathbb G}\cdot Z^{\mathbb P/\mathbb G}$, we can consider the corresponding systems of differential equations \eqref{16.07.17-1}-\eqref{16.07.17-2}. The following results establish the relationship between the solutions of these differential systems, and their Stokes phenomena.

\bsh
\begin{prop}\label{propsolirr}Let $\ell$ be an oriented line in the complex plane, with slope $\phi\in[0;2\pi[$, admissible at both points
\[p:=t^2\sigma_1\in H^2(\mathbb G,\mathbb C),\quad \hat p:=(t^2+\pi i(r-1))\sigma\in H^2(\mathbb P,\mathbb C).
\]Let us denote by $Y^{\mathbb P/\mathbb G}_{\rm formal}(z,u)$ the formal solutions of the differential systems \eqref{16.07.17-1}-\eqref{16.07.17-2} associated with the quantum cohomology of $\mathbb P$ and $\mathbb G$, respectively. If $Y^{(k),\mathbb P/\mathbb G}_{\rm left/right}(z,u)$ denote the solutions of these systems, uniquely characterized by the asymptotic expansion
\[Y^{(k),\mathbb P/\mathbb G}_{\rm left/right}(z,u)\sim Y^{\mathbb P/\mathbb G}_{\rm formal}(z,u),\quad |z|\to \infty,\quad z\in e^{2\pi i k}\Pi_{\rm left/right}(\phi),
\]uniformly in $u$, then we have the following identifications:
\[Y^\mathbb G_{\rm formal}(z,u(p))=\bigwedge\nolimits^r Y^\mathbb P_{\rm formal}(z,u(\hat p)),
\]
\[Y^\mathbb G_{\rm left/right}(z,u(p))=\bigwedge\nolimits^r Y^\mathbb P_{\rm left/right}(z,u(\hat p)).
\]
\end{prop}
\esh

\proof The claim immediately follows from identity \eqref{29.07.17-5}, Proposition \ref{fundsol1}, and from the simple observation that
\[\bigwedge\nolimits^r\left[\exp(zU^\mathbb P(\hat p))\left(\mathbbm 1+ \sum_{h\geq 1}\frac{1}{z^h}A_h\right)\right]=\exp(zU^\mathbb G(p))\left(\mathbbm 1+  \sum_{h\geq 1}\frac{1}{z^h}A'_h\right).
\]
\endproof

\bsh
\begin{cor}\label{08.02.18-1}If $\ell$ is an oriented line in the complex plane, with slope $\phi\in[0;2\pi[$, admissible at both points
\[p:=t^2\sigma_1\in H^2(\mathbb G,\mathbb C),\quad \hat p:=(t^2+\pi i(r-1))\sigma\in H^2(\mathbb P,\mathbb C),\]
and if 
\[S^{\mathbb P}(\hat p;\ell),\ C^{(k),\mathbb P}(\hat p;\ell),\quad S^{\mathbb G}( p;\ell),\ C^{(k),\mathbb G}( p;\ell)
\] denote the Stokes and Central connection matrices of $\mathbb P$, and $\mathbb G$ respectively, computed at a point $\hat p$, and $p$ respectively, with respect to the oriented line $\ell$, then the following identities hold true:
\[S^{\mathbb G}( p;\ell)=\bigwedge\nolimits^rS^{\mathbb P}(\hat p;\ell),\quad C^{(k),\mathbb G}( p;\ell)=i^{-\binom{r}{2}}\left(\bigwedge\nolimits^rC^{(k),\mathbb P}(\hat p;\ell)\right)\cdot e^{\pi i(r-1)\sigma_1\cup (-)}.
\]
\end{cor}
\esh

\subsection{Reduction to (twisted) Kapranov Form} $\quad$

\bsh
\begin{prop}\label{mutatwedgemukai}
Let $(V,\langle\cdot,\cdot\rangle)$ be a Mukai lattice of rank $k$, and define a Mukai structure on the free $\mathbb Z$-module $\wedge^rV$ by setting
\[\langle\alpha_I,\alpha_J\rangle^{\wedge r}:=\det\left(\langle\alpha_{i_h},\alpha_{j_\ell}\rangle\right)_{1\leq h,\ell\leq r},
\]where $\alpha_I:=\alpha_{i_1}\wedge\dots\wedge\alpha_{i_r}$, with $1\leq i_1<\dots<i_r\leq k$, and analogously $\alpha_J$ denote two decomposable elements. If $(\varepsilon_i)_i$ and $(\tilde{\varepsilon}_i)_i$ are two exceptional bases of $V$ related by the action of a braid in $\mathcal B_{k}$, then the exceptional bases $(\varepsilon_I)_I$ and $(\tilde{\varepsilon}_I)_I$ of $\wedge^rV$, obtained by the lexicographical ordering, are in the same orbit with respect to the action of braids in $\mathcal B_{\binom{k}{r}}$ and $(\mathbb Z/2\mathbb Z)^{\times \binom{k}{r}}$.
\end{prop}
\esh

\proof
It clearly suffices to prove the statement for two exceptional bases of $V$ related by the action of an elementary braid. Let us assume, for example, that the exceptional bases of $V$
\[(\varepsilon_1,\dots,\varepsilon_i,\varepsilon_{i+1},\dots, \varepsilon_k),\quad (\varepsilon_1,\dots,\varepsilon_{i+1},\tilde\varepsilon_i,\varepsilon_{i+2},\dots,\varepsilon_k),\quad \tilde\varepsilon_i:=\mathbb R_{\varepsilon_{i+1}}\varepsilon_i,
\]
are related by the action of the elementary braid $\sigma_{i,i+1}$. Let us now consider the exceptional bases of $\wedge^r V$ obtained by the lexicographical ordering. The elements of the second basis can be classified into three different types:
\begin{enumerate}
\item those of the form $\varepsilon_J$ with $\varepsilon_{j_h}\notin\left\{\varepsilon_{i+1},\tilde\varepsilon_i\right\}$ for all $h=1,\dots, r$,
\item those of the form $(\dots\wedge\varepsilon_{i+1}\wedge\tilde\varepsilon_i\wedge\dots)$,
\item and those of the form
\begin{equation}\label{05.07.17-1}\left(\bigwedge_{a=1}^{\ell-1}\varepsilon_{j_a}\right)\wedge\tilde\varepsilon_i\wedge\left(\bigwedge_{a=\ell+1}^r\varepsilon_{j_a}\right),
\end{equation}for some $\ell$.
\end{enumerate}
Using the definition $\mathbb R_{\varepsilon_{i+1}}\varepsilon_i:=\varepsilon_{i}-\langle\varepsilon_i,\varepsilon_{i+1}\rangle\varepsilon_{i+1}$, it is evident that for the elements of the class (2) the following identity holds:
\[\dots\wedge\varepsilon_{i+1}\wedge\tilde\varepsilon_i\wedge\dots=\dots\wedge\varepsilon_{i+1}\wedge\varepsilon_i\wedge\dots=-(\dots\wedge\varepsilon_{i}\wedge\varepsilon_{i+1}\wedge\dots).
\]Consequently, they are the opposites of elements of the first exceptional basis. For the elements of the class (3), notice that all the elements between the \emph{first} one of the type \eqref{05.07.17-1}, and the corresponding one obtained by replacing $\tilde\varepsilon_i$ with $\varepsilon_{i+1}$, are of the form
\[\left(\bigwedge_{a=1}^{\ell-1}\varepsilon_{h_a}\right)\wedge\varepsilon_{i+1}\wedge\left(\bigwedge_{a=\ell+1}^r\varepsilon_{h_a}\right),
\]with $j_a=h_a$ for $a\in\left\{1,\dots,\ell-1\right\}\cup\left\{\ell+1,\dots, n\right\}$, and $j_{n+1}<h_{n+1}$, for some $n\in\left\{\ell+1,\dots,r\right\}$. The scalar product of these elements with the first element \eqref{05.07.17-1} is given by the determinant
\[\det\left(\begin{array}{c|c}
D_1&D_2\\
\hline
0& D_3
\end{array}\right)=0,
\]since the matrices $D_1,D_3$ are upper triangular, $\operatorname{diag}(D_1)=(1,\dots,1,\langle\varepsilon_{i+1},\tilde\varepsilon_i\rangle)$ and $D_3$ has at least one zero element on the diagonal (at least $(D_3)_{n+1,n+1}=0$). Hence, we can successively mutate the first element \eqref{05.07.17-1} on the left, till we obtain the following configuration of an exceptional basis:
\[\left(\dots,\underbrace{\left(\bigwedge_{a=1}^{\ell-1}\varepsilon_{j_a}\right)\wedge\varepsilon_{i+1}\wedge\left(\bigwedge_{a=\ell+1}^r\varepsilon_{j_a}\right)}_{A_{i+1}},\underbrace{\left(\bigwedge_{a=1}^{\ell-1}\varepsilon_{j_a}\right)\wedge\tilde\varepsilon_i\wedge\left(\bigwedge_{a=\ell+1}^r\varepsilon_{j_a}\right)}_{\tilde A_i},\dots\right).
\]At this point, notice that
\[\tilde A_i=A_{i+1}-\langle A_i,A_{i+1}\rangle^{\wedge r}A_{i+1}, \quad A_i:=\left(\bigwedge_{a=1}^{\ell-1}\varepsilon_{j_a}\right)\wedge\varepsilon_i\wedge\left(\bigwedge_{a=\ell+1}^r\varepsilon_{j_a}\right),
\]since $\langle A_i,A_{i+1}\rangle^{\wedge r}=\langle\varepsilon_i,\varepsilon_{i+1}\rangle$. The procedure continues and iterates with the new first term of the type \eqref{05.07.17-1}. At the end of the procedure, one obtaines a factor decomposition of the braids taking the second exceptional basis into the first one (modulo signs for elements of the class (2)). Notice that elements of the class (1) do not mutate.
\endproof

\begin{es}
An example will clarify the procedure. Let us consider the case $(r,k)=(3,6)$ and let $(\varepsilon_1,\dots,\varepsilon_6)$ be an exceptional basis of $V$. Through the action of the braid $\sigma_{23}$ we obtain a new exceptional collection
\[(\varepsilon_1,\varepsilon_3,\tilde\varepsilon_2,\varepsilon_4,\varepsilon_5,\varepsilon_6).
\]By the lexicographical ordering, from the first basis we obtain the exceptional basis 
\begin{equation}\label{05.07.17-2}
\varepsilon_{123},\ \varepsilon_{124},\ \varepsilon_{125},\ \varepsilon_{126},\ \varepsilon_{134},\ \varepsilon_{135},\ \varepsilon_{136},\ \varepsilon_{145},\ \varepsilon_{146},\ \varepsilon_{156},\ \varepsilon_{234},\ \varepsilon_{235},\ \varepsilon_{236},\ \varepsilon_{245},\end{equation}
\[ \varepsilon_{246},\ \varepsilon_{256},\ \varepsilon_{345},\ \varepsilon_{346},\ \varepsilon_{356},\ \varepsilon_{456}.
\]
Analogously, from the second basis we obtain the exceptional one
\begin{equation}\label{05.07.17-3}
{\color{red}\varepsilon_{13\tilde{2}}},\ \varepsilon_{134},\ \varepsilon_{135},\ \varepsilon_{136},\ {\color{blue}\varepsilon_{1\tilde{2}4}},\ {\color{blue}\varepsilon_{1\tilde{2}5}},\ {\color{blue}\varepsilon_{1\tilde{2}6}},\ \varepsilon_{145},\ \varepsilon_{146},\ \varepsilon_{156},\ {\color{red}\varepsilon_{3\tilde{2}4}},\ {\color{red}\varepsilon_{3\tilde{2}5}},\ {\color{red}\varepsilon_{3\tilde{2}6}},\ \varepsilon_{345},
\end{equation}
\[ \varepsilon_{346},\ \varepsilon_{356},\ {\color{blue}\varepsilon_{\tilde{2}45}},\ {\color{blue}\varepsilon_{\tilde{2}46}},\ {\color{blue}\varepsilon_{\tilde{2}56}},\ \varepsilon_{456}.
\]
We want to determine the transformation which transform \eqref{05.07.17-3} into \eqref{05.07.17-2}. In red we have colored elements of the class (2), in blue the elements of the class (3). Black elements are in class (1). Notice that red elements are just the opposite of the corresponding elements in \eqref{05.07.17-2} obtained by the exchange $(3\to 2, \tilde 2\to 3)$. Let us now start  with the first blue element, i.e. $\varepsilon_{1\tilde{2}4}$: we have that
\[\langle\varepsilon_{135},\varepsilon_{1\tilde{2}4}\rangle=0,\quad \langle\varepsilon_{136},\varepsilon_{1\tilde{2}4}\rangle=0.
\]Hence, by acting on \eqref{05.07.17-3} with the braid $\beta_{45}\beta_{34}$, we obtain
\[-\varepsilon_{123},\ \varepsilon_{134},\ {\color{blue}\varepsilon_{1\tilde{2}4}},\ \varepsilon_{135},\ \varepsilon_{136},\dots.
\]Acting now with the braid $\beta_{23}$, we obtain
\[-\varepsilon_{123},\ \varepsilon_{124},\ \varepsilon_{134},\ \varepsilon_{135},\ \varepsilon_{136},\dots.
\] We can continue with the next blue element, i.e. $\varepsilon_{1\tilde{2}5}$, till we obtain the sequence
\[-\varepsilon_{123},\ \varepsilon_{124},\ \varepsilon_{125},\ \varepsilon_{134},\ \varepsilon_{135},\ \varepsilon_{136},\ \varepsilon_{1\tilde{2}6},\dots.
\]By iterating the mutation procedure of the next blue elements, we arrive at the exceptional basis \eqref{05.07.17-2} (modulo signs of the red elements).
\end{es}

\bsh
\begin{lemma}[\cite{gamma1}]\label{31.07.17-1}
The following identity holds true:
\[(j\circ\vartheta)\left[\widehat{\Gamma}^\pm_\mathbb G\cup{\rm Ch}(\mathbb S^\mu\mathcal S^\vee)\right]=(2\pi i)^{-\binom{r}{2}}e^{-\pi i(r-1)\sigma_1}\bigwedge_{h=1}^r\widehat{\Gamma}^\pm_\mathbb P\cup{\rm Ch}(\mathcal O(\mu_h+r-h)).
\]
\end{lemma}
\esh

\proof
As in Section \ref{secclasscohgrass}, denote by $x_1,\dots, x_r$ the Chern roots of the bundle $\mathcal S^\vee$ on $\mathbb G$. Starting from the generalized Euler sequence 
\[0\to\mathcal S\to \mathcal O_{\mathbb G}^{\oplus k}\to\mathcal Q\to 0,
\]and applying to it the tensor product $\mathcal S^\vee\otimes -$, in the Grothendieck group $K_0(\mathbb G)$ we obtain the identity
\[[T\mathbb G]=[\mathcal S^\vee\otimes\mathcal Q]=k[\mathcal S^\vee]-[\mathcal S^\vee\otimes\mathcal S].
\]
Hence, by the multiplicative property of the $\widehat{\Gamma}^\pm$-classes, we obtain
\[\widehat{\Gamma}^\pm_{\mathbb G}=\prod_{i,h=1}^r\frac{\Gamma(1\pm x_i)^k}{\Gamma(1\pm x_i\mp x_h)}.
\]Notice that
\begin{align*}\prod_{i,h=1}^r\Gamma(1\pm x_i\mp x_h)&=\prod_{i< h}\Gamma(1\pm x_i\mp x_h)\Gamma(1\mp x_i\pm x_h)\\
&=\prod_{i<h}\frac{2\pi i (x_i-x_h)}{e^{\pi i (x_i-x_h)}-e^{\pi i (x_h-x_i)}}\\
&=(2\pi i)^{\binom{r}{2}}\prod_{i< h}(x_i-x_h)\prod_{i<h}\frac{e^{\pi i (x_i+x_h)}}{e^{2\pi i x_i}-e^{2\pi i x_h}}\\
&=(2\pi i)^{\binom{r}{2}}\left(\prod_{i< h}\frac{x_i-x_h}{e^{2\pi i x_i}-e^{2\pi i x_h}}\right)e^{(r-1)\pi i\sigma_1},
\end{align*}
where for the last equality we used the fact that $\left\{x_i+x_h\right\}_{i<h}$ are the Chern roots of $\bigwedge\nolimits^2\mathcal S^\vee$, so that
\begin{align*}\prod_{i<h}e^{\pi i (x_i+x_h)}&=\exp\left(\pi i \sum_{i<h}x_i+x_h\right)\\
&=\exp\left(\pi i c_1\left(\bigwedge\nolimits^2\mathcal S^\vee\right)\right)\\
&=\exp\left(\pi i (r-1)c_1(\mathcal S^\vee)\right).
\end{align*}
We have thus obtained the formula
\begin{equation}\label{29.07.17-6}\widehat{\Gamma}^\pm_\mathbb G=(2\pi i)^{-\binom{r}{2}}e^{-\pi i (r-1)\sigma_1}\prod_{i< h}\frac{e^{2\pi i x_i}-e^{2\pi i x_h}}{x_i-x_h}\prod_{i=1}^r\Gamma(1\pm x_i)^k.
\end{equation}
At this point, if we recall that the Chern character defines a morphism of rings, from the definition of Schur polynomials, we obtain the identity
\begin{equation}\label{29.07.17-7}{\rm Ch}(\mathbb S^\mu\mathcal S^\vee)=\frac{\det(e^{2\pi i x_i(\mu_h+r-h)})_{i,h}}{\prod_{i<h}e^{2\pi i x_i}-e^{2\pi i x_h}}.
\end{equation}
The claim follows from equations \eqref{29.07.17-6} and \eqref{29.07.17-7}.
\endproof

\bsh
\begin{teorema}\label{teoexcollgrass}
The central connection matrix, in the lexicographical order, of $QH^\bullet(\mathbb G)$, computed at $t=0$ with respect to an admissible oriented line $\ell$ is the matrix associated with the morphism \textnormal{\textcyr{D}}$^-_\mathbb G\colon K_0(\mathbb G)\otimes\mathbb C\to H^\bullet(\mathbb G,\mathbb C)$ with respect to an exceptional basis of the Grothendieck group $K_0(\mathbb G)$, related by suitable mutations and elements of $\mathbb Z^{\binom{k}{r}}$ to the twisted Kapranov basis
\[\left([\mathbb S^\mu\mathcal S^\vee\otimes\mathscr L]\right)_{\mu},\quad \mathscr L:=\det\left(\bigwedge\nolimits^2\mathcal S^\vee\right).
\]
In particular, the Conjecture \ref{congettura} holds true.
\end{teorema}
\esh

\proof
If $C$ is the matrix associated with the morphism \textcyr{D}$^-_{\mathbb P}$ with respect to
\begin{itemize}
\item the Beilinson basis $([\mathcal O],\dots,[\mathcal O(k-1)])$ of $K_0(\mathbb P)\otimes\mathbb C$,
\item the basis $(1,\sigma,\dots,\sigma^{k-1})$ of $H^\bullet(\mathbb P,\mathbb C)$,
\end{itemize}then by Corollary \ref{corc1} and Lemma \ref{31.07.17-1} it follows that the matrix
\begin{equation}\label{eqcenconngrass}i^{-\binom{r}{2}}\left(\bigwedge\nolimits^r C\right)e^{\pi i (r-1)\sigma_1\cup(-)}
\end{equation}is the matrix associated\footnote{Note that the numerical factor $i^{-\binom{r}{2}}$ in \eqref{eqcenconngrass} can be exactly identified with $i^{\bar d}$ where $d=r(k-r)$, since $r\equiv r^2\ ({\rm mod}\ 2)$.} with \textcyr{D}$^-_\mathbb G$ with respect to
\begin{itemize}
\item the twisted Kapranov basis $\left([\mathbb S^\mu\mathcal S^\vee\otimes\mathscr L]\right)_{\mu}$,
\item the induced Schubert basis $(\sigma_\mu)_{\mu}$.
\end{itemize}
The line bundle $\mathscr L$ is uniquely determined by its first Chern class $c_1(\mathscr L)=(r-1)\sigma_1$, by point (4) of Corollary \ref{corollifexce} (or even because $\mathbb G$ is Fano). Thus, by Corollary \ref{corcongsempl}, it follows that the association
\[\left(\bigwedge\nolimits^rK_0(\mathbb P),\wedge^r\chi^{\mathbb P}\right)\to(K_0(\mathbb G),\chi^{\mathbb G})\colon \bigwedge_{h=1}^r[\mathcal O(\mu_h+r-h)]\mapsto [\mathbb S^\mu\mathcal S^\vee\otimes\mathscr L],
\]defines an isomorphism of Mukai lattices. By Proposition \ref{mutatwedgemukai}, the claim follows.
\endproof

\subsection{Geometry of the Affine Grassmannian and classical/quantum Satake correspondence}
In the quantum setting, the abelian/non-abelian correspondence for Grassmannians admits a further interpretation, discussed by V. Golyshev and L. Manivel in \cite{golymaniv}, as the simplest manifestation of \emph{quantum corrections} to the geometric Satake correspondence due to A. Beilinson, V. Drinfeld, G. Lusztig, V. Ginzburg, I. Mirkovich and K. Vilonen (see \cite{beildrin}, \cite{lusztig}, \cite{ginz1, ginz2} and \cite{mirkvil1, mirkvil2}). 

\subsubsection{Langlands Duality} For any connected split reductive algebraic group $G$ over $\mathbb C$, it is well defined a companion group $G^\vee$, called \emph{Langlands dual} of $G$, whose geometry controls the representation theory of $G$. The definition of such a dual group $G^\vee$ is based on the bijection between isomorphisms classes of reductive groups and their \emph{root data}, a concept introduced by M. Demazure and slightly generalizing the one of \emph{root system} (we refer the reader to \cite{chrissginz}, \cite{springereductive, springerlinalggr}, \cite{ginz1,ginz2}, \cite{SGA3} for more details).
Let us introduce the following group-theoretical data associated with $G$: let $\mathbb T$ be  a maximal torus\footnote{One can also work with the \emph{abstract Cartan} $\mathbb T$ of $G$, defined as the quotient of a Borel subgroup $B\subseteq G$ by its unipotent radical (for different choices of $B$ such quotients are canonically isomorphic). In order to identify $\mathbb T$ with a maximal torus in $G$, an explicit embedding must be given.}
 in $G$, and let 
\[
X^*(\mathbb T):=\Hom(\mathbb T,\mathbb C^*),\quad X_*(\mathbb T):=\Hom(\mathbb C^*,\mathbb T)
\]be respectively the \emph{weight} and \emph{coweight lattices}, and let
\[\Phi\subseteq X^*(\mathbb T),\quad \Phi^\vee\subseteq X_*(\mathbb T)
\]be respectively the finite sets of \emph{roots} and \emph{coroots} (see e.g. \cite{springereductive} for the precise definition). 
The quadruple $(X^*(\mathbb T),X_*(\mathbb T),\Phi,\Phi^\vee)$ will be called the \emph{root datum of $G$}. The choice of a Borel subgroup $B$ containing $\mathbb T$ determines a set of \emph{positive roots} $\Phi^+\subseteq \Phi$ (and consequently of \emph{positive coroots, dominant weights and dominant coweights}). The Langlands dual group $G^\vee$ is defined (up to isomorphism) as the reductive group whose root datum is the dual one $(X_*(\mathbb T),X^*(\mathbb T),\Phi^\vee,\Phi)$. 

This definition of the Langlands dual group $G^\vee$ is completely based on the \emph{combinatorial content} of the root datum, and it relies on the theorem of classification of reductive groups. Following V. Ginzburg (\cite{ginz1,ginz2}), who based on an idea of V. Drinfeld and previous results of G. Lusztig, we can look for Â«an intrinsic new construction of $G^\vee$ which does not appeal to root systems, maximal tori, etc.Â», that we briefly summarize in the subsequent paragraphs.

\subsubsection{The affine Grassmannian $\mathscr Gr_G$} One of the main objects in the construction of Ginzburg is the affine Grassmannian $\mathscr Gr_G$ associated with the group $G$, a space which admits a natural structure of an ind-scheme, i.e. a direct limit of closed embeddings of  schemes of increasing dimension. Such an object can be considered as an algebraic analogue of the loop group defined in the topological setting (see e.g. \cite{segaloop}). Here we are going to describe three different ways for defining this object, and we refer the reader to \cite{ginz1, ginz2}, \cite{beildrin},  for  more detailed descriptions and proofs. See also \cite{beaulaszlo1}.  If $G$ is an algebraic group as above and $R$ is a $\mathbb C$-algebra, as usual we denote by $G(R)$ the group (over $\mathbb C$) of $R$-valued points of $G$.
\begin{enumerate}
\item The affine Grassmannian $\mathscr Gr_G$ is defined as the coset space
\begin{equation}\label{affgr1}\mathscr Gr_G:=G(\mathcal K)/G(\mathcal O), 
\end{equation}where for brevity we set $\mathcal K:=\mathbb C((z)),$ and $\mathcal O:=\mathbb C[\![z]\!]$. Notice that the set $G(\mathcal O)$ admits a structure of a group scheme, and $G(\mathcal K)$ an ind-scheme structure. From this definition, it is clear that a natural left action of $G(\mathcal O)$ is defined on $\mathscr Gr_G$.

\item A second definition of $\mathscr Gr_G$ is a \emph{polynomial} analogue of the previous one. If we define
\[LG:=G(\mathbb C[z^{-1},z]),\quad L^+G:=G(\mathbb C[z]),
\]then the affine Grassmannian can be defined as
\begin{equation}\label{affgr2}\mathscr Gr_G:=LG/L^+G. 
\end{equation}
The natural inclusion \[
 LG/L^+G\to G(\mathcal K)/G(\mathcal O)
\] is indeed not only injective, but actually an $LG$-equivariant isomorphism. From this definition it is clear we have a left action of $L^+G$ on $\mathscr Gr_G$. Any $G(\mathcal O)$-orbit of $\mathscr Gr_G$ is the image of a single $L^+G$-orbit in $LG/L^+G$ (see \cite{ginz2}, Proposition 1.2.4). 

\item Thirdly, the affine Grassmannian $\mathscr Gr_G$ can be defined in a topological setting as the group of based polynomials loops
\[\Omega G_c:=\left\{f\colon\mathbb S^1\to G_c,\ f\text{ polynomial },\ f(1)=1\right\},
\]where $G_c$ denotes the maximal compact subgroup in $G$. Given such a map $f$, it extends uniquely to a polynomial map $f\colon\mathbb C^*\to G$ such that $f(\bar z)=\overline{f(z)}$, where in the rhs the conjugation denotes the involutive automorphism of $G$ whose differential is the Cartan involution $\sigma\colon\frak g\to\frak g$ with $\frak g_c$ as $(+1)$-eigenspace. In this way we obtain an inclusion $\Omega G_c\hookrightarrow LG$, together with an Iwasawa decomposition
\[LG=\Omega\cdot L^+G,\quad \Omega G_c\cap L^+G=\left\{1\right\},
\]which allows to identify $\mathscr Gr_G$ with $\Omega G_c$ (see e.g. Chapter 8 of \cite{segaloop}, and also \cite{nadler}). In particular, such an identification induce a topological group structure on $\mathscr Gr_G$: we will denote by $m\colon\mathscr Gr_G\times\mathscr Gr_G\to\mathscr Gr_G$ the multiplication induced from $\Omega G_c$.
\end{enumerate}

The affine Grassmannian $\mathscr Gr_G$ admits a filtration 
\[\mathscr Gr_1\subseteq \mathscr Gr_2 \subseteq \mathscr Gr_3\subseteq \dots,\quad \mathscr Gr_G=\varinjlim_{i}\mathscr Gr_i,
\]where each $\mathscr Gr_i$ is a finite-dimensional projective variety and the inclusions are projective embeddings. Each variety $\mathscr Gr_i$ is $G(\mathcal O)$-stable, and the action of $G(\mathcal O)$ on $\mathscr Gr_i$ actually factors through a finite-dimensional algebraic group.

Since any coweight $\lambda\in X_*(\mathbb T)$ can be tautologically seen as a $\mathbb C[z^{-1},z]$-point of $G$,
 it determines a coset $ L^+G\cdot\lambda\subseteq LG$, and hence a point\footnote{In the third topological definition, $\lambda$ induces a map $\mathbb S^1\to\mathbb T\hookrightarrow K$, where $K$ is the maximal compact subgroup containing $\mathbb T$. Hence, $\lambda$ can be seen as a point of $\Omega$.} of $\mathscr Gr_G$. The $L^+G$-orbit of this point in $\mathscr Gr_G$ will be denoted by $O_\lambda$. Let us summarize some results about the nature of these orbits.
\begin{itemize}
\item All $L^+G$-orbits in $\mathscr Gr_G$ are of the form $O_\lambda$ for a coweight $\lambda\in X_*(\mathbb T)$.
\item Two orbits $O_\lambda, O_\mu$ are equal if and only if $\lambda,\mu$ are in the same orbit with respect to the action of $W$, the Weyl group of the pair $(G,\mathbb T)$. Hence, we can parametrize the $L^+G$-orbits in $\mathscr Gr_G$ with the \emph{dominant} coweights $\lambda\in X_*(\mathbb T)_+$. The dimension of $O_\lambda$ is given by
\[\dim O_\lambda=\langle2\rho,\lambda\rangle,\quad \rho:=\frac{1}{2}\sum_{\alpha\in\Phi_+}\alpha.
\]
\item The closure of $O_\lambda$ is equal to
\[\overline{O_\lambda}=\bigcup_{\mu\leq \lambda} O_\mu,
\]where $\leq$ denotes the partial ordering on the set of  dominant coweights defined as follows: $\mu\leq \lambda$ if and only if $\lambda-\mu$ is a sum of simple coroots with non-negative coefficients. Typically $\overline{ O_\lambda}$ is singular, and its smooth locus is $O_\lambda$. 

\item Notice that the minimal elements with respect to the partial order $\leq$ are the \emph{minuscule} coweights, i.e. the coweights $\mu\neq 0$ such that $\langle\alpha,\mu\rangle\leq 1$ for every positive root $\alpha\in \Phi^+$. Hence the orbit $O_\lambda$ is closed if and only if $\lambda$ is a minuscule weight. Consequently, $\overline{O_\lambda}$ is smooth if $\lambda$ is minuscule. In such a case, we have that
\[\overline{O_\lambda}=O_{\lambda}=G/P_\lambda,
\]where $P_\lambda$ is the parabolic subgroup of $G$ associated with $\lambda$. Notice that this class of minuscule varieties consists of the classical Grassmannians $\mathbb G(r,k)$, of orthogonal Grassmannians $OG(n,2n)$, even dimensional quadrics $Q^{2n}$, the Cayley plane $\mathbb {OP}^2=E_6/Q_1$ and the Freudenthal variety $E_7/Q_7$. A basic reference of this subject is \cite{lakmusses}.
\end{itemize}
Usually, the orbits $O_\lambda$ are called \emph{Schubert cells}, and their closure $\overline{O_\lambda}$ are called \emph{(spherical) Schubert varieties} of the affine Grassmannian $\mathscr Gr_G$.

\subsubsection{Classical Geometric Satake Correspondence} Given a complex algebraic variety $X$, recall that a (Whitney) stratification of $X$ is given by a filtration $$X=X_n\supseteq X_{n-1}\supseteq\dots\supseteq X_0,$$ where $X_i$ are closed subvarieties such that for each $j$ the locally closed subvariety $X_j\setminus X_{j-1}$ is either empty or non-singular and of complex dimension $j$. The connected components of the loci $X_j\setminus X_{j-1}$ are the \emph{strata} of the stratification and are required to satisfy two conditions, known as Whitney conditions. See \cite{whitney} for precise formulations of these conditions. We refer the reader also to \cite{intersecoh} for more properties and details on stratified spaces. We will denote by $\mathcal S$ the set of strata of a given Whitney stratification of $X$. 

\begin{defi}
Let $X$ be a complex algebraic variety equipped with a (Whitney) stratification $\mathcal S$, defined by the filtration $$X=X_n\supseteq X_{n-1}\supseteq\dots\supseteq X_0,\quad X_i\setminus X_{i-1}\in\mathcal S.$$ We will say that a sheaf (of $\mathbb C$-vector spaces) $\mathscr F$ is \emph{$\mathcal S$-constructible} if $\mathscr F|_{X_i\setminus X_{i-1}}$ are locally constant sheaves of finite rank.
 A complex of sheaves $\mathscr F^\bullet$ will be said to be \emph{cohomologically $\mathcal S$-constructible} if it is bounded and its cohomology sheaves $\mathcal H^i(\mathscr F^\bullet)$ are constructible. The \emph{derived $\mathcal S$-constructible category}, denoted by $\mathcal D^b_{\mathcal S}(X)$, is defined as the full subcategory of complexes in $\mathcal D^b(Sh_{\mathbb C}(X))$ which are cohomologically $\mathcal S$-constructible.  If $Y\subseteq X$ is a locally closed subset which is the union of strata in $\mathcal S$ then, by abuse of notation, we denote by $\mathcal D^b_{\mathcal S}(Y)$ the category $\mathcal D^b_{\mathcal T}(Y)$, where $\mathcal T:=\left\{A\in\mathcal S\colon A\subseteq Y\right\}$.
 \end{defi}

Let us denote by $\mathcal S$ the stratification induced on $\mathscr Gr_G$ by the $G(\mathcal O)$-orbits, and let us denote by $\mathcal D^b_\mathcal S(\mathscr Gr_i)$ the derived $\mathcal S$-constructible categories of $\mathbb C$-sheaves of the sets $\mathscr Gr_i$. The closed embedding $\mathscr Gr_i\hookrightarrow\mathscr Gr_j$ with $i\leq j$ induces an embedding of categories $\mathcal D^b_\mathcal S(\mathscr Gr_i)\hookrightarrow\mathcal D^b_\mathcal S(\mathscr Gr_j)$. We define the derived category of the affine Grassmannian $\mathscr Gr_G$ through a direct limit
\[\mathcal D^b_\mathcal S(\mathscr Gr_G):=\varinjlim\mathcal D^b_\mathcal S(\mathscr Gr_i).
\] 
Particular objects of $\mathcal D^b_\mathcal S(\mathscr Gr_G)$ are the \emph{intersection cohomology sheaves} of the Schubert varieties $\overline{O_\lambda}$ of $\mathscr Gr_G$, whose hypercohomology give the intersection cohomology groups of $\overline{O_\lambda}$: for any dominant coweight $\lambda\in X_*(\mathbb T)_+$ let us denote by $IC^\bullet(\overline{O_\lambda})$ the intersection cohomology sheaf of the Schubert variety $\overline{O_\lambda}$, extended to $0$ to the whole $\mathscr Gr_G$. We refer the reader to \cite{intersecoh} for complete definitions.

\begin{defi}
The category of \emph{perverse sheaves} on $\mathscr Gr_G$, denoted by $\mathcal P(\mathscr Gr_G)$, is defined as the full subcategory of $\mathcal D^b_\mathcal S(\mathscr Gr_G)$ generated by the objects isomorphic to finite sums of intersection cohomology sheaves $IC^\bullet(\overline{O_\lambda})$ of Schubert varieties. 
\end{defi}

Using the topological group structure $m\colon\mathscr Gr_G\times\mathscr Gr_G\to\mathscr Gr_G$, we can define on $\mathcal D^b_\mathcal S(\mathscr Gr_G)$ a \emph{convolution product} by setting
\[A\oast B:=m_*(A\boxtimes B), \quad A,B\in{\rm Obj}\left(\mathcal D^b_\mathcal S(\mathscr Gr_G)\right).
\]

\begin{teorema}[\cite{ginz2}]$\quad$
\begin{enumerate}
\item The category $\mathcal P(\mathscr Gr_G)$ is closed with respect to convolution product. In particular $(\mathcal P(\mathscr Gr_G),\oast)$ is a semisimple rigid tensor category. 
\item The cohomology functor $H^\bullet\colon\mathcal P(\mathscr Gr_G)\to {\rm Vect}_{\mathbb C}$ is exact and fully faithfull.
\end{enumerate}
\end{teorema}

Conditions (1) and (2) above allow us to apply one of the main results of P. Deligne and J.S. Milne (Theorem 2.11 of \cite{delignemilne}), which shows that the category $\mathcal P(\mathscr Gr_G)$ of perverse sheaves on $\mathscr Gr_G$ is a \emph{neutral Tannakian category}, i.e. isomorphic to the category of finite dimensional representations of an affine group scheme $G^*$ (see also \cite{saav}). This realizes the intrinsic characterization of the Langlands dual group $G^\vee$ mentioned at the beginning of this Section.

\begin{cor}[\cite{ginz2}, \cite{mirkvil1, mirkvil2}]\label{corsatequiv}There exists a reductive group $G^*$ whose category ${\rm Rep}_{\mathbb C}(G^*)$ of finite-dimensional $\mathbb C$-representations is equivalent (as tensor category) to the category of perverse sheaves $\mathcal P(\mathscr Gr_G)$:
$$(\mathcal P(\mathscr Gr_G),\oast)\to ({\rm Rep}_{\mathbb C}(G^*),\otimes).$$
Such an equivalence identifies $IC^\bullet(\overline{O_\lambda})$ with the irreducible representation $V_\lambda$ of $G^*$ with extreme weight $\lambda\in X_*(\mathbb T)_+$. The (hyper)-cohomology functor $H^\bullet\colon \mathcal P(\mathscr Gr_G)\to {\rm Vect}_{\mathbb C}$, under this equivalence, goes to the forgetful functor ${\rm Rep}_{\mathbb C}(G^*)\to{\rm Vect}_{\mathbb C}$. Furthermore, the group $G^*$ is isomorphic to the Langlands dual group $G^\vee$.
\end{cor}

In the case $G$ is a semisimple and simply-connected complex Lie group, assumption that will be valid in the rest of this Section, a further more explicit description of the singular cohomology of Schubert varieties is available. In what follows we denote by $U[\frak a]$ the universal enveloping algebra of a Lie algebra $\frak a$; if $V$ is an $\frak a$-module, for any subset $S\subseteq V$ we set
\[{\rm Ann}[\frak a; S]:=\left\{u\in U\frak a\colon\ u(s)=0,\text{ for all }s\in S\right\}.
\] 
\begin{teorema}\label{30.04.18-1}
Let $e\in \frak g^\vee$ be a principal nilpotent in the Lie algebra of the Langlands dual group of $G$, and let $(\frak g^\vee)^e$ be its centralizer. If $\lambda\in X_*(\mathbb T)$ is an anti-dominant coweight, let us denote by $V_\lambda$ an irreducible finite dimensional $\frak g^\vee$-representation with lowest weight $\lambda$, and let $v_\lambda\in V_\lambda$ be a lowest weight vector. Then we have the following isomorphism of graded algebras
\begin{equation}\label{30.04.18-2}
H^\bullet(\overline{O_\lambda},\mathbb C)\cong \frac{U[\frak (g^\vee)^e]}{{\rm Ann}[(g^\vee)^e; v_\lambda]}.
\end{equation}
\end{teorema}

\begin{oss}
The gradings on both the $\frak g^\vee$-module $V_\lambda$ of Corollary \ref{corsatequiv} and on the rhs of the isomorphism of Theorem \ref{30.04.18-1} are defined as follows. By Jacobson-Morozov Theorem, the given principal nilpotent element $e\in \frak g^\vee$ can be completed to a $\frak{sl}_2$-triple $(e,h,f)$. Hence, $h\in\frak g^\vee$ is a semisimple regular element that has integral eigenvalues in any finite dimensional $\frak g^\vee$-module. In Corollary \ref{corsatequiv}, the gradation on $V_\lambda$ is then defined by the eigenvalues of $h$, i.e.
\[V_\lambda=\bigoplus_{k\in\mathbb Z}V_{\lambda,k}(h),\quad V_{\lambda,k}(h):=\left\{v\in V\colon h\cdot v=k\cdot v\right\}.
\]In Theorem \ref{30.04.18-1}, notice that the centralizer $(\frak g^\vee)^e\subseteq \frak g^\vee$ is ${\rm ad}(h)$-stable. The eigenvalue-gradation of $(\frak g^\vee)^e$ is then induced on $U[(\frak g^\vee)^e]$: since the ideal ${\rm Ann}[(g^\vee)^e; v_\lambda]$ is graded, the quotient in the rhs of \eqref{30.04.18-2} is a finite dimensional graded algebra.
\end{oss}

\subsubsection{Quantum cohomology and the affine Grassmannian} Relations between the (small) quantum cohomology ring of homogeneous varieties $G/P$, with $G$ simply-connected and semisimple and $P$ a parabolic subgroup, and the (co)homology of the affine Grassmannian $\mathscr Gr_G$ have been established in literature from several points of view. 

A first description of these relationships can be found in \cite{kostant1, kostant2}. In these papers, basing on his previous work on the Toda lattice \cite{kostant-toda}, B. Kostant described the (small) quantum cohomology of complete flag manifolds associated with $G$ (i.e. $G/B$ with $B$ a Borel subgroup) as rings of rational functions on a unipotent algebraic group, whose Lie algebra is exactly the centralizer $(\frak g^\vee)^e$ discussed in the previous paragraph. Crucial, in \cite{kostant2}, is the role played by the contemporary (and mostly unpublished) work of D. Peterson. Firstly, in his theory of geometric realization of the (small) quantum cohomology of homogeneous spaces $G/P$, Peterson recognized a precise relationship with the affine Schubert calculus. Namely, Peterson claimed the possibility of identifying the small quantum ring $QH^\bullet(G/P)$ with a quotient, after localization, of the homology ring $H_*(\mathscr Gr_G)$ of the affine Grassmannian, giving also an explicit map of Schubert classes. In this way, all three-points genus $0$ Gromov--Witten invariants of $G/P$ can be identified with structural constant of the affine Schubert calculus. Such an identification was proved to hold by T. Lam and M. Shimozono in the equivariant setting:  see \cite{lamshimo} and \cite{lambook}. See also \cite{rietsch} for precise statements and proofs of Peterson description of quantum cohomology of partial flag manifolds $G/P$ (in type $A$) as coordinate rings of strata of a single variety (the so-called \emph{Peterson variety}).

In the paper \cite{golymaniv}, V. Golyshev and L. Manivel addressed the problem of a quantum counterpart of the classical geometric Satake correspondence described in the previous paragraph (Corollary \ref{corsatequiv}). Focusing on Dynkin types $A$ and $D$ (the types preserved by the Langlands duality), the authors showed that a result analogous to the Ginzburg's one (Theorem \ref{30.04.18-1}) holds true for quantum cohomology of minuscule Grassmannians, where the previous role of the principal nilpotent element $e$ is now played by a \emph{cyclic element} (defined by adding a quantum correction to $e$; see \emph{loc. cit.} for precise definitions and details). 

More precisely, by identifying the minuscule Grassmannians $G/P$ with a Schubert variety in $\mathscr Gr_G$, it is shown that the (small) quantum cohomology of $G/P$ admits a module structure over the Lie algebra of the Langlands dual group $G^\vee$ (called \emph{Satake structure}), interacting with a second module structure over the algebra of symmetric functions (called \emph{Schubert structure}). It is shown that the class of \emph{primitive elements} in $H^\bullet(G/P)$, i.e. the classes $x$ such that the $\cup$-product operator $x\cup(-)\colon H^\bullet(G/P)\to H^\bullet(G/P)$ can be expressed as the action of elements of $(\frak g^\vee)^e$, is preserved by the quantum corrections. This means that the small quantum product $x\circ(-)\colon QH^\bullet(G/P)\to QH^\bullet(G/P)$ can be described by the action of elements of a Cartan subalgebra of $\frak g^\vee$, defined as the centralizer of a cyclic element. Furthermore, very explicit formulae of the quantum corrections for multiplication by some \emph{special classes} are given in types $A,D$ and also for the Cayley plane and the Freudenthal variety, the exceptional minuscule spaces of type $E_6$ and $E_7$ respectively (see Theorem 1bis of \cite{golymaniv}). Equation \eqref{satakeprod} furnishes an example of the quantum corrections of Golyshev and Manivel for the primitive elements $p_\ell$ for the Grassmannians $\mathbb G(r,k)$ of type $A$. As V. Golyshev and L. Manivel noticed, an extension of the geometric Satake correspondence for more general Schubert varieties in $\mathscr Gr_G$ presents some foundational problems: namely, a definition of a good quantum analogue of the intersection cohomology for non-smooth varieties is missing.

Proposition \ref{fundsol1}, Corollary \ref{corsoltop}, Proposition \ref{propsolirr} and Corollary \ref{08.02.18-1} should be regarded as a manifestation of the quantum Satake correspondence at the level of solutions of the isomonodromic systems attached to the Frobenius structures of $QH^\bullet(\mathbb P)$ and $QH^\bullet(\mathbb G)$, and also of the corresponding monodromy invariants $(S,C)$. As underlined by Corollary \ref{08.02.18-1} and Theorem 4.8 and Lemma 4.10 of \cite{dubro2}, the monodromy invariants of semisimple Frobenius manifolds are particularly well-behaved with respect to their tensorial operations introduced in \cite{kauf} and \cite{KS}. A more systematic study of the properties of the monodromy invariants with respect to more general tensorial operations of Frobenius manifolds should be the object of future investigations. To the best knowledge of the authors, the study of manifestations of the geometric Satake correspondence at the level of derived categories with emphasis on their (full) exceptional collections  (and not just of their projections on the $K_0$-groups, as done e.g. in Proposition \ref{mutatwedgemukai})  is still missing from the literature.

\subsection{Reinterpretation of the results for $\mathbb G(2,4)$}

In this section we use the results obtained above in order to re-obtain the results of the computations of \cite{CDG}, developed in a different and more straightforward way, for the Grassmannian $\mathbb G(2,4)$. This allows us to understand the geometrical meaning of the numerical values of the entries of the central connection matrix $C$ as well as for the (at that time) ``mysterious'' matrix $A$ of Theorem 6.2 of \cite{CDG}.

According to Section \ref{secbeilinson}, the central connection matrix computed at the point $p=0$ of $QH^\bullet(\mathbb P^3_{\mathbb C})$, with respect to an admissible line $\ell$ of slope $0<\phi<\frac{\pi}{4}$ (and already put in the $\ell$-lexicographical order), and with respect to the topological solution of Proposition \ref{19.07.17-2} has the following columns
\[\left(C^{\mathbb P^3_{\mathbb C}}_{\rm lex}(p)\right)_1=
\left(
\begin{array}{c}
 -\frac{i}{2 \sqrt{2} \pi ^{3/2}} \\
 -\frac{i \sqrt{2} \gamma }{\pi ^{3/2}} \\
 -\frac{i \left(24 \sqrt{2} \gamma ^2+\sqrt{2} \pi ^2\right)}{12 \pi ^{3/2}} \\
 -\frac{i \left(\sqrt{2} \zeta (3)+8 \sqrt{2} \gamma ^3+\sqrt{2} \gamma  \pi ^2\right)}{3 \pi ^{3/2}} \\
\end{array}
\right),\]\[ \left(C^{\mathbb P^3_{\mathbb C}}_{\rm lex}(p)\right)_2=\left(
\begin{array}{c}
 -\frac{3 i}{2 \sqrt{2} \pi ^{3/2}} \\
 \frac{-6 i \gamma +\pi }{\sqrt{2} \pi ^{3/2}} \\
 \frac{-24 i \sqrt{2} \gamma ^2+8 \sqrt{2} \gamma  \pi -3 i \sqrt{2} \pi ^2}{4 \pi ^{3/2}} \\
 \frac{-48 i \sqrt{2} \gamma ^3+24 \sqrt{2} \gamma ^2 \pi -18 i \sqrt{2} \gamma  \pi ^2-\sqrt{2} \pi ^3-6 i \sqrt{2} \zeta (3)}{6 \pi ^{3/2}} \\
\end{array}
\right),
\]
\[\left(C^{\mathbb P^3_{\mathbb C}}_{\rm lex}(p)\right)_3=\left(
\begin{array}{c}
 \frac{i}{2 \sqrt{2} \pi ^{3/2}} \\
 \frac{i (2 \gamma +i \pi )}{\sqrt{2} \pi ^{3/2}} \\
 \frac{24 i \sqrt{2} \gamma ^2-24 \sqrt{2} \gamma  \pi -5 i \sqrt{2} \pi ^2}{12 \pi ^{3/2}} \\
 \frac{16 i \sqrt{2} \gamma ^3-24 \sqrt{2} \gamma ^2 \pi -10 i \sqrt{2} \gamma  \pi ^2+\sqrt{2} \pi ^3+2 i \sqrt{2} \zeta (3)}{6 \pi ^{3/2}} \\
\end{array}
\right),\]\[
\left(C^{\mathbb P^3_{\mathbb C}}_{\rm lex}(p)\right)_4=\left(
\begin{array}{c}
 \frac{i}{2 \sqrt{2} \pi ^{3/2}} \\
 \frac{i \sqrt{2} \gamma -\sqrt{2} \pi }{\pi ^{3/2}} \\
 \frac{24 i \sqrt{2} \gamma ^2-48 \sqrt{2} \gamma  \pi -23 i \sqrt{2} \pi ^2}{12 \pi ^{3/2}} \\
 \frac{8 i \sqrt{2} \gamma ^3-24 \sqrt{2} \gamma ^2 \pi -23 i \sqrt{2} \gamma  \pi ^2+7 \sqrt{2} \pi ^3+i \sqrt{2} \zeta (3)}{3 \pi ^{3/2}} \\
\end{array}
\right).
\]
The corresponding central connection matrix at the point $\hat{p}=\pi i\sigma$ of the small quantum locus of $\mathbb P^3_\mathbb C$ is obtained by the action of the braid 
\[\omega_{1,4}=\beta_{12}\beta_{34},
\]as described in Section \ref{secstoksqcohcpn}. Its columns are the following:

\[\left(C^{\mathbb P^3_{\mathbb C}}_{\rm lex}(\hat p)\right)_1=\left(
\begin{array}{c}
 \frac{i}{2 \sqrt{2} \pi ^{3/2}} \\
 \frac{2 i \gamma +\pi }{\sqrt{2} \pi ^{3/2}} \\
 \frac{24 i \sqrt{2} \gamma ^2+24 \sqrt{2} \gamma  \pi -5 i \sqrt{2} \pi ^2}{12 \pi ^{3/2}} \\
 \frac{16 i \sqrt{2} \gamma ^3+24 \sqrt{2} \gamma ^2 \pi -10 i \sqrt{2} \gamma  \pi ^2-\sqrt{2} \pi ^3+2 i \sqrt{2} \zeta (3)}{6 \pi ^{3/2}} \\
\end{array}
\right),\]\[
\left(C^{\mathbb P^3_{\mathbb C}}_{\rm lex}(\hat p)\right)_2=\left(
\begin{array}{c}
 -\frac{i}{2 \sqrt{2} \pi ^{3/2}} \\
 -\frac{i \sqrt{2} \gamma }{\pi ^{3/2}} \\
 -\frac{i \left(24 \sqrt{2} \gamma ^2+\sqrt{2} \pi ^2\right)}{12 \pi ^{3/2}} \\
 -\frac{i \left(\sqrt{2} \zeta (3)+8 \sqrt{2} \gamma ^3+\sqrt{2} \gamma  \pi ^2\right)}{3 \pi ^{3/2}} \\
\end{array}
\right),
\]
\[\left(C^{\mathbb P^3_{\mathbb C}}_{\rm lex}(\hat p)\right)_3=\left(
\begin{array}{c}
 -\frac{3 i}{2 \sqrt{2} \pi ^{3/2}} \\
 \frac{-3 i \sqrt{2} \gamma +\sqrt{2} \pi }{\pi ^{3/2}} \\
 \frac{-24 i \sqrt{2} \gamma ^2+16 \sqrt{2} \gamma  \pi -i \sqrt{2} \pi ^2}{4 \pi ^{3/2}} \\
 \frac{-24 i \sqrt{2} \gamma ^3+24 \sqrt{2} \gamma ^2 \pi -3 i \sqrt{2} \gamma  \pi ^2+5 \sqrt{2} \pi ^3-3 i \sqrt{2} \zeta (3)}{3 \pi ^{3/2}} \\
\end{array}
\right),
\]
\[\left(C^{\mathbb P^3_{\mathbb C}}_{\rm lex}(\hat p)\right)_4=\left(
\begin{array}{c}
 \frac{i}{2 \sqrt{2} \pi ^{3/2}} \\
 \frac{i (2 \gamma +i \pi )}{\sqrt{2} \pi ^{3/2}} \\
 \frac{24 i \sqrt{2} \gamma ^2-24 \sqrt{2} \gamma  \pi -5 i \sqrt{2} \pi ^2}{12 \pi ^{3/2}} \\
 \frac{16 i \sqrt{2} \gamma ^3-24 \sqrt{2} \gamma ^2 \pi -10 i \sqrt{2} \gamma  \pi ^2+\sqrt{2} \pi ^3+2 i \sqrt{2} \zeta (3)}{6 \pi ^{3/2}} \\
\end{array}
\right).
\]
According to Corollary \ref{08.02.18-1}, the central connection matrix for the Grassmannian $\mathbb G(2,4)$ at the point $t^2=0$ of $QH^\bullet(\mathbb G,\mathbb C)$ with respect to the same line $\ell$ (and already in $\ell$-lexicographical order) is given by
\begin{equation}\label{08.02.18-2}C^{\mathbb G(2,4)}_{\rm lex}(0)=-i\cdot C'\cdot \left(\bigwedge\nolimits^2C^{\mathbb P^3_{\mathbb C}}_{\rm lex}(\hat p)\right),
\end{equation}where 
\[C'=\left(
\begin{array}{cccccc}
 1 & 0 & 0 & 0 & 0 & 0 \\
 \pi  i & 1 & 0 & 0 & 0 & 0 \\
 -\frac{\pi ^2}{2} & \pi  i & 1 & 0 & 0 & 0 \\
 -\frac{\pi ^2}{2} & \pi  i & 0 & 1 & 0 & 0 \\
 -\frac{i \pi ^3}{3} & -\pi ^2 & \pi  i & \pi  i & 1 & 0 \\
 \frac{\pi ^4}{12} & -\frac{i \pi ^3}{3} & -\frac{\pi ^2}{2} & -\frac{\pi ^2}{2} & \pi  i & 1 \\
\end{array}
\right)\]
is the matrix representing the endomorphism
\[H^\bullet(\mathbb G(2,4),\mathbb C)\to H^\bullet(\mathbb G(2,4),\mathbb C) \colon v\mapsto e^{\pi i\sigma_1}\cup v.
\]We explicitly show the result of the multiplication \eqref{08.02.18-2} by columns:
\[\left(C^{\mathbb G(2,4)}_{\rm lex}(0)\right)_1=\left(
\begin{array}{c}
 \frac{1}{4 \pi ^2} \\
 \frac{\gamma }{\pi ^2} \\
 \frac{48 \gamma ^2+\pi ^2}{24 \pi ^2} \\
 \frac{48 \gamma ^2+\pi ^2}{24 \pi ^2} \\
 \frac{-\zeta (3)+16 \gamma ^3+\gamma  \pi ^2}{3 \pi ^2} \\
 \frac{-192 \gamma  \zeta (3)+768 \gamma ^4-\pi ^4+96 \gamma ^2 \pi ^2}{144 \pi ^2} \\
\end{array}
\right),\]

\[ \left(C^{\mathbb G(2,4)}_{\rm lex}(0)\right)_2=\left(
\begin{array}{c}
 \frac{5}{4 \pi ^2} \\
 \frac{10 \gamma +i \pi }{2 \pi ^2} \\
 \frac{240 \gamma ^2+48 i \gamma  \pi +17 \pi ^2}{24 \pi ^2} \\
 \frac{240 \gamma ^2+48 i \gamma  \pi +17 \pi ^2}{24 \pi ^2} \\
 \frac{160 \gamma ^3+48 i \gamma ^2 \pi +34 \gamma  \pi ^2-3 i \pi ^3-10 \zeta (3)}{6 \pi ^2} \\
 \frac{3840 \gamma ^4+1536 i \gamma ^3 \pi +1632 \gamma ^2 \pi ^2-288 i \gamma  \pi ^3-29 \pi ^4-960 \gamma  \zeta (3)-96 i \pi  \zeta (3)}{144 \pi ^2} \\
\end{array}
\right),
\]

\[\left(C^{\mathbb G(2,4)}_{\rm lex}(0)\right)_3=\left(
\begin{array}{c}
 -\frac{1}{2 \pi ^2} \\
 \frac{-4 \gamma -i \pi }{2 \pi ^2} \\
 \frac{-48 \gamma ^2-24 i \gamma  \pi +5 \pi ^2}{12 \pi ^2} \\
 \frac{-48 \gamma ^2-24 i \gamma  \pi -7 \pi ^2}{12 \pi ^2} \\
 \frac{-64 \gamma ^3-48 i \gamma ^2 \pi -4 \gamma  \pi ^2-3 i \pi ^3+4 \zeta (3)}{6 \pi ^2} \\
 \frac{-768 \gamma ^4-768 i \gamma ^3 \pi -96 \gamma ^2 \pi ^2-144 i \gamma  \pi ^3+\pi ^4+192 \gamma  \zeta (3)+48 i \pi  \zeta (3)}{72 \pi ^2} \\
\end{array}
\right),
\]

\[\left(C^{\mathbb G(2,4)}_{\rm lex}(0)\right)_4=\left(
\begin{array}{c}
 -\frac{1}{2 \pi ^2} \\
 \frac{-4 \gamma -i \pi }{2 \pi ^2} \\
 \frac{-48 \gamma ^2-24 i \gamma  \pi -7 \pi ^2}{12 \pi ^2} \\
 \frac{-48 \gamma ^2-24 i \gamma  \pi +5 \pi ^2}{12 \pi ^2} \\
 \frac{-64 \gamma ^3-48 i \gamma ^2 \pi -4 \gamma  \pi ^2-3 i \pi ^3+4 \zeta (3)}{6 \pi ^2} \\
 \frac{-768 \gamma ^4-768 i \gamma ^3 \pi -96 \gamma ^2 \pi ^2-144 i \gamma  \pi ^3+\pi ^4+192 \gamma  \zeta (3)+48 i \pi  \zeta (3)}{72 \pi ^2} \\
\end{array}
\right),
\]
\[\left(C^{\mathbb G(2,4)}_{\rm lex}(0)\right)_5=\left(
\begin{array}{c}
 \frac{1}{4 \pi ^2} \\
 \frac{2 \gamma +i \pi }{2 \pi ^2} \\
 \frac{48 \gamma ^2+48 i \gamma  \pi -11 \pi ^2}{24 \pi ^2} \\
 \frac{48 \gamma ^2+48 i \gamma  \pi -11 \pi ^2}{24 \pi ^2} \\
 \frac{32 \gamma ^3+48 i \gamma ^2 \pi -22 \gamma  \pi ^2-3 i \pi ^3-2 \zeta (3)}{6 \pi ^2} \\
 \frac{768 \gamma ^4+1536 i \gamma ^3 \pi -1056 \gamma ^2 \pi ^2-288 i \gamma  \pi ^3+23 \pi ^4-192 \gamma  \zeta (3)-96 i \pi  \zeta (3)}{144 \pi ^2} \\
\end{array}
\right),
\]

\[\left(C^{\mathbb G(2,4)}_{\rm lex}(0)\right)_6=\left(
\begin{array}{c}
 \frac{1}{4 \pi ^2} \\
 \frac{\gamma +i \pi }{\pi ^2} \\
 \frac{48 \gamma ^2+96 i \gamma  \pi -47 \pi ^2}{24 \pi ^2} \\
 \frac{48 \gamma ^2+96 i \gamma  \pi -47 \pi ^2}{24 \pi ^2} \\
 \frac{16 \gamma ^3+48 i \gamma ^2 \pi -47 \gamma  \pi ^2-15 i \pi ^3-\zeta (3)}{3 \pi ^2} \\
 \frac{768 \gamma ^4+3072 i \gamma ^3 \pi -4512 \gamma ^2 \pi ^2-2880 i \gamma  \pi ^3+671 \pi ^4-192 \gamma  \zeta (3)-192 i \pi  \zeta (3)}{144 \pi ^2} \\
\end{array}
\right).
\]
The reader can recognize, up to irrelevant signs of three columns (due to different choices of branches of the matrix $\Psi$), a \emph{perfect matching} with the entries of the central connection matrix exhibited in Appendix A of \cite{CDG}.  

In Section 6 of \cite{CDG} it was underlined the difference between the computed central connection matrix $C^{\mathbb G(2,4)}_{\rm lex}(0)$ and the one originally predicted in \cite{dubro4}: it was shown that the central connection matrix $C^{\mathbb G(2,4)}_{\rm lex}(0)$ above has the form
\[C^{\mathbb G(2,4)}_{\rm lex}(0)=A^{-1}X,\]where the matrix $X$ is such that \[\sum_\lambda X^\lambda_\ell\sigma_\lambda=\frac{1}{4\pi^2}\widehat{\Gamma}^-_{\mathbb G(2,4)}\cup {\rm Ch}(E_\ell),
\]for an explicit mutation $(E_\ell)_\ell$ of the Kapranov exceptional collection. By explicit computation, it was found that the difference $A$ is given by the matrix associated with the $\cup$-multiplication by $\exp(2\pi i\sigma_1)$:
\[A\colon H^\bullet(\mathbb G(2,4),\mathbb C)\to H^\bullet(\mathbb G(2,4),\mathbb C) \colon v\mapsto e^{2\pi i\sigma_1}\cup v.
\]
 At this point, we are also able to understand what is the geometrical meaning of the entries of the matrix $A$ (or, better to say, of $A^{-1}$) appearing in Theorem 6.2 of \cite{CDG}: indeed, since
 \[{\rm Ch}(\det\mathcal S^\vee)=\exp(2\pi i\sigma_1),\quad c_1(\mathbb G(2,4))=4\sigma_1,
 \]we can identity the operator $A^{-1}$ with
 \[H^\bullet(\mathbb G(2,4),\mathbb C)\to H^\bullet(\mathbb G(2,4),\mathbb C) \colon v\mapsto e^{-\pi ic_1(\mathbb G(2,4))}\cup {\rm Ch}(\det\mathcal S^\vee)\cup v.
 \]In other words, the contribution of the characteristic classes ${\rm Ch}(\det\mathcal S^\vee)$ and $e^{-\pi ic_1(\mathbb G(2,4))}$, prescribed by the Conjecture \ref{congettura}, are hidden in the entries of the matrix $A$. 
 
 Finally, as a further verification of our computations, let us show how to re-obtain the braids, which put the monodromy data $C^{\mathbb G(2,4)}_{\rm lex}(0), S^{\mathbb G(2,4)}_{\rm lex}(0)$ into the (twisted) Kapranov form, by applying the results of Proposition \ref{mutatwedgemukai}. For doing this, let us notice that the braid which transforms the data associated with the Beilison collection $\frak B:=(\mathcal O,\mathcal O(1),\mathcal O(2),\mathcal O(3))$ on $\mathbb P^3_\mathbb C$ into the monodromy data $C^{\mathbb P^3_\mathbb C}_{\rm lex}(\hat p), S^{\mathbb P^3_\mathbb C}_{\rm lex}(\hat p)$ is
 \[\underbrace{(\beta^{-1}_{12}\beta^{-1}_{23}\beta^{-1}_{34})\beta^{-1}_{12}}_{{\rm Lemma\ }\ref{bellatreccia}}\underbrace{\beta_{12}\beta_{34}}_{\omega_{1,4}}=\beta^{-1}_{12}\beta^{-1}_{23}.
 \]
 Let us now consider a Mukai lattice $V$ of rank 4 and let $(\varepsilon_1,\varepsilon_2,\varepsilon_3,\varepsilon_4)$ be an exceptional basis. The braid $\beta^{-1}_{12}$ transforms this basis into 
\[(\varepsilon_2,\varepsilon_{\tilde 1},\varepsilon_3,\varepsilon_4).
\]If we consider on the Mukai lattice $\bigwedge\nolimits^2 V$ the exceptional basis 
\[(\varepsilon_{12},\varepsilon_{13},\varepsilon_{14},\varepsilon_{23},\varepsilon_{24},\varepsilon_{34}),\quad\text{and}\quad (\varepsilon_{2\tilde 1},\varepsilon_{23},\varepsilon_{24},\varepsilon_{\tilde 13},\varepsilon_{\tilde 14},\varepsilon_{34}),
\]
by applying the argument of Proposition \ref{mutatwedgemukai} we immediately see that the second collection can be transformed into the first one (up to a sign) by the braid 
\begin{equation}\label{12.02.18-1}\beta_{34}\beta_{23}\beta_{45}.
\end{equation}
Analogously, starting from the two exceptional bases $(\varepsilon_1,\varepsilon_2,\varepsilon_3,\varepsilon_4)$ and $(\varepsilon_1,\varepsilon_3,\varepsilon_{\tilde 2},\varepsilon_4)$, we obtain the two exceptional bases 
\[(\varepsilon_{12},\varepsilon_{13},\varepsilon_{14},\varepsilon_{23},\varepsilon_{24},\varepsilon_{34}),\quad\text{and}\quad (\varepsilon_{13},\varepsilon_{1\tilde 2},\varepsilon_{14},\varepsilon_{3\tilde 2},\varepsilon_{34},\varepsilon_{\tilde 24}).
\]The second basis is transformed into the first one (up to a sign) by the braid
\begin{equation}\label{12.02.18-2}\beta_{12}\beta_{56}.
\end{equation}
By taking the product of \eqref{12.02.18-2} and \eqref{12.02.18-1}, we obtain the braid $\beta_{12}\beta_{56}\beta_{34}\beta_{23}\beta_{45}$.
 
This braid differs from the one of the paper \cite{CDG} just by an irrelevant factor $\beta_{34}$, which coincides with a mere permutation of the central objects of the 5-block. This shows the complete agreement between   the previous results of \cite{CDG} and  those presented in this paper.

\subsection[Symmetries and Quasi-Periodicity of the Stokes matrices]{Symmetries and Quasi-Periodicity of the Stokes matrices along the small quantum locus}\label{symgrass} We conclude this Section with the following result, concerning the symmetries and quasi-periodicity properties of the Stokes matrix $S$ of $QH^\bullet(\mathbb G)$ computed at points of the small quantum cohomology. It is an immediate consequence of the analogous properties of the Stokes matrix for $QH^\bullet (\mathbb P)$ and of Corollary \ref{08.02.18-1}.

\bsh
\begin{teorema}
\label{quasiperiodgrass}
The Stokes matrix $S_{\mathbb G(r,k)}(p,\phi)$, computed at a point $p\ni H^2(\mathbb G,\mathbb C)$ with respect to an admissible line $\ell$ of slope $\phi\in\mathbb R$ and in the $\ell$-lexicographical order, satisfies the following conditions:
\begin{enumerate}
\item it has the following functional form
\[S_{\mathbb G(r,k)}(t\sigma_1,\phi)=S({\rm Im }t+k\phi);
\]
\item it is \emph{quasi-periodic} along the small quantum locus, in the sense that
\[S_{\mathbb G(r,k)}(p,\phi)\sim S_{\mathbb G(r,k)}\left(p,\phi+\frac{2\pi i}{k}\right),
\]where $A\sim B$ means that the matrices $A$ and $B$ are in the same orbit under the action of $(\mathbb Z/2\mathbb Z)^{\binom{k}{r}}$. Moreover, we have that
\[S_{\mathbb G(r,k)}(p,\phi)= S_{\mathbb G(r,k)}\left(p,\phi+{2\pi i}\right);
\]
\item the upper-diagonal entries 
\[S_{\mathbb G(r,k)}(p,\phi)_{j,j+1},\quad S_{\mathbb G(r,k)}\left(p,\phi+\frac{\pi i}{k}\right)_{j,j+1}
\]differ for some signs, and we have that
\[|S_{\mathbb G(r,k)}(p,\phi)_{j,j+1}|\in\left\{\binom{k}{1},\dots,\binom{k}{k-1}\right\}\cup\left\{0\right\}.
\]
\end{enumerate}
\end{teorema}
\esh

From this Theorem, Corollary \ref{corbeilip12}, Proposition \ref{18.07.17-1} and from Lemma \ref{31.07.17-1}, we finally deduce the following result.

\bsh
\begin{cor}\label{corkapranovgrass}
The Kapranov exceptional collection $(\mathbb S^\lambda\mathcal S^\vee)_\lambda$, twisted by a suitable line bundle, is associated with the monodromy data of $\mathbb G(r,k)$ at points of the small quantum locus if and only if $(r,k)=(1,2),(1,3),(2,3)$. In this cases, the line bundle is trivial, and the Kapranov collection coincides with the Beilinson one\footnotemark.
\end{cor}
\esh
\footnotetext[35]{Notice that $\mathbb G(2,3)\cong\mathbb P((\mathbb C^3)^\vee)\cong\mathbb P^2_\mathbb C$ by duality. }

\newpage
\appendix
\section{Basic notions on Pure Motives}\label{chowmot}
 In this Appendix we briefly recall basic notions and properties of Chow motives, referring the interested reader to \cite{scholl}, \cite{andre}, \cite{murre} for more detailed and complete introductions to this vast and fascinating topic. 
Let $\mathbb {F,K}$ be two fields. Let us denote by $\mathcal V_{\mathbb K}$ the category of smooth projective varieties\footnote{Here by \emph{variety over $\mathbb K$} we mean a reduced $\mathbb K$-scheme. In particular, we do not assume it to be irreducible.} over $\mathbb K$. For any $X\in {\rm Ob}(\mathcal V_{\mathbb K})$ and $d\in\mathbb Z$, let us denote 
\begin{enumerate}
\item by $Z^d(X)$ the \emph{group of $d$-cycles} of $X$, i.e. the free abelian group generated by irreducible subvarieties of $X$ of codimension $d$;
\item by ${ CH}^d(X):=Z^d(X)/\sim_{rat}$ the $d$-codimensional Chow group of $X$ (where $\sim_{\rm rat}$ denotes the \emph{rational equivalence of cycles}),
\item by ${CH}^d(X)_{\mathbb F}:={CH}^d(X)\otimes_{\mathbb Z}\mathbb F$ and ${CH}^\bullet (X)_\mathbb F:=\bigoplus_d {CH}^d(X)_{\mathbb F}$.
\end{enumerate}
If $X,Y\in {\rm Ob}(\mathcal V_{\mathbb K})$ with $X$ irreducible and of pure dimension $d$, and $r\in\mathbb Z$, we define the \emph{group of correspondences of degree $r$ from $X$ to $Y$} as
\[{\rm Corr}^{r}(X,Y):={CH}^{r+d}(X\times Y)_{\mathbb F}.
\]For $X:=\coprod X_i$, with $X_i$ irreducible, we set
\[{\rm Corr}^{r}(X,Y):=\bigoplus {\rm Corr}^{r}(X_i,Y).
\]Given a third object $Z\in{\rm Ob}(\mathcal V_{\mathbb K})$, and $s\in \mathbb Z$ we can define the composition of correspondences as follows
\[\xymatrix@R-2pc{{\rm Corr}^r(X,Y)\otimes {\rm Corr}^s(Y,Z)\ar[r]& {\rm Corr}^{r+s}(X,Z)\\
f\otimes g\ar@{|->}[r]&(\pi_{XZ})_*\left(\pi_{XY}^*f\cdot \pi_{YZ}^*g\right)
}\] where $\pi_{XY/YZ/XZ}$ denote the projections from the triple product $X\times Y\times Z$ to the product of two spaces, and where the intersection product is performed in ${CH}^\bullet(X\times Y\times Z)_{\mathbb F}$.
At this point we can define the category ${\rm CHM}(\mathbb K)_{\mathbb F}$ of the rational Chow motives over $\mathbb K$ as the category 
\begin{itemize}
\item whose objects are triples $(X,p,m)$ where $X\in {\rm Ob}(\mathcal V_{\mathbb K})$, $m\in\mathbb Z$ and $p\in {\rm Corr}^0(X,X)$ is a \emph{projector}, i.e. an idempotent wrt the composition product of correspondences;
\item morphisms from $(X,p,m)$ to $(Y,q,n)$ are defined as elements of the set
\[q\circ {\rm Corr}^{n-m}(X,Y)\circ p,
\]and composition of morphisms comes from composition of correspondences.
\end{itemize}
There is a naturally defined contravariant functor $\frak h(\cdot)_{\mathbb F}\colon \mathcal V_{\mathbb K}^{\rm op}\to {\rm CHM}(\mathbb K)_{\mathbb F}$, which on the objects is defined as
\[\frak h(X)_{\mathbb F}:=(X,{\rm id}_X, 0),
\]and which associates to a morphism $f\colon X\to Y$ the correspondence $\frak h(f)_{\mathbb F}:=[\Gamma_f^T]\in {\rm Corr}^0(Y,X)$ given by (the rational equivalence class of) the \emph{transpose }of its graph $\Gamma_f\subseteq X\times Y$. 

Three natural operations are defined on ${\rm CHM}(\mathbb K)_{\mathbb F}$:
\begin{itemize}
\item given two objects $(X,p,m),(Y,q,m)$, we define their \emph{direct sum} as 
\[(X,p,m)\oplus (Y,q,m):=\left(X\coprod Y, p\amalg q,m\right),
\]where $\amalg$ denotes the disjoint union. For the definition of the general case $(X,p,m),(Y,q,n)$ with $m\neq n$ see \cite{scholl}.
\item given two objects $(X,p,m),(Y,q,m)$, we define their \emph{tensor product} as 
\[(X,p,m)\otimes(Y,q,n):=(X\times Y, p\times q,m+n).
\]The motive of a point $\frak h({\rm pt })_{\mathbb F}:=({\rm Spec}(\mathbb K), {\rm id}, 0)$ coincides with the \emph{unit} motive, denoted by $\mathbbm 1$.
\item Given an object $(X,p,m)$ we define its \emph{dual} as 
\[(X,p,m)^\vee:=(X,p^T,d-m),
\]in the case $X$ is of pure dimension $d$, and which acts on morphisms as the transposition of correspondences. 
\end{itemize}
These operations make ${\rm CHM}(\mathbb K)_{\mathbb F}$ a $\mathbb F$-linear, pseudo-abelian\footnote{An additive category $\mathcal C$ is called \emph{pseudo-abelian} (or \emph{Karoubian}) if for all all objects $X\in {\rm Ob}(\mathcal C)$ all projectors $p\in Hom_{\mathcal C}(X,X)$, $p\circ p=p$, have a kernel.}, rigid tensor category (see \cite{andre}).

The motive $\mathbb L:=({\rm Spec}(\mathbb K),{\rm id},-1)$ is called \emph{Lefschetz motive}. From the above definitions, it immediately follows that for all smooth projective varieties $X\in{\rm Ob}\left(\mathcal V_{\mathbb K}\right)$ the following canonical isomorphism holds true
\begin{equation}\label{mot1}CH^r(X)_{\mathbb F}\cong {\rm Hom_{CHM(\mathbb K)_{\mathbb F}}}(\mathbb L^{\otimes r},\frak h(X)_{\mathbb F}).
\end{equation}

\subsubsection{Universality of Chow Motives} Let $(\mathcal T, {\bf L}, H, (tr_X)_X, (c^r_X)_{r,X})$ be the datum of
\begin{enumerate}
\item a $\mathbb F$-linear, pseudo-abelian, rigid tensor category $\mathcal T$;
\item a $\otimes$-invertible object ${\bf L}\in{\rm Ob}(\mathcal T)$;
\item a monoidal functor $H\colon\mathcal V_{\mathbb K}^{op}\to\mathcal T$ such that the diagram
\[\mathbb P^1\to{\rm Spec}(\mathbb K)=\left\{\infty\right\}\hookrightarrow\mathbb P^1
\]induces the decomposition
\[H(\mathbb P^1)=\mathbbm 1\oplus{\bf L};
\]
\item for all $X\in{\rm Ob}(\mathcal V_{\mathbb K})$ of pure dimension $d$, a morphism 
\[tr_X\colon H(X)\to{\bf L}^{\otimes d},
\]such that
\begin{enumerate}
\item $tr_{X\times Y}=tr_X\otimes tr_Y$,
\item it identifies the dual object $H(X)^\vee$ with $H(X)\otimes {\bf L}^{\otimes(-d)}$, in the sense that the evaluation\footnote{Since any object of a \emph{rigid} tensor category is reflexive, i.e. $X\cong X^{\vee\vee}$, the coevaluation map is simply given by the transpose of the evaluation map.} is given by
\[\xymatrix{H(X)\otimes H(X)\otimes{\bf L}^{\otimes(-d)}\ar[rrr]^{\ (tr_X\circ id)\circ(\mu_X\circ id)}&&&{\bf L}^{\otimes d}\otimes{\bf L}^{ \otimes(-d)}\cong\mathbbm 1,\\
}
\]where the morphism $\mu_X$ is induced by the diagonal map $\delta_X\colon X\to X\times X$;
\end{enumerate}
\item for all $X\in{\rm Ob}(\mathcal V_{\mathbb K})$, a family of $\mathbb F$-linear morphisms
\[c^r_X\colon CH^r(X)_{\mathbb F}\to {\rm Hom}_{\mathcal T}(\mathbbm 1, H(X)\otimes {\bf L}^{\otimes (-r)}),
\]such that
\begin{enumerate}
\item they are contravariant in $X$,
\item they satisfy the identity $c_{X\times Y}^n=\sum_{r+s=n}c^r_X\otimes c^s_Y$,
\item they are normalized so that, when $X$ is of pure dimension $d$, the morphism
\[CH^d(X)_{\mathbb F}\to \End(\mathbbm 1),
\]obtained by composition of $c^d_X$ and $tr_X$ coincides with the \emph{degree morphism} on the $0$-cycles of $X$.
\end{enumerate}
\end{enumerate}
The prototypical example is given by taking 
\begin{itemize}
\item the category ${\rm CHM}(\mathbb K)_{\mathbb F}$, 
\item as invertible object the Lefschetz motive $\mathbb L$, 
\item the functor $\frak h(\cdot)_{\mathbb F}\colon \mathcal V_{\mathbb K}^{op}\to {\rm CHM}(\mathbb K)_{\mathbb F}$,
\item as morphisms $tr_X$ the correspondences, denoted $Tr_X$, defined by the graph $$[\Gamma_X]\in{\rm Corr}^0(X,{\rm Spec}(\mathbb K))$$ of the structural moprhisms $X\to{\rm Spec}(\mathbb K)$;
\item as morphisms $c^r_X$ the isomorphisms of equation \eqref{mot1}, denoted $\gamma^r_X$.
\end{itemize}

\begin{teorema}[Universality of Chow Motives] 
The 5-tuple $({\rm CHM}(\mathbb K)_{\mathbb F}, \mathbb L, \frak h(\cdot)_{\mathbb F}, (Tr_X)_X, (\gamma^r_X)_{r,X})$ is \emph{universal} among all other 5-tuples $(\mathcal T, {\bf L}, H, (tr_X)_X, (c^r_X)_{r,X})$ verifying points (1)-(5) above. This means that any functor $H\colon \mathcal V_{\mathbb K}^{op}\to\mathcal T$ as above factorizes through $\frak h(\cdot)_{\mathbb F}$: there exists a \emph{realization} $\omega_H\colon {\rm CHM}(\mathbb K)_{\mathbb F}\to\mathcal T$ which makes commutative the diagram
\[\xymatrix{\mathcal V_{\mathbb K}^{op}\ar[rr]^{\frak h(\cdot)_{\mathbb F}}\ar@/_1pc/[rrrr]_{H}&&{\rm CHM}(\mathbb K)_{\mathbb F}\ar[rr]^{\omega_H}&&\mathcal T}
\]Moreover, we have that ${\bf L}=\omega_H(\mathbb L)$, $\omega_H\circ Tr_X=tr_X\circ \omega_H$ and the morphisms $c^r_X$ are induced by $\omega_H$.
\end{teorema}

If $\mathcal T={\rm GrVec}_{\mathbb F}^{<\infty,+}$ (category of finite dimensional graded $\mathbb F$-vector spaces with only non-negative degrees), the functors $H\colon \mathcal V_{\mathbb K}^{op}\to \mathcal T$ as above, such that
\[ Gr^i\left(\omega_H(\mathbb L)\right)=0, \quad i\neq 0,2,
\]are called \emph{Weil cohomology functors}. Examples are Betti and (algebraic) de Rham realizations (if ${\rm char}(\mathbb K)=0$), crystalline cohomology (if ${\rm char}(\mathbb K)=p$) or étale cohomology (both ${\rm char}(\mathbb K)=0$ and ${\rm char}(\mathbb K)=p$). See \cite{andre} for further details.

If $\mathbb F=\mathbb Q,\mathbb K=\mathbb C$ we can consider $\mathcal T= HS_{\mathbb Q}$ the category of $\mathbb Q$-Hodge structures, whose objects are finite dimensional $\mathbb Q$-vector spaces $V$ whose complexification $V_{\mathbb C}$ admits a decomposition 
\[V_{\mathbb C}=\bigoplus_{p,q\in\mathbb Z}V^{p,q},\quad \text{such that }V^{p,q}=\overline{V^{q,p}}.
\]Such a category is a neutral Tannakian category over $\mathbb Q$ (see \cite{saav}). Hodge theory defines a canonical rational Hodge structure on the rational Betti cohomology $H_{B}(X,\mathbb Q)$ of any complex smooth projective variety $X$. Hence, the Betti cohomology functor can be enriched to a functor with values in $HS_{\mathbb Q}$, and the corresponding realization on ${\rm CHM}(\mathbb C)_{\mathbb Q}$ is called the \emph{Hodge realization}.

\newpage
\section{Tabulation of Stokes matrices for $\mathbb G(r,k)$ for small $k$}
In this appendix we tabulate all the Stokes matrices computed along the small quantum cohomology of Grassmannians $\mathbb G(r,k)$ for $k\leq 5$, with respect to an oriented line of slope $\phi\in\mathbb R$, and for a suitable choice of the branch of the $\Psi$-matrix. From these tables, the quasi-periodicity properties proved in Section \ref{symcpn} and Section \ref{symgrass} are evident. The matrices are obtained in the following way: the matrix $S$ for $\mathbb P^{k-1}_\mathbb C$ with $0<{\rm Im}(t)+k\phi<\pi$ is the one computed by the third author in \cite{guzzetti1}. The other Stokes matrices of $\mathbb P^{k-1}_\mathbb C$ are obtained through an action of the braids $\omega_{1,k},\omega_{2,k}$ described in Section \ref{secstoksqcohcpn}. The Stokes matrices for $\mathbb G(r,k)$ are obtained by applying Corollary \ref{08.02.18-1}. Colors keep track of the shifts of the quantum Satake identification: a matrix in the $r$-th column is the $r$-th exterior power of the matrix in the first column and of the same color.

\begin{table}[ht!]
\centering
\caption{Case $k=2$}
\begin{tabular}{|l|c|l|}
\hline
                                                & \multicolumn{2}{c|}{$\mathbb P^1_\mathbb C$}                                        \\ \hline
\multicolumn{1}{|c|}{$0<{\rm Im}(t)+2\phi<\pi$} & \multicolumn{2}{c|}{$\left(\begin{array}{cc} 1 & 2 \\ 0 & 1 \\\end{array}\right)$}  \\ \hline
$\pi<{\rm Im}(t)+2\phi<2\pi$                    & \multicolumn{2}{c|}{$\left(\begin{array}{cc} 1 & -2 \\ 0 & 1 \\\end{array}\right)$} \\ \hline
\end{tabular}
\end{table}

\begin{table}[ht!]
\centering
\caption{Case $k=3$}
\begin{tabular}{|c|c|l|l|c|l|l|}
\hline
\multicolumn{1}{|l|}{}                               & \multicolumn{3}{c|}{$\mathbb P^2_\mathbb C$}                                                                                      & \multicolumn{3}{c|}{$\mathbb G(2,3)$}                                                                                             \\ \hline
$0<{\rm Im}(t)+3\phi<\pi$                            & \multicolumn{3}{c|}{$\left(\begin{array}{ccc} 1 & 3 & -3 \\ 0 & 1 & -3 \\ 0 & 0 & 1 \\\end{array}\right)$}                        & \multicolumn{3}{c|}{{\color[HTML]{FE0000} $\left(\begin{array}{ccc} 1 & 3 & -3 \\ 0 & 1 & -3 \\ 0 & 0 & 1 \\\end{array}\right)$}} \\ \hline
$\pi<{\rm Im}(t)+3\phi<2\pi$                         & \multicolumn{3}{c|}{{\color[HTML]{FE0000} $\left(\begin{array}{ccc} 1 & -3 & -6 \\ 0 & 1 & 3 \\ 0 & 0 & 1 \\\end{array}\right)$}} & \multicolumn{3}{c|}{{\color[HTML]{009901} $\left(\begin{array}{ccc} 1 & 3 & 6 \\ 0 & 1 & 3 \\ 0 & 0 & 1 \\\end{array}\right)$}}   \\ \hline
$2\pi<{\rm Im}(t)+3\phi<3\pi$                        & \multicolumn{3}{c|}{{\color[HTML]{009901} $\left(\begin{array}{ccc} 1 & 3 & 3 \\ 0 & 1 & 3 \\ 0 & 0 & 1 \\\end{array}\right)$}}   & \multicolumn{3}{c|}{{\color[HTML]{00D2CB} $\left(\begin{array}{ccc} 1 & -3 & -3 \\ 0 & 1 & 3 \\ 0 & 0 & 1 \\\end{array}\right)$}} \\ \hline
$3\pi<{\rm Im}(t)+3\phi<4\pi$                        & \multicolumn{3}{c|}{{\color[HTML]{00D2CB} $\left(\begin{array}{ccc} 1 & 3 & -6 \\ 0 & 1 & -3 \\ 0 & 0 & 1 \\\end{array}\right)$}} & \multicolumn{3}{c|}{{\color[HTML]{6665CD} $\left(\begin{array}{ccc} 1 & 3 & -6 \\ 0 & 1 & -3 \\ 0 & 0 & 1 \\\end{array}\right)$}} \\ \hline
$4\pi<{\rm Im}(t)+3\phi<5\pi$                        & \multicolumn{3}{c|}{{\color[HTML]{6665CD} $\left(\begin{array}{ccc} 1 & -3 & -3 \\ 0 & 1 & 3 \\ 0 & 0 & 1 \\\end{array}\right)$}} & \multicolumn{3}{c|}{{\color[HTML]{FFC702} $\left(\begin{array}{ccc} 1 & -3 & 3 \\ 0 & 1 & -3 \\ 0 & 0 & 1 \\\end{array}\right)$}} \\ \hline
{\color[HTML]{000000} $5\pi<{\rm Im}(t)+3\phi<6\pi$} & \multicolumn{3}{c|}{{\color[HTML]{FFC702} $\left(\begin{array}{ccc} 1 & -3 & 6 \\ 0 & 1 & -3 \\ 0 & 0 & 1 \\\end{array}\right)$}} & \multicolumn{3}{c|}{$\left(\begin{array}{ccc} 1 & -3 & -6 \\ 0 & 1 & 3 \\ 0 & 0 & 1 \\\end{array}\right)$}                        \\ \hline
\end{tabular}
\end{table}

\begin{table}[]
\centering
\caption{Case $k=4$}
\begin{tabular}{|c|c|l|l|c|l|l|c|l|l|}
\hline
${\rm Im}(t)+4\phi$                  & \multicolumn{3}{c|}{$\mathbb P^3_\mathbb C$}                                                                                                                        & \multicolumn{3}{c|}{$\mathbb G(2,4)$}                                                                                                                                                                                                                          & \multicolumn{3}{c|}{$\mathbb G(3,4)$}                                                                                                                               \\ \hline
$]0;\pi[$                            & \multicolumn{3}{c|}{$\left(\begin{array}{cccc} 1 & -4 & -20 & 10 \\ 0 & 1 & 6 & -4 \\ 0 & 0 & 1 & -4 \\ 0 & 0 & 0 & 1 \\\end{array}\right)$}                        & \multicolumn{3}{c|}{{\color[HTML]{FE0000} $\left(\begin{array}{cccccc} 1 & -6 & -20 & -20 & -70 & 20 \\ 0 & 1 & 4 & 4 & 16 & -6 \\ 0 & 0 & 1 & 0 & 4 & -4 \\ 0 & 0 & 0 & 1 & 4 & -4 \\ 0 & 0 & 0 & 0 & 1 & -6 \\ 0 & 0 & 0 & 0 & 0 & 1 \\\end{array}\right)$}} & \multicolumn{3}{c|}{{\color[HTML]{009901} $\left(\begin{array}{cccc} 1 & 4 & 20 & -10 \\ 0 & 1 & 6 & -4 \\ 0 & 0 & 1 & -4 \\ 0 & 0 & 0 & 1 \\\end{array}\right)$}}  \\ \hline
$]\pi;2\pi[$                         & \multicolumn{3}{c|}{{\color[HTML]{FE0000} $\left(\begin{array}{cccc} 1 & 4 & -4 & -10 \\ 0 & 1 & -6 & -20 \\ 0 & 0 & 1 & 4 \\ 0 & 0 & 0 & 1 \\\end{array}\right)$}} & \multicolumn{3}{c|}{{\color[HTML]{009901} $\left(\begin{array}{cccccc} 1 & 6 & 4 & -4 & -6 & -20 \\ 0 & 1 & 4 & -4 & -16 & -70 \\ 0 & 0 & 1 & 0 & -4 & -20 \\ 0 & 0 & 0 & 1 & 4 & 20 \\ 0 & 0 & 0 & 0 & 1 & 6 \\ 0 & 0 & 0 & 0 & 0 & 1 \\\end{array}\right)$}} & \multicolumn{3}{c|}{{\color[HTML]{00D2CB} $\left(\begin{array}{cccc} 1 & -4 & -4 & -10 \\ 0 & 1 & 6 & 20 \\ 0 & 0 & 1 & 4 \\ 0 & 0 & 0 & 1 \\\end{array}\right)$}}  \\ \hline
$]2\pi;3\pi[$                        & \multicolumn{3}{c|}{{\color[HTML]{009901} $\left(\begin{array}{cccc} 1 & -4 & -20 & -10 \\ 0 & 1 & 6 & 4 \\ 0 & 0 & 1 & 4 \\ 0 & 0 & 0 & 1 \\\end{array}\right)$}}  & \multicolumn{3}{c|}{{\color[HTML]{00D2CB} $\left(\begin{array}{cccccc} 1 & 6 & -20 & 20 & -70 & -20 \\ 0 & 1 & -4 & 4 & -16 & -6 \\ 0 & 0 & 1 & 0 & 4 & 4 \\ 0 & 0 & 0 & 1 & -4 & -4 \\ 0 & 0 & 0 & 0 & 1 & 6 \\ 0 & 0 & 0 & 0 & 0 & 1 \\\end{array}\right)$}} & \multicolumn{3}{c|}{{\color[HTML]{6665CD} $\left(\begin{array}{cccc} 1 & 4 & -20 & -10 \\ 0 & 1 & -6 & -4 \\ 0 & 0 & 1 & 4 \\ 0 & 0 & 0 & 1 \\\end{array}\right)$}} \\ \hline
$]3\pi;4\pi[$                        & \multicolumn{3}{c|}{{\color[HTML]{00D2CB} $\left(\begin{array}{cccc} 1 & 4 & 4 & -10 \\ 0 & 1 & 6 & -20 \\ 0 & 0 & 1 & -4 \\ 0 & 0 & 0 & 1 \\\end{array}\right)$}}  & \multicolumn{3}{c|}{{\color[HTML]{6665CD} $\left(\begin{array}{cccccc} 1 & -6 & -4 & -4 & -6 & 20 \\ 0 & 1 & 4 & 4 & 16 & -70 \\ 0 & 0 & 1 & 0 & 4 & -20 \\ 0 & 0 & 0 & 1 & 4 & -20 \\ 0 & 0 & 0 & 0 & 1 & -6 \\ 0 & 0 & 0 & 0 & 0 & 1 \\\end{array}\right)$}} & \multicolumn{3}{c|}{{\color[HTML]{FFC702} $\left(\begin{array}{cccc} 1 & -4 & -4 & 10 \\ 0 & 1 & 6 & -20 \\ 0 & 0 & 1 & -4 \\ 0 & 0 & 0 & 1 \\\end{array}\right)$}} \\ \hline
$]4\pi;5\pi[$                        & \multicolumn{3}{c|}{{\color[HTML]{6665CD} $\left(\begin{array}{cccc} 1 & 4 & -20 & -10 \\ 0 & 1 & -6 & -4 \\ 0 & 0 & 1 & 4 \\ 0 & 0 & 0 & 1 \\\end{array}\right)$}} & \multicolumn{3}{c|}{{\color[HTML]{FFC702} $\left(\begin{array}{cccccc} 1 & 6 & -20 & -20 & 70 & 20 \\ 0 & 1 & -4 & -4 & 16 & 6 \\ 0 & 0 & 1 & 0 & -4 & -4 \\ 0 & 0 & 0 & 1 & -4 & -4 \\ 0 & 0 & 0 & 0 & 1 & 6 \\ 0 & 0 & 0 & 0 & 0 & 1 \\\end{array}\right)$}} & \multicolumn{3}{c|}{{\color[HTML]{2914E3} $\left(\begin{array}{cccc} 1 & 4 & -20 & 10 \\ 0 & 1 & -6 & 4 \\ 0 & 0 & 1 & -4 \\ 0 & 0 & 0 & 1 \\\end{array}\right)$}}  \\ \hline
{\color[HTML]{000000} $]5\pi;6\pi[$} & \multicolumn{3}{c|}{{\color[HTML]{FFC702} $\left(\begin{array}{cccc} 1 & -4 & -4 & 10 \\ 0 & 1 & 6 & -20 \\ 0 & 0 & 1 & -4 \\ 0 & 0 & 0 & 1 \\\end{array}\right)$}} & \multicolumn{3}{c|}{{\color[HTML]{2914E3} $\left(\begin{array}{cccccc} 1 & -6 & -4 & 4 & 6 & -20 \\ 0 & 1 & 4 & -4 & -16 & 70 \\ 0 & 0 & 1 & 0 & -4 & 20 \\ 0 & 0 & 0 & 1 & 4 & -20 \\ 0 & 0 & 0 & 0 & 1 & -6 \\ 0 & 0 & 0 & 0 & 0 & 1 \\\end{array}\right)$}} & \multicolumn{3}{c|}{{\color[HTML]{CB0000} $\left(\begin{array}{cccc} 1 & -4 & 4 & 10 \\ 0 & 1 & -6 & -20 \\ 0 & 0 & 1 & 4 \\ 0 & 0 & 0 & 1 \\\end{array}\right)$}}  \\ \hline
$]6\pi;7\pi[$                        & \multicolumn{3}{c|}{{\color[HTML]{2914E3} $\left(\begin{array}{cccc} 1 & -4 & 20 & 10 \\ 0 & 1 & -6 & -4 \\ 0 & 0 & 1 & 4 \\ 0 & 0 & 0 & 1 \\\end{array}\right)$}}  & \multicolumn{3}{c|}{{\color[HTML]{CB0000} $\left(\begin{array}{cccccc} 1 & -6 & 20 & -20 & 70 & -20 \\ 0 & 1 & -4 & 4 & -16 & 6 \\ 0 & 0 & 1 & 0 & 4 & -4 \\ 0 & 0 & 0 & 1 & -4 & 4 \\ 0 & 0 & 0 & 0 & 1 & -6 \\ 0 & 0 & 0 & 0 & 0 & 1 \\\end{array}\right)$}} & \multicolumn{3}{c|}{$\left(\begin{array}{cccc} 1 & -4 & -20 & 10 \\ 0 & 1 & 6 & -4 \\ 0 & 0 & 1 & -4 \\ 0 & 0 & 0 & 1 \\\end{array}\right)$}                        \\ \hline
$]7\pi;8\pi[$                        & \multicolumn{3}{c|}{{\color[HTML]{CB0000} $\left(\begin{array}{cccc} 1 & 4 & -4 & 10 \\ 0 & 1 & -6 & 20 \\ 0 & 0 & 1 & -4 \\ 0 & 0 & 0 & 1 \\\end{array}\right)$}}  & \multicolumn{3}{c|}{$\left(\begin{array}{cccccc} 1 & 6 & -4 & -4 & 6 & 20 \\ 0 & 1 & -4 & -4 & 16 & 70 \\ 0 & 0 & 1 & 0 & -4 & -20 \\ 0 & 0 & 0 & 1 & -4 & -20 \\ 0 & 0 & 0 & 0 & 1 & 6 \\ 0 & 0 & 0 & 0 & 0 & 1 \\\end{array}\right)$}                        & \multicolumn{3}{c|}{{\color[HTML]{FE0000} $\left(\begin{array}{cccc} 1 & 4 & -4 & -10 \\ 0 & 1 & -6 & -20 \\ 0 & 0 & 1 & 4 \\ 0 & 0 & 0 & 1 \\\end{array}\right)$}} \\ \hline
\end{tabular}
\end{table}

\begin{table}[]
\begin{adjustwidth}{-2cm}{-2cm}  
\begin{center}
\begin{tabular}{|c|c|l|l|c|l|l|}
\hline
${\rm Im}(t)+5\phi$ & \multicolumn{3}{c|}{$\mathbb P^4_\mathbb C$}                                                                                                                                              & \multicolumn{3}{c|}{$\mathbb G(2,5)$}                                                                                                                                                                                                                                                                                                                                                                                                                                                                                                                          \\ \hline
$]0;\pi[$           & \multicolumn{3}{c|}{$\left(\begin{array}{ccccc} 1 & 5 & -5 & -40 & 15 \\ 0 & 1 & -10 & -95 & 40 \\ 0 & 0 & 1 & 10 & -5 \\ 0 & 0 & 0 & 1 & -5 \\ 0 & 0 & 0 & 0 & 1 \\\end{array}\right)$}  & \multicolumn{3}{c|}{{\color{red}$\left(\begin{array}{cccccccccc} 1 & 10 & -5 & -15 & -5 & 10 & 40 & 75 & 325 & -50 \\ 0 & 1 & -10 & -45 & -5 & 50 & 225 & 435 & 1990 & -325 \\ 0 & 0 & 1 & 5 & 0 & -5 & -25 & -45 & -225 & 40 \\ 0 & 0 & 0 & 1 & 0 & 0 & -5 & 0 & -45 & 15 \\ 0 & 0 & 0 & 0 & 1 & -10 & -45 & -95 & -435 & 75 \\ 0 & 0 & 0 & 0 & 0 & 1 & 5 & 10 & 50 & -10 \\ 0 & 0 & 0 & 0 & 0 & 0 & 1 & 0 & 10 & -5 \\ 0 & 0 & 0 & 0 & 0 & 0 & 0 & 1 & 5 & -5 \\ 0 & 0 & 0 & 0 & 0 & 0 & 0 & 0 & 1 & -10 \\ 0 & 0 & 0 & 0 & 0 & 0 & 0 & 0 & 0 & 1 \\\end{array}\right)$}} \\ \hline
$]\pi;2\pi[$        & \multicolumn{3}{c|}{{\color{red}$\left(\begin{array}{ccccc} 1 & -5 & -45 & 15 & 35 \\ 0 & 1 & 10 & -5 & -15 \\ 0 & 0 & 1 & -10 & -45 \\ 0 & 0 & 0 & 1 & 5 \\ 0 & 0 & 0 & 0 & 1 \\\end{array}\right)$}} & \multicolumn{3}{c|}{{\color{blue}$\left(\begin{array}{cccccccccc} 1 & -10 & -95 & -40 & -45 & -435 & -185 & 75 & 50 & 175 \\ 0 & 1 & 10 & 5 & 5 & 50 & 25 & -10 & -10 & -50 \\ 0 & 0 & 1 & 5 & 0 & 5 & 25 & -5 & -25 & -185 \\ 0 & 0 & 0 & 1 & 0 & 0 & 5 & 0 & -5 & -40 \\ 0 & 0 & 0 & 0 & 1 & 10 & 5 & -5 & -10 & -75 \\ 0 & 0 & 0 & 0 & 0 & 1 & 5 & -10 & -50 & -435 \\ 0 & 0 & 0 & 0 & 0 & 0 & 1 & 0 & -10 & -95 \\ 0 & 0 & 0 & 0 & 0 & 0 & 0 & 1 & 5 & 45 \\ 0 & 0 & 0 & 0 & 0 & 0 & 0 & 0 & 1 & 10 \\ 0 & 0 & 0 & 0 & 0 & 0 & 0 & 0 & 0 & 1 \\\end{array}\right)$}}     \\ \hline
$]2\pi;3\pi[$       & \multicolumn{3}{c|}{{\color{blue}$\left(\begin{array}{ccccc} 1 & 5 & -5 & -40 & -15 \\ 0 & 1 & -10 & -95 & -40 \\ 0 & 0 & 1 & 10 & 5 \\ 0 & 0 & 0 & 1 & 5 \\ 0 & 0 & 0 & 0 & 1 \\\end{array}\right)$}}  & \multicolumn{3}{c|}{{\color{auburn}$\left(\begin{array}{cccccccccc} 1 & 10 & 5 & -15 & -5 & -10 & 40 & -75 & 325 & 50 \\ 0 & 1 & 10 & -45 & -5 & -50 & 225 & -435 & 1990 & 325 \\ 0 & 0 & 1 & -5 & 0 & -5 & 25 & -45 & 225 & 40 \\ 0 & 0 & 0 & 1 & 0 & 0 & -5 & 0 & -45 & -15 \\ 0 & 0 & 0 & 0 & 1 & 10 & -45 & 95 & -435 & -75 \\ 0 & 0 & 0 & 0 & 0 & 1 & -5 & 10 & -50 & -10 \\ 0 & 0 & 0 & 0 & 0 & 0 & 1 & 0 & 10 & 5 \\ 0 & 0 & 0 & 0 & 0 & 0 & 0 & 1 & -5 & -5 \\ 0 & 0 & 0 & 0 & 0 & 0 & 0 & 0 & 1 & 10 \\ 0 & 0 & 0 & 0 & 0 & 0 & 0 & 0 & 0 & 1 \\\end{array}\right)$}} \\ \hline
$]3\pi;4\pi[$       & \multicolumn{3}{c|}{{\color{auburn}$\left(\begin{array}{ccccc} 1 & -5 & -45 & -15 & 35 \\ 0 & 1 & 10 & 5 & -15 \\ 0 & 0 & 1 & 10 & -45 \\ 0 & 0 & 0 & 1 & -5 \\ 0 & 0 & 0 & 0 & 1 \\\end{array}\right)$}} & \multicolumn{3}{c|}{{\color{ambra}$\left(\begin{array}{cccccccccc} 1 & 10 & -95 & -40 & 45 & -435 & -185 & -75 & -50 & 175 \\ 0 & 1 & -10 & -5 & 5 & -50 & -25 & -10 & -10 & 50 \\ 0 & 0 & 1 & 5 & 0 & 5 & 25 & 5 & 25 & -185 \\ 0 & 0 & 0 & 1 & 0 & 0 & 5 & 0 & 5 & -40 \\ 0 & 0 & 0 & 0 & 1 & -10 & -5 & -5 & -10 & 75 \\ 0 & 0 & 0 & 0 & 0 & 1 & 5 & 10 & 50 & -435 \\ 0 & 0 & 0 & 0 & 0 & 0 & 1 & 0 & 10 & -95 \\ 0 & 0 & 0 & 0 & 0 & 0 & 0 & 1 & 5 & -45 \\ 0 & 0 & 0 & 0 & 0 & 0 & 0 & 0 & 1 & -10 \\ 0 & 0 & 0 & 0 & 0 & 0 & 0 & 0 & 0 & 1 \\\end{array}\right)$}}     \\ \hline
$]4\pi;5\pi[$       & \multicolumn{3}{c|}{{\color{ambra}$\left(\begin{array}{ccccc} 1 & 5 & 5 & -40 & -15 \\ 0 & 1 & 10 & -95 & -40 \\ 0 & 0 & 1 & -10 & -5 \\ 0 & 0 & 0 & 1 & 5 \\ 0 & 0 & 0 & 0 & 1 \\\end{array}\right)$}}  & \multicolumn{3}{c|}{{\color{ametista}$\left(\begin{array}{cccccccccc} 1 & -10 & -5 & 15 & -5 & -10 & 40 & 75 & -325 & -50 \\ 0 & 1 & 10 & -45 & 5 & 50 & -225 & -435 & 1990 & 325 \\ 0 & 0 & 1 & -5 & 0 & 5 & -25 & -45 & 225 & 40 \\ 0 & 0 & 0 & 1 & 0 & 0 & 5 & 0 & -45 & -15 \\ 0 & 0 & 0 & 0 & 1 & 10 & -45 & -95 & 435 & 75 \\ 0 & 0 & 0 & 0 & 0 & 1 & -5 & -10 & 50 & 10 \\ 0 & 0 & 0 & 0 & 0 & 0 & 1 & 0 & -10 & -5 \\ 0 & 0 & 0 & 0 & 0 & 0 & 0 & 1 & -5 & -5 \\ 0 & 0 & 0 & 0 & 0 & 0 & 0 & 0 & 1 & 10 \\ 0 & 0 & 0 & 0 & 0 & 0 & 0 & 0 & 0 & 1 \\\end{array}\right)$}} \\ \hline
\end{tabular}
\caption{Case $k=5$ (first part)}
\end{center}
\end{adjustwidth}
\end{table}

\begin{table}[]
\begin{adjustwidth}{-2cm}{-2cm}  
\begin{center}
\begin{tabular}{|c|c|l|l|c|l|l|}
\hline
${\rm Im}(t)+5\phi$ & \multicolumn{3}{c|}{$\mathbb G(3,5)$}                                                                                                                                                                                                                                                                                                                                                                                                                                                                                                                          & \multicolumn{3}{c|}{$\mathbb G(4,5)$}                                                                                                                                                     \\ \hline
$]0;\pi[$           & \multicolumn{3}{c|}{{\color{blue}$\left(\begin{array}{cccccccccc} 1 & 10 & 5 & -5 & -10 & -75 & -15 & -40 & -325 & 50 \\ 0 & 1 & 5 & -10 & -50 & -435 & -45 & -225 & -1990 & 325 \\ 0 & 0 & 1 & 0 & -10 & -95 & 0 & -45 & -435 & 75 \\ 0 & 0 & 0 & 1 & 5 & 45 & 5 & 25 & 225 & -40 \\ 0 & 0 & 0 & 0 & 1 & 10 & 0 & 5 & 50 & -10 \\ 0 & 0 & 0 & 0 & 0 & 1 & 0 & 0 & 5 & -5 \\ 0 & 0 & 0 & 0 & 0 & 0 & 1 & 5 & 45 & -15 \\ 0 & 0 & 0 & 0 & 0 & 0 & 0 & 1 & 10 & -5 \\ 0 & 0 & 0 & 0 & 0 & 0 & 0 & 0 & 1 & -10 \\ 0 & 0 & 0 & 0 & 0 & 0 & 0 & 0 & 0 & 1 \\\end{array}\right)$}} & \multicolumn{3}{c|}{{\color{auburn}$\left(\begin{array}{ccccc} 1 & -5 & -5 & -40 & 15 \\ 0 & 1 & 10 & 95 & -40 \\ 0 & 0 & 1 & 10 & -5 \\ 0 & 0 & 0 & 1 & -5 \\ 0 & 0 & 0 & 0 & 1 \\\end{array}\right)$}}  \\ \hline
$]\pi;2\pi[$        & \multicolumn{3}{c|}{{\color{auburn}$\left(\begin{array}{cccccccccc} 1 & 10 & -45 & 95 & -435 & -75 & -40 & 185 & 50 & 175 \\ 0 & 1 & -5 & 10 & -50 & -10 & -5 & 25 & 10 & 50 \\ 0 & 0 & 1 & 0 & 10 & 5 & 0 & -5 & -10 & -75 \\ 0 & 0 & 0 & 1 & -5 & -5 & -5 & 25 & 25 & 185 \\ 0 & 0 & 0 & 0 & 1 & 10 & 0 & -5 & -50 & -435 \\ 0 & 0 & 0 & 0 & 0 & 1 & 0 & 0 & -5 & -45 \\ 0 & 0 & 0 & 0 & 0 & 0 & 1 & -5 & -5 & -40 \\ 0 & 0 & 0 & 0 & 0 & 0 & 0 & 1 & 10 & 95 \\ 0 & 0 & 0 & 0 & 0 & 0 & 0 & 0 & 1 & 10 \\ 0 & 0 & 0 & 0 & 0 & 0 & 0 & 0 & 0 & 1 \\\end{array}\right)$}}     & \multicolumn{3}{c|}{{\color{ambra}$\left(\begin{array}{ccccc} 1 & 5 & -45 & -15 & -35 \\ 0 & 1 & -10 & -5 & -15 \\ 0 & 0 & 1 & 10 & 45 \\ 0 & 0 & 0 & 1 & 5 \\ 0 & 0 & 0 & 0 & 1 \\\end{array}\right)$}} \\ \hline
$]2\pi;3\pi[$       & \multicolumn{3}{c|}{{\color{ambra}$\left(\begin{array}{cccccccccc} 1 & -10 & -5 & -5 & -10 & 75 & -15 & -40 & 325 & 50 \\ 0 & 1 & 5 & 10 & 50 & -435 & 45 & 225 & -1990 & -325 \\ 0 & 0 & 1 & 0 & 10 & -95 & 0 & 45 & -435 & -75 \\ 0 & 0 & 0 & 1 & 5 & -45 & 5 & 25 & -225 & -40 \\ 0 & 0 & 0 & 0 & 1 & -10 & 0 & 5 & -50 & -10 \\ 0 & 0 & 0 & 0 & 0 & 1 & 0 & 0 & 5 & 5 \\ 0 & 0 & 0 & 0 & 0 & 0 & 1 & 5 & -45 & -15 \\ 0 & 0 & 0 & 0 & 0 & 0 & 0 & 1 & -10 & -5 \\ 0 & 0 & 0 & 0 & 0 & 0 & 0 & 0 & 1 & 10 \\ 0 & 0 & 0 & 0 & 0 & 0 & 0 & 0 & 0 & 1 \\\end{array}\right)$}} & \multicolumn{3}{c|}{{\color{ametista}$\left(\begin{array}{ccccc} 1 & -5 & -5 & 40 & 15 \\ 0 & 1 & 10 & -95 & -40 \\ 0 & 0 & 1 & -10 & -5 \\ 0 & 0 & 0 & 1 & 5 \\ 0 & 0 & 0 & 0 & 1 \\\end{array}\right)$}}  \\ \hline
$]3\pi;4\pi[$       & \multicolumn{3}{c|}{{\color{ametista}$\left(\begin{array}{cccccccccc} 1 & 10 & -45 & -95 & 435 & 75 & -40 & 185 & 50 & -175 \\ 0 & 1 & -5 & -10 & 50 & 10 & -5 & 25 & 10 & -50 \\ 0 & 0 & 1 & 0 & -10 & -5 & 0 & -5 & -10 & 75 \\ 0 & 0 & 0 & 1 & -5 & -5 & 5 & -25 & -25 & 185 \\ 0 & 0 & 0 & 0 & 1 & 10 & 0 & 5 & 50 & -435 \\ 0 & 0 & 0 & 0 & 0 & 1 & 0 & 0 & 5 & -45 \\ 0 & 0 & 0 & 0 & 0 & 0 & 1 & -5 & -5 & 40 \\ 0 & 0 & 0 & 0 & 0 & 0 & 0 & 1 & 10 & -95 \\ 0 & 0 & 0 & 0 & 0 & 0 & 0 & 0 & 1 & -10 \\ 0 & 0 & 0 & 0 & 0 & 0 & 0 & 0 & 0 & 1 \\\end{array}\right)$}}     & \multicolumn{3}{c|}{{\color{ballblue}$\left(\begin{array}{ccccc} 1 & 5 & -45 & -15 & 35 \\ 0 & 1 & -10 & -5 & 15 \\ 0 & 0 & 1 & 10 & -45 \\ 0 & 0 & 0 & 1 & -5 \\ 0 & 0 & 0 & 0 & 1 \\\end{array}\right)$}} \\ \hline
$]4\pi;5\pi[$       & \multicolumn{3}{c|}{{\color{ballblue}$\left(\begin{array}{cccccccccc} 1 & -10 & -5 & -5 & -10 & 75 & 15 & 40 & -325 & -50 \\ 0 & 1 & 5 & 10 & 50 & -435 & -45 & -225 & 1990 & 325 \\ 0 & 0 & 1 & 0 & 10 & -95 & 0 & -45 & 435 & 75 \\ 0 & 0 & 0 & 1 & 5 & -45 & -5 & -25 & 225 & 40 \\ 0 & 0 & 0 & 0 & 1 & -10 & 0 & -5 & 50 & 10 \\ 0 & 0 & 0 & 0 & 0 & 1 & 0 & 0 & -5 & -5 \\ 0 & 0 & 0 & 0 & 0 & 0 & 1 & 5 & -45 & -15 \\ 0 & 0 & 0 & 0 & 0 & 0 & 0 & 1 & -10 & -5 \\ 0 & 0 & 0 & 0 & 0 & 0 & 0 & 0 & 1 & 10 \\ 0 & 0 & 0 & 0 & 0 & 0 & 0 & 0 & 0 & 1 \\\end{array}\right)$}} & \multicolumn{3}{c|}{{\color{cadmiumgreen}$\left(\begin{array}{ccccc} 1 & -5 & -5 & 40 & -15 \\ 0 & 1 & 10 & -95 & 40 \\ 0 & 0 & 1 & -10 & 5 \\ 0 & 0 & 0 & 1 & -5 \\ 0 & 0 & 0 & 0 & 1 \\\end{array}\right)$}}  \\ \hline
\end{tabular}
\caption{Case $k=5$ (second part)}
\end{center}
\end{adjustwidth}
\end{table}

\begin{table}[]
\begin{adjustwidth}{-2cm}{-2cm} 
\begin{center}
\begin{tabular}{|c|c|l|l|c|l|l|}
\hline
${\rm Im}(t)+5\phi$ & \multicolumn{3}{c|}{$\mathbb P^4_\mathbb C$}                                                                                                                                              & \multicolumn{3}{c|}{$\mathbb G(2,5)$}                                                                                                                                                                                                                                                                                                                                                                                                                                                                                                                          \\ \hline
$]5\pi;6\pi[$       & \multicolumn{3}{c|}{{\color{ametista}$\left(\begin{array}{ccccc} 1 & 5 & -45 & -15 & 35 \\ 0 & 1 & -10 & -5 & 15 \\ 0 & 0 & 1 & 10 & -45 \\ 0 & 0 & 0 & 1 & -5 \\ 0 & 0 & 0 & 0 & 1 \\\end{array}\right)$}} & \multicolumn{3}{c|}{{\color{ballblue}$\left(\begin{array}{cccccccccc} 1 & 10 & -95 & -40 & -45 & 435 & 185 & 75 & 50 & -175 \\ 0 & 1 & -10 & -5 & -5 & 50 & 25 & 10 & 10 & -50 \\ 0 & 0 & 1 & 5 & 0 & -5 & -25 & -5 & -25 & 185 \\ 0 & 0 & 0 & 1 & 0 & 0 & -5 & 0 & -5 & 40 \\ 0 & 0 & 0 & 0 & 1 & -10 & -5 & -5 & -10 & 75 \\ 0 & 0 & 0 & 0 & 0 & 1 & 5 & 10 & 50 & -435 \\ 0 & 0 & 0 & 0 & 0 & 0 & 1 & 0 & 10 & -95 \\ 0 & 0 & 0 & 0 & 0 & 0 & 0 & 1 & 5 & -45 \\ 0 & 0 & 0 & 0 & 0 & 0 & 0 & 0 & 1 & -10 \\ 0 & 0 & 0 & 0 & 0 & 0 & 0 & 0 & 0 & 1 \\\end{array}\right)$}}     \\ \hline
$]6\pi;7\pi[$       & \multicolumn{3}{c|}{{\color{ballblue}$\left(\begin{array}{ccccc} 1 & -5 & -5 & 40 & 15 \\ 0 & 1 & 10 & -95 & -40 \\ 0 & 0 & 1 & -10 & -5 \\ 0 & 0 & 0 & 1 & 5 \\ 0 & 0 & 0 & 0 & 1 \\\end{array}\right)$}}  & \multicolumn{3}{c|}{{\color{cadmiumgreen}$\left(\begin{array}{cccccccccc} 1 & -10 & -5 & 15 & 5 & 10 & -40 & -75 & 325 & 50 \\ 0 & 1 & 10 & -45 & -5 & -50 & 225 & 435 & -1990 & -325 \\ 0 & 0 & 1 & -5 & 0 & -5 & 25 & 45 & -225 & -40 \\ 0 & 0 & 0 & 1 & 0 & 0 & -5 & 0 & 45 & 15 \\ 0 & 0 & 0 & 0 & 1 & 10 & -45 & -95 & 435 & 75 \\ 0 & 0 & 0 & 0 & 0 & 1 & -5 & -10 & 50 & 10 \\ 0 & 0 & 0 & 0 & 0 & 0 & 1 & 0 & -10 & -5 \\ 0 & 0 & 0 & 0 & 0 & 0 & 0 & 1 & -5 & -5 \\ 0 & 0 & 0 & 0 & 0 & 0 & 0 & 0 & 1 & 10 \\ 0 & 0 & 0 & 0 & 0 & 0 & 0 & 0 & 0 & 1 \\\end{array}\right)$}} \\ \hline
$]7\pi;8\pi[$       & \multicolumn{3}{c|}{{\color{cadmiumgreen}$\left(\begin{array}{ccccc} 1 & -5 & 45 & 15 & -35 \\ 0 & 1 & -10 & -5 & 15 \\ 0 & 0 & 1 & 10 & -45 \\ 0 & 0 & 0 & 1 & -5 \\ 0 & 0 & 0 & 0 & 1 \\\end{array}\right)$}} & \multicolumn{3}{c|}{{\color{candypink}$\left(\begin{array}{cccccccccc} 1 & -10 & 95 & 40 & -45 & 435 & 185 & -75 & -50 & 175 \\ 0 & 1 & -10 & -5 & 5 & -50 & -25 & 10 & 10 & -50 \\ 0 & 0 & 1 & 5 & 0 & 5 & 25 & -5 & -25 & 185 \\ 0 & 0 & 0 & 1 & 0 & 0 & 5 & 0 & -5 & 40 \\ 0 & 0 & 0 & 0 & 1 & -10 & -5 & 5 & 10 & -75 \\ 0 & 0 & 0 & 0 & 0 & 1 & 5 & -10 & -50 & 435 \\ 0 & 0 & 0 & 0 & 0 & 0 & 1 & 0 & -10 & 95 \\ 0 & 0 & 0 & 0 & 0 & 0 & 0 & 1 & 5 & -45 \\ 0 & 0 & 0 & 0 & 0 & 0 & 0 & 0 & 1 & -10 \\ 0 & 0 & 0 & 0 & 0 & 0 & 0 & 0 & 0 & 1 \\\end{array}\right)$}}       \\ \hline
$]8\pi;9\pi[$       & \multicolumn{3}{c|}{{\color{candypink}$\left(\begin{array}{ccccc} 1 & 5 & -5 & 40 & 15 \\ 0 & 1 & -10 & 95 & 40 \\ 0 & 0 & 1 & -10 & -5 \\ 0 & 0 & 0 & 1 & 5 \\ 0 & 0 & 0 & 0 & 1 \\\end{array}\right)$}}    & \multicolumn{3}{c|}{{\color{caribbeangreen}$\left(\begin{array}{cccccccccc} 1 & 10 & -5 & 15 & -5 & 10 & -40 & 75 & -325 & 50 \\ 0 & 1 & -10 & 45 & -5 & 50 & -225 & 435 & -1990 & 325 \\ 0 & 0 & 1 & -5 & 0 & -5 & 25 & -45 & 225 & -40 \\ 0 & 0 & 0 & 1 & 0 & 0 & -5 & 0 & -45 & 15 \\ 0 & 0 & 0 & 0 & 1 & -10 & 45 & -95 & 435 & -75 \\ 0 & 0 & 0 & 0 & 0 & 1 & -5 & 10 & -50 & 10 \\ 0 & 0 & 0 & 0 & 0 & 0 & 1 & 0 & 10 & -5 \\ 0 & 0 & 0 & 0 & 0 & 0 & 0 & 1 & -5 & 5 \\ 0 & 0 & 0 & 0 & 0 & 0 & 0 & 0 & 1 & -10 \\ 0 & 0 & 0 & 0 & 0 & 0 & 0 & 0 & 0 & 1 \\\end{array}\right)$}} \\ \hline
$]9\pi;10\pi[$      & \multicolumn{3}{c|}{{\color{caribbeangreen}$\left(\begin{array}{ccccc} 1 & -5 & -45 & 15 & -35 \\ 0 & 1 & 10 & -5 & 15 \\ 0 & 0 & 1 & -10 & 45 \\ 0 & 0 & 0 & 1 & -5 \\ 0 & 0 & 0 & 0 & 1 \\\end{array}\right)$}} & \multicolumn{3}{c|}{$\left(\begin{array}{cccccccccc} 1 & -10 & -95 & 40 & -45 & -435 & 185 & 75 & -50 & -175 \\ 0 & 1 & 10 & -5 & 5 & 50 & -25 & -10 & 10 & 50 \\ 0 & 0 & 1 & -5 & 0 & 5 & -25 & -5 & 25 & 185 \\ 0 & 0 & 0 & 1 & 0 & 0 & 5 & 0 & -5 & -40 \\ 0 & 0 & 0 & 0 & 1 & 10 & -5 & -5 & 10 & 75 \\ 0 & 0 & 0 & 0 & 0 & 1 & -5 & -10 & 50 & 435 \\ 0 & 0 & 0 & 0 & 0 & 0 & 1 & 0 & -10 & -95 \\ 0 & 0 & 0 & 0 & 0 & 0 & 0 & 1 & -5 & -45 \\ 0 & 0 & 0 & 0 & 0 & 0 & 0 & 0 & 1 & 10 \\ 0 & 0 & 0 & 0 & 0 & 0 & 0 & 0 & 0 & 1 \\\end{array}\right)$}     \\ \hline
\end{tabular}
\caption{Case $k=5$ (third part)}
\end{center}
\end{adjustwidth}
\end{table}


\begin{table}[]
\begin{adjustwidth}{-2cm}{-2cm} 
\begin{center}
\begin{tabular}{|c|c|l|l|c|l|l|}
\hline
${\rm Im}(t)+5\phi$ & \multicolumn{3}{c|}{$\mathbb G(3,5)$}                                                                                                                                                                                                                                                                                                                                                                                                                                                                                                                          & \multicolumn{3}{c|}{$\mathbb G(4,5)$}                                                                                                                                                     \\ \hline
$]5\pi;6\pi[$       & \multicolumn{3}{c|}{{\color{cadmiumgreen}$\left(\begin{array}{cccccccccc} 1 & 10 & -45 & -95 & 435 & 75 & 40 & -185 & -50 & 175 \\ 0 & 1 & -5 & -10 & 50 & 10 & 5 & -25 & -10 & 50 \\ 0 & 0 & 1 & 0 & -10 & -5 & 0 & 5 & 10 & -75 \\ 0 & 0 & 0 & 1 & -5 & -5 & -5 & 25 & 25 & -185 \\ 0 & 0 & 0 & 0 & 1 & 10 & 0 & -5 & -50 & 435 \\ 0 & 0 & 0 & 0 & 0 & 1 & 0 & 0 & -5 & 45 \\ 0 & 0 & 0 & 0 & 0 & 0 & 1 & -5 & -5 & 40 \\ 0 & 0 & 0 & 0 & 0 & 0 & 0 & 1 & 10 & -95 \\ 0 & 0 & 0 & 0 & 0 & 0 & 0 & 0 & 1 & -10 \\ 0 & 0 & 0 & 0 & 0 & 0 & 0 & 0 & 0 & 1 \\\end{array}\right)$}}     & \multicolumn{3}{c|}{{\color{candypink}$\left(\begin{array}{ccccc} 1 & 5 & -45 & 15 & 35 \\ 0 & 1 & -10 & 5 & 15 \\ 0 & 0 & 1 & -10 & -45 \\ 0 & 0 & 0 & 1 & 5 \\ 0 & 0 & 0 & 0 & 1 \\\end{array}\right)$}}   \\ \hline
$]6\pi;7\pi[$       & \multicolumn{3}{c|}{{\color{candypink}$\left(\begin{array}{cccccccccc} 1 & -10 & -5 & 5 & 10 & -75 & 15 & 40 & -325 & 50 \\ 0 & 1 & 5 & -10 & -50 & 435 & -45 & -225 & 1990 & -325 \\ 0 & 0 & 1 & 0 & -10 & 95 & 0 & -45 & 435 & -75 \\ 0 & 0 & 0 & 1 & 5 & -45 & 5 & 25 & -225 & 40 \\ 0 & 0 & 0 & 0 & 1 & -10 & 0 & 5 & -50 & 10 \\ 0 & 0 & 0 & 0 & 0 & 1 & 0 & 0 & 5 & -5 \\ 0 & 0 & 0 & 0 & 0 & 0 & 1 & 5 & -45 & 15 \\ 0 & 0 & 0 & 0 & 0 & 0 & 0 & 1 & -10 & 5 \\ 0 & 0 & 0 & 0 & 0 & 0 & 0 & 0 & 1 & -10 \\ 0 & 0 & 0 & 0 & 0 & 0 & 0 & 0 & 0 & 1 \\\end{array}\right)$}}   & \multicolumn{3}{c|}{{\color{caribbeangreen}$\left(\begin{array}{ccccc} 1 & -5 & 5 & 40 & -15 \\ 0 & 1 & -10 & -95 & 40 \\ 0 & 0 & 1 & 10 & -5 \\ 0 & 0 & 0 & 1 & -5 \\ 0 & 0 & 0 & 0 & 1 \\\end{array}\right)$}}  \\ \hline
$]7\pi;8\pi[$       & \multicolumn{3}{c|}{{\color{caribbeangreen}$\left(\begin{array}{cccccccccc} 1 & -10 & 45 & -95 & 435 & -75 & 40 & -185 & 50 & 175 \\ 0 & 1 & -5 & 10 & -50 & 10 & -5 & 25 & -10 & -50 \\ 0 & 0 & 1 & 0 & 10 & -5 & 0 & -5 & 10 & 75 \\ 0 & 0 & 0 & 1 & -5 & 5 & -5 & 25 & -25 & -185 \\ 0 & 0 & 0 & 0 & 1 & -10 & 0 & -5 & 50 & 435 \\ 0 & 0 & 0 & 0 & 0 & 1 & 0 & 0 & -5 & -45 \\ 0 & 0 & 0 & 0 & 0 & 0 & 1 & -5 & 5 & 40 \\ 0 & 0 & 0 & 0 & 0 & 0 & 0 & 1 & -10 & -95 \\ 0 & 0 & 0 & 0 & 0 & 0 & 0 & 0 & 1 & 10 \\ 0 & 0 & 0 & 0 & 0 & 0 & 0 & 0 & 0 & 1 \\\end{array}\right)$}}     & \multicolumn{3}{c|}{$\left(\begin{array}{ccccc} 1 & -5 & -45 & 15 & 35 \\ 0 & 1 & 10 & -5 & -15 \\ 0 & 0 & 1 & -10 & -45 \\ 0 & 0 & 0 & 1 & 5 \\ 0 & 0 & 0 & 0 & 1 \\\end{array}\right)$} \\ \hline
$]8\pi;9\pi[$       & \multicolumn{3}{c|}{$\left(\begin{array}{cccccccccc} 1 & 10 & -5 & -5 & 10 & 75 & -15 & 40 & 325 & -50 \\ 0 & 1 & -5 & -10 & 50 & 435 & -45 & 225 & 1990 & -325 \\ 0 & 0 & 1 & 0 & -10 & -95 & 0 & -45 & -435 & 75 \\ 0 & 0 & 0 & 1 & -5 & -45 & 5 & -25 & -225 & 40 \\ 0 & 0 & 0 & 0 & 1 & 10 & 0 & 5 & 50 & -10 \\ 0 & 0 & 0 & 0 & 0 & 1 & 0 & 0 & 5 & -5 \\ 0 & 0 & 0 & 0 & 0 & 0 & 1 & -5 & -45 & 15 \\ 0 & 0 & 0 & 0 & 0 & 0 & 0 & 1 & 10 & -5 \\ 0 & 0 & 0 & 0 & 0 & 0 & 0 & 0 & 1 & -10 \\ 0 & 0 & 0 & 0 & 0 & 0 & 0 & 0 & 0 & 1 \\\end{array}\right)$} & \multicolumn{3}{c|}{{\color{red}$\left(\begin{array}{ccccc} 1 & 5 & -5 & -40 & 15 \\ 0 & 1 & -10 & -95 & 40 \\ 0 & 0 & 1 & 10 & -5 \\ 0 & 0 & 0 & 1 & -5 \\ 0 & 0 & 0 & 0 & 1 \\\end{array}\right)$}}  \\ \hline
$]9\pi;10\pi[$      & \multicolumn{3}{c|}{{\color{red}$\left(\begin{array}{cccccccccc} 1 & -10 & -45 & -95 & -435 & 75 & 40 & 185 & -50 & -175 \\ 0 & 1 & 5 & 10 & 50 & -10 & -5 & -25 & 10 & 50 \\ 0 & 0 & 1 & 0 & 10 & -5 & 0 & -5 & 10 & 75 \\ 0 & 0 & 0 & 1 & 5 & -5 & -5 & -25 & 25 & 185 \\ 0 & 0 & 0 & 0 & 1 & -10 & 0 & -5 & 50 & 435 \\ 0 & 0 & 0 & 0 & 0 & 1 & 0 & 0 & -5 & -45 \\ 0 & 0 & 0 & 0 & 0 & 0 & 1 & 5 & -5 & -40 \\ 0 & 0 & 0 & 0 & 0 & 0 & 0 & 1 & -10 & -95 \\ 0 & 0 & 0 & 0 & 0 & 0 & 0 & 0 & 1 & 10 \\ 0 & 0 & 0 & 0 & 0 & 0 & 0 & 0 & 0 & 1 \\\end{array}\right)$}}     & \multicolumn{3}{c|}{{\color{blue}$\left(\begin{array}{ccccc} 1 & 5 & 45 & -15 & -35 \\ 0 & 1 & 10 & -5 & -15 \\ 0 & 0 & 1 & -10 & -45 \\ 0 & 0 & 0 & 1 & 5 \\ 0 & 0 & 0 & 0 & 1 \\\end{array}\right)$}} \\ \hline
\end{tabular}
\caption{Case $k=5$ (fourth part)}
\end{center}
\end{adjustwidth}
\end{table}

\newpage

\bibliographystyle{alpha}
\bibliography{biblio}

\end{document}